\documentclass[reqno]{amsart}

\usepackage{
amssymb,
amsmath,
amsthm,
eucal,
dsfont,
multicol,
multirow,
mathrsfs,
tikz,
graphicx,
hyperref,
makecell
}
\usepackage{lipsum}
\makeatletter\renewcommand*{\eqref}[1]{%
\hyperref[{#1}]{\textup{\tagform@{\!\!\ref*{#1}}}}%
}\makeatother %add when putting this file to arXiv

\makeatletter
 
  \@addtoreset{equation}{section}
 \makeatother
\theoremstyle{plain}
\newtheorem{theorem}{Theorem}[section]
\newtheorem{lemma}[theorem]{Lemma}
\newtheorem{proposition}[theorem]{Proposition}

\newtheorem{corollary}[theorem]{Corollary}
\theoremstyle{definition}
\newtheorem{definition}[theorem]{Definition}
\newtheorem{remark}[theorem]{Remark}

\newtheorem{assumption}[theorem]{Assumption}

\def\supp{\mathop{\mathrm{supp}}\nolimits}

\def\Ker{\mathop{\mathrm{Ker}}\nolimits}

\def\Re{\mathop{\mathrm{Re}}\nolimits}
\def\Im{\mathop{\mathrm{Im}}\nolimits}

\def\sgn{\mathop{\mathrm{sgn}}\nolimits}

\def\R{{\mathbb{R}}}
\def\Z{{\mathbb{Z}}}
\def\N{{\mathbb{N}}}
\def\C{{\mathbb{C}}}

\def\k{{ \textbf{k}}}
\def\i{{\rm i}}
\def\d{{\rm d}}
\def\<{{\langle}}
\def\>{{\rangle}}

\title[$L^p$-boundedness of wave operators for higher order Schr\"odinger  operator ]
{The $L^p$-boundedness of wave operators for higher order  Schr\"odinger  operator with zero singularities in low  even dimensions }
\author{Han Cheng,  Avy Soffer, Zhao Wu, Xiaohua Yao}

\address {Han Cheng, Institute of Applied Physics and Computational Mathematics, Beijing, 100088, People's Republic of China}
\email{chmathh@163.com}

\address{Avy Soffer, Mathematics Department, Rutgers University, New Brunswick, NJ, 08903, USA}
\email{soffer@math.rutgers.edu}

\address {Zhao Wu, Faculty of Mathematics and Statistics, Hubei University, Wuhan, 430062, People's Republic of China}
\email{wuzhao@hubu.edu.cn}
\address {Xiaohua Yao, School of Mathematics and Statistics,   Key Laboratory of Nonlinear Analysis  and Applications (Ministry of Education), and Hubei Key Laboratory of Mathematical Sciences, Central China Normal University, Wuhan, 430079, P.R. China}
\email{yaoxiaohua@ccnu.edu.cn}
\date{\today}
\keywords{$L^p$-boundedness, Wave operator, Resonance, Schr\"odinger operator}
\begin{document}
\begin{abstract}

Let $H = (-\Delta)^m + V$ be  Schr\"odinger type operator, where $m \geq 1$ and $V$ is a decaying real-valued potential. We are devoted to establishing the  $L^p$-boundedness of the wave operators $W_\pm(H; (-\Delta)^m)$ with all types of zero-energy singularities in low dimensions $1 \leq n \leq 4m$. 

 In view of the recent work \cite{CSWY24} addressing the {\it odd-dimensional cases}, this paper extends the analysis to even dimensions, proving that $W_\pm(H; (-\Delta)^m) \in \mathbb{B}(L^p)$ with sharp ranges of $p$ determined by the nature of the zero-energy singularities for {\it all even dimension}  \( n \) satisfying $1\le n\le 4m$. 

More specifically, as in the odd-dimensional cases, there exist two similar but
different indexes:  $m_n$ ({\it the total types of zero resonances}) and  {\it the  critical type} $k_c$ for all even dimensions  $1 \le n \le 4m$, such that if zero is a regular point of $H$ or a $\mathbf{k}$-th kind resonance with $1 \le \mathbf{k} \le k_c$,  the wave operators
 $W_\pm(H; (-\Delta)^m) \in \mathbb{B}(L^p)$ for all $1 < p < \infty$. If zero is a $\mathbf{k}$-th kind resonance with $k_c < \mathbf{k} \le m_n-1$, we show that $W_\pm(H; (-\Delta)^m) \in \mathbb{B}(L^p)$ for $1 < p < p_{\mathbf{k}}$, where
$$p_{\mathbf{k}} = \frac{n}{\max\{n-2m-1, 0\} + \mathbf{k} - k_c};$$
Additionally, if zero is a $m_n$-th kind resonance of $H$ or an eigenvalue of $H$ (i.e., $\mathbf{k}=m_n+1$), then
 $W_\pm(H; (-\Delta)^m) \in \mathbb{B}(L^p)$
for all
$1 < p \le 2$.

In particular, it was shown that except for some endpoint cases, these admissible ranges of $p$ are sharp for all zero resonance types of $H$, including the case of zero eigenvalues.  The key ideas of the $L^p$-sharpness is to establish the optimality of local-decay ({\it  Kato-Jensen type}) estimates for $e^{itH}$ in weighted $L^2$ spaces.  As a consequence, we establish the sharp $\mathbb{B}(L^{p}, L^{p^{\prime}})$-decay estimates of Schr\"odinger  group $e^{itH}$ with all kinds of zero-energy singularities. 

It is worth noting that the core aspect of this work is the complete classification of zero resonances and the full derivation of asymptotic expansions for the resolvent of $H=(-\Delta)^m + V$ near the zero threshold, valid for all {\it even dimensions} $n \leq 4m$ and  $m\ge 1$. 

%It is worthy to noting that  $(-\Delta)^m$ is \emph{degenerate} at zero energy when $m \geq 2$, and asymptotic expansions of free resolvent $R_0(z)$ contain several \emph{logarithmic terms} in even dimensions. These expansions for $R(z)$ at zero-energy are  more intricate than those in the odd-dimensional case.

\end{abstract}
\maketitle
\tableofcontents

\section{Introduction}
\setcounter{equation}{0}
\subsection{Backgrounds}
In this paper we study the $L^p$-boundedness of the wave operators $W_\pm(H;(-\Delta)^m)$ associated with the higher-order Sch\"{o}dinger operator $H=(-\Delta)^m+V(x)$ where $m\geq 1$, which are defined through the strong limits in $L^2(\mathbb{R}^n)$:
\begin{equation}\label{eq-wave operator-defi}
W_\pm(H;(-\Delta)^m) := \text{s-}\lim_{t\to\pm\infty} e^{itH}e^{-it(-\Delta)^m},
\end{equation}
where $e^{itH}$ and $e^{-it(-\Delta)^m}$ denote the evolution groups generated by $H$ and $(-\Delta)^m$ respectively.

As $V(x)$ is a real-valued bounded potential satisfying  $|V(x)|\lesssim\langle x\rangle^{-\beta}$ for $\beta>1,$ it was known that $W_\pm (H;(-\Delta)^m)$ exist on $L^2(\mathbb{R}^n)$ as partial isometric operators and are asymptotically complete ( c.f.  e.g.  Kuroda \cite{Kuroda}). In particular,  the wave operators $W_\pm$ satisfy the following identities
\[W_\pm W_\pm^*=P_{ac}(H),\ \ W_\pm^*W_\pm=I,\]
and the intertwining formula  (see e.g. Reed and Simon \cite{Reed-Simon-III}):
\begin{equation}\label{eq-intertwing identity}
f(H)P_{ac}(H)=W_\pm f((-\Delta)^m)W^*_\pm,
\end{equation}
where $f$  is  any Borel measurable function of  $\mathbb{R}$ and $P_{ac}(H)$ denotes the projection onto the absolutely continuous spectral subspace of $H.$ 

By the identity \eqref{eq-intertwing identity} and the $L^p-L^q$ estimates for the free operator $f((-\Delta)^m),$ the $L^p$-boundedness of $W_\pm(H;(-\Delta)^m)$ can immediately be used to establish the $L^p-L^q$ estimates for the perturbed operator $f(H)P_{ac}(H)$ by
\begin{equation}\label{eq-disperesive estimate}
  \|f(H)P_{ac}(H)\|_{L^p\rightarrow L^q}\le \|W_\pm\|_{L^q\rightarrow L^q}\|f((-\Delta)^m)\|_{L^p\rightarrow L^q}\|W_\pm^*\|_{L^p\rightarrow L^p},
\end{equation}
 In fact, under a suitable condition on $f,$ it is comparably easy to establish the $L^p-L^q$ estimates for the free operator $f((-\Delta)^m)$ by Fourier multiplier methods. Thus, in order to obtain the estimates \eqref{eq-disperesive estimate}, it suffices to prove all possible $L^p$-boundedness of $W_\pm(H;(-\Delta)^m)$ and $W_\pm^*(H;(-\Delta)^m)$.

Due to this characteristic, the \( L^p \)-boundedness of wave operators for the Schr\"odinger operator \( -\Delta + V \) (i.e., \( m = 1 \)) has been extensively studied and serves as a fundamental tool in the analysis of various nonlinear dispersive equations. As a result, it is both natural and intriguing to extend the investigation of \( L^p \)-boundedness to wave operators associated with more general higher-order elliptic operators of the form \( (-\Delta)^m + V(x) \) for $m\ge1$. While a significant body of work has been devoted to this topic in recent studies (see, e.g., \cite{EG22}, \cite{GG21}, \cite{GY23}, \cite{MWY22}, \cite{MWY23-1}, \cite{MWY23-2}), these results are still less developed compared to the case of the Schr\"odinger operator \( -\Delta + V(x) \). For further background and a survey of existing literature, we refer to Subsection \ref{subsec:1.3} below.

Note that Erdo\u{g}an-Green in \cite{EG22,EG23} have established the $L^p$-boundedness of  wave operators $W_\pm(H;(-\Delta)^m)$ for all $1 \leq p \leq \infty$ when $n > 2m$ in the regular case at zero energy. More recently, Erdo\u{g}an-Green-LaMarster \cite{EGL24,EGL25} further investigated the corresponding $L^p$-bounds of $W_\pm(H;(-\Delta)^m)$ with zero eigenvalue of $H$ when $n > 4m$. Our goal is to study the remaining cases with all possible types of zero energy resonances of $H$ for $1 \leq n \leq 4m$.

In the previous work \cite{CSWY24} for {\bf any odd dimension} $ 1\le n\le 4m$, we have demonstrated that there exist \emph{  $\bar{m}_n $ distinct types of zero resonances for the operator $H$, along with a critical resonance type $\bar{k}_c$}, defined by 
$$\bar{m}_n:=\min\Big\{m, \ 2m-\frac{n-1}{2} \Big\},\ \ \bar{k}_c=\max\Big\{m-\frac {n-1}{2},\  0\Big\}.$$
such that the following conclusions hold:
\begin{itemize}
    \item If zero is a regular point of $H$ (i.e., $\mathbf{k}=0$) or a $\mathbf{k}$-th kind resonance of $H$ with $1 \leq \mathbf{k} \leq \bar{k}_c$, then $W_\pm\big(H;(-\Delta)^m\big) \in \mathbb{B}(L^p(\mathbb{R}^n))$ for all $1 < p < \infty$;
    
    \item If zero is a $\mathbf{k}$-th kind resonance of $H$ with $\bar{k}_c < \mathbf{k} \leq \bar{m}_n$, then $W_\pm\big(H;(-\Delta)^m\big) \in \mathbb{B}(L^p(\mathbb{R}^n))$ for all
    \begin{equation}\label{finite-range}
1 < p < \frac{n}{n - 2m + \mathbf{k} + \bar{k}_c - 1}. \end{equation}
    \item If  an eigenvalue of $H$ (i.e., $\mathbf{k}=\bar{m}_n+1$), then $W_\pm\big(H;(-\Delta)^m\big) \in \mathbb{B}(L^p(\mathbb{R}^n))$ for all \( 1 < p < \frac{2n}{n - 1} \).
\end{itemize}
 In particular, we point out that these $L^p$-bounds of $ W_\pm(H;(-\Delta)^m)$ are sharp except for endpoint cases, and are also consistent with the existing results for $m=1,2$.

In this paper, we extend the analysis to even dimensions, proving that $W_\pm(H; (-\Delta)^m)$ are bounded on $L^p(\mathbb{R}^n)$, where the permissible range of $p$ is determined by the nature of the zero-energy singularities for \textbf{all even dimensions} $n$ satisfying $1 \leq n \leq 4m$ (see Theorems~\ref{thm-main result-0} and~\ref{thm-main result-1} below).

More specifically, as in the odd-dimensional cases, there exist the following two similar but different indexes $m_n$ ({\it the total number of zero resonance types}) and  $k_c$ ({\it the critical type}) for  {\bf any even dimension}  \( n \) satisfying $1\le n\le 4m$:
$$m_n=\min\Big\{m+1, \ 2m-\frac{n}{2}+1\Big\}, \ \ \ k_c=\max\Big\{m-\frac {n}{2}+1,\  0\Big\},$$
such that the following conclusions hold:
\begin{itemize}
    \item If zero is a regular point of $H$ (i.e., $\mathbf{k}=0$) or a $\mathbf{k}$-th kind resonance of $H$ with $1 \leq \mathbf{k} \leq k_c$, then $W_\pm\big(H;(-\Delta)^m\big) \in \mathbb{B}(L^p(\mathbb{R}^n))$ for all $1 < p < \infty$;
    
    \item If zero is a $\mathbf{k}$-th kind resonance of $H$ with $k_c < \mathbf{k} \leq m_n-1$, then $W_\pm\big(H;(-\Delta)^m\big) \in \mathbb{B}(L^p(\mathbb{R}^n))$ for all $1 < p < p_{\mathbf{k}}$, where
    \begin{equation}\label{p_k}
        p_{\mathbf{k}} = \frac{n}{\max\{n-2m-1, 0\} + \mathbf{k} - k_c};
    \end{equation}
    
    \item If zero is a $m_n$-th kind resonance of $H$ or an eigenvalue of $H$ (i.e., $\mathbf{k}=m_n+1$), then $W_\pm\big(H;(-\Delta)^m\big) \in \mathbb{B}(L^p(\mathbb{R}^n))$ for all $1 < p \leq 2$.
\end{itemize}
Again, when $k_c<\k\le m_m+1$, we also proved that, except for some endpoints, {\it the permissible range of $p$ above is sharp} for all zero resonance types, including the case of zero eigenvalue (see Theorem \ref{thm-unbound} below). The key ideas of the unboundedness involve reducing the problem to establishing the optimality of time-decay estimates for $e^{itH}$ in weighted $L^2$ spaces (referred to as {\it  Kato-Jensen-type} estimates). 

We emphasize that the core aspects of these results above are the complete classification of zero resonances and the full derivation of asymptotic expansions for the resolvent \( R(z) \) of \( H = (-\Delta)^m + V \) near the zero threshold, valid for all even dimensions \( n \leq 4m \) and $m\ge 1$, see Theorems \ref{thm-M inverse-even} and \ref{thm-M inverse-even-1} below. 

Notably, the symbol polynomial \( P(\xi) = |\xi|^{2m} \) of \( (-\Delta)^m \) is {\it degenerate} at \( \xi = 0 \) when \( m \geq 2 \), and the asymptotic expansions of the free resolvent \( R_0(z) \)   contain {\it logarithmic terms} in even dimensions. The asymptotic expansions of \( R(z) \) are more complex than in the odd-dimensional case, as addressed in \cite{CHHZ-2024}. Further details are provided in Section \ref{sec-proof-for-sec2.2}.

\subsection{The main results}
%\subsubsection{Zero resonances  and the conditions on $V$ }
 In order to state our main results, 
we first introduce the concept of the zero resonances for $H=(-\Delta)^m + V$ on $\R^n$, which are divided into two cases: $1\le n\le 2m$ and $2m<n\le 4m$ when $n$ is an even integer and $1\le n\le 4m$.  

Recall that zero is an eigenvalue of $H$ if there exists a nonzero $\phi \in L^2(\mathbb{R}^n)$ such that
\begin{equation}\label{equ1.2}
    ((-\Delta)^m + V)\phi = 0
\end{equation}
in the distributional sense. However, for some potentials, there may exist non-trivial solutions to \eqref{equ1.2} in the weighted spaces
\[
L_s^2(\mathbb{R}^n) := \{f \in L_{\mathrm{loc}}^2(\mathbb{R}^n) : \langle x \rangle^s f \in L^2(\mathbb{R}^n)\}
\]
for some $s < 0$. This is typically referred to as zero resonance of $H$, which presents challenges in establishing dispersive estimates or $L^p$-boundedness of the wave operators associated with $H$.

For dimensions $1 \leq n \leq 2m$, $H$ always has \textit{$m+1$ distinct types of zero resonance}, while for $2m < n \leq 4m$, there are \textit{$2m - \frac{n}{2} + 1$ different types of zero resonance}, decreasing as $n$ increases. In particular, when $n > 4m$, $H$ has \textit{no zero resonance} and may only have a zero eigenvalue (see, e.g., Feng et al. \cite{FSWY20} or Erdogan-Green-LaMarster \cite{EGL24}).

In what follows, we provide the precise definition of zero resonance types for $H = (-\Delta)^m + V$ when $n$ is an even integer satisfying $1 \leq n \leq 4m$.

\begin{definition}\label{defi-resonance-1}
Let $n$ be an even integer satisfying $1 \leq n \leq 2m$, and consider the higher-order Schr\"{o}dinger operator $H = (-\Delta)^m + V$, where $V \in L^\infty(\mathbb{R}^n)$ is a real-valued potential.  Then
\vskip0.1cm
(\romannumeral1) 
We say that zero is { \bf the first kind resonance} of $H$ if equation 
$((-\Delta)^{m}+V)\phi=0$ has a non-trivial solution $\phi$ in $\bigcap\limits_{s<-m}L_{s}^2(\mathbb{R}^n)$, but only the zero solution in $ \bigcap\limits_{s<1-m}L_{s}^2(\mathbb{R}^n)$.  Furthermore, all such solutions $\phi$ in $\bigcap\limits_{s<-m}L_{s}^2(\mathbb{R}^n)$ must have the following form:
\begin{equation}\label{eq-phi-minus-polynomial-1}
	\phi(x)=\phi_e(x)+\sum_{|\alpha|= m-\frac n2}C_\alpha x^\alpha,	
\end{equation}
where $C_\alpha\in \C$  and $\phi_e(x) \in \bigcap\limits_{s<1-m}L_{s}^2(\mathbb{R}^n)$.
\vskip0.1cm
(\romannumeral2) 
Let $2 \leq \mathbf{k} \leq m+1$. We say that zero is a \textbf{$\mathbf{k}$-th kind resonance} of $H$ if the equation $H\phi = 0$ has a non-trivial solution $\phi \in \bigcap_{s < \mathbf{k}-1-m} L_s^2(\mathbb{R}^n)$, but only the zero solution in $\bigcap_{s < \mathbf{k}-m} L_s^2(\mathbb{R}^n)$. 

(\romannumeral3)  We say that zero is {\bf the eigenvalue} of $H$, if $H\phi=0$  has a non-trivial solution in $L^2(\R^n)$. 
\vskip0.1cm
(\romannumeral4)  We say that zero is {\bf the regular point} of $H$, if zero is neither a resonance nor an eigenvalue. 
\end{definition}

\begin{definition}\label{defi-resonance-2}
Let $n$ be an even integer satisfying $2m< n \leq 4m$, and consider the higher-order Schr\"{o}dinger operator $H = (-\Delta)^m + V$, where $V \in L^\infty(\mathbb{R}^n)$ is a real-valued potential. 
	\vskip0.1cm
(\romannumeral1) Let $1\le \mathbf{k}\le 2m-\frac{n}{2}+1$. Then we say that zero is {\bf the $\k$-th kind resonance} of $H$  if  equation $H\phi=0$
has a nontrivial solution $\phi$ in $\bigcap\limits_{s<\k-2m+\frac n2-1}L_{s}^2(\mathbb{R}^n)$, but has only the zero solution to $H\phi=0$ in $\bigcap\limits_{s<\k-2m+\frac n2}L_{s}^2(\mathbb{R}^n)$. 
\vskip0.1cm
(\romannumeral2)  We say that zero is {\bf the eigenvalue} of $H$, if $H\phi=0$ has a non-trivial solution in $L^2(\R^n)$. 
\vskip0.1cm
(\romannumeral3)  We say that zero is {\bf the regular point} of $H$, if zero is neither a resonance nor an eigenvalue. 
\end{definition}
\begin{remark}\mbox{}

\begin{itemize}

\item The classifications in Definition \ref{defi-resonance-1} and \ref{defi-resonance-2} of zero energy are essentially consistent with those in the cases  $m=1$ and $m=2$ considered in \cite{GY23, JN01, LSY23, Yajima-2022-4}, where especially, they use the $L^p$ and the weighted space $\langle x\rangle^{s} L^{\infty}$  to characterize the decay of the solution $\phi$ to equation $H\phi=0$. 
Moreover, we note that the additional property \eqref{eq-phi-minus-polynomial-1} also appears in Lemma 6.4 of \cite{JN01} for $(m,n)=(1,2)$ (i.e. the case  of 
 $-\Delta +V$ on $\R^2$).

\vskip0.2cm
\item As mentioned above, zero resonance classifications are closely  to the asymptotic expansions of the resolvent $R(z)$ of $H=(-\Delta)^m+V$ near the zero threshold.    
 Compared with the Schr\"odinger operator $-\Delta+V$  ( i.e.,  $m=1$ ) ( see, e.g. \cite{JK, J80,JN01} ), the zero resonance classifications and asymptotic expansions of  $(-\Delta)^m+V$   are complicated  since  $(-\Delta)^m$ is degenerate at zero energy when $m\ge 2$, especially, for low even dimension $1\le n\le 4m$. 
We refer  to  Theorems \ref{thm-M inverse-even} and \ref{thm-M inverse-even-1} for more details. 

\end{itemize}
\end{remark}

Let $m\ge 1$ and $n$ be an even integer,  and set 
\begin{equation}\label{equ0.1}
	m_n:=\min\Big\{m+1,\ 2m-\frac n2 +1 \Big\}=
	\begin{cases}
		m+1,&\  \mbox{if}\ \  1\le n\le 2m, \ {\rm even};\\
		\mbox{$2m-\frac{n}{2}+1$},&\ \mbox{if}\ \ 2m<n\le 4m, \ {\rm even}.
		%	0,\quad\qquad\,\, \,\,\mbox{if}\, n\ge 4m+1,
	\end{cases}
\end{equation}
The integer \( m_n \) represents {\bf the total number of zero resonance types} for the operator \( H = (-\Delta)^m + V \) on \(\mathbb{R}^n\), as defined above. Furthermore, we note that for \( n > 4m \) (regardless of whether \( n \) is even or odd), \( (-\Delta)^m + V \) admits no zero resonance but may still possess a zero eigenvalue (cf. Feng et al. \cite{FSWY20}).  

For convenience, in the following,  we also regard  {\it  the regular zero energy as  the $0$\text{-}th resonance of $H$ and  the zero eigenvalue as the $(m_n+1)$\text{-}th kind resonance of $H$}.

\vskip0.2cm

Throughout the paper, unless otherwise specified, we make the following assumption on $V$

\begin{assumption}\label{Assumption}
Let $m\ge 1$, $n$ be an even integer with $1\le n\le 4m$,  and  $0\le {\bf k}\le m_n+1$. Assume that $H=(-\Delta)^m+V$ where  $V$ is a real-valued potential. Furthermore,  assume that  the following  three conditions hold:
\vskip0.1cm
 (\romannumeral1) Zero is the $\mathbf{k}$\text{-}th kind resonance of $H$  and $H$ has no positive embedded eigenvalues.
\vskip0.1cm
 (\romannumeral2) $\lvert V(x)\rvert\lesssim \langle x \rangle^{-\beta}$ for some $\beta$ satisfying
\begin{equation}\label{decayonV}
\beta>
\begin{cases}
4m-n+4\mathbf{k}+6, &\quad\text{if\ \ $1\le n\le 2m$},\\
n+4\mathbf{k}+5, &\quad\text{if\ \ $2m< n\le 4m$}.
\end{cases}
\end{equation}
\vskip0.1cm
 (\romannumeral3)  $\|\langle\mathscr{F}\left(\langle \cdot\rangle^{1+\sigma}V(\cdot)\right) \|_{\frac{4m-1-\delta}{2m-\delta}}<C$ for some $\sigma>\frac{4m}{4m-1-\delta}+\delta$ and a sufficiently small $\delta>0$ when $n=4m$.
\end{assumption}

\begin{remark}\mbox{}
\begin{itemize}
 \item We note that the assumption about the absence of embedded positive eigenvalues is indispensable. The classical Kato theorem indicates that $-\Delta +V$ has no positive eigenvalues if $V(x)=o(|x|^{-1})$ as $|x|\to \infty$. Compared to the second-order case,  there exists a $V$ in $C_0^{\infty}(\mathbb{R}^n)$ such that $H=(-\Delta)^m+V$ has a positive eigenvalue when $m$ is even (cf. Feng et al. \cite[Section 7.1]{FSWY20}).  On the other hand, the authors in \cite{FSWY20} have proved that $H=(-\Delta)^m +V$ has no positive eigenvalues if $V\in C^1(\R^n)$ is bounded and satisfies the repulsive condition $x\cdot \nabla V(x) \leq 0$ for $x\in \R^n.$  
\vskip 0.1cm
\item Moreover, Assumption (iii) is a  condition of smoothness  used to estimate the high-energy part, as in Erdo\u{g}an-Green \cite{EG22} for $n=4m$. In particular, we note  that the smoothness assumption of $V$ for $H$ is necessary when $n\ge4m$, as shown by Erdo\u{g}an-Goldberg-Green\cite{EGG-JFA-2023}.
\end{itemize}
\end{remark}

We now introduce the main results of this paper.
Let $\mathbb{B}(X,Y)$ be the space of bounded operators from Banach space $X$ to Banach space $Y,$ and  $\mathbb{B}(X)=\mathbb{B}(X,X).$  Denote 
\begin{equation}\label{defkc}
	k_c:=\max\Big\{ m-\frac n2+1,\ 0\Big\}=\begin{cases}
		m-\frac{n}{2}+1,\ \ &\ \text{if}\ \ 1\le n\le2m, \ {\rm even},\\
		0,\ \ & \ \text{if}\ \ 2m< n\le 4m, \ {\rm even}.
	\end{cases}
\end{equation}
As shown in Theorems \ref{thm-main result-0} and  \ref{thm-main result-1}  below, $k_c$ is the critical type of $H$ in the sense that when ${\bf k}>k_c$, the range of the value $p$ of the $L^p$-bound of wave operators $W_\pm\big(H;(-\Delta)^m\big)$ shrinks into  a finite interval depending on $\bf{k}$.

\begin{theorem}\label{thm-main result-0}
Let $1\le n\le 2m$ be even, $m_n=m+1$, $0\le \mathbf{k}\le m_n+1$, and  $k_c=m-\frac n2 +1$. Suppose that $H=(-\Delta)^m+V$ satisfies  Assumption \ref{Assumption} (depending on the type $\bf{k}$ of zero resonance).    Then the following statements hold:

%\vspace{7pt}
%\indent\emph{(\romannumeral1)} If  zero is the regular point of $H$  ( i.e. $\mathbf{k}=0$ ), then  $W_\pm\big(H;(-\Delta)^m\big)\in \mathbb{B}(L^p(\mathbb{R}^n))$ for all $1<p<\infty$;

\vspace{4pt}
\indent\emph{(\romannumeral1)} If  zero is the regular point of $H$  ( i.e. $\mathbf{k}=0$ ) or zero is the $\bf{k}$-th kind resonance of $H$ with $1\le \mathbf{k}\le k_c$, then $ W_\pm\big(H;(-\Delta)^m\big)\in \mathbb{B}(L^p(\mathbb{R}^n))$ for all $1<p<\infty;$

\vspace{4pt}
\indent\emph{(\romannumeral2)} If  zero is the $\bf{k}$-th kind resonance of  $H$ with $k_c<\mathbf{k}\le m_n-1$, then  $W_\pm\big(H;(-\Delta)^m\big)\in \mathbb{B}(L^p(\mathbb{R}^n))$ for all $1<p<\frac{2n}{n-2m+2\mathbf{k}-2}$;

\vspace{4pt}
\indent\emph{(\romannumeral3)} If zero is the $m_n$-th kind resonance of $H$  or zero is an eigenvalue of $H$ (i.e. $\mathbf{k}=m_n+1$), then $ W_\pm\big(H;(-\Delta)^m\big)\in \mathbb{B}(L^p(\mathbb{R}^n))$ for all $1<p\le 2.$     
\end{theorem}

%\vskip0.2cm
\begin{theorem}\label{thm-main result-1}
	Let $2m< n\le 4m$ be even, $m_n=2m-\frac n2 +1$, $0\le \mathbf{k}\le m_n +1$, and $k_c=0$. Assume that $H=(-\Delta)^m+V$ satisfies  Assumption \ref{Assumption} (depending on the type $\bf{k}$ of zero resonance).    Then the following statements hold:
	
	\vspace{4pt}
	\indent\emph{(\romannumeral1)} If  zero is the regular point of $H$  ( i.e. $\mathbf{k}=0$ ), then  $W_\pm\big(H;(-\Delta)^m\big)\in \mathbb{B}(L^p(\mathbb{R}^n))$ for all $1<p<\infty$;
	
	\vspace{4pt}
	\indent\emph{(\romannumeral2)} If  zero is the $\bf{k}$-th kind resonance of $H$ and  $1\le \mathbf{k}\le m_n-1$, then  $W_\pm\big(H;(-\Delta)^m\big)\in \mathbb{B}(L^p(\mathbb{R}^n))$  for all $1<p<\frac{n}{n-2m+\mathbf{k}-1}$.
	
	\vspace{4pt}
	\indent\emph{(\romannumeral 3)} If zero is the $m_n$-th kind resonance of $H$ or zero is an eigenvalue of $H$ (i.e. $\mathbf{k}=m_n+1$), then $ W_\pm\big(H;(-\Delta)^m\big)\in \mathbb{B}(L^p(\mathbb{R}^n))$ for all $1<p\le 2$.     
\end{theorem}
%\vskip 0.2cm

\begin{remark} We give some comments on Theorems \ref{thm-main result-0} and  \ref{thm-main result-1} above as follows:

\begin{itemize}
\item From Theorem~\ref{thm-main result-0} (i), in even dimensions $1 \leq n \leq 2m$, there exist some mild zero-energy resonances (i.e., $1 \leq \mathbf{k} \leq k_c = m - \frac{n}{2} + 1$) for which the wave operators $W_\pm\big(H;(-\Delta)^m\big)$ remain bounded on $L^p$ for all $1 < p < \infty$. 
However, for even dimensions $2m < n \leq 4m$, Theorem~\ref{thm-unbound} below shows that the range of $p$ for which $W_\pm$ are $L^p$-bounded must be a finite interval when any zero resonance occurs (i.e., $1 \leq \mathbf{k} \leq m_n$).
\vskip0.1cm
\item In the regular case, when $n > 2m$ is even, the wave operators $W_\pm\big(H;(-\Delta)^m\big)$ are bounded on $L^p(\mathbb{R}^n)$ for all $1 \leq p \leq \infty$, including the endpoints $p=1$ and $p=\infty$ (cf. Erdogan and Green \cite{EG23}). For the case when $n \leq 2m$ is even, the endpoint boundedness of $W_\pm\big(H;(-\Delta)^m\big)$ at both $p=1$ and $p=\infty$ are generally expected to fails, even in the simplest case of $m=1$ (cf. Yajima \cite{Yajima-99-CMP}). For odd dimensions, specifically when $n=1,3$ and $m=2$, explicit counterexamples demonstrating this failure can be found in the works of Mizutani-Wan-Yao \cite{MWY22,MWY23-2}. 
\vskip0.1cm
 \item In the following Table \ref{table}, we list the known results of the $L^p$-boundedness of the wave operators $W_\pm\big(H;(-\Delta)^m\big)$ when $1\le n\le 4m$ is even. In particular, our results  are consistent with the existing results for $m=1,2$.

\end{itemize}
\end{remark}

In the following, we show that the finite intervals of $p$ in both theorems above are sharp except for some endpoints if zero is the $\k$-th kind resonance of $H$ with $k_c<\mathbf{k}\le m_n+1$ and $1\le n\le 4m$. 

\begin{theorem}\label{thm-unbound}
Assume that the same assumptions of Theorem \ref{thm-main result-0} and Theorem \ref{thm-main result-1}  hold.  Then 

\emph{(\romannumeral1)} If zero is the $\k$-th kind resonance of $H$ with $k_c < \mathbf{k}\le  m_n-1$, the wave operators $W_\pm(H; (-\Delta)^m) \notin \mathbb{B}(L^p(\mathbb{R}^n))$ for all $p_{\k}<p\le \infty$, where $p_\k$ is defined in \eqref{p_k}, equivalently,
\begin{equation}\label{critical p}
p_{\k}=\begin{cases}
		\frac{2n}{n-2m+2\mathbf{k}-2},\ \ &\ \ \ 1\le n\le2m, \ {\rm even},\\
		\frac{n}{n-2m+\mathbf{k}-1},\ \ & \ \ \ 2m< n\le 4m, \ {\rm even}.
	\end{cases}
\end{equation}

\emph{(\romannumeral2)}  If zero is the $m_n$-th kind resonance of $H$ or if zero is an eigenvalue of $H$ and there exists a non-zero solution to $H\phi=0$ in $\displaystyle\bigcap_{s<0}L^2_{s}(\mathbb{R}^n) \setminus L^2(\mathbb{R}^n)$,  then  $W_\pm(H; (-\Delta)^m) \notin \mathbb{B}(L^p(\mathbb{R}^n))$ for all $2<p\le \infty$.  
\end{theorem}
%\vskip 0.2cm
%We make 
\begin{remark} Some remarks on Theorem \ref{thm-unbound} are given as follows:
\begin{itemize}
\item The non-zero solutions $\phi$ of $H\phi=0$ in $\bigcap_{s<0} L^2_{s}(\mathbb{R}^n) \setminus L^2(\mathbb{R}^n)$ are referred to as {\it p-wave resonances} in the case  $(m,n) = (1,2)$ by  Yajima \cite{Yajima-2022}, and as {\it d-wave resonances} in the case $(m,n) = (2,4)$ by Galtbayar-Yajima\cite{GY23}. In fact, if zero is an eigenvalue of $H$ but the equation $H\phi = 0$ admits no nonzero solution in $\bigcap_{s<0} L^2_s(\mathbb{R}^n) \setminus L^2(\mathbb{R}^n)$, then the wave operators \(W_\pm\big(H;(-\Delta)^m\big)\) may be bounded on \(L^p\) with an interval of \(p\) values strictly containing \((1, 2]\). Indeed,  when zero is an eigenvalue but $H\phi=0$ has no solution in $\bigcap_{s<-2m+\frac{n}{2}} L^2_s(\mathbb{R}^n) \setminus L^2(\mathbb{R}^n)$ for even dimensions $2m<n\le 4m$,  Erdo\u{g}an-Green-Lamaster \cite{EGL25} recently have shown that the wave operators $W_\pm(H; (-\Delta)^m)$ are bounded on $L^p$ for $1 < p < \frac{2n}{n-2}$.
 \vskip0.1cm
 \item

We emphasize that the results in Theorem \ref{thm-unbound} (ii) apply to the classical Schrödinger operator (i.e., \( m = 1 \)). In particular, for the critical case \( n = 2 \) and \( m = 1 \), the presence of a {\it p-wave resonance} implies that the wave operators \( W_\pm(H; -\Delta) \) are unbounded for all \( p > 2 \). This was recently established by Yajima \cite{Yajima-2022} using methods different from ours.

 \vskip0.1cm
 \item Theorem \ref{thm-unbound} (ii)  have resolved  the problem  concerning unboundedness of wave operators posed by Galtbayar-Yajima \cite {GY23}  for the case of $(m, n)=(2,4)$. Precisely, when there exists a {\it p-wave resonance} ( i.e. a nonzero solution to $H\phi=0$ in $\bigcap_{s<-1} L^2_{s}(\mathbb{R}^n) \setminus L^2(\mathbb{R}^n)$),  the wave operators \( W_\pm(H; \Delta^2)\) are unbounded for $p>4$ ( except for the endpoint $p=4$ ), and in the presence of {\it d-wave resonance}  (i.e. a nonzero solution to $H\phi=0$ in $\bigcap_{s<0} L^2_{s}(\mathbb{R}^n) \setminus L^2(\mathbb{R}^n)$), the wave operators \( W_\pm(H; \Delta^2)\) become unbounded for $p>2$, 
\vskip0.1cm
 \item 
  The proof of Theorem \ref{thm-unbound} was given by contradiction in Section \ref{sec-unb}. The key idea involves reducing the problem to establishing the optimality of time-decay estimates for $e^{itH}$ in weighted $L^2$ spaces (referred to as {\it  Kato-Jensen-type} estimates, cf. \cite{FSWY20}). These optimal decay estimates, in turn, depend on a detailed asymptotic expansion of the resolvent $R^\pm(\lambda^{2m})$ of $H = (-\Delta)^m + V$ near the zero threshold. Following the same approach, one could prove that $W_\pm(H; (-\Delta)^m) \notin \mathbb{B}(L^p(\mathbb{R}^n))$ for all $2 < p \leq \infty$ when the equation $H\phi = 0$ admits a non-zero solution in $L^2 \setminus L^2_1$. We omit it for brevity.
\end{itemize}
\end{remark}

\subsubsection{Some applications}
 In this subsection, as a consequence of the $L^p$ bound of the wave operator,  we obtain the $L^p$-decay estimate of Schr\"odinger group $e^{itH}$. 
 
 It is well-known that the $L^p$-decay estimate  for the Schr\"odinger operator $-\Delta + V$  has been an active research topic over the past three decades and has found broad applications in nonlinear Schr\"odinger equations. 
For $n \geq 3$, Journé, Soffer, and Sogge \cite{JSS} first established the $L^1-L^\infty$ estimate of Schr\"odinger operator in the regular case using methods distinct from wave operator techniques. Subsequent contributions for lower dimensions were made by Weder \cite{Weder1}, Rodnianski and Schlag \cite{RodSchl}, Goldberg and Schlag \cite{Goldberg-Schlag}, and Schlag \cite{Schlag-CMP}. Further developments and references can be found in Schlag’s survey work \cite{Schlag}.  

In particular, there exist many recent advances in dispersive estimates for higher-order Schrödinger operators $H = (-\Delta)^m + V$ with $m \geq 2$. For the case $m = 2$, Kato-Jensen estimates or $ L^1-L^\infty$ estimates for $e^{itH}$ have been established in works such as \cite{EGT21, FSY, Green-Toprak-4D, SWY22, LSY23}. For general $m \geq 2$, further results can be found in \cite{CHHZ-2024, EGG_ArXiv, FSWY20}.

Recall that the following $L^p$ estimates of the free higher-order Schr\"{o}dinger group $e^{it(-\Delta)^m}$ always hold  (see e.g. \cite{ZYF}) 
\begin{equation}\label{eq: for free}
\left\|e^{i t(-\Delta)^m}\right\|_{L^p-L^{p'}}\lesssim |t|^{-\frac{n}{m}(\frac{1}{p}-\frac {1}{2})}, \quad t\neq0,
\end{equation}
for all  $1\le p\le2$ and $\frac{1}{p'}+\frac{1}{p}=1$.  Hence by the $L^p$-boundedness of  $W_\pm(H;-\Delta)^m$  in Theorems \ref{thm-main result-0} and   \ref{thm-main result-1}, we 
find immediately that the following corollary holds from the  intertwining formula (also see \eqref{eq-intertwing identity}):
\begin{equation*}
e^{itH}P_{ac}(H)=W_\pm\ e^{it(-\Delta)^m}\ W_\pm^*.
\end{equation*}

\begin{corollary}
Let $1\le n\le 4m$ be even. Under Assumption \ref{Assumption} with resonance type $\mathbf{k}$,  one has 
\begin{equation}\label{decayest}
\Big\|e^{i tH}P_{ac}(H)\Big\|_{L^p-L^{p'}}\lesssim |t|^{-\frac{n}{m}(\frac{1}{p}-\frac 12)}, \quad t\neq0,
\end{equation}
where  $\frac{1}{p'}+\frac{1}{p}=1$ and the exponent $p'$  satisfies that
\begin{equation}\label{range of p}
\begin{cases}
2\le p'< \infty, &\quad \text{if}\ \ ~ 0\le \mathbf{k} \le k_c,\\
2\le p'< p_{\k}, &\quad \text{if}\ \ ~  k_c+1\le \mathbf{k} \le m_n-1,\\
p'=2, &\quad \text{if}\ \ ~\mathbf{k}= m_n,\ m_n+1.
\end{cases}
\end{equation}
\end{corollary}
%\vskip0.2cm
\begin{remark} 

If $0 \leq \k \leq \kappa_c$, then the decay estimate \eqref{decayest} remains valid for case $p' = \infty$. This has been established in \cite{EG23} for even dimensions $n \geq 2m$, in \cite{LSY23} for the case $(n, m) = (2, 2)$,  or in \cite{Green-Toprak-4D} for $(n, m) = (4, 2)$.
However, the estimate \eqref{decayest} does not hold for any $p' > p_{\k}$ when $k_c < \k \le m_n-1$, and also fails to hold for any $p'>2$ when $\k=m_n$, or $\k=m_n+1$ with a {\it p-wave resoance} of $H$. This is due to the optimality of the \emph{Kato-Jensen estimates} for $e^{itH}P_{ac}(H)$, as demonstrated in Proposition~\ref{Pro: lower bound of the time decay}. 
\end{remark}

\begin{table}[h]
\renewcommand\arraystretch{2}
\centering
\begin{tabular}{|p{2cm}<{\centering}|p{4.4cm}<{\centering}|p{9cm}<{\centering}|}
 \hline 
  {\bf Operators}      &  {$\bf n$,\  $\bf m_n$ {\bf and} $\bf{k}$}     & {\bf $L^p$-bounds of  $W_\pm$ with type ${\bf k }$} \\
  \hline
  \multirow{2}{*}{\makecell{$-\Delta+V$\\ $(m=1)$}}      &  \makecell{$n=2,\ m_n=2$  \vspace{0.05cm}\\ $0\le \mathbf{k} \le 3$  }   \vspace{0.05cm}&  \makecell{$\mathbf{k}=0:$\ \ $1<p<\infty$, Yajima \cite{Yajima-99-CMP};  \vspace{0.05cm}\\ $\mathbf{k}=1:$\ \ $1<p<\infty$, \\Erdog\u{a}n-Goldberg-Green \cite{Erdogan-Goldberg-Green-JFA-2018}; \vspace{0.05cm}
\\ $\mathbf{k}=2,3:$\ \ $1<p\le 2$, Yajima \cite{Yajima-2022}. \vspace{0.05cm}}\\
  \cline{2-3}
  & \makecell{$n=4,\ m_n=2$  \vspace{0.05cm}\\ $0\le {\bf k}\le 3$}   \vspace{0.05cm}& \makecell{$\mathbf{k}=0:$\ \ $1\le p\le\infty$, Yajima \cite{Yajima-1995-even};  \vspace{0.05cm}\\ $\mathbf{k}=1,2,3:$\ \ $1<p\le 2$, Yajima \cite{Yajima-2022-4};\vspace{0.05cm}
\\ only zero eigen. but not reso.:\\ $1\le p<4$,\vspace{0.05cm}
Goldberg-Green \cite{Goldberg-Green-2017}. \vspace{0.05cm}}\\
  \hline
  \multirow{3}{*}{\makecell{$\Delta^2+V$\\ $(m=2)$}}      & \makecell{$n=2,\ m_n=3$  \vspace{0.05cm}\\ $0\le \k\le 4$ }   \vspace{0.05cm}&  \makecell{$k=0,1,2:\ 1<p<\infty,$\ \  Thm \ref{thm-main result-0} (i),\ (ii);  \vspace{0.05cm}\\ $k=3,4:\ 1<p\le2,$\ \ Thm \ref{thm-main result-0} (iv).  \vspace{0.05cm}}\\ 
  \cline{2-3}
  &\makecell{ $n=4,\ m_n=3$  \vspace{0.05cm}\\ $0\le \k\le 4$  
  } & \makecell{$\mathbf{k}=0,1$:\ \ $1<p<\infty$;\ \ $\mathbf{k}=2$:\ \ $1<p<4$; \vspace{0.05cm}\\
  and $\mathbf{k}=3,4$:\ \ $1<p<2$;\\ Galtbayar-Yajima \cite{GY23}. \vspace{0.05cm}}\\
  \cline{2-3}
  & \makecell{$n=6,8$  \vspace{0.05cm}\\  $0\le \k \le m_n+1$  \vspace{0.05cm}} & \makecell{$\mathbf{k}=0$:\ \ $1\le p\le\infty$, Erdo\u{g}an-Green \cite{EG22};  \vspace{0.05cm}\\   $1\le \k \le m_n+1$,  Thm \ref{thm-main result-1}.
  }  \vspace{0.05cm}\\
  \hline
   \multirow{2}{*}{\makecell{$(-\Delta)^m+V$\\ $(m\ge3)$}}     & \makecell{$1\le n\le 2m$ even  \vspace{0.05cm}\\ $k_c=m-\frac{n}{2}+1$\\ $m_n=m+1$   \vspace{0.05cm}
   \\ $0\le \k \le m_n+1$} & \makecell{$0\le \k\le k_c:\ 1<p<\infty,$\ \  Thm \ref{thm-main result-0} (i),\ (ii);  \vspace{0.05cm}\\ $k_c< \k< m_n:\ 1<p<\frac{2n}{2\mathbf{k}-2m-2+n},$\\ Theorem \ref{thm-main result-0} (iii);  \vspace{0.05cm}\\ $k=m_n,m_n+1:\ 1<p\le2$,\ \  Thm \ref{thm-main result-0} (iv).  }\vspace{0.05cm}\\
  \cline{2-3}
  & \makecell{$2m< n\le 4m$ even \vspace{0.05cm}\\ $k_c=0$\\$m_n=2m-\frac{n}{2}+1$   \vspace{0.05cm}\\ $0\le \k \le m_n+1$}&  \makecell{$\mathbf{k}=0$:\ \ $1\le p\le\infty$, Erdo\u{g}an-Green \cite{EG22};  \vspace{0.05cm}\\
$1\le \k\le m_n-1:\ 1<p<\frac{n}{\mathbf{k}-2m-1+n},$\\ Theorem \ref{thm-main result-1} (ii);  \vspace{0.05cm}\\ $\k=m_n,m_n+1:\ 1<p\le2$,\ \  Thm \ref{thm-main result-1}  (iii).}   \vspace{0.05cm}\\
  \cline{2-3}
  \hline
\end{tabular}
\vskip 0.2cm    
\caption{
The $L^p$-boundedness for $W_\pm\big(H;(-\Delta)^m\big)$ when $1\le n\le 4m$ is even}\label{table}
\end{table}

\subsection{Further related backgrounds.}\label{subsec:1.3}
 In this subsection, we recall some known results on the $L^p$ boundedness of the wave operators 
 associated with classical Schr\"odinger operator and  higher-order elliptic operators.

For the classical Schr\"odinger operator \( H = -\Delta + V(x) \), much of the research on the \( L^p \)-boundedness of wave operators \( W_\pm(H; -\Delta) \) has been motivated by their applications in linear and nonlinear dispersive equations. If zero energy is regular, Yajima \cite{Yajima-JMSJ-95, Yajima-1995-even} proved that for certain smoothness and decay conditions on \( V \), the wave operators are bounded on \( L^p(\mathbb{R}^n) \) for all \( 1 \leq p \leq \infty \) when \( n \geq 3 \). Beceanu \cite{Beceanu} and Beceanu-Schlag \cite{Beceanu-Schlag-20} further refined these results for \( n = 3 \). In \( n = 1, 2 \), the wave operators are bounded on \( L^p(\mathbb{R}^n) \) for \( 1 < p < \infty \) when zero energy is regular (cf. \cite{Weder-99-CMP, Gal-Yaj-00, AF06, Yajima-99-CMP, Jensen-Yajima-02-CMP}).

When zero resonance or eigenvalue occurs, the range of \( p \) may shrink. For \( n > 4 \), the wave operators are bounded on \( L^p(\mathbb{R}^n) \) for \( \frac{n}{n-2} < p < \frac{n}{2} \) \cite{Finco-Yajima-2006-even, Yajima-JMSJ-2006}, and this was extended to \( 1 < p < n \) by Yajima \cite{Yajima-2016} and Goldberg-Green \cite{Goldberg-Green-2016}. For \( n = 4 \), the range is \( 1 \leq p < 4 \) if zero is an eigenvalue but not a resonance \cite{Goldberg-Green-2017}, and \( 1 < p < 4 \) if zero is a resonance \cite{Yajima-2022-4}. For \( n = 3 \), Yajima \cite{Yajima-2016-3d} found the range \( 1 \leq p < 3 \) for eigenvalues, and \( 1 < p < 3 \) for resonances. Finally, for \( n = 2 \), Erdoğan-Goldberg-Green \cite{Erdogan-Goldberg-Green-JFA-2018} showed boundedness for \( 1 < p < \infty \) if zero is a resonance or eigenvalue, while Yajima \cite{Yajima-2022} proved boundedness for \( 1 < p \leq 2 \) and unboundedness for \( 2 < p < \infty \) if zero is a \( p \)-wave resonance.

For \( n = 4 \), Goldberg-Green \cite{Goldberg-Green-2017} proved that the wave operators are bounded on \( L^p(\mathbb{R}^4) \) for \( 1 \leq p < 4 \) if zero is an eigenvalue but not a resonance, and for \( 1 < p < 4 \) if zero is a resonance. For \( n = 3 \), Yajima \cite{Yajima-2016-3d} showed that \( 1 \leq p < 3 \) if zero energy is an eigenvalue, and \( 1 < p < 3 \) if zero energy is a resonance. In \( n = 2 \), Erdoğan-Goldberg-Green \cite{Erdogan-Goldberg-Green-JFA-2018} showed that the wave operators are bounded on \( L^p(\mathbb{R}^2) \) for \( 1 < p < \infty \) if there is an \( s \)-wave resonance or an eigenvalue at zero. Yajima \cite{Yajima-2022} proved that for a \( p \)-wave resonance at zero, the wave operator is bounded on \( L^p(\mathbb{R}^2) \) for \( 1 < p \leq 2 \) and unbounded for \( 2 < p < \infty \).

For higher-order Schr\"odinger operators, \( H = (-\Delta)^m + V \) (with \( m \geq 2 \)), Goldberg-Green \cite{GG21} first proved that for \( n = 3 \) and \( m = 2 \), the wave operators \( W_\pm(H; \Delta^2) \) are bounded on \( L^p(\mathbb{R}^3) \) for \( 1 < p < \infty \) under the regularity assumption of zero energy. This was generalized to higher-order Schr\"odinger operators on \( \mathbb{R}^n \) for \( n > 2m \) in the regular case by Erdoğan-Green \cite{EG22, EG23}.

Recent progress includes the study of the fourth-order Schr\"odinger operator \( \Delta^2 + V(x) \). For \( n = 1 \), Mizutani-Wan-Yao \cite{MWY22} showed that \( W_\pm(H; \Delta^2) \) are bounded on \( L^p(\mathbb{R}) \) for \( 1 < p < \infty \), and unbounded at \( p = 1 \) and \( p = \infty \), regardless of whether zero energy is regular or a resonance. For \( n = 3 \), Mizutani-Wan-Yao \cite{MWY23-2} showed that the wave operators are unbounded at \( p = 1 \) and \( p = \infty \) and established weak \( (1,1) \) estimates in the regular case. In \cite{MWY23-1}, they proved that the range of \( p \) is \( 1 < p < \infty \) if zero is a first-kind resonance, and \( 1 < p < 3 \) if zero is a second- or third-kind resonance. For \( n = 4 \), Galtbayar-Yajima \cite{GY23} showed that the wave operators are bounded on \( L^p(\mathbb{R}^4) \) for \( 1 < p < \infty \) if zero energy is regular or corresponds to an \( s \)-wave resonance. If zero energy corresponds to \( s \)- or \( p \)-wave resonances or eigenfunctions, the range of \( p \) is \( 1 < p < 4 \), and for \( d \)-wave resonance, the range is \( 1 < p \leq 2 \).

In the present paper, and jointly with the work \cite{CSWY24}, we have shown that the wave operators \( W_\pm(H; (-\Delta)^m) \) are bounded on \( L^p(\mathbb{R}^n) \), with  the sharp range of \( p \) depending on all possible types of zero-energy resonances of \( H \) for \( 1 \leq n \leq 4m \). These results are also consistent with the existing ones for $m=1,2$.

\subsection{The outline of the proof}
In this subsection, we outline the proofs of Theorems~\ref{thm-main result-0} and~\ref{thm-main result-1}. Our proof begins with the stationary representation of the wave operators (see, e.g., Kuroda \cite{Kuroda})
\begin{equation}\label{eq-stationary formula}
W_{\pm}(H;(-\Delta)^m) = I - \frac{m}{\pi i}\int_0^{+\infty}\lambda^{2m-1}R^{\mp}(\lambda^{2m})V\big(R_0^+(\lambda^{2m})-R_0^-(\lambda^{2m})\big)d\lambda,
\end{equation}
where $R^{\pm}(\lambda^{2m})$ and $R_0^\pm(\lambda^{2m})$ denote the limiting resolvent operators of $H$ and $H_0$, respectively. We focus on proving the $L^p$-boundedness of $W_{-}$ since the relation $\overline{W_+f} = W_-\overline{f}$ allows us to deduce the result for $W_+$ from that of $W_-$.

As the identity operator $I$ is bounded on $L^p(\mathbb{R}^n)$ for all $1 \leq p \leq \infty$, we need only consider the integral term in \eqref{eq-stationary formula}, which we denote by $W$. Let $\chi \in C_0^\infty(\mathbb{R}^n)$ be a cutoff function satisfying $\chi = 1$ on $(-\lambda_0/2, \lambda_0/2)$ with $\text{supp}\,\chi \subset [-\lambda_0, \lambda_0]$, where $\lambda_0$ is the constant appearing in Theorems~\ref{thm-M inverse-even} and~\ref{thm-M inverse-even-1}. Using $\chi$, we decompose $W$ into low and high energy components as
\begin{equation}\label{eq-wave operator-low-0}
W^L = \frac{m}{\pi i}\int_0^\infty \lambda^{2m-1}\chi(\lambda)R^+(\lambda^{2m})V\big(R_0^+(\lambda^{2m})-R_0^-(\lambda^{2m})\big)d\lambda,
\end{equation}
\begin{equation}\label{eq-wave operator-high-0}
W^H = \frac{m}{\pi i}\int_0^\infty \lambda^{2m-1}(1-\chi(\lambda))R^+(\lambda^{2m})V\big(R_0^+(\lambda^{2m})-R_0^-(\lambda^{2m})\big)d\lambda.
\end{equation}

\vskip0.3cm
\subsubsection{The low energy part $W^{L}$} 
We begin by analyzing the low-energy component $W^{L}$. Let $v(x) = \sqrt{|V(x)|}$ and $U(x) = \text{sgn}\,V(x)$, and define
\begin{equation} \label{defM}
M^\pm(\lambda) = U + vR_0^\pm(\lambda^{2m})v.
\end{equation}
Since $H = (-\Delta)^m + V$ has no embedded positive eigenvalues, it is well-known that $M^\pm(\lambda)$ is invertible on $L^2(\mathbb{R}^n)$ for all $\lambda > 0$ (see, e.g., Agmon \cite{Agmon} and Kuroda \cite{Kuroda}). Moreover, we have the symmetric second resolvent identity
\[
R^\pm(\lambda^{2m})V = R_0^\pm(\lambda^{2m})v (M^\pm(\lambda))^{-1}v.
\]
Using this identity, we can rewrite $W^L$ as
\begin{equation}\label{eq-wave operator-low-1}
W^L = \frac{m}{\pi i}\int_0^\infty \lambda^{2m-1}\chi(\lambda) R_0^+(\lambda^{2m})v (M^+(\lambda))^{-1}v \big(R_0^+(\lambda^{2m}) - R_0^-(\lambda^{2m})\big) d\lambda.
\end{equation}
The key technical challenge lies in understanding the refined integral kernel structure of the following operator composition as $\lambda \to 0^+$:
\[
R_0^+(\lambda^{2m})v (M^+(\lambda))^{-1}v \big(R_0^+(\lambda^{2m}) - R_0^-(\lambda^{2m})\big).
\]
To establish this, a key step involves deriving the asymptotic behavior of $(M^\pm(\lambda))^{-1}$ as $\lambda \to 0^+$.

\vskip0.3cm
\textbf{Step 1. Proof of asymptotic expansions of $(M^\pm(\lambda))^{-1}$ around zero.}
We derive the asymptotic expansion (see Theorems~\ref{thm-M inverse-even} and~\ref{thm-M inverse-even-1}):
\begin{equation*}
    \big(M^\pm(\lambda)\big)^{-1} = \sum_{i,j \in J_\mathbf{k}} \lambda^{2m-n} \rho_i(\lambda)\rho_j(\lambda) Q_i \big(M_{i,j}^\pm + \Gamma_{i,j}^\pm(\lambda)\big) Q_j,
\end{equation*}
where $M_{i,j}^\pm \in \mathbb{B}(L^2)$ are $\lambda$-independent bounded operators, $\Gamma_{i,j}^\pm(\lambda)$ are $\lambda$-dependent operators, and the index set $J_{\mathbf{k}}$ depends on the zero-energy resonance type $\mathbf{k}$. We remark that the functions $\rho_i(\lambda)$ differ between the cases: $n \leq 2m$ ($n$ even) and $2m < n \leq 4m$ ($n$ even).

To demonstrate our approach, we examine the case where $n \le  2m$ ($n$ even) and zero is an eigenvalue. In this case, $M^\pm(\lambda)$ admit the following expansions  from \eqref{eq-expansion-Mpm-3}, expressed in different orders:
\begin{equation}\label{eq-expansion-Mpm-when-m+1}
\begin{aligned}
M^\pm(\lambda)&=\sum_{j=0}^{m -\frac n2-1}
a_j^\pm\lambda^{n-2m+2j}vG_{2j}v+b_0\log\lambda vG_{2m-n}v+T_0+a_{m-\frac n2}^\pm vG_{2m-n}v\\
&+\sum_{j=m -\frac n2+1}^{2m -\frac n2-1}
a_j^\pm\lambda^{n-2m+2j}vG_{2j}v+b_1\lambda^{2m}\log\lambda vG_{4m-n}v+
b_1\lambda^{2m}vG^0_{4m-n}v \\
&+a_{2m-\frac n2}^\pm\lambda^{2m}vG_{4m-n}v+\lambda^{n-2m}vr_{4m-n+1}^\pm(\lambda)v,
\end{aligned}
\end{equation}
where all terms are $L^2$-bounded  under our decay assumption on $V$.
Furthermore, the remainder terms $vr_{4m-n+1}^\pm(\lambda)v=O(\lambda^{2m+1})$  (see \eqref{eq-est-for-rlam}).

When zero is the eigenvalue of $H$ (i.e., $\k = m_n+1$), we  introduce a set  $J_{m_n+1}$ of indices associated with the expansion \eqref{eq-expansion-Mpm-when-m+1}\footnote{The indices in $J_{m_n+1}$ correspond to terms from the first $a_0^\pm\lambda^{n-2m}vG_0v$ to $b_1\lambda^{2m}vG_{4m-n}^0v$ in the expansion \eqref{eq-expansion-Mpm-when-m+1}. The fractional exponents were chosen for computational convenience to characterize the cancellation order of the projections $Q_j$ specified in \eqref{eq-cancellation-Q_j}.} as follows:
\begin{equation*}
	J_{m_n+1} = \left\{0, 1, \ldots, m-\frac{n}{2},\ m-\frac{n}{2}+\frac 12,\ m-\frac{n}{2}+\frac 34,\ m-\frac{n}{2}+1, \ldots, 2m-\frac{n}{2},\ 2m-\frac{n}{2}+\frac 12\right\}. 
\end{equation*} Correspondingly, we define a family of orthogonal projection operators $\{Q_j\}_{j \in J_{m_n+1}}$ in \eqref{eq-Q_j-even}, which satisfy $\sum\limits_{j\in J_{m_n+1}}Q_j=I$ by Proposition \ref{prop-resonance-equivalence}, and also introduce the following weight functions:
\begin{equation*}
 	\rho_j(\lambda)=\begin{cases}
 		\lambda^{-j}(\log\lambda)^{-\frac{1}{2}},\ \ &\  j= m-\frac{n}{2},\ 2m-\frac{n}{2},\\
 		\lambda^{-j+\frac12},\ \ &\  j= m-\frac{n}{2}+\frac{1}{2},\ 2m-\frac{n}{2}+\frac{1}{2},\\
 		\lambda^{-(m-\frac{n}{2})}(\log\lambda)^{\frac{1}{2}},\ \ &\  j= m-\frac{n}{2}+\frac{3}{4},\\
 		\lambda^{-j}, &\ \text{otherwise}.\\
 	\end{cases}
 \end{equation*}
For each $\lambda>0$,  we define the isomorphism:
\begin{equation*}\
B_{\lambda} = (\lambda^{-j}Q_j)_{j\in J_{m_n+1}} : \bigoplus_{j\in J_{m_n+1}} Q_jL^2 \to L^2, \quad (f_j)_{j\in J_{m_n+1}} \mapsto \sum_{j\in J_{m_n+1}} \rho_j(\lambda) Q_jf_j,
\end{equation*}
with its transport operator
$B_{\lambda}^T$. In matrix form, we write 
$$
B_{\lambda}= \left(Q_0, \rho_1(\lambda)Q_1, \ldots, \rho_{2m-\frac n2 +\frac 12}(\lambda)Q_{2m-\frac{n}{2}+\frac 12} \right),\quad B_{\lambda}^T= \begin{pmatrix} Q_0 \\ \rho_1(\lambda)Q_1\\ \vdots \\ \rho_{2m-\frac n2 +\frac 12}(\lambda)Q_{2m-\frac{n}{2}+\frac 12} \end{pmatrix}
\qquad  .
$$
 The isomorphism property of $B_{\lambda}$ yields the following equivalence: 
 \begin{equation*}
     \mbox{ $M^{\pm}(\lambda)$ is invertible on $L^2$ }\iff \mbox{$B_{\lambda}^TM^{\pm}(\lambda)B_{\lambda}$ is invertible on $\bigoplus_{j\in J_{m+1}} Q_jL^2$}
 \end{equation*}
with its inverse operator
%These two isomorphisms allow us to have the following identity
\begin{equation}\label{eq-inverse-A-to-M-0}
	\left(M^{\pm}(\lambda)\right)^{-1}=B_{\lambda} \left(B_{\lambda}^TM^{\pm}(\lambda)B_{\lambda}\right)^{-1}B_{\lambda}^T.
\end{equation}
This equivalence enables us to transform  the original problem  into studying the invertibility  of the block operator matrices $B_{\lambda}^TM^{\pm}(\lambda)B_{\lambda}$  and their asymptotic behaviors around $\lambda=0$.  This framework allows us to systematically analyze  $\left(M^{\pm}(\lambda)\right)^{-1}$ in a unified approach. 
%satisfying $\sum_{j \in J_{m+1}} Q_jL^2 = L^2$. Then, for each $\lambda>0$, denote two isomorphisms:
%\begin{equation*}
%	B_{\lambda} = (\lambda^{-j}Q_j)_{j\in J_{m+1}} : \bigoplus_{j\in J_{m+1}} Q_jL^2 \to L^2, \quad (f_j)_{j\in J_{m+1}} \mapsto \sum_{j\in J_{m+1}} \lambda^{-j}Q_jf_j,
%\end{equation*}
%and its adjoint operator $B_{\lambda}^*$. With these two isomorphisms, we can write
%\begin{equation*}
%	\left(M^{\pm}(\lambda)\right)^{-1}=B_{\lambda} \left(B_{\lambda}^*M^{\pm}(\lambda)B_{\lambda}\right)^{-1}B_{\lambda}^*
%\end{equation*}
%provided that $B_{\lambda}^*M^{\pm}(\lambda)B_{\lambda}$ is invertible on $\bigoplus_{j\in J_{m+1}} Q_jL^2$. Through this identity, we transform the original problem into studying the invertibility and asymptotic behaviors of the block operator matrices $B_{\lambda}^*M^{\pm}(\lambda)B_{\lambda}$ when $\lambda$ is around zero energy. This framework allows us to handle all cases in a unified approach. 
Based on the cancellation property of $Q_j$ in \eqref{eq:cancell of Q_j-1},  one get
$$\rho_i(\lambda)\rho_j(\lambda)Q_iM^\pm(\lambda)Q_j=\lambda^{n-2m}(d_{i,j}^\pm+ r_{i,j}^\pm(\lambda)).$$
Therefore, we can write the matrices $B_{\lambda}^TM^{\pm}(\lambda)B_{\lambda}$ as
\begin{equation}\label{eq:decomposition-intro}
	\mathbb{A^\pm} = \lambda^{2m-n}B_{\lambda}^TM^{\pm}(\lambda)B_{\lambda} = D^{\pm} + (r_{i,j}^\pm(\lambda))_{i,j\in J_{m_n+1}},
\end{equation}
where $D^{\pm}$ represent the leading order terms, and the remainder terms satisfy the better estimates $r_{i,j}^\pm(\lambda)=O((\log\lambda)^{-1/2})$ when $\lambda \in (0,1/2)$ (see \eqref{eq: expansion for r_ij^pm}).  The cancellation properties of the projections $Q_j$ induce a structure in $D^{\pm}$:
\begin{equation*}
D_{m_n+1}^\pm=
	\begin{pmatrix}
		D_{0,0}^{\pm} & 0 & D_{0,2}^{\pm}&0\\
		0 & D_{1,1}^{\pm} & 0&0\\
	    D_{2,0}^{\pm}  & 0 & D_{2,2}^{\pm}&0\\
	     0 & 0 & 0&D_{3,3}^{\pm}\\
	\end{pmatrix},
\end{equation*}
This structure reduces the problem to proving the invertibility of the following matrices in even dimensions $1 \leq n \leq 2m$:
\begin{equation}\label{eq:key-matrices}
\begin{aligned}
A_0^\pm &= \begin{pmatrix}
D_{0,0}^{\pm} & D_{0,2}^{\pm} \\
D_{2,0}^{\pm} & D_{2,2}^{\pm}
\end{pmatrix}, 
\quad D_{1,1}^{\pm}, &\quad D_{3,3}^{\pm}.
\end{aligned}
\end{equation}
It is not difficult to obtain the invertibility of $D_{1,1}^{\pm}$ and $ D_{3,3}^{\pm}$. The main challenge is transformed into proving the invertibility of $A_0^\pm$ on $\bigoplus_{j\in \tilde{J}_{0}\cup \tilde{J}_{2}} Q_jL^2$, where 
$$\tilde{J}_{0}=\left\{0, 1, \ldots, m-\frac{n}{2}-1\}, \quad  \tilde{J}_{2}=\{m-\frac{n}{2}+1, \ldots, 2m-\frac{n}{2}-1\right\} $$

Through the resolvent relation
$$R_0(-\lambda^{2m})(x-y)=R_0^{\pm}((e^{\frac{\pm i\pi}{2m}}\lambda)^{2m})(x-y),$$
the matrices $A_0^\pm$ are equivalent to another matrix, i.e.,
\begin{align*}
\mathbb{A}_0^\pm :=&
\text{diag}\left((-1)^iQ_i \right)_{i\in\tilde{J}_{0}\cup\tilde{J}_{2} }\,\text{diag}\left( U_l^\pm\right)_{l=0,2}\cdot A_0^{\pm}\cdot 
\text{diag}\left( U_l^\pm\right)_{l=0,2} \\
=&\left( b_0\log\lambda+a_{2m-n}^\pm \right)A_0^{s}+A_0^{r}.
\end{align*}
Here, the matrix $A_0^r$ has the block structure
\begin{equation*}
A_0^{r} = \begin{pmatrix}
\mathbf{D}_{00} & \mathbf{D}_{0,2}^{r} \\
\mathbf{D}_{2,0}^{r} & \mathbf{D}_{2,2}
\end{pmatrix}.
\end{equation*}
To prove the invertibility of the $\mathbb{A}_0^\pm$,  reapplying the preceding framework yields
\begin{equation*}
	\begin{split}
		 &\begin{pmatrix}
			(\log\lambda)^{-\frac{1}{2}} K_0 \\
			K_1
		\end{pmatrix} 
		\cdot  \mathbb{A}_0^\pm \cdot \begin{pmatrix}
			(\log\lambda)^{-\frac{1}{2}} K_0, \, K_1 
		\end{pmatrix}\\
		=&\begin{pmatrix}
			b_0K_0 A_0^{s}K_0 & 0 \\
			0 & K_1 A_0^{r}K_1
		\end{pmatrix}+\begin{pmatrix}
			(\log\lambda)^{-1} a_{2m-n}^{\pm}K_0 A_0^{s}K_0 & (\log\lambda)^{-\frac{1}{2}}K_0A_0^{r}K_1 \\
			(\log\lambda)^{-\frac{1}{2}}K_1A_0^{r}K_0 & 0
		\end{pmatrix},
	\end{split}
\end{equation*}
where $K_1$ is the projection onto the kernel of $A_0^{s}$ on $\bigoplus_{j\in \tilde{J}_0\cup\tilde{J}_2}Q_jL^2$, and $K_0=I-K_1$. By the perturbation arguments, the problem reduces to establishing the invertibility of the block operator matrix $K_1 A_0^{r} K_1$ on the subspace $K_0\left(\bigoplus_{j \in \tilde{J}_0 \cup \tilde{J}_2} Q_j L^2\right)$. This follows from the invertibility of $A_0^r$ on $\bigoplus_{j \in \tilde{J}_0 \cup \tilde{J}_2} Q_j L^2$.

The proof of the invertibility of $A_0^r$ on the direct sum space $\bigoplus_{j\in \tilde{J}_{0}\cup \tilde{J}_{2}} Q_jL^2$ is mainly based on  \textbf{two key observations (see Proposition \ref{pro-proerty-of-D00-D22})}:

\begin{itemize}

\item[($O_1$)] The block $\mathbf{D}_{00}$ is strictly positive definite on $\bigoplus_{j\in \tilde{J}_0}Q_jL^2$, corresponding to the Gram matrix of the system 
\[
\left\{\frac{\xi^\alpha}{(2\pi)^{n/2}\alpha!(1+|\xi|^{2m})^{1/2}}\right\}_{|\alpha|\in  \tilde{J}_0}.
\]
	
    \item[($O_2$)] The block $-\textbf{D}_{2,2}$ is strictly positive on $\bigoplus_{j\in \tilde{J}_2}Q_{j}^1 L^2$,  corresponding to  the Gram matrix of the system 
	\[
	\left\{\frac{\xi^\alpha}{(2\pi)^{n/2}\alpha!(1+|\xi|^{2m})^{1/2}|\xi|^m}\right\}_{|\alpha|\in \tilde{J}_2},
	\]
\end{itemize}
where 
$$
Q_j^1=Q_j-Q_j^0,\quad Q_j^0L^2=Q_jL^2\cap \{x^\alpha v(x):\, |\alpha|\le j\}^{\perp}.
$$
Through the abstract Feshbach formula,  it immediately yelds the invertibility of both $A_0^r$ is consistent with the invertibility of 
$$\textbf{d}=\textbf{D}_{2,2}- \textbf{D}_{2,0}^{r}\big(\textbf{D}_{0,0}\big)^{-1}\textbf{D}_{0,2}^{r},$$
which will be proved by observations  $(O_1)$ and $(O_2)$. Therefore, we conclude the invertbility of $D_{m_n+1}^\pm$.
Applying a standard Neumann series expansion and \eqref{eq-inverse-A-to-M-0} then yields the desired expansion of $\left(M^{\pm}(\lambda)\right)^{-1}$.

\vskip0.3cm
\textbf{Step 2. Proof of the $L^p$-boundedness of $W^L$.}
Based on the asymptotic expansion of $\big(M^\pm(\lambda)\big)^{-1}$ above, for even dimensions $n \leq 4m-2$, the integral kernel $W^L(x,y)$ can be expressed as a finite sum of terms of the form
\begin{equation*}
\int_0^\infty \lambda^{4m-n-1} \rho_i(\lambda)\rho_j(\lambda)\chi(\lambda)
\Big\langle \big(M_{i,j}^+ + \Gamma_{i,j}^+(\lambda)\big)Q_jv\big(R_0^{+}-R_0^{-}\big)(\lambda^{2m})(\cdot,y),\, Q_ivR_0^{-}(\lambda^{2m})(\cdot,x) \Big\rangle d\lambda,
\end{equation*}
where $i, j \in J_{\mathbf{k}}$. Using the properties of each $\Gamma_{i,j}^\pm(\lambda)$, we can further simplify the analysis by considering the modified integral
\begin{equation*}
\int_0^\infty \lambda^{4m-n-1} \widetilde{\rho}_i(\lambda)\widetilde{\rho}_j(\lambda)\chi(\lambda)
\Big\langle \big(M_{i,j}^+ + \Gamma_{i,j}^{+,r}(\lambda)\big)Q_jv\big(R_0^{+}-R_0^{-}\big)(\lambda^{2m})(\cdot,y),\, Q_ivR_0^{-}(\lambda^{2m})(\cdot,x) \Big\rangle d\lambda,
\end{equation*}
which we denote by $W_{i,j}^L(x,y)$. 
While for dimension $n = 4m$, we reduce the problem to studying the $L^p$-boundedness of the following integral
\begin{align*}
\int_0^\infty \lambda^{4m-n-1}&\widetilde{\rho}_i(\lambda)\widetilde{\rho}_j(\lambda)\chi(\lambda) \\
\times &\Big\langle \Big(M_{i,j}^++\Gamma_{i,j}^{+,r}(\lambda)\Big)Q_jv\Big(R_0^{+}-R_0^{-}\Big)(\lambda^{2m})(\cdot,y),\, Q_ivR_0^{-}(\lambda^{2m})VR_0^{-}(\lambda^{2m})(\cdot,x)\Big\rangle d\lambda.
\end{align*}

We shall demonstrate the $L^p$-boundedness of $W_{i,j}^L$ by  analyzing its integral kernel in three cases:
$$|x|> 2|y|,\,\, |x|< \frac12|y|\, \,\,  {\rm and}\,\,\,  \frac12|y|\le |x|\le 2|y|.$$ The corresponding operators for these cases are $W_{i,j,1}^L$, $W_{i,j,2}^L$ and $W_{i,j,3}^L$, respectively.

For the first two cases, we write $W_{i,j,1}^L(x,y)$ and $W_{i,j,2}^L(x,y)$ , respectively as follows:
$$\mathbf{1}_{\{|x|> 2|y|\}}\int_0^\infty\chi(\lambda) \Big(e^{i\lambda(|x|+|y|)}T_{i,j}^+(\lambda,x,y)+e^{i\lambda(|x|-|y|)}T_{i,j}^-(\lambda,x,y)\Big)d\lambda,$$
$$\mathbf{1}_{\{|x|< \frac12|y|\}}\int_0^\infty\chi(\lambda) \Big(e^{i\lambda(|x|+|y|)}T_{i,j}^+(\lambda,x,y)+e^{i\lambda(|x|-|y|)}T_{i,j}^-(\lambda,x,y)\Big)d\lambda.$$
Here $T_{i,j}^\pm(\lambda;x,y)$ satisfy some point-wise estimates (see Lemmas \ref{lemma-T-even-1}, \ref{lemma-T-even-2} and \ref{lemma-T-even-2-1}).  Hence, by using  the oscillatory integral estimate, we will prove that $W_{i,j,1}^L$ is in $\mathbb{B}(L^p(\mathbb{R}^n))$ for $1< p \le \infty$ ( see Proposition \ref{propo-low energy-even-1} ), and $W_{i,j,2}^L$ is in $\mathbb{B}(L^p(\mathbb{R}^n))$ for $1\le p < \sigma_i$ (see Proposition \ref{propo-low energy-even-2}  for the specific value $\sigma_i$). The upper bounds on the range of $p$ in Theorems \ref{thm-main result-0} and  \ref{thm-main result-1}  are derived from the results established in this case.

For the last term $W_{i,j,3}^L$, it suffices to establish the $L^p$-boundedness of both the good parts $I_{i,j,1,t}^{\pm,g}$ (for $t=2,3$) and $I_{i,j,s}^{\pm,g}$ (for $s=2,3$), as well as the bad part $I_{i,j,1,1}^{\pm,b}$, for all $1<p<\infty$, which are precisely defined in equations \eqref{eq-def-Iijb} - \eqref{eq-def-Iij-g2}. The \textbf{good parts} consist of admissible kernels that can be handled using the oscillatory integral estimates developed previously. 
To analyze the \textbf{bad part}, we employ dyadic decomposition to split the operator into admissible kernels and a singular kernel. The singular kernel takes the essential form
\begin{equation*}
\frac{y^\alpha}{|y|^{\frac{n-1}{2}+\vartheta(j)}} \frac{x^{\tilde{\alpha}}}{|x|^{\frac{n-1}{2}+\vartheta(i)}} \frac{1}{|x|\pm|y|},
\end{equation*}
where $|\alpha|=|\tilde{\alpha}|=\vartheta(j)$. 
The $L^p$-boundedness of the corresponding operators for $1<p<\infty$ follows from the $L^p$-boundedness of the \textbf{truncated Hilbert transform} on the line (see Proposition~\ref{lem.estforwij3b1} for details). This completes the proof for the low energy part.

\subsubsection{The high energy part $W^H$}

We now outline the proof for the high-energy part. Since Erdoğan-Green \cite{EG22, EG23} have established that \( W^H \) is bounded on \( L^p(\mathbb{R}^n) \) for all \( 1 \leq p \leq \infty \) when \( n > 2m \) , we focus on the case \( 1 \leq n \leq 2m \) and $m\ge 1$. We break this into two cases: \( n = 2m \) and \( n < 2m \).  The proof can be carried out using the methods from \cite{Yajima-JMSJ-95}, \cite{Yajima-99-CMP} and \cite{EG22}. Notably, Galtbayar-Yajima \cite{GY23, GY25} have shown that \( W^H \) is bounded on \( L^p(\mathbb{R}^n) \) for \( 1 < p < \infty \) using frequency-space methods for \( (m,n) = (2,4) \) and $(m,n) = (2,2)$, respectively. But our  proof gives the full range $1\le p\le \infty$
 for $m\ge2$ and all even dimensions $n\le2m$.

%We now sketch the proof of the high energy part.
%Since Erdo\u{g}an-Green \cite{EG22,EG23} have  proved that $W^H$ is bounded on $L^p(\R^n)$ for all $1\le p\le\infty$ when $n>2m$, we shall focus on the case $1\le n\le 2m$ here.  We decompose the cases into $n=2m$ and $n< 2m$. It should be mentioned that Yajima \cite{Yajima-99-CMP} and Galtbayar-Yajima \cite{GY23} proved that the $L^p$-boundedness of  $W^H$  for $1<p<\infty$ by using the analysis method in the frequency space when $(m,n)=(1,2)$ and $(2,4)$. Our method can extend this result to the endpoints.

When $n=2m$, applying the second resolvent identity, $W^H$ can be divided as follows:
\begin{equation*}
	W^H=W^H_1+W^H_2,
\end{equation*}
where
\begin{equation*}
	\begin{split}
		W_1^H=& \int_0^\infty\lambda^{2m-1}(1-\chi(\lambda))R_0^+(\lambda^{2m})V(R_0^+(\lambda^{2m})-R_0^-(\lambda^{2m}))d\lambda,\\
		W_2^H=&\int_0^\infty\lambda^{2m-1}(1-\chi(\lambda))R_0^+(\lambda^{2m})VR^+(\lambda^{2m})V(R_0^+(\lambda^{2m})-R_0^-(\lambda^{2m}))d\lambda.
	\end{split}
\end{equation*}
We rewrite $W_1^H$ as $$W_1^H=W_{1,1}^H\circ \frac{(1-\chi(D))}{|D|^{1+\delta}},$$ 
for a small number $0<\delta<1/10$. Here, $W_{1,1}^Hf(x)$  can be represented by the following integral:
\begin{align*}
	W_{1,1}^H f(x)=
    \int_{S^{n-1}} \int_{\R^n}\int_{\omega(y-x)}^{\infty} K(t,\omega,y)  f (x-2\omega(\omega\cdot x)-y+2\omega(\omega\cdot y)-2t\omega)dt dy d\omega ,\end{align*}
where
$$K(t,\omega,y)=\int_0^{\infty}  e^{-ist}s^{1+\delta}\, \widehat{V}(s \omega) h_{s,\omega}(y) ds.$$
Then, by Minkowski inequality  we obtain the $L^p$-boundedness of $W_{1,1}^H$ by showing that 
$K(t,\omega,y)\in L^{1}(\R \times S^{n-1} \times \R^n ).$
While $W_2^H$ can be estimated into an admissible operator by using oscillatory integral estimates due to better decay on $\lambda$.

When $n<2m$, a similar decomposition of $W^H$ is provided, and it can be shown that $W^H$ is an admissible operator by applying  oscillatory integral estimates again.
\subsection{ Notations}
 We first list some common notations and conventions as follows:
\begin{itemize}
\item  $A\lesssim B$ means $A \leq  CB$, where $C>0$ is an absolute constant. $A\lesssim_l B$ means $A \leq  C_lB$ with a constant $C_l$ dependent on $l$.
 $\mathbb{N}_0=\{0,1,\cdots\}$,  $\mathbb{N}^+=\{1,2,\cdots\}$,  $\mathbb{Z}=\{0,\pm1, \pm2,\cdots\}$.
\vspace{0.1cm}
\item $L^p=L^p(\R^n;\mathbb{C})$ for $1\le p\le \infty$, and $L^{2}_\sigma(\mathbb{R}^n)=\{f; \langle x\rangle^\sigma f(x)\in L^2 \}$.  $L^{p,q}(\R^n)$  denotes the Lorentz space for $1\le p\le \infty$ and  $1\le q\le \infty$ (cf. \cite{gra}).
\vspace{0.1cm}
\item $\lfloor l\rfloor$ denotes the greatest integer less than and equal to  $l$. 
\vspace{0.1cm}
\item $\mathbb{B}(X,Y)$ denotes the space of bounded operators from Banach space $X$ to Banach space $Y$, and $\mathbb{B}(X)=\mathbb{B}(X,X)$.
\vspace{0.1cm}
\item $\mathbf{1}_{\Omega}$ denotes the characteristic function of a set $\Omega\subset \mathbb{R}^n$.
\vspace{0.1cm}
\item $T(x,y)$  denotes  the integral kernel of an operator $T$, namely, 
$$(Tf)(x)=\int T(x,y) f(y) dy.$$

\end{itemize}
%\vspace{0.2cm}
\begin{itemize}
\item $m_n$  in \eqref{equ0.1},   denotes {\bf  the total number of zero resonance types} of $H=(-\Delta)^m+V$ on $\mathbb{R}^n$ when $n$ is even and $1\le n\le 4m.$
\vspace{0.1cm}
\item $k_c$  in \eqref{defkc},   denotes {\bf the critical resonance type} of $H$ in the sense that the range of the value $p$ of $L^p$-bound of wave operators will shrink into  a finite interval when  $\k>k_c$.
\vspace{0.1cm}
\item $J_\k$ is the index set with $0\le \k\le m_n$, which is defined in \eqref{eq-Jk-even} and \eqref{eq-Jk-even-2}.
\vspace{0.1cm}
\item $S_j$ is the  orthogonal projection on $L^2(\R^n)$ for each $j\in J_\k,$ which satisfies the following cancellations: 
\[S_j(x^\alpha v)=0, \quad \text{for}~\alpha\in {\mathbb{N}^n_0},~|\alpha|\leq\max\{0,\lfloor j\rfloor\}.\]
See \eqref{eq-S_j-even} and \eqref{eq-S_j-even-1}.
\vspace{0.1cm}
\item $Q_j$    denotes another orthogonal projection on $L^2(\R^n)$ for each $j\in J_\k,$ which satisfies that
\[Q_j(x^\alpha v)=0, \quad \text{for ~$ \alpha\in {\mathbb{N}_0^n}, ~~|\alpha|\leq\max\{0,\lfloor j+1/2\rfloor\}-1$}.\]
See  \eqref{eq-Q_j-even} and \eqref{eq-cancellation-Q_j}.

\end{itemize}
 
%\subsection{The organization of paper}
%In section \ref{sec:3},  we are devoted to express the properties of the free resolvent, the asymptotic behavior of $\big(M^\pm(\lambda)\big)^{-1}$ as $\lambda \to 0^+$ and some integral kernel estimates, which will be used to prove the $L^p$-boundedness of $W_\pm$. 
%In section \ref{sec:low energy}, we analyze the low energy part of $W_\pm$.
%The high energy part of $W_\pm$ is studied in section \ref{sec:high energy}.
%Section \ref{sec-unb} is devoted to proving the $L^p$-unboundedness of $W_\pm$ when $k_c<\mathbf{k}\le m_n+1.$
 %In section \ref{sec-proof-for-sec2.2}, we present the proof of the asymptotic expansion of $\big(M^\pm(\lambda)\big)^{-1}$. 
%Finally, three key lemmas are proved in section \ref{sec:proof}.

%\vspace{1cm}

%\section{Preliminaries}\label{sec:3}
%This section is dedicated to  studying the properties of the free resolvent of $(-\Delta)^m$, the asymptotic behavior of $\big(M^\pm(\lambda)\big)^{-1}$ as $\lambda \to 0^+$,  and some integral kernel estimates, which will be used to prove the $L^p$-boundedness of $W_\pm\big(H;(-\Delta)^m\big)$. 

\vspace{1.2em}
\section{Properties of the free resolvent of $(-\Delta)^m$}\label{sec-2.1}
In this section, we recall some basic properties of the free resolvent of $(-\Delta)^m$ with $m\ge1$. 
For $z\in \mathbb{C}\setminus[0,+\infty),$  set
\begin{equation*}
R_0(z)=((-\Delta)^m-z)^{-1},\ \ \ R(-\Delta;z)=(-\Delta-z)^{-1}.
\end{equation*}
Note that $R_0(z)$ has the following expression (see e.g. Feng et al.\cite{FSWY20}):
\begin{equation}\label{eq-R0z}
R_0(z)=\frac{1}{mz^{1-\frac{1}{m}}}\sum_{k=0}^{m-1}e^{i\frac{2\pi k}{m}}R(-\Delta;e^{i\frac{2\pi k}{m}}z^{\frac{1}{m}}),
\end{equation}
where  $\arg(z^{\frac{1}{m}})\in (0, \frac{2\pi}{m})$.
Hence by the formula \eqref{eq-R0z} and the limiting absorption principle for the free resolvent $R(-\Delta;\,z)$ (see e.g. Agmon \cite{Agmon}), it follows that 
\begin{equation}\label{defR0}
R_0^\pm(\lambda^{2m})=\big((-\Delta)^m-(\lambda^{2m}\pm i0)\big)^{-1}=\frac{1}{m\lambda^{2m}}\sum_{k\in I^\pm}\lambda_k^2R^\pm(-\Delta;\lambda_k^2),\quad \lambda>0,
\end{equation}
where $\lambda_k=\lambda e^{i\frac{k\pi}{m}}$, $R^\pm(-\Delta;\lambda_k^2)=(-\Delta-(\lambda_k^2\pm i0))^{-1}$ and
\begin{equation}\label{eq-I}
	I^+=\{0,1,\ldots,m-1\},\ \ \ \  I^-=\{1,\ldots,m\}.
\end{equation}
Thus, the integral kernels of $R_0^\pm(\lambda^{2m})$ can be obtained by
\begin{equation}\label{eq-free resovelt-LAP}
	R_0^\pm(\lambda^{2m})(x,y)=\frac{1}{m\lambda^{2m}}\sum_{k\in I^\pm}\lambda_k^2R^\pm(-\Delta;\lambda_k^2)(x,y),
\end{equation}
where
\begin{equation*}
	R(-\Delta;z)(x,y)=\frac{i}{4}\Big(\frac{z^{\frac{1}{2}}}{2\pi|x-y|}\Big)^{\frac{n}{2}-1}H^{(1)}_{\frac{n}{2}-1}(z^{\frac{1}{2}}|x-y|).
\end{equation*}
Here $\Im z^{\frac{1}{2}}\ge 0$ and $H^{(1)}_{\frac{n}{2}-1}=J_{\frac n2-1}(z) + iY_{\frac n2-1}(z)$ is the first Hankel function of the order $\frac n2 -1$ (c.f. \cite[(3.1)]{J80}). Then, we have 
\begin{equation} \label{eq: expression of R0-1}
\begin{split}
R_0^+(\lambda^{2m})(x,y)=& \frac{i}{ 4m\lambda^{2m}}\sum_{k=0}^{m-1}\lambda_k^2\Big(\frac{\lambda_k}{2\pi|x-y|}\Big)^{\frac{n}{2}-1}H^{(1)}_{\frac{n}{2}-1}(\lambda_k|x-y|), \\	
R_0^-(\lambda^{2m})(x,y)=& \frac{i}{ 4m\lambda^{2m}}\sum_{k=1}^{m}\lambda_k^2\Big(\frac{\lambda_k}{2\pi|x-y|}\Big)^{\frac{n}{2}-1}H^{(1)}_{\frac{n}{2}-1}(\lambda_k|x-y|).
\end{split}
\end{equation}
It is well known that the Hankel functions are holomorphic in $z$ throughout the $z$-plane cut along the negative real axis. Below, we list some properties of Hankel functions of a fixed integer order $\nu\in \N_0$, which can be found in \cite[Chapter~9]{AS}.
\begin{itemize}
	\item The expansion of the Hankel function: 
	\begin{equation}\label{eq:expansion of Hankel}
		\begin{split}
			H_{\nu}^{(1)}(z) &= J_{\nu}(z) + iY_{\nu}(z) \\
			&= (z/2)^{\nu} \sum_{k=0}^{\infty} \frac{(-z^{2}/4)^{k}}{k!\,\Gamma(\nu+k+1)} - \frac{i}{\pi}(z/2)^{-\nu} \sum_{k=0}^{\nu-1} \frac{(\nu-k-1)!}{k!} (z^{2}/4)^{k} \\
			&\quad + \frac{2i}{\pi} \ln(z/2)(z/2)^{\nu} \sum_{k=0}^{\infty} \frac{(-z^{2}/4)^{k}}{k!\,\Gamma(\nu+k+1)} \\
			&\quad - \frac{i}{\pi}(z/2)^{\nu} \sum_{k=0}^{\infty} \bigl( \psi(k+1) + \psi(\nu+k+1) \bigr) \frac{(-z^{2}/4)^{k}}{k!\,\Gamma(\nu+k+1)}.
		\end{split}	
	\end{equation} 
	Here $\psi(k)$ denotes the digamma function, with $\psi(k) = \sum_{j=1}^{k-1} \frac{1}{j} - \gamma$ for $k\ge 1$, where $\gamma$ is the Euler constant.
	
	\item The analytic continuation of Hankel functions: 
	\begin{equation}\label{eq: relationship of Hankel functions}
	H^{(1)}_{\nu}(z e^{\pi i})=-e^{-\nu \pi i} H^{(2)}_{\nu}(z),\qquad H^{(2)}_{\nu}(z e^{\pi i})=-e^{\nu \pi i} H^{(1)}_{\nu}(z),  
	\end{equation} 
     where  $H^{(2)}_{\nu}(z)=J_{\nu}(z)-iY_{\nu}(z)$.
     
    \item Asymptotic expansions for large arguments:
    \begin{equation}\label{eq: Asymptotic expansions for large arguments}
    \begin{split}
    &H^{(1)}_{\nu}(z)\sim \sqrt{2/(\pi z)} e^{i(z-\frac12 \nu \pi-\frac{1}{4}\pi) },\ -\pi<\arg\, z<2\pi, \\
    &H^{(2)}_{\nu}(z)\sim \sqrt{2/(\pi z)} e^{-i(z-\frac12 \nu \pi-\frac{1}{4}\pi)},\ -2\pi<\arg\, z<\pi. 
    \end{split}
    \end{equation}
\end{itemize}

Next, we provide some expansions of the integral kernel of $R_0^\pm(\lambda^{2m})$ by the \eqref{eq-free resovelt-LAP}.
\begin{lemma}\label{lemma-free resolvent expansio-smalln}
Let $1\le n\le 2m$ be even and $\theta\in\mathbb{N}_0$. Then for any $\lambda>0$, we have the following expansions:

\noindent{\rm(i)}~ If $0\le \theta\le 2m-n-1$, one has
\begin{equation}\label{eq-free expansion-even-1}
\begin{split}
R_0^\pm(\lambda^{2m})(x,y)
=&\sum_{0\le j\le\lfloor \frac{\theta-1}{2}\rfloor}a_j^\pm\lambda^{n-2m+2j}|x-y|^{2j}+\lambda^{n-2m}r_{\theta}^\pm(\lambda|x-y|);
\end{split}
\end{equation}
 {\rm(ii)}~   ~ If $2m-n+1\le \theta\le 4m-n-1$, one has
\begin{equation}\label{eq-free expansion-even-2}
\begin{split}
R_0^\pm(\lambda^{2m})(x,y)
=&\sum_{0\le j\le\lfloor \frac{\theta-1}{2}\rfloor}a_j^\pm\lambda^{n-2m+2j}|x-y|^{2j}+b_0|x-y|^{2m-n}\log(\lambda|x-y|)+\lambda^{n-2m}r_{\theta}^\pm(\lambda|x-y|);
\end{split}
\end{equation}
{\rm(iii)}~  If $\theta=4m-n+1$, one has
\begin{equation}\label{eq-free expansion-even-3}
	\begin{split}
		R_0^\pm &(\lambda^{2m})(x,y)
		=\sum_{0\le j\le\lfloor \frac{\theta-1}{2}\rfloor}a_j^\pm \lambda^{n-2m+2j}|x-y|^{2j}+\\
		&\ \ \ \ \ b_0|x-y|^{2m-n}\log(\lambda|x-y|)+b_1\lambda^{2m}|x-y|^{4m-n}\log(\lambda|x-y|) +\lambda^{n-2m}r_{\theta}^\pm(\lambda|x-y|).
	\end{split}
\end{equation}
Here, $\lfloor l \rfloor$ denotes the greatest integer less than or equal to $l$.  In the  \eqref{eq-free expansion-even-1}--\eqref{eq-free expansion-even-3}, $a_j^\pm\in\mathbb{C}\setminus\mathbb{R}$ ,  $b_0,~b_1\in \R
\setminus\{0\}$ and all the remainder terms $r_{\theta}^\pm(z)\in C^\infty\big((0,\infty)\big)$ satisfy that 
\begin{equation}\label{eq-r estimate-even}
	\left|\frac{d^l}{dz^l}\Big( r_{\theta}^{\pm}(z)\Big)  \right| \lesssim_l
\begin{cases}
	|z|^{\theta-l}, &\quad |z|\leq 1,\\
	|z|^{\max\{\theta-l,\,-\frac{n-1}{2}\}}, &\quad |z|\geq 1.
\end{cases}
\end{equation}
\end{lemma}

\begin{lemma}\label{lemma-free resolvent expansio-largen}
	Let $2m< n\le 4m$ be even and $\theta\in\mathbb{N}_0$. Then for any $\lambda>0$, we have the following expansions:
	
\noindent{\rm(i)}   ~ If $2m-n+1\le \theta\le 4m-n-1$, one has
	\begin{equation}\label{eq-free expansion-even-4}
		\begin{split}
			R_0^\pm(\lambda^{2m})(x,y)
			=&b_0|x-y|^{2m-n}+\sum_{0\le j\le\lfloor \frac{\theta-1}{2}\rfloor}a_j^\pm\lambda^{n-2m+2j}|x-y|^{2j}+\lambda^{n-2m}r_{\theta}^\pm(\lambda|x-y|);
		\end{split}
	\end{equation}
	{\rm(ii)}~  If $\theta=4m-n+1$, one has
	\begin{equation}\label{eq-free expansion-even-5}
		\begin{split}
			R_0^\pm(\lambda^{2m})(x,y)
			=&b_0|x-y|^{2m-n}+\sum_{0\le j\le\lfloor \frac{\theta-1}{2}\rfloor}a_j^\pm \lambda^{n-2m+2j}|x-y|^{2j}+\\
			&b_1\lambda^{2m}|x-y|^{4m-n}\log(\lambda|x-y|) +\lambda^{n-2m}r_{\theta}^\pm(\lambda|x-y|).
		\end{split}
	\end{equation}
	Here in the  \eqref{eq-free expansion-even-4}--\eqref{eq-free expansion-even-5}, $a_j^\pm\in\mathbb{C}\setminus\mathbb{R}$ ,  $b_0,~b_1\in\R\setminus\{0\}$ and all the remainder terms $r_{\theta}^\pm(z)\in C^\infty\big((0,\infty)\big)$  also satisfy \eqref{eq-r estimate-even}.
\end{lemma}

\begin{lemma}\label{lemma-free resolvent expansion-difference}
In all dimensions $n\ge 1$, for $\theta\in \N_0$,  we have 
	\begin{equation}\label{eq-free expansion-difference}
		{\small\begin{split}
				\big(R_0^+- R_0^-\big)(\lambda^{2m})(x,y)
				=&\sum_{0\le j\le \lfloor\frac{\theta-1}{2}\rfloor}(a_j^+-a_j^-)\lambda^{n-2m+2j}|x-y|^{2j}+\lambda^{n-2m}\big(r_{\theta}^+-r_{\theta}^-\big)(\lambda|x-y|)),
		\end{split}}
	\end{equation}
where $a_j^\pm\in\mathbb{C}\setminus\mathbb{R}$ are the same as those in corresponding dimensions in Lemma \ref{lemma-free resolvent expansio-smalln} and Lemma \ref{lemma-free resolvent expansio-largen},   and all the remainder terms $r_{\theta}^\pm(z)\in C^\infty\big((0,\infty)\big)$  satisfy \eqref{eq-r estimate-even}.
\end{lemma}

Furthermore, we also provide an alternative representation of the integral kernels $ R_0^\pm(\lambda^{2m})(x, y)$ (see \cite[Lemma~6.2]{EG22} for the corresponding result in even dimensions $n \ge 2m$).
\begin{lemma}\label{lemma-free kernel}
Let $1\le n\le 4m-1$ be even, then we have
\begin{align}\label{eq-free kernel-F}
R_0^\pm(\lambda^{2m})(x,y)&=\frac{\lambda^{\frac{n+1}{2}-2m}}{|x-y|^{\frac{n-1}{2}}}e^{\pm i\lambda|x-y|}F_n^\pm(\lambda|x-y|). 
\end{align}
Here, for any $l\in\mathbb{N}_0,$  the functions $F_n^\pm(z)\in C^\infty((0,\,\infty))$ satisfy that
\begin{equation}\label{eq-F decay-odd}
 \big|\frac{d^l}{dz^l} F^\pm_n(z)\big|\lesssim_l |z| ^{-l}.
\end{equation}

\end{lemma}

\begin{lemma}\label{lemma-free kernel-1-1}
	For $n=4m$, one has
	\begin{align}\label{eq-free kernel-F-01}
R_0^\pm(\lambda^{2m})(x,y)&
=\frac{e^{\pm i\lambda|x-y|}}{|x-y|^{\frac n2}}F_{n,0}^\pm(\lambda|x-y|)
+\lambda^{\frac{1}{2}}\frac{e^{\pm i\lambda|x-y|}}{|x-y|^{\frac{n-1}{2}}}F_{n,1}^\pm(\lambda|x-y|), \\
\label{eq-free kernel-F-02}
R_0^{+}(\lambda^{2m})(x,y)-R_0^{-}(\lambda^{2m})(x,y)&=\lambda^{\frac{1}{2}}
\frac{e^{i\lambda|x-y|}}{|x-y|^{\frac{n-1}{2}}}\widetilde{F}_{n}^+(\lambda|x-y|)-\lambda^{\frac{1}{2}}\frac{e^{-i\lambda|x-y|}}{|x-y|^{\frac{n-1}{2}}}\widetilde{F}_{n}^-(\lambda|x-y|).
	\end{align}
	Here, for any $l\in\mathbb{N}_0,$  the functions  $F_{n,0}^\pm(z)$, $F_{n,1}^\pm(z)$ and $\widetilde{F}_{n}^{\pm}(z)$ satisfy \eqref{eq-F decay-odd}.	
\end{lemma}

We will now briefly provide the proof for the lemmas \ref{lemma-free resolvent expansio-smalln} to \ref{lemma-free kernel}.
\begin{proof}[{\bf Proof of Lemma \ref{lemma-free resolvent expansio-smalln}--\ref{lemma-free kernel}.}]
If $\lambda|x-y|<1$ is small enough and $n$ is even, it follows from the expansion of Hankel functions in \eqref{eq:expansion of Hankel} that
\begin{equation}\label{eq-epansion-small-argument}
\begin{cases}
R^+(-\Delta;\lambda_k^2)(x,y)&=|x-y|^{2-n}\sum_{j=0}^{\infty}B_j(\lambda_k|x-y|),\quad  \text{for}\ k=0,1,\ldots,m-1, \vspace{0.4cm}\\
R^-(-\Delta;\lambda_k^2)(x,y)&=|x-y|^{2-n}\sum_{j=0}^{\infty}B_j(\lambda_k|x-y|),\quad  \text{for}\ k=1,\ldots,m,
\end{cases}
\end{equation}
where  $\lambda_k=\lambda e^{i\frac{k\pi}{m}}$ and 

\begin{align*}
	B_j(z)=
	\begin{cases}
		\mu_j z^{2j},&\text{if}\,\, 0\le j\le \frac{n}{2}-2,\\[4pt]
		\upsilon_j z^{2j}+\mu_j z^{2j}\log(z),&\text{if}\,\,j\ge \frac{n}{2}-1,\\[4pt]
	\end{cases}
\end{align*}
with  $\upsilon_j\in \mathbb{C}\setminus{\mathbb{R}}$, $\mu_j\in \R $ and $\mu_j\neq 0$. Also, note that 
\begin{equation}\label{eq-sum-lambda_k}
\sum_{k=0}^{m-1}\lambda_k^{2j}=\sum_{k=1}^{m}\lambda_k^{2j}=0, \quad \text{for}\ j\in \Z \setminus (m\Z),	
\end{equation}
and  $$\log(\lambda_k|x-y|)=\frac{\i k\pi}{m}+\log(\lambda|x-y|).$$
Based on the identities in \eqref{eq-free resovelt-LAP},  the following expansion holds:
\begin{equation}\label{eq-expansion-2m-small-argument}
	{R}_0^{\pm}(\lambda^{2m})(x,y)=\lambda^{n-2m}\sum_{l=1}^{[\frac{n-1}{2m}]}b_{l-1}(\lambda |x-y|)^{2ml-n}+\lambda^{n-2m}\sum_{j=0}^{\infty}B_j^{1,\pm}(\lambda|x-y|),
\end{equation}
where
\begin{align*}
	B_j^{1,\pm}(z)=
	\begin{cases}
		a_{j}^\pm z^{2j},&\text{if}\  2j+n\in \Z \setminus 2m\Z ,\\[4pt]
		a_{j}^\pm z^{2j}+b_{\frac{2j+n}{2m}-1}z^{2j}\log(z),&\text{if}\,\,  2j+n\in 2m\Z.\\[4pt]
	\end{cases}
\end{align*}
In addition, we have  $b_j=\mu_{m(j+1) -1}$ and 
\begin{equation}\label{eq： def for a_j^pm}
a_j^\pm=
\begin{cases}
\frac{1}{m}\mu_{j+\frac n2 -1}\sum_{k\in I^\pm}e^{ i\frac{(2j+n)k\pi}{m}}, \quad & \text{if}\  2j+n\in \Z \setminus 2m\Z ,\\[4pt]
\upsilon_{j+\frac n2 -1}+\frac{i}{m}\mu_{j+\frac n2 -1}\sum_{k\in I^\pm}\frac{k\pi}{m}, \quad & \text{if}\  2j+n\in  2m\Z. \\[4pt]
\end{cases}
\end{equation}

Meanwhile, the asymptotic behaviors in \eqref{eq: relationship of Hankel functions} of Hankel functions for large arguments imply that for $\lambda|x-y|>c_0$ with a fixed $c_0>0$, ${R}_0^{\pm}(\lambda^{2m})(x,y)$ can be expressed as
\begin{equation}\label{eq-resolvent-large-argument-1}
{R}_0^{\pm}(\lambda^{2m})(x,y)=\frac{\lambda^{\frac{n+1}{2}-2m}}{|x-y|^{\frac{n-1}{2}}}e^{\pm i\lambda|x-y|}\Phi^\pm(\lambda|x-y|),	
\end{equation}
where $\Phi^\pm(z)\in C^\infty{(0, \infty)}$ satisfy
\begin{equation}\label{eq-resolvent-large-argument-2}
	\left| \frac{d^l}{dz^l}\Phi^\pm(z)\right| \lesssim_l \langle z\rangle^{-l}, \quad \text{for}\ l\in \N_0,\ |z|>c_0.
\end{equation}
At this point, the results in Lemmas~\ref{lemma-free resolvent expansio-smalln}--\ref{lemma-free kernel} can be derived from \eqref{eq-expansion-2m-small-argument}--\eqref{eq-resolvent-large-argument-2} through direct computations. For clarity, we only demonstrate the proof of \eqref{eq-free expansion-even-1} as an illustrative example in the following. The other arguments follow similarly.

Firstly, noting that $\lambda^{2m-n}R_0^\pm(\lambda^{2m})(x,y)$ can be expressed as a function of $\lambda|x-y|$, for $0 \leq \theta \leq 2m - n - 1$ with $n \leq 2m$, we define 
\begin{equation}\label{r_theta}
	\begin{split}
		r_{\theta}^\pm(\lambda|x-y|) = \lambda^{2m-n}R_0^\pm(\lambda^{2m})(x,y) 
	 - \sum_{0 \leq j \leq \lfloor \frac{\theta-1}{2} \rfloor} a_j^\pm \lambda^{2j} |x-y|^{2j}.
	\end{split}
\end{equation}
From \eqref{eq-expansion-2m-small-argument}, we obtain
\begin{equation*}
	r_{\theta}^\pm(z) = \sum_{j=\lfloor \frac{\theta+1}{2} \rfloor}^{\infty} B_j^{1,\pm}(z),
\end{equation*}
which implies that $r_{\theta}^\pm(z)$ satisfies \eqref{eq-r estimate-even} when $|z| \ll 1$. 

Next to consider the arguments for $|z|>c_0$.  Combining the definition \eqref{r_theta}  of $r_{\theta}^\pm(z)$ with \eqref{eq-resolvent-large-argument-1}, we have
\begin{equation*}
	r_{\theta}^\pm(z) = \frac{1}{|z|^{\frac{n-1}{2}}} e^{\pm iz} \Phi^\pm(z) - \sum_{0 \leq j \leq \lfloor \frac{\theta-1}{2} \rfloor} a_j^\pm |z|^{2j}.
\end{equation*}
Together with \eqref{eq-resolvent-large-argument-2}, this shows that $r_{\theta}^\pm(z)$ satisfies \eqref{eq-r estimate-even} when $|z| > c_0 > 0$. Thus, we complete the proof of \eqref{eq-free expansion-even-1}. 
\end{proof}

\begin{remark}
	By the definition in \eqref{eq： def for a_j^pm}, we have that 
	\begin{equation}\label{eq-equi-aj+=aj-}
		a_j^+e^{\frac{ i\pi(n-2m+2j)}{2m}}=a_j^-e^{\frac{ -i\pi(n-2m+2j)}{2m}}, \quad \text{for $2j+m \in \Z\setminus 2m\Z$},
	\end{equation}

	In even dimensions $n\le 4m$, for $j=0,\, 1,\ldots,\,2m-\frac{n}{2}-1$ and $j\neq m-\frac n2$, we have used the following algebra expansions in the asymptotic expansion for $M^{\pm}(\lambda^2)$:
	\begin{equation}\label{eq-def-Aalbet}
		a_j^\pm |x-y|^{2j}=\sum_{\substack{|\alpha|+|\beta|=2j}} A_{\alpha,\beta}^\pm x^\alpha y^\beta,
	\end{equation}
	where $a_j^\pm$ are given in Lemma \ref{lemma-free resolvent expansio-smalln} and Lemma \ref{lemma-free resolvent expansio-largen}, and $A_{\alpha,\beta}^{\pm}$ are constants in $\C$, which play important roles in the proofs of asymptotical expansions of $\big(M^{\pm}(\lambda)\big)^{-1}$ in Section \ref{sec-proof-for-sec2.2}, and can be calculated specifically in the following way.
\end{remark}

\begin{lemma}\label{lem-Aalbet}

	The constants $A_{\alpha,\beta}^\pm$ in \eqref{eq-def-Aalbet} can be expressed as
	\begin{equation}\label{eq-def-Aalbet-1}
		A_{\alpha,\beta}^\pm=A_{\alpha,\beta} \, e^{\frac{\mp i\pi(n-2m+|\alpha|+|\beta|)}{2m}},
	\end{equation}
	and
	\begin{equation}\label{eq-def-Aalbet-2}
		A_{\alpha,\beta}=
		\begin{cases}
			\frac{(-1)^{|\beta|}i^{|\alpha|+|\beta|}}{(2\pi)^n\alpha!\beta!}\int_{\R^n}\frac{\xi^{\alpha+\beta}}{1+|\xi|^{2m}}\,\d\xi, &\quad |\alpha|+|\beta|<2m-n, \vspace{0.5cm}\\
			-\frac{(-1)^{|\beta|}i^{|\alpha|+|\beta|}}{(2\pi)^n\alpha!\beta!}\int_{\R^n}\frac{\xi^{\alpha+\beta}}{|\xi|^{2m}(1+|\xi|^{2m})}\,\d\xi, &\quad 2m-n<|\alpha|+|\beta|<4m-n.	
		\end{cases}	
	\end{equation}	
\end{lemma}

\begin{proof}
	By the argument similar to \eqref{eq-expansion-2m-small-argument}, it follows that
	\begin{equation}\label{eq-expansion-2m-small-argument-1}
		{R}_0^{+}(-\lambda^{2m})(x,y)={R}_0^{-}(-\lambda^{2m})(x,y)=\lambda^{n-2m}\sum_{l=1}^{[\frac{n-1}{2m}]}b_{l-1}(\lambda |x-y|)^{2ml-n}+\lambda^{n-2m}\sum_{j=0}^{\infty}B_j^{1}(\lambda|x-y|),
	\end{equation}
	where
	\begin{align*}
		B_j^{1}(z)=
		\begin{cases}
			a_{j}z^{2j},&\text{if}\  2j+n\in \Z \setminus 2m\Z ,\\[4pt]
			a_{j}z^{2j}+b_{\frac{2j+n}{2m}-1}z^{2j}\log(z),&\text{if}\,\,  2j+n\in 2m\Z,\\[4pt]
		\end{cases}
	\end{align*}
	and
	\begin{equation}\label{eq: def for a_j}
		\begin{aligned}
		a_j=
		\begin{cases}
			 \frac{i}{m}\mu_{j+\frac n2 -1}\sum_{k\in I^\pm}exp{\left( i\frac{(2j+n-2m)(k\pm1/2)\pi}{m}\right) }, &\quad  \text{if}\  2j+n\in \Z \setminus 2m\Z, \\[4pt]
			exp\left(\frac{i(2j+n-2m)\pi}{2m}\right)\left( \upsilon_j+\frac{i}{m}\mu_{j+\frac n2 -1}\sum_{k\in I^\pm }\frac{(k\pm\frac12)\pi}{m}\right), &\quad  \text{if}\  2j+n\in 2m\Z.
		\end{cases}	
	\end{aligned}
\end{equation}
Here, we have used the fact \eqref{eq-sum-lambda_k}.
We remark that by the identity \eqref{eq: def for a_j}, one has 
\begin{equation}\label{eq:def of aj-1}
 a_j=a_j^+e^{\frac{ i\pi(n-2m+2j)}{2m}}=a_j^-e^{\frac{ -i\pi(n-2m+2j)}{2m}}, \quad \text{for $2j+n \in \Z\setminus 2m\Z$},   
\end{equation}
Thus, in even dimensions $n\le 4m$, for $j=0,\, 1,\ldots,\,2m-n/2-1$ and $j\neq m-n/2$, we have the following expression
	$$a_j|x-y|^{2j}=\sum_{\substack{|\alpha|+|\beta|=2j}} A_{\alpha,\beta} x^\alpha y^\beta,$$
	where $A_{\alpha,\beta}$ is given in \eqref{eq-def-Aalbet-1}. We consider \eqref{eq-expansion-2m-small-argument-1} with $\lambda=1$. Then when $n\le 2m$ and $|\alpha|+|\beta|< 2m-n$, it follows from \eqref{eq-expansion-2m-small-argument-1} and the dominated convergence theorem that
	\begin{equation}\label{eq-limits-xy-0}
		\begin{split}
			A_{\alpha,\beta}\alpha!\beta!
			=\lim\limits_{x,y\to0}\partial_x^\alpha\partial_y^\beta R_0^+(-1)(x-y)
			&=\lim\limits_{x,y\to0}\frac{1}{(2\pi)^n}\partial_x^\alpha\partial_y^\beta\left(\int_{\R^n}\frac{e^{\i(x-y)\cdot\xi}}{|\xi|^{2m}+1}\,\d\xi\right)
			\\&=\frac{\i^{|\alpha|+|\beta|}(-1)^{|\beta|}}{(2\pi)^n}\int_{\R^n}\frac{\xi^{\alpha+\beta}}{|\xi|^{2m}+1}\, \d\xi.
		\end{split}
	\end{equation}
	Note that $b_0|x|^{2m-n}\log|x|$ is the fundamental solution of $(-\Delta)^m$ in even dimensions $n\le 2m$ and $b_0|x|^{2m-n}$ is the fundamental solution of $(-\Delta)^m$ in even dimensions $n> 2m$. Therefore, in the distributional sense, we have that \text{if the dimension $n\le 2m$ is even},
	\begin{equation*}
		(-\Delta)^m\partial^{\alpha+\beta}(b_0|\cdot|^{2m-n}\log|\cdot|)=\partial^{\alpha+\beta}(-\Delta)^m(b_0|\cdot|^{2m-n}\log|\cdot|)=\partial^{\alpha+\beta}\delta_0, 
	\end{equation*}
	and \text{if the dimension $n> 2m$ is even},
	\begin{equation*}
		(-\Delta)^m\partial^{\alpha+\beta}(b_0|\cdot|^{2m-n})=\partial^{\alpha+\beta}(-\Delta)^m(b_0|\cdot|^{2m-n})=\partial^{\alpha+\beta}\delta_0.
	\end{equation*}
	Taking Fourier transform in both sides of the above two identities, we then have
	\begin{equation*}
		\begin{split}
			\mathscr{F}(\partial^{\alpha+\beta}(b_0|\cdot|^{2m-n}\log|\cdot|))(\xi)&=i^{|\alpha|+|\beta|}\ \xi^{\alpha+\beta}\ |\xi|^{-2m}, \quad n\le 2m,\\
			\mathscr{F}(\partial^{\alpha+\beta}(b_0|\cdot|^{2m-n}))(\xi)&=i^{|\alpha|+|\beta|}\ \xi^{\alpha+\beta}\ |\xi|^{-2m}, \quad n> 2m,	
		\end{split}	
	\end{equation*}
	where $i^{|\alpha|+|\beta|}\ \xi^{\alpha+\beta}\ |\xi|^{-2m}$ is a tempered distribution.
	Applying \eqref{eq-expansion-2m-small-argument-1} and the dominated convergence theorem, 
	when the dimension $n\le 2m$ is even and $2m-n<|\alpha|+|\beta|< 4m-n$, the following identity can be derived:
	\begin{equation}\label{eq-Aalpbeta-limit-1}
		\begin{split}
A_{\alpha,\beta}\alpha!\beta!&=\lim\limits_{x,y\to0}\partial_x^\alpha\partial_y^\beta \left( R_0^+(-1)(x-y)-b_0|x-y|^{2m-n}\log|x-y|\right) \\	&=\lim\limits_{x,y\to0}\lim\limits_{\varepsilon\to0^+}\frac{1}{(2\pi)^n}\left(\partial_x^\alpha\partial_y^\beta\int_{\R^n}\frac{e^{\i(x-y)\cdot\xi-\varepsilon|\xi|^2}}{|\xi|^{2m}+1}\,\d\xi-(-1)^{|\beta|}i^{|\alpha|+|\beta|}\int_{\R^n}\frac{\xi^{\alpha+\beta}e^{\i(x-y)\cdot\xi-\varepsilon|\xi|^2}}{|\xi|^{2m}}\,\d\xi\right)\\
			&=-\frac{\i^{|\alpha|+|\beta|}(-1)^{|\beta|}}{(2\pi)^n}\int_{\R^n}\frac{\xi^{\alpha+\beta}}{|\xi|^{2m}(|\xi|^{2m}+1)}\, \d\xi.	
		\end{split}
	\end{equation}
	Similarly, when $2m <n\le 4m$ is even and $2m-n<|\alpha|+|\beta|< 4m-n$, we have 
	\begin{equation}\label{eq-Aalpbeta-limit-2}
		\begin{split}
			A_{\alpha,\beta}\alpha!\beta!&=\lim\limits_{x,y\to0}\partial_x^\alpha\partial_y^\beta \left( R_0^+(-1)(x-y)-b_0|x-y|^{2m-n}\right) \\
			&=-\frac{\i^{|\alpha|+|\beta|}(-1)^{|\beta|}}{(2\pi)^n}\int_{\R^n}\frac{\xi^{\alpha+\beta}}{|\xi|^{2m}(|\xi|^{2m}+1)}\, \d\xi.	
		\end{split}
	\end{equation}
	Now, \eqref{eq-def-Aalbet-2} can be obtained by \eqref{eq-limits-xy-0}, \eqref{eq-Aalpbeta-limit-1} and \eqref{eq-Aalpbeta-limit-2}.
\end{proof}

%\vspace{1.2em}
\section{The asymptotic expansion for $\big(M^{\pm}(\lambda)\big)^{-1}$ with zero  resonances.}\label{sec-expansion}
As seen in \eqref{eq-wave operator-low-1}, a key step in establishing the $L^p$-boundedness of the wave operators $W_\pm(H; (-\Delta)^m)$ in the low-energy regime involves deriving the asymptotic expansion of $\big(M^\pm(\lambda)\big)^{-1}$ near the zero-energy threshold. In this section, we complete this analysis.

We first define some notations used throughout the paper.  Let $m\ge 1$ and $n$ be an even integer. Recall that 
\begin{equation*}\
	m_n=\max\big\{m+1,\ 2m-\frac n2+1 \big\}=
	\begin{cases}
		m+1,&\ \ \mbox{if}\ \  1\le n\le 2m,\\
		2m-\frac n2+1,&\ \ \mbox{if}\ \ 2m<n\le 4m,
		%	0,\quad\qquad\,\, \,\,\mbox{if}\, n\ge 4m+1,
	\end{cases}
\end{equation*}
 which is the total number of zero-resonance types of $H=(-\Delta)^{m}+V$ on $\R^n$.

Let $\k$ be an integer and   $0\le \k\le m_n+1$.  In the even dimension $n\leq 2m$, we use  $J_{\k}$ to denote the following sets of real number:

\begin{equation}\label{eq-Jk-even}
	J_{\k}:=\begin{cases}
		 \big\{j\in\mathbb{N}_0;\,j\le m-\frac{n}{2}\big\}\cup\ \big\{m-\frac{n}{2}+\frac{1}{2}\big\},\quad&\mbox{if}\   \k=0,\\[0.3cm]
		\{m-\frac{n}{2}+\frac{1}{2},\ m-\frac{n}{2}+\frac{3}{4}\}\cup\{j\in\mathbb{N}_0;\,j\le m-\frac{n}{2}+\mathbf{k}-1\},\quad&\mbox{if}\ 1\le  \k\le m_n,\\[0.3cm]
		J_{m_n}\cup \{2m-\frac{n}{2}+\frac{1}{2}\},&\mbox{if}\  \k=m_n+1.
	\end{cases}
\end{equation}
For $2m< n\leq 4m$ and $n$ is even,  we denote
\begin{equation}\label{eq-Jk-even-2}
	J_{\k}:=\begin{cases}\{m-\frac{n}{2}\},\quad&\mbox{if}\, \, \k=0,\\[0.3cm]
		\{m-\frac{n}{2}\}\cup\{j\in\mathbb{N}_0;\,0\leq j\le \k-1 \},\quad&\mbox{if}\, \, 1\le \k\le m_n,\\[0.3cm]
		J_{m_n}\cup\{2m-\frac{n}{2}+\frac12\},&\mbox{if}\,\, \k=m_n+1.
	\end{cases}
\end{equation}
Note that the set $J_{\mathbf{k}}$ may contain fractional indices in certain cases. From the definition of $J_{\mathbf{k}}$, it follows directly that
\begin{equation}\label{eq-max-Jk}
    \max J_{\mathbf{k}} \leq \max\left\{m - \frac{n}{2}, 0\right\} + \mathbf{k} + 1.
\end{equation}

To show the main results of this section,  we define some operators as follow
\begin{equation}\label{eq-G-even}
	G_{2j}f=\int_{\R^n}|x-y|^{2j}f(y)\,\d y,\ \ \ \ G^{0}_{2j}f=\int_{\R^n}|x-y|^{2j}\log|x-y|f(y)\,\d y,\quad j\in J_{m_n+1},
\end{equation}
and
\begin{equation}\label{eq-T_0-even}
	T_0=
\begin{cases}
U+b_{0}vG^{0}_{2m-n}v, &\quad n\le 2m, \\[0.2cm]
U+b_{0}vG_{2m-n}v, &\quad 2m< n\le 4m. 
\end{cases}	
\end{equation}
where $v(x)=\sqrt{|V(x)|}, U(x)=\text{sgn} V(x)$. Here,  $b_0|x|^{2m-n}\log|x|$ for $n\le 2m$ and $b_0|x|^{2m-n}$ for $2m<n\le 4m$ are the fundamental solutions  of $(-\Delta)^m$  where the constant  $b_0$ is given in Lemma \ref{lemma-free resolvent expansio-smalln} and Lemma \ref{lemma-free resolvent expansio-largen}. 
\vskip0.3cm
We next introduce some subspaces of $L^2(\R^n)$ and their associated orthogonal projections $S_j$.  We define $S_{j}=I$ when $j\le -1$.  In the next,  $X^\perp$ denotes the the orthogonal complement to  a subspace  $X$ of $L^2(\R^n)$,  and  $\rm{Ker}(T)$  denotes the kernel  space of some operator $T$ on $L^2(\R^n)$.  

For $j\in J_{m_n}$ and a fixed integer   $0\le \k\le m_n$, if $n$ is even and $1\le n \le2m$, we define
{\small
\begin{equation}\label{eq-S_j-even}
{\small
S_j L^2=\begin{cases}
\{x^\alpha v;\,\,|\alpha|\le j\}^\perp,&\mbox{if}\,0\le j\le m-\frac{n}{2},\\[0.2cm]
\mbox{ker}\,(S_{m-\frac{n}{2}} T_0 S_{m-\frac{n}{2}})\bigcap S_{m-\frac{n}{2}} L^2,\,&\mbox{if}\,j=m-\frac{n}{2}+\frac{1}{2},\\[0.2cm]
\mbox{ker}\,(S_{m-\frac{n}{2}-1} T_0 S_{m-\frac{n}{2}})\bigcap S_{m-\frac{n}{2}+\frac12} L^2,\,&\mbox{if}\,j=m-\frac{n}{2}+\frac{3}{4},\\[0.2cm]
\{x^\alpha v;\,\,\,|\alpha|\le j\}^\perp\bigcap\mbox{ker}(S_{2m-n-j-1}T_0S_{m-\frac{n}{2}}),
&\mbox{if}\,m-\frac{n}{2}+1\le j\le\max J_{\k},\\[0.2cm]
\{0\},&\mbox{if}\,j> \max J_{\k};
\end{cases}
}
\end{equation}
}
For $j\in J_{m_n}$, if $n$ is even and $2m<n\leq 4m$, we define
\begin{equation}\label{eq-S_j-even-1}
S_j L^2=\begin{cases}
\mbox{ker}\,(T_0) ,\,&\mbox{if}\,\ j=m-\frac{n}{2},\\[0.2cm]
\{x^\alpha v;\,\,\,|\alpha|\le j\}^\perp
\bigcap S_{m-\frac{n}{2}} L^2,
&\mbox{if}\,\ 0\le j\le \max J_{\k},\\[0.2cm]
\{0\},&\mbox{if}\,\ j>\max J_{\k}.
		\end{cases}
\end{equation}
And we also define $$S_{2m-\frac{n}{2}+\frac 12}L^2=\{0\}.$$

\begin{remark}
\mbox{}
\begin{itemize}

\item   In the  case of  a fixed resonance type  $\mathbf{k}$ of  $H=(-\Delta)^m+V$,   we emphasize that these projections $S_j$ are well-defined on $L^2$ by the decay condition on $V$ in Assumption \ref{Assumption},  and  one has $S_j=0$ if  $j\in J_{m_n+1}\setminus J_{\k}$. 
\item  From the definition of $S_jL^2$ above,  
 the following  inclusion relations hold for  $1\le n\le 2m$,
\begin{equation*}
	S_0\ge S_1\ge \ldots \ge  S_{m-\frac{n}{2}+\frac12}\ge S_{m-\frac{n}{2}+\frac 34}\ge \ldots\ge S_{2m-\frac{n}{2}+\frac 12}=0,
\end{equation*}
and  for $2m< n\le 4m$, 
\begin{equation*}
	S_{m-\frac{n}{2}}\ge S_0\ge S_1\ge \ldots\ge S_{2m-\frac{n}{2}+\frac 12}=0.
\end{equation*}
 In particular,   zero is the regular point of $H$ ( i.e. $\k=0$ )  is equivalent that $S_{m-\frac{n}{2}+\frac12}L^2=\{0\}$ when $n\le 2m$ and $S_{m-\frac{n}{2}}L^2=\{0\}$ when $2m<n\le 4m$.
\item Furthermore, for each $j\in J_{m_n+1}$,  one always has the cancellation relations:
\begin{equation}\label{CforQ_j}
	S_j(x^\alpha v)=0, \quad \text{for}~\alpha\in {\mathbb{N}^n_0},~|\alpha|\leq\max\{0,\lfloor j\rfloor\},
\end{equation} where $\lfloor j\rfloor$ denotes the greatest integer at most $j$. In particular, the integer function $\lfloor j\rfloor$ (not $j$ itself) exactly captures the cancellation order of all projections $S_j$ with $j\in J_{m_n+1}$. This actually explains the motivations of choosing three fractional values of $j=m-\frac{n}{2}+\frac{1}{2},\ m-\frac{n}{2}+\frac{3}{4},\  2m-\frac{n}{2}+\frac 12$, as the index elements of $J_{m_n+1}$. 
\end{itemize}
\end {remark}
\vskip0.3cm
\begin{remark}
In the remark, we will
 give some motivation for defining the index sets $J_{\k}$ and these particular subspaces $\big\{S_jL^2\big\}_{j\in J_{\k}}$. Applying the expansions of $R_0^{\pm}(x,y)$ for $\theta=4m-n+1$ in lemma \ref{lemma-free resolvent expansio-smalln} and Lemma \ref{lemma-free resolvent expansio-largen}, $M^\pm(\lambda)$ can be expanded in the  increasing power order of $\lambda$ as follows. In even dimensions $n\le 2m$, we have
\begin{equation}\label{eq-expansion-Mpm-3-0}
	\begin{aligned}
		M^\pm(\lambda)&=\sum\limits_{0\le j\le 2m -\frac n2}
		a_j^\pm\lambda^{n-2m+2j}vG_{2j}v+b_0\log\lambda vG_{2m-n}v\\
		&+T_0+b_1\lambda^{2m}\log\lambda vG_{4m-n}v+b_1\lambda^{2m}vG^0_{4m-n}v+\lambda^{n-2m}vr_{4m-n+1}^\pm(\lambda)v,
	\end{aligned}
\end{equation}
and in even dimensions $2m<n\le 4m$, we have
\begin{equation}\label{eq-expansion-Mpm-4-0}
	\begin{aligned}
		M^\pm(\lambda)=&T_0+\sum\limits_{0\le j\le  2m-\frac n2}a_j^\pm\lambda^{n-2m+2j}vG_{2j}v+\\
		&b_1\lambda^{2m}\log\lambda vG_{4m-n}v+b_1\lambda^{2m}vG^0_{4m-n}v+\lambda^{n-2m}vr_{4m-n+1}^\pm(\lambda)v.
	\end{aligned}
\end{equation}
If the potential $V$ has a rapid decay ( $|V(x)| \lesssim \langle x \rangle^{-(8m-n+2)}$ )\footnote{
	The decay requirement for $V$ can be weakened when zero energy corresponds to a $\mathbf{k}$-th kind resonance of $H$ with $\mathbf{k} < m_n$. 
     When $\mathbf{k} < m_n$, a truncated expansion suffices, eliminating the need for the complete expansions \eqref{eq-expansion-Mpm-3-0} or \eqref{eq-expansion-Mpm-4-0}. 
	For details, we refer to Section \ref{sec-proof-for-sec2.2}.
}, then by virtue of the Hilbert-Schmidt property, each term in both expansions is $L^2$-bounded. 

To prove the invertibility of $M^\pm(\lambda)$ near zero, we need to establish their injectivity on $L^2$. As a representative case, we consider dimensions $n \leq 2m$. For any $f \in L^2(\mathbb{R}^n)$, when $n \leq 2m$ is even, $M^{\pm}(\lambda)f$ admits the asymptotic expansion based on \eqref{eq-expansion-Mpm-3-0}.

From perturbation theory arguments, if the most singular term $a_0^\pm\lambda^{n-2m}vG_{0}vf$ in the expansion of $M^{\pm}(\lambda)f$ doesn't vanish (which occurs when $f \notin \ker(vG_{0}v)$ since all coefficients $a_j \neq 0$), then $M^{\pm}(\lambda)f \neq 0$ in a neighborhood of zero. 

We observe that $\ker(vG_{0}v)$ coincides exactly with the space $S_0L^2 = \{v\}^\perp$. When $f \in \ker(vG_{0}v)$, the leading singular term vanishes, which requires us to examine whether the next singular term $a_1^\pm\lambda^{n-2m+2}vG_{2}vf$ disappears. This leads us to consider whether $f$ belongs to $\ker(vG_{2}v)$ restricted to $S_0L^2$, which is precisely the space $S_1L^2 = \{x^\alpha v : |\alpha|\leq 1\}^\perp$. 
Iterating this argument motivates the introduction of the subspaces $S_jL^2$ for $j \leq m - \frac{n}{2}$. 

For $f \in S_{m-\frac{n}{2}}L^2$, we must subsequently examine whether $T_0 f = 0$ holds on $S_{m-\frac{n}{2}}L^2$. This leads to the definition of the subspace $S_{m-\frac{n}{2}+\frac{1}{2}}L^2$. The fractional index $m-\frac{n}{2}+\frac{1}{2}$ appears because that
 the moment cancellation order for elements in $S_{m-\frac{n}{2}+\frac{1}{2}}L^2$ matches that of $S_{\lfloor m-\frac{n}{2}+\frac{1}{2} \rfloor}L^2=S_{ m-\frac{n}{2}}L^2$ by \eqref{CforQ_j}.
Other fractional indices arise from analogous considerations of cancellation orders. Note that $T_0=(I-S_{m-\frac n2})T_0+ S_{m-\frac n2}T_0$. The rest subspace are introduced by the same argument.
\end{remark}

In the following, we introduce another group of orthogonal projections $Q_j$ for each $j\in J_{m_n+1}$ by
\begin{equation}\label{eq-Q_j-even}
Q_j:=
\begin{cases}   
I-S_0,  &\quad \text{if $n\le 2m$ and $j=0$}, \\
I-S_{m-\frac n2},  &\quad \text{if $2m< n \le 4m$ and $j=m-\frac n2$}, \\
S_{j'}-S_j,   &\quad \text{else},     
\end{cases}
\end{equation}
where $$ j'=\max\Big\{l\in J_{m_n+1};\,\, \ \ l<j \Big\}.$$  When zero energy is the $\k$-th kind of resonance,   it then follows from the definition of $Q_j$ and \eqref{CforQ_j}  and  the following cancellation relations hold:
\begin{equation}\label{eq-cancellation-Q_j}
	Q_j(x^\alpha v)=0, \quad \text{for ~$ \alpha\in {\mathbb{N}_0^n}, ~~|\alpha|\leq\max\{0,\lfloor j+1/2\rfloor\}-1$}.	
\end{equation}

For convenience,  let 
\begin{equation}\label{thetaj}
\vartheta(j)=\max\{0,\lfloor j+1/2 \rfloor\}.
\end{equation}

\vskip0.3cm
By the definition of these orthogonal projections $\{Q_j\}_{j\in J_{\k}}$,  we can give   the following equivalent characterizations of zero-energy resonances of $H$. It should be mentioned that similar situations have been treated in \cite{CHHZ-2024,JK,JN01,SWY22}.

\begin{proposition}\label{prop-resonance-equivalence}
Let $1\le n \le 4m$ be even, $0 \le \mathbf{k} \le m_n + 1$, and assume that $V(x)$ satisfies \eqref{decayonV} in Assumption \ref{Assumption}. Then the following statements are equivalent:
\vskip0.1cm
\emph{(\romannumeral1)} Zero energy is a regular point of $H$ (i.e., $\mathbf{k} = 0$) if and only if $\sum_{j \in J_0} Q_j = I$;
\vskip0.1cm
\emph{(\romannumeral2)} Zero energy is the $\mathbf{k}$-th kind resonance with $1\le \mathbf{k} \le m_n + 1$ if and only if $\sum_{j \in J_{\mathbf{k}-1}} Q_j \neq I$ and $\sum_{j \in J_{\mathbf{k}}} Q_j = I$.
\end{proposition}
\begin{proof}
Here, we only provide the proof in dimensions $n\le 2m$ since the proof in dimensions $2m<n\le 4m$ follows similarly.

Let $j_\k = \max J_{\mathbf{k}}$. Observe that 
$$\sum_{j \in J_{\k}} Q_j = I-S_{j_{\k}}.$$ 
 Then the condition $\sum_{j \in J_{\mathbf{k}-1}} Q_j \neq I$ and $\sum_{j \in J_{\mathbf{k}}} Q_j = I$ is equivalent to $S_{j_{\k-1}} \neq 0$ and $S_{j_{\k}} = 0$ on $L^2$, respectively. Therefore, in view of Definitions \ref{defi-resonance-1} and \ref{defi-resonance-2}, to prove this proposition, it suffices to prove that for $\k=1,2,\ldots,m_n+1$,  
$$\psi \in S_{j_{\k-1}} L^2 \Leftrightarrow (-\Delta)^m \phi+V\phi=0,\quad \text{for some }~  0\neq \phi(x) \in \bigcap_{s < \mathbf{k} - m_n} L_s^2(\mathbb{R}^n), \quad  \psi(x) = Uv \phi(x).$$  
In particular, when $\k=1$, we have
\begin{equation}\label{decom of phi}
\phi(x)=\phi_e(x)+\sum_{|\alpha|= m-\frac{n}{2}}C_\alpha x^\alpha,	
\end{equation}
where $C_\alpha\in \C$, $\phi_e(x) \in \bigcap\limits_{s<-m+1}L_{s}^2(\mathbb{R}^n)$ and $\sum_{|\alpha|= m-\frac{n}{2}}C_\alpha x^\alpha\neq 0$.

We first prove the sufficiency. In even dimensions $n\le 2m$,  if $\psi\in S_{j_{\k-1}}L^2$ for $\k=1,2,\ldots,m_n+1$, then we have $\psi\in S_{m-\frac{n}{2}+\frac12}L^2$ since $j_{\k-1}\ge m-\frac{n}{2}+\frac12$. By the definition of $S_{m-\frac{n}{2}+\frac12}L^2$, it follows that
\begin{equation*}
U\psi(x)=-b_0vG^0_{2m-n}v\psi(x)+T_0\psi(x),
\end{equation*}
Note that $\psi\in S_{j_{\k-1}}L^2$ also implies that $S_{m-\frac{n}{2}-\k+1}T_0\psi(x)=0$, thus
\begin{align}\label{equ.psi}
\psi(x)=-Ub_0vG_{2m-n}^0v\psi(x)+U(I-S_{m-\frac{n}{2}-\k+1})T_0\psi(x). 
\end{align}
	The definition of $S_{m-\frac{n}{2}-\k+1}$ indicates that
	\begin{equation}(I-S_{m-\frac{n}{2}-\k+1})T_0\psi(x)=\begin{cases}
 \sum_{|\alpha|\le  m-\frac{n}{2}-\k+1}C_\alpha x^{\alpha}v, \quad  &\text{if}\   1\le\k \le m-\frac{n}{2}+1 ,\\
 0, \quad  &\text{if}\  m-\frac{n}{2}+1< \k \le m_n+1.\\
\end{cases}
 \end{equation}   
Set
\begin{equation}\label{eq:def for phi function}
\phi(x):=
\begin{cases}
-b_0G_{2m-n}^0v\psi(x)+\sum_{|\alpha|\le  m-\frac{n}{2}-\k+1}C_\alpha  x^{\alpha}, &\quad \text{if}~ 1\le\k \le m-\frac{n}{2}+1, \\
-b_0G_{2m-n}^0v\psi(x),
&\quad \text{if}~ m-\frac{n}{2}+1<\k.
\end{cases}
\end{equation}
Then one has $\psi=Uv\phi$ and
$$\phi(x)=-b_0G_{2m-n}^0V\phi(x)+\sum_{|\alpha|\le  m-\frac{n}{2}-\k+1}C_\alpha  x^{\alpha}, \quad \ \text{if} \ 1\le\k \le m-\frac{n}{2}+1,$$
as well as
$$\phi(x)=-b_0G_{2m-n}^0V\phi(x),
\quad \text{if}~ m-\frac{n}{2}-\k+1<0.$$
Since $b_0|x|^{2m-n}\log|x|$ is the fundamental solutions of $(-\Delta)^m$, one has 
$$(-\Delta)^{m}b_0G_{2m-n}^0V\phi(x)=V\phi.$$
Thus, we conclude that 
$$(-\Delta)^{m}\phi(x)+ V\phi(x)=0,$$
in the sense of distribution.
 
Note that 
\begin{equation}\label{eq: cancellation property-phi}
\int_{\R^n}v(x)x^\alpha \psi(x)dx=\int_{\R^n}V(x)x^\alpha \phi(x) dx=0,\quad \text{for all}\  \alpha\in \N_0^n,\  |\alpha|\le \lfloor j_{\k-1} \rfloor. 	
\end{equation}
It follows from Lemma \ref{lemma-kernel estimate-for-vanish-2} below that
\begin{align*}
	|G_{2m-n}^0v\psi(x)|&=\left| \int|x-y|^{2m-n}(\log|x-y|)v(y)\psi(y) dy \right|\\
    &\lesssim\langle x\rangle^{2m-n-\lfloor j_{\k-1}\rfloor-1}(1+\log\langle x\rangle),
\end{align*}
where by the definition,
\begin{equation*}
\lfloor j_{\k-1} \rfloor=
\begin{cases}
 m-\frac n2, \quad  &\text{for}\   \k =1,\\
 m-\frac n2+\mathbf{k}-2, \quad  &\text{for}\  2\le \k \le m_n+1.\\
\end{cases}
\end{equation*}
Therefore, we can conclude from the definition of $\phi$ in \eqref{eq:def for phi function} that $\phi(x)\in \bigcap_{s<\mathbf{k}-m_n}L_{s}^2(\mathbb{R}^n)=\bigcap_{s<\k-m-1}L_{s}^2(\mathbb{R}^n)$. In particular, if $\k=1$, and set
$$\phi_e=-b_0G_{2m-n}^0v\psi(x)+\sum_{|\alpha|\le  m-\frac{n}{2}-1}C_\alpha  x^{\alpha},$$
then the \eqref{decom of phi} holds with $\phi_e\in \bigcap_{s<-m+1}L_{s}^2(\mathbb{R}^n) $.

Conversely, in even dimensions $n\le 2m$, assume that there  exists  $\phi(x)\in \bigcap_{s<\mathbf{k}-m_n}L_{s}^2(\mathbb{R}^n)$ such that
$\big((-\Delta)^m+V\big)\phi=0$ for $\k=1,\ldots,m_n+1$, as well as \eqref{decom of phi} holds if $\k=1$. Choose $\eta(x)\in C_0^{\infty}(\R^n)$  such that $\eta(x)=1$ for $|x|\le 1$ and $\eta(x)=0$ for $|x|>2$. Let $\psi=Uv\phi$. Then for any $0<\delta<1$ and $\alpha\in \N_0^n$,  we obtain that  
\begin{equation}\label{equA.6}
		\begin{aligned}
			\left|\int_{\R^n}x^{\alpha}v(x)\psi(x)\eta(\delta x)dx \right|
                =&\left|\int_{\R^n}x^{\alpha}V(x)\phi(x)\eta(\delta x)dx \right|\\
			=& \left| \int_{\R^n}x^{\alpha}\eta(\delta x)(-\Delta)^m\phi(x)\,\d x\right|\\
		    =&\Big|\sum_{\alpha\geq\beta\   \atop |\beta|+|\gamma|=2m, } C_{\beta,\gamma}\delta^{|\gamma|}\int_{\R^n}x^{\alpha-\beta}\phi(x)(\partial^{\gamma}\eta)(\delta x)\,\d x\Big| \\
			\lesssim& \sum_{\alpha\geq\beta\   \atop |\beta|+|\gamma|=2m } \delta^{2m-\frac n2-|\alpha|-\theta} \left\|\langle x\rangle^{|\alpha|-|\beta|+\theta}\partial^{\gamma}\eta(x) \right\|_{L^2}  \left\|\langle x\rangle^{-\theta} \phi(x)\right\|_{L^2},
		\end{aligned}
\end{equation}
	where the last line follows from scaling and $\theta$ is defined as $m_n-\k+\varepsilon$ for some $0<\varepsilon<< 1$. Specially, when $\k=1$ and $n\le 2m$, we let $\theta=m_n-2+\varepsilon$ and $\phi$ is replaced by $\phi_e$ in the inequality \eqref{equA.6}.
Hence, let $\delta$ tend to $0$, then it follows that
	\begin{equation}\label{equA.7}
		\int_{\R^n}x^{\alpha}v(x)\psi(x)dx=0, \quad \text{ for $|\alpha|\le  {m}-\frac n2$ if $\k=1$ and $|\alpha|\le  {m}-\frac n2+\k -2$ if $\k>1$}.
	\end{equation}
This yields that $\psi(x)\in S_{m-\frac{n}{2}}L^2$ for $1\le \k \le m_n$. When $\k=1$, by the definition of $S_{m-\frac n2}$ and \eqref{equA.7}, we have 
	$$S_{m-\frac n2}vG_jv\psi(x)=S_{m-\frac n2}vG_jV\phi(x)=0,\quad \text{for}~ j=0,1,\ldots,m-\frac n2. $$
We conclude by the dominated convergence theorem and the expansion \eqref{eq-expansion-2m-small-argument-1} that for any $f\in C_0^\infty(\R^n)$, 
\begin{equation*}
	\begin{aligned}
		\langle S_{m-\frac{n}{2}}(v b_0G_{2m-n}^0v)S_{m-\frac{n}{2}}\psi, f\rangle&=\langle  S_{m-\frac{n}{2}} (vb_0G_{2m-n}^0)V\phi(x), f\rangle \\
		&=\lim_{\lambda\rightarrow0} \langle S_{m-\frac{n}{2}} vR_0(-\lambda^{2m})V\phi(x), f\rangle \\
		&=\lim_{\lambda\rightarrow0} \langle  S_{m-\frac{n}{2}} v\big((-\Delta)^{m}+\lambda^{2m}\big)^{-1}V\phi(x), f\rangle \\
		&= \langle v(-\Delta)^{-m}V\phi(x), S_{m-\frac{n}{2}}f\rangle.
	\end{aligned}
\end{equation*}
 Since $((-\Delta)^m+V)\phi(x)=0$,  it follows that in the distributional sense,
	\begin{equation*}
		\begin{aligned}
			S_{m-\frac{n}{2}}T_0S_{m-\frac{n}{2}}\psi=S_{m-\frac{n}{2}}T_0(Uv\phi) 
			&=S_{m-\frac{n}{2}}(U+b_0vG_{2m-n}^0v)Uv\phi\\
			&=S_{m-\frac{n}{2}}v(I+(-\Delta)^{-m}V)\phi \\
			&=S_{m-\frac{n}{2}}v(-\Delta)^{-m}((-\Delta)^{m}+V)\phi 
		=0,
		\end{aligned}
	\end{equation*}
	which implies that $\psi(x)\in S_{m-\frac{n}{2}+\frac12}L^2$.  If $\k=2$, repeating the above argument and noting \eqref{equA.7}, for any $f\in C_0^\infty(\R^n)$, 
	\begin{equation*}
		\begin{aligned}
		    \langle S_{m-\frac{n}{2}-1}(v b_0G_{2m-n}^0v)\psi, f\rangle 
			&=\langle  S_{m-\frac{n}{2}-1} (vb_0G_{2m-n}^0)V\phi(x), f\rangle \\
			&=\lim_{\lambda\rightarrow0} \langle S_{m-\frac{n}{2}-1} vR_0(-\lambda^{2m})V\phi(x), f\rangle \\
			&=\lim_{\lambda\rightarrow0} \langle  S_{m-\frac{n}{2}-1} v\big((-\Delta)^{m}+\lambda^{2m}\big)^{-1}V\phi(x), f\rangle \\
			&= \langle v(-\Delta)^{-m}V\phi(x), S_{m-\frac{n}{2}-1}f\rangle.
		\end{aligned}
	\end{equation*}
	Thus we have 
	$S_{m-\frac{n}{2}-1}T_0S_{m-\frac{n}{2}}\psi(x)=0,$
	which implies that $\psi(x)\in S_{m-\frac{n}{2}+\frac34}L^2$. Similarly, if $\k>2$,  we can obtain  
	$$S_{2m-n-j_{\k-1}-1}T_0S_{m-\frac{n}{2}}\psi(x)=0,$$
	which concludes that $\psi(x)\in S_{j_{\k-1}}L^2$.
	Therefore the proof is completed.	
\end{proof}

We now present the results on asymptotic expansion of $(M^\pm(\lambda))^{-1}$ near zero energy threshold, whose proof is provided in Section \ref{sec-proof-for-sec2.2}.

\begin{theorem} \label{thm-M inverse-even}
  Let $ n\le 2m$ be even and zero energy is the $\mathbf{k}$-th kind of resonance with $0\le \mathbf{k}\le m_n+1$.  Assume that $V(x)$ satisfies \eqref{decayonV} in Assumption \ref{Assumption}. Then there exists a small $\lambda_0\in(0,1)$ depending on $V, m, n$ such that $M^\pm(\lambda)$ is invertible on $L^2$ for all $0<\lambda<\lambda_0$, in addition, we have
 \begin{equation}\label{eq-M expansion-even}
 	(M^{\pm}(\lambda))^{-1}=\sum_{i,j\in J_{\mathbf{k}}}\lambda^{2m-n}\rho_i(\lambda)\rho_j(\lambda)Q_i(M_{i,j}^\pm+\Gamma_{i,j}^\pm(\lambda))Q_j,
 \end{equation}
 where $M_{i,j}^\pm$, $\Gamma_{i,j}^\pm(\lambda)$ are bounded operators on $L^2(\mathbb{R}^n)$, and $\rho_j(\lambda)$ is defined as 
 \begin{equation}\label{eq-def-for-rhoj-1}
 	\rho_j(\lambda)=\begin{cases}
 		\lambda^{-j}(\log\lambda)^{-\frac{1}{2}},\ \ &\  j= m-\frac{n}{2},\ 2m-\frac{n}{2},\\
 		\lambda^{-j+\frac12},\ \ &\  j= m-\frac{n}{2}+\frac{1}{2},\ 2m-\frac{n}{2}+\frac{1}{2},\\
 		\lambda^{-(m-\frac{n}{2})}(\log\lambda)^{\frac{1}{2}},\ \ &\  j= m-\frac{n}{2}+\frac{3}{4},\\
 		\lambda^{-j}, &\text{otherwise}.\\
 	\end{cases}
 \end{equation}
Meanwhile, $M_{i,j}^\pm=M_{j,i}^\pm=0$ for the following index pairs $(i,j)$:
\begin{equation}\label{eq-Mij=0}
\begin{cases}
i = m - \frac{n}{2} + \frac{3}{4}, & j \neq m - \frac{n}{2}, m - \frac{n}{2} + \frac{3}{4}, \\
i = 2m - \frac{n}{2}, & j \neq 2m - \frac{n}{2}, \\
i = 2m - \frac{n}{2} + \frac{1}{2}, & j \neq 2m - \frac{n}{2} + \frac{1}{2}, \\
\end{cases}
\end{equation}
 and 
$$M_{j,j}^+=M_{j,j}^- \quad  \text{ for $j=2m-\frac n2,\, 2m-\frac n2+\frac 12$}.$$  
 In addition, $\Gamma_{i,j}^\pm(\lambda)$  are $\mathbb{B}(L^2)$-valued $C^{\frac{n}{2}+1}$ functions of $\lambda$ which can be expressed as follows:
\begin{equation}\label{eq-expansion-for-Gamma-i,j}
 \Gamma_{i,j}^\pm(\lambda)=(\log\lambda)^{-\frac 12}\Gamma_{i,j}^{\pm, 0}+(\log\lambda)^{-1}\Gamma_{i,j}^{\pm, 1}+(\log\lambda)^{-\frac32}\Gamma_{i,j}^{\pm, 2}+(\log\lambda)^{-2}\Gamma_{i,j}^{\pm, 3}+\Gamma_{i,j}^{\pm, r}(\lambda),
\end{equation}
where $\Gamma_{i,j}^{\pm, s} \in \mathbb{B}(L^2)$ ($s=0,1,2,3$) are independent of $\lambda$ and $\Gamma_{i,j}^{\pm, r}(\lambda)$ satisfy
 \begin{equation}\label{eq-est-gamma-ij-1}
\left\|\frac{d^l}{d\lambda^l}\Gamma_{i,j}^{\pm,r}(\lambda)\right\|_{\mathbb{B}(L^2)}\lesssim \left|\log \lambda \right|^{-\frac 52}\lambda^{-l}, \ \ l=0,\,1,\,\cdots,\,\frac n2 +1.
 \end{equation}
Especially, when $i=j=2m-\frac{n}{2}$ or $2m-\frac{n}{2}+\frac{1}{2}$, for any $l=0,\,1,\,\cdots,\,\frac n2 +1$, one has
 \begin{equation}\label{eq-est-gamma-ij-2}
 	\begin{split}
 		&\left\|\frac{d^l}{d\lambda^l}\Gamma_{ 2m-\frac{n}{2},  2m-\frac{n}{2}}^\pm(\lambda)\right\|_{\mathbb{B}(L^2)} \lesssim \left|\log \lambda \right|^{-1}\lambda^{-l},\\
 		&\left\|\frac{d^l}{d\lambda^l}\Gamma_{2m-\frac{n}{2}+\frac{1}{2},2m-\frac{n}{2}+\frac{1}{2}}^\pm(\lambda)
 		\right\|_{\mathbb{B}(L^2)} \lesssim \lambda^{1-l}.
 	\end{split}
 \end{equation}
 \end{theorem}

\vskip0.3cm

 In the even dimensions $2m<n\le 4m$, there is a similar result.
 \begin{theorem} \label{thm-M inverse-even-1}
 	Let $2m <n\le 4m$ be even, zero energy is the $\mathbf{k}$-th kind of resonance with $0\le \mathbf{k}\le m_n+1$. Assume that $V(x)$ satisfies \eqref{decayonV} in Assumption \ref{Assumption}. Then there exists a small $\lambda_0\in(0,1)$ depending on $V, m, n$ such that $M^\pm(\lambda)$ is invertible on $L^2$ for all $0<\lambda<\lambda_0$, in addition, we have
 	\begin{equation}\label{eq-M expansion-even-1}
 		(M^{\pm}(\lambda))^{-1}=\sum_{i,j\in J_{\mathbf{k}}}\lambda^{2m-n}\rho_i(\lambda)\rho_j(\lambda)Q_i(M_{i,j}^\pm+\Gamma_{i,j}^\pm(\lambda))Q_j,
 	\end{equation}
 	where $M_{i,j}^\pm$, $\Gamma_{i,j}^\pm(\lambda)$ are bounded operators on $L^2(\mathbb{R}^n)$, and $\rho_j(\lambda)$ is defined as
 	\begin{equation}\label{eq-def-for-rhoj-2}
 		\rho_j(\lambda)=\begin{cases}
 			\lambda^{-j}(\log\lambda)^{-\frac{1}{2}},\ \ &\  j=2m-\frac{n}{2},\\
 			\lambda^{-(2m-\frac{n}{2})},\ \ &\  j= 2m-\frac{n}{2}+\frac{1}{2},\\
 			\lambda^{-j}, &\text{otherwise}.\\
 		\end{cases}
 	\end{equation}
 	Meanwhile, $M_{i,j}^\pm=0$ for 
 	\begin{equation}\label{eq-Mij=0-1}
 	\begin{cases}
 	i=m-\frac n2, j\neq m-\frac n2\ \text{or}\ \ j=m-\frac n2, i\neq m-\frac n2, &\\
 	i\neq2m-\frac{n}{2}+\frac{1}{2},j= 2m-\frac{n}{2}+\frac{1}{2}, &
	\\
 	\end{cases}
\end{equation}
and 
$$M_{j,j}^+=M_{j,j}^- \quad  \text{ for $j=2m-\frac n2,\, 2m-\frac n2+\frac 12$}.$$ 
 In addition, $\Gamma_{i,j}^\pm(\lambda)$  are $\mathbb{B}(L^2)$-valued $C^{\frac{n}{2}+1}$ functions of $\lambda$,  and satisfy \eqref{eq-expansion-for-Gamma-i,j}—\eqref{eq-est-gamma-ij-2}. Specially,  
    $$\Gamma_{i, m-\frac{n}{2}}^{\pm,s}(\lambda)=\Gamma_{m-\frac{n}{2}, i}^{\pm,s}(\lambda)=0\  \text{in \eqref {eq-expansion-for-Gamma-i,j} for s=0,1,2,3,}$$ 
    and for $l=0,\,1,\, \ldots, \frac n2+1$, one has that when $i\neq 2m-\frac n2$ or $2m<n<4m$, 
 \begin{equation}\label{eq-est-gamma-ij-2-1}
 		\begin{split}
\left\|\frac{d^l}{d\lambda^l}\Gamma_{i, m-\frac{n}{2}}^{\pm,r}(\lambda)\right\|_{\mathbb{B}(L^2)}
+\left\|\frac{d^l}{d\lambda^l}\Gamma_{ m-\frac{n}{2},  i}^{\pm,r}(\lambda)\right\|_{\mathbb{B}(L^2)} \lesssim  \lambda^{\frac n2-m-l}, 
 		\end{split}
 	\end{equation}
 	and when $i=2m-\frac n2$ and $n=4m$,
  \begin{equation}\label{eq-est-gamma-ij-2-2}
 	\begin{split}
 	\left\|\frac{d^l}{d\lambda^l}\Gamma_{i, m-\frac{n}{2}}^{\pm,r}(\lambda)\right\|_{\mathbb{B}(L^2)}
 	+\left\|\frac{d^l}{d\lambda^l}\Gamma_{ m-\frac{n}{2},  i}^{\pm,r}(\lambda)\right\|_{\mathbb{B}(L^2)} \lesssim\lambda^{\frac n2-m-l}|\log\lambda|^{\frac 12}.
 	\end{split}
 \end{equation}

 \end{theorem}
\vskip0.3cm

\begin{remark}
The similar results in odd dimensions $1\le n\le 4m$  appear in \cite{CHHZ-2024}. However, the presence of the logarithmic terms in the expansion of $R_0^{\pm}(\lambda^2)(x,y)$ makes obtaining the expansion of $(M^{\pm}(\lambda))^{-1}$ more difficult. It should be noted that there are alternative forms of expansions for $(M^{\pm}(\lambda))^{-1}$ for the fourth-order Schr\"odinger operator in dimensions $n=2$ and $n=4$, as demonstrated in \cite{LSY23} and  \cite{GY23}, respectively. 
\end{remark}
%\begin{remark}

%\vspace{1.2em}
\section{Three cancellation lemmas}
In this section, we will state the following three important lemmas, which take advantage of the cancellation property of $Q_j$ to study the structure of  kernels  $ \Big(Q_jvR_0^\pm(\lambda^{2m})(\cdot, y)\Big)(x)$ corresponding to  $Q_jvR_0^\pm(\lambda^{2m})$ by the following integrals
$$\Big(Q_jvR^\pm(\lambda^{2m})f\Big)(x)=Q_jv\int_{\R^n}R_0^\pm(\lambda^{2m})(x, y)\, f(y)dy.$$
 Their proofs  are given in Section \ref{sec:proof}.  
%\vskip0.4cm
\begin{lemma}\label{lemma-QjvR-even-1}
  Let $1\le n\le 4m-2$ be even, and $0\le \mathbf{k}\le m_n+1$. Assume that $V(x)$ satisfies \eqref{decayonV} in Assumption \ref{Assumption}.    Then for any  $j\in J_\mathbf{k}$ and $\lambda\in(0,1/2)$, we have
  \begin{equation}\label{eq-QjR kernel-even}
  \Big(Q_jvR_0^\pm(\lambda^{2m})(\cdot,\,y)\Big)(x)=e^{\pm i\lambda |y|}\omega_{j}^{\pm}(\lambda,x,y)
  \end{equation}
where for any $l=0,1,\ldots,\frac{n}{2}+1$, the functions $\omega_{j}^{\pm}(\lambda,x,y)$ satisfy
\begin{align}\label{eq-wj-estimate-even-1}
\big\|\partial_\lambda^l \omega_{j}^{\pm}(\lambda,x,y)\big\|_{L_x^2}&\lesssim \lambda^{\min\{n-2m+\vartheta(j),\, 0\}-l}, \quad \text{if}\  j\neq 2m- n,
\\ \label{eq-wj-estimate-even-2}
\big\|\partial_\lambda^l \omega_{j}^{\pm}(\lambda,x,y)\big\|_{L_x^2}&\lesssim_\varepsilon \lambda^{-\varepsilon-l}\langle y\rangle^{-\varepsilon} , \quad \text{for}\ j= 2m-n,\ n\le 2m, \  
\text{for all} \   \varepsilon \in (0,1/2), 
\\ \label{eq-wj-estimate-even-3}
\big\|\partial_\lambda^l \omega_{j}^{\pm}(\lambda,x,y)\big\|_{L_x^2}&\lesssim \lambda^{\min\{\frac{n+1}{2}-2m+\vartheta(j),\, 0\}-l}\langle y\rangle^{-\frac{n-1}{2}}.
\end{align}
Furthermore, when $j\in J_\mathbf{k}$ and $j> \max\{2m-n,\, 0\}$, then for any $l=0,1,\ldots,\frac{n}{2}+1$,  we have
\begin{equation}\label{eq-wj-estimate-even-4}
\big\|\partial_\lambda^l \omega_{j}^{\pm}(\lambda,x,y)\big\|_{L_x^2}\lesssim \lambda^{-l}\langle y\rangle^{-\min\{n-2m+\vartheta(j), \,\frac{n-1}{2}\}}.\\
\end{equation}
In the inequalities above,  $\vartheta(j)=\max\{0,\lfloor j+1/2 \rfloor\}$  as given in \eqref{thetaj}.
\end{lemma}

%\vskip0.4cm
\begin{lemma}\label{lemma-QjvR-even-1-1}
	Let $n= 4m$, and $0\le \mathbf{k}\le m_n+1$. Assume that $V(x)$ satisfies \eqref{decayonV} in Assumption \ref{Assumption}.    Then for any  $j\in J_\mathbf{k}$ and $\lambda\in(0,1/2)$, we have
	\begin{equation}\label{eq-QjR kernel-even-1-1}
		\Big(Q_jvR_0^\pm(\lambda^{2m})VR_0^\pm(\lambda^{2m})(\cdot,\,y)\Big)(x)=e^{\pm i\lambda |y|}\omega_{j}^{\pm}(\lambda,x,y)
	\end{equation}
	where for any $l=0,1,\ldots,\frac{n}{2}+1$, the functions $\omega_{j}^{\pm}(\lambda,x,y)$ satisfy
	\begin{align}\label{eq-wj-estimate-even-1-2}
		\big\|\partial_\lambda^l \omega_{j}^{\pm}(\lambda,x,y)\big\|_{L_x^2}\lesssim \lambda^{-l}\langle y\rangle^{-\frac{n-1}{2}}.
	\end{align}
\end{lemma}
%\vskip0.4cm
\begin{remark}
In dimension $n = 4m$, we need to analyze the terms $Q_j v R_0^\pm(\lambda^{2m}) V R_0^\pm(\lambda^{2m})(\cdot,y)$ instead of $Q_j v R_0^\pm(\lambda^{2m})(\cdot,y)$. This modification is necessary because the integral kernels $R_0^\pm(\lambda^{2m})(\cdot,y)$ fail to belong to $L^2_{\mathrm{loc}}(\mathbb{R}^n)$, which follows from their singular behavior $|\cdot - y|^{-n+2m}$ as established in \eqref{eq-free expansion-even-5}.
\end{remark}

%\vskip0.4cm
\begin{lemma}\label{lemma-QjvR+-R-even-1}
Let $1\le n\le 4m$ be even, and $0\le \mathbf{k}\le m_n+1$. Assume that $V(x)$ satisfies \eqref{decayonV} in Assumption \ref{Assumption}. Then for any  $j\in J_\mathbf{k}$ and $\lambda\in(0,1)$, we have
\begin{equation}\label{eq3.900}
  \Big(Q_jv(R_0^+-R_0^-)(\lambda^{2m})(\cdot,\,y)\Big)(x)=e^{ i\lambda |y|}\tilde{\omega}_{j}^{+}(\lambda,x,y)+e^{-i\lambda |y|}\tilde{\omega}_{j}^{-}(\lambda,x,y),
  \end{equation}
where the functions $\tilde{\omega}_{j}^{\pm}(\lambda,x,y)$ satisfy
\begin{align}\label{eq-k2-even-1}
\big\|\partial_\lambda^l \tilde{\omega}_{j}^{\pm}(\lambda,x,y)\big\|_{L_x^2}&\lesssim\lambda^{n-2m+\vartheta(j)-l},\\ \label{eq-k2-even-2}
\big\|\partial_\lambda^l \tilde{\omega}_{j}^{\pm}(\lambda,x,y)\big\|_{L_x^2}&\lesssim\lambda^{\frac{n+1}{2}-2m+\vartheta(j)-l}\langle y\rangle^{-\frac{n-1}{2}}.
\end{align}
\end{lemma}

%\vspace{0.4em}

We also need the following result on a more refined analysis of $Q_jvR_0^\pm(\lambda^{2m})(\cdot, x)$.

\begin{lemma}\label{lemma-QjvR-even-2}
Let $1\le n\le 4m$ be even, and $0\le \mathbf{k}\le m_n+1$. Assume that $V(x)$ satisfies \eqref{decayonV} in Assumption \ref{Assumption}.    Then for any  $j\in J_\mathbf{k}$ and $\lambda\in(0,1)$, we have
\begin{equation}\label{eq3.1000}
\omega_{j}^\pm(\lambda,x,y)=\omega_{j,1}^\pm(\lambda,x,y)+\omega_{j,2}^\pm(\lambda,x,y)+\omega_{j,3}^\pm(\lambda,x,y),
\end{equation}
where $\omega_{j}^\pm$ are given in Lemmas \ref{lemma-QjvR-even-1} and \ref{lemma-QjvR-even-1-1}.
Here, 
\begin{equation}\label{eq3.100001}
\omega_{j,1}^\pm(\lambda,x,y)=\mathbf{1}_{\{|y|>1\}}\, \lambda^{\frac{n+1}{2}-2m+\vartheta(j)}\cdot\sum_{|\alpha|=|\beta|=\vartheta(j)}A_{\alpha,\beta}^\pm \frac{y^\alpha \Big(Q_j(v z^\beta)\Big)(x)}{|y|^{\frac{n-1}{2}+\vartheta(j)}},
\end{equation}
and for any $l=0,\ldots,\frac{n}{2}+1$, each $\varepsilon\in [0,1]$, we have that the functions $\omega_{j,2}^\pm(\lambda,x,y)$ and $ \omega_{j,3}^\pm(\lambda,x,y)$ satisfy
\begin{equation}\label{eq3.1001}
\Big\|\partial^l_\lambda \omega_{j,2}^\pm(\lambda,x,y)\Big\|_{L_x^2}\lesssim 
\begin{cases}
 \lambda^{\frac{n+1}{2}-2m+\vartheta(j)-\varepsilon-l}\langle y\rangle^{-\frac{n-1}{2}-\varepsilon} \quad & j< 2m-\frac n2 ,\\	
 \lambda^{-l-\frac{\varepsilon}{2}}\langle y\rangle^{-\frac{n}{2}-\frac{\varepsilon}{2}}, \quad & j=2m-\frac n2, \, 2m-\frac n2+\frac 12,
\end{cases}
\end{equation}
 and 
\begin{equation}\label{eq3.10011}
	\Big\|\partial^l_\lambda \omega_{j,3}^\pm(\lambda,x,y)\Big\|_{L_x^2}\lesssim \lambda^{\frac{n+3}{2}-2m+\vartheta(j)}\langle y\rangle^{-\frac{n-1}{2}}. 
\end{equation}
In the inequalities above,  $\vartheta(j)=\max\{0,\lfloor j+1/2 \rfloor\}$ is given in \eqref{thetaj}.
\end{lemma}
%\vskip0.3cm
\begin{remark}
When $n\le 4m-2$, let $\{\psi_{s}^j\}_{s \in \Lambda_{j}}$ be an orthonormal basis of $Q_{j}L^2$.\footnote{The subspace $Q_{m-\frac{n}{2}+\frac{1}{2}}L^2$ is infinite-dimensional, whereas for $j \neq m-\frac{n}{2}+\frac{1}{2}$, the subspace $Q_{j}L^2$ is finite-dimensional.} We can represent $\omega_{j}^{\pm}(\lambda,x,y)$ as
	$$
	\omega_{j}^{\pm}(\lambda,x,y)=\sum_{s\in \Lambda_{j}} e^{\mp i\lambda|y|} \big\langle  R_0^\pm(\lambda^{2m})(\cdot,y), \psi_{s}^j \big\rangle\,  \psi_{s}^j(x) = \sum_{s\in \Lambda_{j}} \alpha_{s}^j(\lambda, y) \psi_{s}^j(x).
	$$
	In fact, we know that 
	$$
	\big\|\partial_\lambda^l \omega_{j}^{\pm}(\lambda,x,y)\big\|_{L_x^2} = \Big(\sum_{s\in \Lambda_{j}} |\lambda^l\alpha_{s}^j(\lambda, y)|^{2}\Big)^{1/2}.
	$$
\end{remark}

\section{Low energy estimates }\label{sec:low energy}
In this section, we are devoted to studying the $L^p$-boundedness of $W^L$ which is defined in \eqref{eq-wave operator-low-0}. 
Recall that $p_{\k}$ is a number defined as 
\begin{equation}\label{eq: def of p_k}
p_{\mathbf{k}}=
\begin{cases}
	\frac{2n}{2\mathbf{k}+n-2m-2}, &\quad k_c<\mathbf{k}\le m_n-1,\ 1\le n\le 2m,\\
	\frac{n}{\mathbf{k}+n-2m-1}, &\quad k_c<\mathbf{k}\le m_n-1,\ 2m< n\le 4m,\\
	2, &\quad \mathbf{k}=m_n,\ m_n+1.
\end{cases}	
\end{equation}
We will prove the following theorem.
\begin{theorem}\label{thm-low energy-regular-even}
Let $1\le n\le 4m$ be even. Assume that Assumption \ref{Assumption}, which depends on the resonance type $\mathbf{k}$, holds. Then the following statements hold:

\vspace{7pt}
\indent\emph{(\romannumeral1)} If $\mathbf{k}=0$ (zero is the regular point of $H$), then $W^L\in \mathbb{B}(L^p(\mathbb{R}^n))$ for all $1<p<\infty$;

\vspace{7pt}
\indent\emph{(\romannumeral2)} If  $1\le \mathbf{k}\le k_c$, then $W^L\in \mathbb{B}(L^p(\mathbb{R}^n))$ for all $1<p<\infty;$

\vspace{7pt}
\indent\emph{(\romannumeral3)} If  $k_c+1\le \mathbf{k}\le m_n-1$, then $W^L\in \mathbb{B}(L^p(\mathbb{R}^n))$ for all $1<p<p_{\k}$;

\vspace{7pt}
\indent\emph{(\romannumeral4)} If $\mathbf{k}=m_n,\, m_n+1$, then $W^L\in \mathbb{B}(L^p(\mathbb{R}^n))$ for all $1<p\le 2.$
\end{theorem}
The proof of the above theorem are presented in the next four subsections.

\subsection {The reductions of the integral kernel for \texorpdfstring{$W^L$}{WL}  }\label{sec-reprensentation}
In this subsection, we decompose  $W^L(x,y)$ into appropriate integral representations by using the asymptotic expansions of $\big(M^{\pm}(\lambda)\big)^{-1}$  in Theorems \ref{thm-M inverse-even}  and   \ref{thm-M inverse-even-1},  and  the formulas of $Q_jv R^\pm_0(\lambda^{2m})$ from Lemma \ref{lemma-QjvR-even-1} to Lemma \ref{lemma-QjvR-even-2}. Our studies naturally are divided into the following two cases: $n \leq 4m-2$ and $n = 4m$.

\subsubsection{Reduction in even dimensions $n\le 4m-2$.}
First, recall that  
\begin{equation*}
\begin{split}
W^L=&\frac{m}{\pi i}\int_0^\infty \lambda^{2m-1}\chi(\lambda)R_0^+(\lambda^{2m})v(M^+(\lambda))^{-1}(\lambda)v(R_0^+(\lambda^{2m})-R_0^+(\lambda^{2m}))d\lambda.
%:=&\frac{m}{\pi i}\int_0^\infty H(\lambda)\chi(\lambda)d\lambda.
\end{split}
\end{equation*}
If zero is a resonance of the $\mathbf{k}$-th kind,   then it  from the expansions of $(M^\pm(\lambda))^{-1}$ given in Theorem \ref{thm-M inverse-even} (see \eqref{eq-M expansion-even}) and Theorem \ref{thm-M inverse-even-1} (see \eqref{eq-M expansion-even-1}) follows that 
\begin{equation*}
	\begin{split}
		W^L=&\frac{m}{\pi i}\sum_{i,j\in J_\k}\int_0^\infty \chi(\lambda)\lambda^{4m-n-1}\rho_i(\lambda) \rho_j(\lambda)  R_0^{+}(\lambda^{2m})vQ_i\big(M_{i,j}^++\Gamma_{i,j}^+(\lambda)\big)Q_jv\big(R_0^{+}-R_0^{-}\big)(\lambda^{2m})d\lambda, 
		%:=&\frac{m}{\pi i}\sum_{i,j\in J_\k} \int_0^\infty  \chi(\lambda)H_{i,j}(\lambda)d\lambda,
	\end{split}
\end{equation*}
where the function $\rho_j(\lambda)$ for $j\in J_\mathbf{k}$ (identically, $\rho_i(\lambda)$ for $i\in J_\mathbf{k}$) is given in \eqref{eq-def-for-rhoj-1} or \eqref{eq-def-for-rhoj-2}. Specifically, when $1\le n\le 2m$,
\begin{equation}\label{eq-def-for-rhoj-1'}
 	\rho_j(\lambda)=\begin{cases}
 		\lambda^{-j}(\log\lambda)^{-\frac{1}{2}},\ \ &\  j= m-\frac{n}{2},\ 2m-\frac{n}{2},\\
 		\lambda^{-j+\frac12},\ \ &\  j= m-\frac{n}{2}+\frac{1}{2},\ 2m-\frac{n}{2}+\frac{1}{2},\\
 		\lambda^{-(m-\frac{n}{2})}(\log\lambda)^{\frac{1}{2}},\ \ &\  j= m-\frac{n}{2}+\frac{3}{4},\\
 		\lambda^{-j}, &\ \text{otherwise}.\\
 	\end{cases}
 \end{equation}
 when $2m< n\le 4m$,
\begin{equation}
\label{eq-def-for-rhoj-2'}
 		\rho_j(\lambda)=\begin{cases}
 			\lambda^{-j}(\log\lambda)^{-\frac{1}{2}},\ \ &\  j=2m-\frac{n}{2},\\
 			\lambda^{-(2m-\frac{n}{2})},\ \ &\  j= 2m-\frac{n}{2}+\frac{1}{2},\\
 			\lambda^{-j}, &\ \text{otherwise}.\\
 		\end{cases}
 \end{equation}
Let
\begin{equation}\label{H_ij}
H_{i,j}(\lambda)=\lambda^{4m-n-1}\rho_i(\lambda) \rho_j(\lambda)  R_0^{+}(\lambda^{2m})vQ_i\big(M_{i,j}^++\Gamma_{i,j}^+(\lambda)\big)Q_jv\big(R_0^{+}-R_0^{-}\big)(\lambda^{2m}).
\end{equation}
Then
$$W^L=\frac{m}{\pi i}\sum_{i,j\in J_\k} \int_0^\infty  \chi(\lambda)H_{i,j}(\lambda)d\lambda.$$
From \eqref{eq-expansion-for-Gamma-i,j} we make the following further decomposition: 
\begin{equation}\label{eq-reduction-to-log}
H_{i,j}(\lambda)=\sum_{s=0}^3H_{i,j,s}(\lambda)+H_{i,j,r}(\lambda)+H_{i,j,M}(\lambda),
\end{equation}
where
\begin{equation*}
\begin{aligned}
H_{i,j,s}(\lambda)&=	(\log\lambda)^{-\frac{1+s}{2}}\lambda^{4m-n-1}\rho_i(\lambda) \rho_j(\lambda)  R_0^{+}(\lambda^{2m})vQ_i\Gamma_{i,j}^{+,s}Q_jv\big(R_0^{+}-R_0^{-}\big)(\lambda^{2m}), \\
H_{i,j,r}(\lambda&)=\lambda^{4m-n-1}\rho_i(\lambda) \rho_j(\lambda)  R_0^{+}(\lambda^{2m})vQ_i\Gamma_{i,j}^{+,r}(\lambda)Q_jv\big(R_0^{+}-R_0^{-}\big)(\lambda^{2m}),
\end{aligned}
\end{equation*}
and
\begin{equation*}
H_{i,j,M}(\lambda)=	\lambda^{4m-n-1}\rho_i(\lambda) \rho_j(\lambda)  R_0^{+}(\lambda^{2m})vQ_iM_{i,j}^+Q_jv\big(R_0^{+}-R_0^{-}\big)(\lambda^{2m}).
\end{equation*}

Note that for a Schwartz function $f$, 
\begin{equation}\label{eq-repren-spectral-measure}
\big(R_0^{+}-R_0^{-}\big)(\lambda^{2m})f=\frac{i\pi}{m(2\pi)^n \lambda^{2m-1}}\int_{S^{n-1}}e^{i\lambda \omega x} \hat{f}(\lambda \omega)d\omega. 
\end{equation}
Then for any symbol $m(\lambda)$, the following commutating formula holds:
\begin{equation}\label{commuting}m(\lambda)\big(R_0^{+}-R_0^{-}\big)(\lambda^{2m})=\big(R_0^{+}-R_0^{-}\big)(\lambda^{2m})m(|D|),\end{equation}
where $D=-i\nabla$. Thus, one can apply the formula \eqref{commuting}  
to transfer some log-type factors in the \eqref{eq-def-for-rhoj-1'}-\eqref{eq-reduction-to-log} into   Fourier multiplier operators on the right side:
\begin{itemize}
\item $(\log\lambda)^{-\frac{1+s}{2}}\chi(\lambda)\big(R_0^{+}-R_0^{-}\big)(\lambda^{2m})=\chi(\lambda)\big(R_0^{+}-R_0^{-}\big)(\lambda^{2m})(\log|D|)^{-\frac{1+s}{2}}\chi(|D|/2);$
\vskip0.2cm
%And by the definition of $\rho_i(\lambda)$, in even dimensions $n\le 2m$ we also have that  
\item $\rho_{m-\frac n2}(\lambda)\chi(\lambda)\big(R_0^{+}-R_0^{-}\big)(\lambda^{2m})=\lambda^{-{m-\frac n2}}\chi(\lambda)\big(R_0^{+}-R_0^{-}\big)(\lambda^{2m})(\log|D|)^{-\frac{1}{2}}\chi(|D|/2);$
\vskip0.2cm
\item 
$\rho_{2m-\frac n2}(\lambda)\chi(\lambda)\big(R_0^{+}-R_0^{-}\big)(\lambda^{2m})=\lambda^{-{2m-\frac n2}}\chi(\lambda)\big(R_0^{+}-R_0^{-}\big)(\lambda^{2m})(\log|D|)^{-\frac{1}{2}}\chi(|D|/2).$
\end{itemize}
where we use the equality $\chi(\lambda)=\chi(\lambda)\chi(\lambda/2)$.

Note that  the Mihlin-H\"omander multiplier theorem (see \cite[Theorem 6.2.7]{gra}) can show that the operators $(\log|D|)^{-\frac{1+s}{2}}\chi(|D|/2)$ for $s=0,1,2,3$ are bounded on $L^p$ for $1<p<\infty$. Therefore, to prove Theorem \ref{thm-low energy-regular-even}, it suffices to prove $L^p$-boundedness of the following  reduced operators: 
\begin{equation}\label{eq-to-prove-part-1}
\int_0^\infty \lambda^{4m-n-1}\widetilde{\rho}_i(\lambda)\widetilde{\rho}_j(\lambda)  R_0^{+}(\lambda^{2m})vQ_i\Gamma_{i,j}^{+,s}Q_jv\big(R_0^{+}-R_0^{-}\big)(\lambda^{2m})\chi(\lambda) d\lambda, \quad s=0,1,2,3, 
\end{equation}
and 
\begin{equation}\label{eq-to-prove-part-1-1}
	\int_0^\infty \lambda^{4m-n-1}\widetilde{\rho}_i(\lambda)\widetilde{\rho}_j(\lambda)  R_0^{+}(\lambda^{2m})vQ_i\big(M_{i,j}^++\Gamma_{i,j}^{+,r}(\lambda)\big)Q_jv\big(R_0^{+}-R_0^{-}\big)(\lambda^{2m})\chi(\lambda)d\lambda,
\end{equation}
where the function $\widetilde{\rho}_j(\lambda)$ is given by 
\begin{equation}\label{eq-widetilde-rho-j}
\widetilde{\rho}_j(\lambda)=
\begin{cases}
\lambda^{-(m-\frac n2)}, \quad &\text{if $j=m-\frac n2$ and $n\le 2m$},\\
\lambda^{-(2m-\frac n2)}, \quad &\text{if $j=2m-\frac n2$},\\
\rho_j(\lambda), \quad &\text{other}.
\end{cases}
\end{equation}
%That is, when $n\le 2m$,
%\begin{equation}\label{eq-widetilde-%rho-j-1}
	%\widetilde{\rho}_j(\lambda)=
%\begin{cases}
	%\lambda^{-j+\frac12},\ \ &\  j= m-\frac{n}{2}+\frac{1}{2},\ 2m-\frac{n}{2}+\frac{1}{2},\\
	%\lambda^{-(m-\frac{n}{2})}(\log\lambda)^{\frac{1}{2}},\ \ &\  j= m-\frac{n}{2}+\frac{3}{4},\\
		%\lambda^{-j}, &\text{other},\\
	%\end{cases}
%\end{equation}
%when $2m<n\le 4m$,
%\begin{equation}\label{eq-widetilde-rho-j-2}
%	\widetilde{\rho}_j(\lambda)=
%\begin{cases}
%	\lambda^{-(2m-\frac{n}{2})},\ \ &\  j= 2m-\frac{n}{2}+\frac{1}{2},\\
%	\lambda^{-j}, &\text{other}.\\
%\end{cases}
%\end{equation}
Since the \eqref{eq-to-prove-part-1}, and \eqref{eq-to-prove-part-1-1} can be estimated in parallel, we only give the analysis for the integral \eqref{eq-to-prove-part-1-1} here.   In the sequel, we may denote   
\begin{equation}\label{G_ij}G_{i,j}(\lambda)=\lambda^{4m-n-1}\widetilde{\rho}_i(\lambda)\widetilde{\rho}_j(\lambda)  R_0^{+}(\lambda^{2m})vQ_i\big(M_{i,j}^++\Gamma_{i,j}^{+,r}(\lambda)\big)v\big(R_0^{+}-R_0^{-}\big)(\lambda^{2m}).
\end{equation}
The integral kernel $G_{i,j}(\lambda)(x,y)$ of $G_{i,j}(\lambda)$ can be represented as
\begin{equation}\label{eqforGij}
	\begin{split}
		G_{i,j}(\lambda)(x,y)=\lambda^{4m-n-1}\widetilde{\rho}_i(\lambda)\widetilde{\rho}_j(\lambda)\Big\langle \big(M_{i,j}^++\Gamma_{i,j}^{+,r}(\lambda)\big) Q_jv\big(R_0^{+}-R_0^{-}\big)(\lambda^{2m})(\cdot, y),\, Q_ivR_0^{-}(\lambda^{2m})(\cdot,x)\Big\rangle.
	\end{split}
\end{equation}
Here, we used the fact $R_0^{-}(\lambda^{2m})(\cdot,x)=R_0^{-}(\lambda^{2m})(x,\cdot)$. 
\begin{remark}
The induction from \eqref{H_ij} to \eqref{G_ij} for $1\le n\le4m$ will be instrumental in estimating $W_{i,j,3}^L$ below. 
\end{remark}
%\vspace{0.8cm}
\subsubsection {Reduction in the dimension $n= 4m$} We use the second resolvent identity 
$$R^{\pm}(\lambda^{2m})=R_0^{\pm}(\lambda^{2m})-R_0^{\pm}(\lambda^{2m})VR^{\pm}(\lambda^{2m}) $$
to write $W^{L}$ as 
\begin{align*}
W^L=&\frac{m}{\pi i}\int_0^\infty \lambda^{2m-1}\chi(\lambda)R_0^+(\lambda^{2m})V(R_0^+(\lambda^{2m})-R_0^+(\lambda^{2m}))d\lambda \\
&-\frac{m}{\pi i}\int_0^\infty \lambda^{2m-1}\chi(\lambda)R_0^+(\lambda^{2m})VR_0^+(\lambda^{2m})v(M^+(\lambda))^{-1}(\lambda)v(R_0^+(\lambda^{2m})-R_0^+(\lambda^{2m}))d\lambda\\
:=&W^{L,0}-\frac{m}{\pi i}\int_0^\infty H(\lambda)\chi(\lambda)d\lambda.
\end{align*}
Under the condition (iii) of Assumption~\ref{Assumption}, it was established by Erdo\u{g}an and Green in \cite{EG23} that $W^{L,0}$ is bounded on $L^p(\mathbb{R}^n)$ for all $1 \leq p \leq \infty$. Thus, it suffices to analyze the remained integral part $\int_0^\infty H(\lambda)\chi(\lambda)d\lambda$, where 
$$H(\lambda)=\lambda^{2m-1}R_0^+(\lambda^{2m})VR_0^+(\lambda^{2m})v(M^+(\lambda))^{-1}(\lambda)v(R_0^+(\lambda^{2m})-R_0^+(\lambda^{2m})).$$

Following the similar reductions of $H(\lambda)$ used for even dimensions $n \leq 4m-2$ above, we only need to prove the $L^p$-boundedness of the following reduced operators for each $i,j\in J_{\k}$ in the case $n=4m$:
\begin{equation}\label{eq-to-prove-part-1-0-1}
	\int_0^\infty \lambda^{4m-n-1}\widetilde{\rho}_i(\lambda)\widetilde{\rho}_j(\lambda)  R_0^{+}(\lambda^{2m})VR_0^{+}(\lambda^{2m})vQ_i\big(M_{i,j}^++\Gamma_{i,j}^{+,r}(\lambda)\big)v\big(R_0^{+}-R_0^{-}\big)(\lambda^{2m})\chi(\lambda) d\lambda,  
\end{equation}

As the same as in \eqref{G_ij} for the case $n\le 4m-2$,  We still use the following notation for $n=4m$: 
\begin{equation}\label{G_ij_1}G_{i,j}(\lambda)=\lambda^{4m-n-1}\widetilde{\rho}_i(\lambda)\widetilde{\rho}_j(\lambda)  R_0^{+}(\lambda^{2m})VR_0^{+}(\lambda^{2m})vQ_i\big(M_{i,j}^++\Gamma_{i,j}^{+,r}(\lambda)\big)v\big(R_0^{+}-R_0^{-}\big)(\lambda^{2m}),
\end{equation}
and  the corresponding integral kernel $G_{i,j}(\lambda)(x,y)$ are represented as
\begin{equation}\label{eqforGij-1}
	\begin{split}
		G_{i,j}(\lambda)(x,y)=&\lambda^{4m-n-1}\widetilde{\rho}_i(\lambda)\widetilde{\rho}_j(\lambda)\\
		\times&\Big\langle \big(M_{i,j}^++\Gamma_{i,j}^{+,r}(\lambda)\big) Q_jv\big(R_0^{+}-R_0^{-}\big)(\lambda^{2m})(\cdot, y),\, Q_ivR_0^{-}(\lambda^{2m})VR_0^{-}(\lambda^{2m})(\cdot,x)\Big\rangle.
	\end{split}
\end{equation}

%\vspace{0.3cm}
\subsubsection {Representation for the Integral Kernels of the Low-Energy Components}
By the above reduction, in order to prove Theorem\ref{thm-low energy-regular-even},  it suffices to  consider the $W_{i,j}^L$ for each $i,j\in J_{\k}$ given by
\begin{equation}\label{eq.rewrifrowl}
		W_{i,j}^L:=\int_0^\infty G_{i,j}(\lambda)\chi(\lambda)d\lambda,
\end{equation}
where $G_{i,j}(\lambda)$ is given by  the \eqref{G_ij} for $n\le 4m-2$ and the \eqref{G_ij_1} for $n=4m$. We then split $W_{i,j}^L$ into three parts as 
\begin{equation}\label{eqsplitwij}
	W_{i,j}^L=W_{i,j,1}^L+W_{i,j,2}^L+W_{i,j,3}^L,
\end{equation}
where their corresponding integral kernels satisfy
\begin{equation*}
\begin{cases}
W_{i,j,1}^L(x,y)=&{\bf1}_{\{|x|>2|y|\}}W_{i,j}^L(x,y),\\
W_{i,j,2}^L(x,y)=&{\bf1}_{\{|x|<\frac 12|y|\}}W_{i,j}^L(x,y),\\
W_{i,j,3}^L(x,y)=&{\bf1}_{\{\frac12|y|\le|x|\le 2|y|\}}W_{i,j}^L(x,y).
\end{cases}	
\end{equation*}
The proof of $L^p$-boundedness of $W_{i,j,1}^L$, $W_{i,j,2}^L$ and $W_{i,j,3}^L$ is given in the three subsections below, correspondingly. And Theorem \ref{thm-low energy-regular-even} can be immediately derived from these results.

\vspace{0.5cm}
{\bf (1) For the kernels $W_{i,j,1}^L(x,y)$ and $W_{i,j,2}^L(x,y)$.} Recall that by \eqref{eq-QjR kernel-even} in Lemma \ref{lemma-QjvR-even-1}, \eqref{eq-QjR kernel-even-1-1} in Lemma \ref{lemma-QjvR-even-1-1} and \eqref{eq3.900} in Lemma \ref{lemma-QjvR+-R-even-1}, we have
\begin{align*}
&\Big(Q_jvR_0^\pm(\lambda^{2m})(\cdot,\,y)\Big)(x)=e^{\pm i\lambda |y|}\omega_{j}^{\pm}(\lambda,x,y), \quad n\le 4m-2, \\
&\Big(Q_jvR_0^\pm(\lambda^{2m})VR_0^\pm(\lambda^{2m})(\cdot,\,y)\Big)(x)=e^{\pm i\lambda |y|}\omega_{j}^{\pm}(\lambda,x,y), \quad  n=4m,\\
&\Big(Q_jv(R_0^+-R_0^-)(\lambda^{2m})(\cdot,\,y)\Big)(x)=e^{ i\lambda |y|}\tilde{\omega}_{j}^{+}(\lambda,x,y)+e^{-i\lambda |y|}\tilde{\omega}_{j}^{-}(\lambda,x,y).
\end{align*}
Thus, we can represent $G_{i,j}(\lambda)(x,y)$ as 
\begin{equation}\label{eqexpfroGij}
	G_{i,j}(\lambda)(x,y)= e^{i\lambda(|x|+|y|)}T_{i,j}^+(\lambda;x,y)+e^{i\lambda(|x|-|y|)}T_{i,j}^-(\lambda;x,y),
\end{equation}
with
\begin{equation}\label{eq-def-Tij}
	\begin{aligned}
		T_{i,j}^{\pm}(\lambda;x,y)=&\lambda^{4m-n-1}\widetilde{\rho}_i(\lambda)\widetilde{\rho}_j(\lambda)\Big\langle\big(M_{i,j}^++\Gamma_{i,j}^{+,r}(\lambda)\big)\widetilde{\omega}_{j}^{\pm}(\lambda,\cdot,y),\, \omega_{i}^{-}(\lambda,\cdot,x)\Big\rangle, \\
	\end{aligned}
\end{equation}
Now, the integral kernels $W_{i,j,1}^L(x,y)$ and $W_{i,j,2}^L(x,y)$ can be written as 
\begin{equation}\label{eq.w-odd-1}
	W_{i,j,1}^L(x,y)
	=\mathbf{1}_{\{|x|>2|y|\}}\int_0^\infty \Big(e^{i\lambda(|x|+|y|)}T_{i,j}^+(\lambda; x,y)+e^{i\lambda(|x|-|y|)}T_{i,j}^-(\lambda;x,y)\Big)\chi(\lambda)d\lambda,
\end{equation}
\begin{equation}\label{eq.w-odd-2}
	W_{i,j,2}^L(x,y)
	=\mathbf{1}_{\{|x|<\frac 12|y|\}}\int_0^\infty \Big(e^{i\lambda(|x|+|y|)}T_{i,j}^+(\lambda;x,y)+e^{i\lambda(|x|-|y|)}T_{i,j}^-(\lambda;x,y)\Big)\chi(\lambda)d\lambda.
\end{equation}
\vskip0.5em
{\bf (2) For the kernel $W_{i,j,3}^L(x,y)$.} By using \eqref{eq3.1000} in Lemma \ref{lemma-QjvR-even-2}, we can rewrite $G_{i,j}(\lambda)(x,y)$ as
\begin{equation*}
	\begin{aligned}
		G_{i,j}(\lambda)(x,y)=&\sum_{s=1}^3
		e^{i\lambda|x|}\lambda^{4m-n-1}\widetilde{\rho}_i(\lambda)\widetilde{\rho}_j(\lambda) \Big\langle \big(M_{i,j}^++\Gamma_{i,j}^{+,r}(\lambda)\big)Q_jv(R_0^+-R_0^-)(\lambda)(\cdot,y),\, \omega_{i,s}^{-}(\lambda,\cdot,x)\Big\rangle\\
		:=&\sum_{s=1}^3G_{i,j,s}(\lambda)(x,y).
	\end{aligned}
\end{equation*}
{For the term $G_{i,j,1}(\lambda)(x,y)$},  applying \eqref{eq3.1000} again, we then have 
\begin{equation}\label{eq4.10101}
	G_{i,j,1}(\lambda)(x,y)=\sum_{s=1}^3\Big(e^{i\lambda(|x|+|y|)}T_{i,j,1,s}^{+}(\lambda;x,y)-e^{i\lambda(|x|-|y|)}T_{i,j,1,s}^{-}(\lambda;x,y)\Big),
\end{equation}
where
\begin{equation}\label{eq4.101012}
	T_{i,j,1,s}^{\pm}(\lambda;x,y)=\lambda^{4m-n-1}\rho_i(\lambda)\rho_j(\lambda)\Big\langle \big(M_{i,j}^++\Gamma_{i,j}^{+,r}(\lambda)\big)\omega_{j,s}^{\pm}(\lambda,\cdot,y),\, \omega_{i,1}^{-}(\lambda, \cdot,x)\Big\rangle, \ \ ~s=1,2,3.	
\end{equation}
{ For the terms $G_{i,j,2}(\lambda)(x,y)$ and $G_{i,j,3}(\lambda)(x,y)$}, using \eqref{eq3.900} in Lemma \ref{lemma-QjvR+-R-even-1}, we have that
\begin{equation}
	{G}_{i,j,s}(\lambda)(x,y)= e^{i\lambda(|x|+|y|)}{T}_{i,j,s}^{+}(\lambda;x,y)+e^{i\lambda(|x|-|y|)}{T}_{i,j,s}^{-}(\lambda;x,y),\quad \ \  ~s=2,\,3,
\end{equation}
where
\begin{equation}\label{eq.4.100123}
{T}_{i,j,s}^{+}(\lambda;x,y)=\lambda^{4m-n-1}\widetilde{\rho}_i(\lambda)\widetilde{\rho}_j(\lambda)
\Big\langle \big(M_{i,j}^++\Gamma_{i,j}^{+,r}(\lambda)\big)\tilde{\omega}_{j}^{\pm}(\lambda,\cdot,y),\, \omega_{i,s}^{-}(\lambda,\cdot,x)\Big\rangle.	
\end{equation}
Thus, the integral kernel  $W_{i,j,3}^{L}(x,y)$  in \eqref{eqsplitwij} can be written as a linear combination of the following terms:
\begin{equation}\label{eq-def-Iijb}
	I_{i,j,1,1}^{\pm,b}(x,y)=\mathbf{1}_{\{\frac12|y|\le|x|\le 2|y|\}}\int_0^\infty e^{i\lambda(|x|\pm|y|)}T_{i,j,1,1}^{\pm}(\lambda;x,y)\chi(\lambda)d\lambda,
\end{equation}
\begin{equation}\label{eq-def-Iij-g}
	I_{i,j,1,s}^{\pm,g}(x,y)=\mathbf{1}_{\big\{\frac12|y|\le|x|\le 2|y|\big\}}\int_0^\infty e^{i\lambda(|x|\pm|y|)}T_{i,j,1,s}^\pm(\lambda;x,y)\chi(\lambda)d\lambda, \quad s=2,3,
\end{equation}
and
\begin{equation}\label{eq-def-Iij-g2}
	I_{i,j,s}^{\pm,g}(x,y)=\mathbf{1}_{\big\{\frac12|y|\le|x|\le 2|y|\big\}}\int_0^\infty e^{i\lambda(|x|\pm|y|)}{T}_{i,j,s}^\pm(\lambda;x,y)\chi(\lambda)d\lambda, \quad s=2,3,
\end{equation}
where the corresponding operators are denoted by  $I_{i,j,1,1}^{\pm,b}$, $I_{i,j,1,t}^{\pm,g}$ and $I_{i,j,s}^{\pm,g}$.

Consequently, to prove the results on $	W_{i,j,3}$, it suffices to establish the  $L^p$-boundedness of the {\bf  good parts} $I_{i,j,1,t}^{\pm,g}$, $I_{i,j,s}^{\pm,g}$, and the {\bf  bad parts} $I_{i,j,1,1}^{\pm,b}$ for all $1<p<\infty$, which are presented in Section \ref{subsec:3.2}.

\subsection {Some criteria to \texorpdfstring{$L^p$}{Lp}-boundedness } \label{sec-pre-for-low}
 Before proving Theorem \ref{thm-low energy-regular-even},  we list several basic criteria to establish $L^p$-boundedness of  operators which also can be found in \cite{CSWY24}.  

Let $T(x,y)$ be the integral kernel of the operator $T$. And we say that the integral kernel $T(x,y)$ is admissible if and only if $T(x,y)$ satisfies
\begin{equation}\label{eq-admissible condition}
\sup_{x\in\mathbb{R}^n}\int |T(x,y)|dy+\sup_{y\in\mathbb{R}^n}\int|T(x,y)|dx<\infty.
\end{equation}
Let us start with the classical Schur's lemma.
\begin{lemma}\label{Schur Lemma}
If $T(x,y)$ is admissible, then $T$ is bounded on $L^p(\mathbb{R}^n)$ for all $1\le p\le\infty$.
\end{lemma}

We also need the following lemmas (see \cite{CSWY24}). 

\begin{lemma}\label{lemma-bounedness-1}
Assume that the integral kernel $T(x,y)$ obeys the following pointwise bounds:
\[|T(x,y)|\lesssim \langle x\rangle^{-n}{\bf1}_{\{|x|>2|y|\}}.\] 
Then the operator $T$ is bounded on $L^p(\mathbb{R}^n)$ for all $1< p\le \infty.$
\end{lemma}

\begin{lemma}\label{lemma-bounedness-11}
Assume that the integral kernel $T(x,y)$ obeys the following pointwise bounds:
\[|T(x,y)|\lesssim \langle x\rangle^{-\alpha}\langle y\rangle^{\alpha-n}{\bf1}_{\{|x|<\frac 12|y|\}}\] 
for some $\alpha\in[0,n)$. Then the operator $T$ is bounded on $L^p(\mathbb{R}^n)$ for all $1\le  p< \frac{n}{\alpha}$.
\end{lemma}

\begin{lemma}\label{lemma-bounedness-2}
Assume that the integral kernel $T(x,y)$ obeys the following pointwise bounds:
\[|T(x,y)|\le C_\varepsilon \langle x\rangle^{-(n-\varepsilon)}\langle |x|-|y|\rangle^{-\varepsilon}{\bf1}_{\{\frac12|y|\le|x|\le 2|y|\}}\] 
for any $\varepsilon\in(0,1).$ Then the operator $T$ is bounded on $L^p(\mathbb{R}^n)$ for all $1< p< \infty.$
\end{lemma}

Besides, we also introduce a result on oscillatory integral estimates in the low energy. 
\begin{lemma}\label{oscillatory estimates}
Let $b>-1$ and $k>b+1$, and assume that $f(\lambda)\in C^k(\mathbb{R})$ such that
\begin{equation*}
\Big|\frac{d^l}{d\lambda^l} f(\lambda)|\Big|\lesssim \lambda^{b-l}, \quad l=0,1,\ldots,k.
\end{equation*}
Then,we have  
\begin{equation}\label{eq.oscillatory est}
\left| \int_0^\infty e^{i\lambda x}f(\lambda)\chi(\lambda)d\lambda\right| \lesssim \langle x\rangle^{-(b+1)},\quad x\in\mathbb{R},
\end{equation}
where $\chi(\lambda)$ is a smooth function such that $\chi\equiv1$ on $(-\lambda_0/2,\,\lambda_0/2)$ and $\supp\chi\subset[-\lambda_0,\,\lambda_0]$ for some fixed $\lambda_0>0$.
\end{lemma}

%In what follows, we will first prove the results in Theorem \ref{thm-low energy-regular-even} for  even dimensions $n\le 2m$.  The proof for even dimensions $2m< n\le 4m$ follows  the same method and will be presented later.  
%\subsection{The $L^p$-boundedness of  $W_{i,j,1}^L$}\label{subsec:3.1}\mbox{}\vspace{8pt}\\
%\indent In this subsection, we will analyze the  $L^p$-boundedness of the operator  $W_{i,j,1}^L$. 

We next establish the  $L^p$-boundedness of $W_{i,j,1}^L$, $W_{i,j,2}^L$ and $W_{i,j,3}^L$ in the following thee subsections, respectively. 
\subsection{The $L^p$-boundedness of $W_{i,j,1}^L$}\label{sec-3.3-wij1}

\begin{proposition}\label{propo-low energy-even-1}
Let $1\le n\le 4m$ be even and $0\le \mathbf{k}\le m_n+1$. Then, for any $i,j\in J_\mathbf{k}$, the operator $W_{i,j,1}^{L}$ is bounded on $L^p(\mathbb{R}^n)$ for all $1<p\le \infty$.
\end{proposition}

By noting \eqref{eq.w-odd-1}, we first study the behavior of $T_{i,j}^\pm(\lambda;x,y)$.

\begin{lemma}\label{lemma-T-even-1}
Let $1\le n\le 4m$ be even and $0\le \mathbf{k}\le m_n+1$. Then
for any $i,j\in J_\mathbf{k}$, \,$\lambda\in(0,\lambda_0)$ and $l=0,\ldots, \frac{n}{2}+1$, where $\lambda_0$ given in Theorem \ref{thm-M inverse-even}, we have that
\begin{equation}\label{eq-estimates- for-T-1}
\left| \partial^l_\lambda T_{i,j}^\pm(\lambda;x,y)\right| \lesssim \lambda^{\frac{n-1}{2}-l}\langle x\rangle^{-\frac{n-1}{2}}, \quad \text{if $i< 2m-\frac n2$}.
\end{equation}
If $i= 2m-\frac n2$ or $2m-\frac n2 +\frac 12$, then we have that 
$$T_{i,j}^\pm(\lambda;x,y)=T_{i,j,1}^\pm(\lambda;x,y)+T_{i,j,2}^\pm(\lambda;x,y)+T_{i,j,3,}^\pm(\lambda;x,y),$$ 
where these functions satisfy that
\begin{align}\label{eq-estimates-for-T-2} 
\left| \partial^l_\lambda T_{i,j,1}^\pm(\lambda;x,y)+T_{i,j,3}^\pm(\lambda;x,y)\right| \lesssim \lambda^{\frac{n-1}{2}-l}\langle x\rangle^{-\frac{n-1}{2}},  \\
\label{eq-estimates-for-T-3}
\left| \partial^l_\lambda T_{i,j,2}^\pm(\lambda;x,y)\right| \lesssim \lambda^{\frac{n}{2}-1-l}\langle x\rangle^{-\frac{n}{2}}. 	
\end{align}
\end{lemma}
  
\begin{proof}
We first prove the case $i<2m-n$.
Recall that 
\[T_{i,j}^{\pm}(\lambda;x,y)=\lambda^{4m-n-1}\widetilde{\rho}_i(\lambda)\widetilde{\rho}_j(\lambda)\Big\langle M_{i,j}^+\widetilde{\omega}_{j}^{\pm}(\lambda,\cdot,y),\, \omega_{i}^{-}(\lambda,\cdot,x)\Big\rangle,\] 
then by Leibniz's rule and Schwartz inequality,  for any $l=0,1,\cdots \frac{n}{2}+1$ and $\lambda\in(0,\lambda_0)$, we have
\begin{equation}\label{eq-leibnitz-expression- for-T}
	\begin{aligned}
		\big|\partial^l_\lambda T_{i,j}^\pm(\lambda;x,y)\big| 	\lesssim& \sum_{l_1+l_2+l_3+l_4=l}\lambda^{4m-n-1-l_1}|\widetilde{\rho}_i(\lambda)||\widetilde{\rho}_j(\lambda)|\\
	 &\times\left\|\frac{d^{l_2}}{d\lambda^{l_2}}\big(M_{i,j}^++\Gamma_{i,j}^{+,r}(\lambda)\big)\right\|_{\mathbb{B}(L^2)} \cdot \left\|\partial_\lambda^{l_3}\widetilde{\omega}_{j}^{\pm}(\lambda,\cdot,y) \right\|_{L^2} \left\|\partial_\lambda^{l_4}\omega_{i}^{-}(\lambda,\cdot,x) \right\|_{L^2},
	\end{aligned}
\end{equation}	
where we have used the expressions of $\widetilde{\rho}_i(\lambda)$ and  $\widetilde{\rho}_j(\lambda)$  in \eqref{eq-widetilde-rho-j}.

By \eqref{eq-est-gamma-ij-1} , we have 
$$\left\|\frac{d^{l_2}}{d\lambda^{l_2}}\big(M_{i,j}^++\Gamma_{i,j}^{+,r}(\lambda)\big)\right\|_{\mathbb{B}(L^2)}\lesssim \lambda^{-l_2}.$$ 
When $i< 2m-\frac n2$,   it follows from the estimates in \eqref{eq-wj-estimate-even-3}, \eqref{eq-wj-estimate-even-1-2} and \eqref{eq-k2-even-1} that
\begin{align*}
&\left\|\partial_\lambda^{l_3}\widetilde{\omega}_{j}^{\pm}(\lambda,\cdot,y)\right\|_{L^2}\lesssim\lambda^{n-2m+\vartheta(j)-l_3},\\ 
&\left\|\partial_\lambda^{l_4}\omega_{i}^{-}(\lambda,\cdot,x) \right\|_{L^2}\lesssim \lambda^{\min\{\frac{n+1}{2}-2m+\vartheta(i),\, 0\}-l_4}\langle x\rangle^{-\frac{n-1}{2}} .
\end{align*}
Combining these with the facts 
$$|\widetilde{\rho}_i(\lambda)|\le \lambda^{-\vartheta(i)} ,  \quad i\neq m-\frac n2 +\frac 34$$
and
$$|\widetilde{\rho}_i(\lambda)|\le \lambda^{-(m-\frac n2+\frac 12)} ,  \quad i= m-\frac n2 +\frac 34,  $$
we conclude that
\begin{equation*}
\big|\partial^l_\lambda T_{i,j}^\pm(\lambda;x,y)\big|\lesssim \lambda^{\frac{n-1}{2}+\vartheta(i)+\vartheta(j)-i-j-l}\langle x\rangle^{-\frac{n-1}{2}}	\lesssim \lambda^{\frac{n-1}{2}-l}\langle x\rangle^{-\frac{n-1}{2}}, \quad i\neq m-\frac n2 +\frac 34,
\end{equation*} 
and 
\begin{equation*}
	\big|\partial^l_\lambda T_{i,j}^\pm(\lambda;x,y)\big|\lesssim \lambda^{\frac {n}{2}+\vartheta(j)-j}\langle x\rangle^{-\frac{n-1}{2}}	\lesssim \lambda^{\frac{n-1}{2}-l}\langle x\rangle^{-\frac{n-1}{2}},\quad i= m-\frac n2 +\frac 34.
\end{equation*} 
In the above, we also used $\vartheta(j)\ge j$ by the definition.

In the case of $i=2m-\frac n2,\,2m-\frac n2+\frac12$, we use \eqref{eq3.1000} to write $T_{i,j}^\pm$ as
\begin{equation}\label{eq-split-for-T}
\begin{aligned}
	T_{i,j}^\pm(\lambda;x,y)&=\sum_{s=1}^3T_{i,j,s}^\pm(\lambda;x,y) \\
	&=\sum_{s=1}^3\lambda^{4m-n-1}\widetilde{\rho}_i(\lambda)\widetilde{\rho}_j(\lambda)\Big\langle \big(M_{i,j}^++\Gamma_{i,j}^{+,r}(\lambda)\big)\widetilde{\omega}_{j}^{\pm}(\lambda,\cdot,y),\, \omega_{i,s}^{-}(\lambda,\cdot,x)\Big\rangle.
\end{aligned}
\end{equation}
Combining  with the estimates \eqref{eq3.100001}--\eqref{eq3.10011}, and following the same method, one can obtain \eqref{eq-estimates-for-T-2}--\eqref{eq-estimates-for-T-3}.
\end{proof}

Now,  we return to prove Proposition \ref{propo-low energy-even-1}.
\begin{proof}[\bf{Proof of Proposition \ref{propo-low energy-even-1}.}]
Recall that \eqref{eq.w-odd-1} states 
\begin{equation*}
	W_{i,j,1}^L(x,y)
	=\mathbf{1}_{\{|x|>2|y|\}}\int_0^\infty \Big(e^{i\lambda(|x|+|y|)}T_{i,j}^+(\lambda; x,y)+e^{i\lambda(|x|-|y|)}T_{i,j}^-(\lambda;x,y)\Big)\chi(\lambda)d\lambda.
\end{equation*}
It is an oscillatory integral with a same form to that in \eqref{eq.oscillatory est} in Lemma \ref{oscillatory estimates}.
Thus, the estimates \eqref{eq-estimates- for-T-1}--\eqref{eq-estimates-for-T-3} and Lemma \ref{oscillatory estimates} indicate that for
all $i,j\in J_{\k}$, we have
 $$ \big|W_{i,j,1}^{L}(x,y)\big| \lesssim  \langle x \rangle^{-n}\mathbf{1}_{\{|x|>2|y|\}}. $$
According to Lemma \ref{lemma-bounedness-1}, $W_{i,j,1}^{L}$ is bounded on $L^p(\mathbb{R}^n)$ for all $ 1 < p \le  \infty $.
\end{proof}

\subsection{The $L^p$-boundedness of $W_{i,j,2}^L$}\label{sec-3.4-wij2}
%The $L^p$-boundedness of $W_{i,j,2}^{L}$ is stated as follows.
\begin{proposition}\label{propo-low energy-even-2}
Let $1\le n\le 4m$ be even and $0\le \mathbf{k}\le m_n+1$. Then, for any $i,j\in J_\mathbf{k}$, we have the following statements:
	\vskip 0.2cm
\indent\emph{(\romannumeral1)} If ~~$i< \max\{2m-n+1,\, 0\}$, then $W_{i,j,2}^{L}\in \mathbb{B}(L^p(\mathbb{R}^n))$ for all  $1\le p<\infty$ ;
		\vskip 0.2cm
\indent\emph{(\romannumeral2)} If~~$\max\{2m-n+1,\, 0\}\le i<2m-\frac n2$, then $W_{i,j,2}^{L}\in \mathbb{B}(L^p(\mathbb{R}^n))$ for all $1\le p<\frac{n}{n-2m+\vartheta(i)}$;
		\vskip 0.2cm
\indent\emph{(\romannumeral3)} If~~$i=2m-\frac n2$ or $2m-\frac n2+\frac 12$, then $W_{i,j,2}^{L}\in \mathbb{B}(L^p(\mathbb{R}^n))$ for all $1\le p< 2$.
	\vskip 0.2cm
 Here $\vartheta(i)=\max\{0,\lfloor i+1/2 \rfloor\}$ is given in \eqref{thetaj}.
\end{proposition}

\begin{remark}\label{rmk-on-p}
We emphasize that the restriction on the range of the exponent $p$ in Theorems~\ref{thm-main result-0}, \ref{thm-main result-1}, and \ref{thm-low energy-regular-even} is entirely determined by this proposition.
\end{remark}

We divide the proof into two cases: $1 \leq n \leq 2m$ and $2m < n \leq 4m$. To proceed, we establish two additional estimates for the integral kernel $T_{i,j}^\pm(\lambda,x,y)$ to prove Proposition~\ref{propo-low energy-even-2}.
\begin{lemma}\label{lemma-T-even-2} 
Let $1\le n\le 2m$ be even and $0\le \mathbf{k}\le m_n+1=m+2$.  Then
for any $i,j\in J_\mathbf{k}$, \,$\lambda\in(0,\lambda_0)$ and $l=0,\ldots, \frac{n}{2}+1$, the following holds.

\emph{(\romannumeral1)} 
 If $i=2m-n$, one has 
\begin{equation} \label{eq-estimates-for-T-2-1-0}
	\big|\partial^l_\lambda T_{i,j}^\pm(\lambda,x,y)\big|\lesssim_\varepsilon \lambda^{\frac{n-1}{2}-\varepsilon-l}\langle y\rangle^{-\frac{n-1}{2}}\langle x\rangle^{-\varepsilon}, \quad \text{for each $\varepsilon\in (0,1/2)$};  
\end{equation}

\emph{(\romannumeral2)} If $n=2m$ and $i= \frac 34$, one has 
\begin{equation} \label{eq-estimates-for-T-2-1-1}
	\big|\partial^l_\lambda T_{i,j}^\pm(\lambda,x,y)\big|\lesssim_\varepsilon \lambda^{\frac{n-1}{2}-\varepsilon-l}\langle y\rangle^{-\frac{n-1}{2}}\langle x\rangle^{-\frac 12}, \quad \text{for each $\varepsilon\in (0,1/2)$};  
\end{equation}

\emph{(\romannumeral3)} If indices $i < 2m - n + 1$ with the exclusion $i= 2m - n$ and $i=\frac{3}{4}$, one has
\begin{equation} \label{eq-estimates-for-T-2-1}
	\big|\partial^l_\lambda T_{i,j}^\pm(\lambda,x,y)\big|\lesssim\lambda^{\frac{n-1}{2}-l}\langle y\rangle^{-\frac{n-1}{2}};  
\end{equation}

\emph{(\romannumeral4)} If $2m-n+1\le i< 2m-\frac n2$, one has 
\begin{equation} \label{eq-estimates-for-T-2-2}
\big|\partial^l_\lambda T_{i,j}^\pm(\lambda,x,y)\big|\lesssim\lambda^{2m-\frac{n+1}{2}-i-l}
\langle x \rangle^{2m-n-\vartheta(i)}\langle y\rangle^{-\frac{n-1}{2}}.  
\end{equation}
Here $\vartheta(i)=\max\{0,\lfloor i+1/2 \rfloor\}$ is given in \eqref{thetaj}.
\end{lemma}
\begin{proof}	
 Using Leibniz's rule  and  the definition of $T_{i,j}^\pm(\lambda;x,y)$, the proof follows from the similar process for the proof of Lemma \ref{lemma-T-even-1}. According \eqref{eq-leibnitz-expression- for-T}, the key is to give the estimate for $\left\|\partial_\lambda^{l_2}\widetilde{\omega}_{j}^{\pm}(\lambda,\cdot,y) \right\|_{L^2}$ and $\left\|\partial_\lambda^{l_3}\omega_{i}^{-}(\lambda,\cdot,x) \right\|_{L^2}$.

\textbf{(i).} The estimate \eqref{eq-estimates-for-T-2-1-0} follows from the following three facts.
\begin{itemize}
\item the $L^2$ bound from \eqref{eq-wj-estimate-even-2} that
$
\|\partial_\lambda^l \omega_i^\pm\|_{L^2} \lesssim_\varepsilon \lambda^{-\varepsilon-l}\langle x\rangle^{-\varepsilon};
$
\vskip0.1cm
\item the $L^2$ bound from \eqref{eq-wj-estimate-even-1-2} that
\begin{equation}\label{eq-est-for-wj-y}
\|\partial_\lambda^l \tilde{\omega}_j^\pm\|_{L^2} \lesssim \lambda^{\frac{n+1}{2}-2m+\vartheta(j)-l}\langle y\rangle^{-\frac{n-1}{2}}; 
\end{equation}

\item the behaviors from \eqref{eq-widetilde-rho-j} that
\[
\widetilde{\rho}_i(\lambda)=\lambda^{2m-n}, \quad |\widetilde{\rho}_j(\lambda)| \leq \lambda^{-\vartheta(i)}.\]
\end{itemize}

\textbf{(ii).} The estimate \eqref{eq-estimates-for-T-2-1-1} follows from \eqref{eq-est-for-wj-y} and
\begin{itemize}
\item the facts
\[
\widetilde{\rho}_i(\lambda)=(\log\lambda)^{1/2}, \quad \vartheta(i)=1,\quad \text{for}~~ i=m-\frac n2+\frac 34;
\]

\item the decay from \eqref{eq-wj-estimate-even-4} that
\[
\|\partial_\lambda^l \omega_i^\pm\|_{L^2} \lesssim \lambda^{-l}\langle x\rangle^{-1}
\]
\end{itemize}

Similarly, \eqref{eq-estimates-for-T-2-1} is immediately obtained by \eqref{eq-wj-estimate-even-1},  \eqref{eq-wj-estimate-even-1-2} and \eqref{eq-est-for-wj-y}, meanwhile, \eqref{eq-estimates-for-T-2-2} is derived from \eqref{eq-wj-estimate-even-4},  \eqref{eq-wj-estimate-even-1-2} and \eqref{eq-est-for-wj-y}.
\end{proof}

By applying parallel analysis, we obtain a similar result when $2m< n\le 4m$.
\begin{lemma}\label{lemma-T-even-2-1} 
	Let $2m< n\le 4m$ be even and $0\le \mathbf{k}\le m_n+1$. Assume that $V$ satisfies \eqref{decayonV} in Assumption \ref{Assumption}. Then
	for any $i,j\in J_\mathbf{k}$,  $\lambda\in(0,\lambda_0)$ and $l=0,\ldots, \frac{n}{2}+1$, the following holds: 
    
\emph{(\romannumeral1)} When $i=m-\frac n2$, and either $j\neq 0$ or $n<4m$, then 
\begin{equation} \label{eq-estimates-for-T-2-2-1}
\big|\partial^l_\lambda T_{i,j}^\pm(\lambda,x,y)\big|\lesssim \lambda^{\frac{n-1}{2}-l}\langle y\rangle^{-\frac{n-1}{2}};    
\end{equation}
 When $i=m-\frac n2$, $j=2m-\frac n2=0$ for $n=4m$, then 
\begin{equation} \label{eq-estimates-for-T-2-2-2}
	\big|\partial^l_\lambda T_{i,j}^\pm(\lambda,x,y)\big|\lesssim \lambda^{\frac{n-1}{2}-l}|\log\lambda|^{\frac 12}\langle y\rangle^{-\frac{n-1}{2}};    
\end{equation}

\emph{(\romannumeral2)} When $0\le i<2m-\frac n2$,  then 
\begin{equation} \label{eq-estimates-for-T-2-2-3}
	\big|\partial^l_\lambda T_{i,j}^\pm(\lambda,x,y)\big|\lesssim\lambda^{2m-\frac{n+1}{2}-i-l}
	\langle x \rangle^{2m-n-\vartheta(i)}\langle y\rangle^{-\frac{n-1}{2}}.  
\end{equation}
\end{lemma}
\begin{proof}
Note that $M_{i,j}^+=0$ when $i=m-\frac n2, j\neq m-\frac n2\ \text{or}\ \ j=m-\frac n2, i\neq m-\frac n2$. Also,  observe the improved estimates for $\Gamma_{i,j}^{\pm,r}(\lambda)$ in \eqref{eq-est-gamma-ij-2-1} and \eqref{eq-est-gamma-ij-2-2}. 
The estimate \eqref{eq-estimates-for-T-2-2-1}--\eqref{eq-estimates-for-T-2-2-3} is immediately derived from \eqref{eq-wj-estimate-even-1}, \eqref{eq-wj-estimate-even-3}, \eqref{eq-wj-estimate-even-4} and the inequality \eqref{eq-leibnitz-expression- for-T}.
\end{proof}

We next prove Proposition \ref{propo-low energy-even-2}.
\begin{proof}[\bf{Proof of Proposition \ref{propo-low energy-even-2}.}]
Recall that
\[W_{i,j,2}^L(x,y)
	=\mathbf{1}_{\{|x|<\frac 12|y|\}}\int_0^\infty \Big(e^{i\lambda(|x|+|y|)}T_{i,j}^+(\lambda,x,y)+e^{i\lambda(|x|-|y|)}T_{i,j}^-(\lambda,x,y)\Big)\chi(\lambda)d\lambda.\]
Using Lemma \ref{oscillatory estimates}, we prove the results for even dimensions in two cases: $ n \le 2m$ and $ 2m < n \le 4m$.

We first consider the case $n\le 2m$. 

When $i,\, j \in J_\mathbf{k}$ and $i=2m-n$ 
$$
 \eqref{eq-estimates-for-T-2-1-0}~\implies~ |W_{i,j,2}^{L}(x,y)|  \lesssim \langle y \rangle^{-n-\varepsilon}\langle x\rangle^{-\varepsilon}\mathbf{1}_{\{|x|<\frac 12|y|\}}, \quad \text{for each $\varepsilon\in (0,1/2)$}.
$$
When  $i= \frac 34$ and $n=2m$,  
$$\eqref{eq-estimates-for-T-2-1-1}~\implies~|W_{i,j,2}^{L}(x,y)|  \lesssim \langle y \rangle^{-n-\varepsilon}\langle x\rangle^{-\varepsilon}\mathbf{1}_{\{|x|<\frac 12|y|\}}, \quad \text{for each $\varepsilon\in (0,1/2)$}.$$
 When $i<2m-n+1$,  $i\neq \frac 34$ and $i\neq 2m-n$, 
$$
\eqref{eq-estimates-for-T-2-1}~\implies~|W_{i,j,2}^{L}(x,y)|  \lesssim \langle y \rangle^{-n}\mathbf{1}_{\{|x|<\frac 12|y|\}}.
$$
The above three estimates, together with  Lemma \ref{lemma-bounedness-11}, show that $W_{i,j,2}^{L}\in \mathbb{B}(L^p)$ for all $1\le p<\infty$ if $i\le 2m-n+1$. 

When $2m-n+1\le i< 2m-\frac n2$
$$
\eqref{eq-estimates-for-T-2-2}~\implies~|W_{i,j,2}^{L}(x,y)|  \lesssim \langle x \rangle^{-(n-2m+\vartheta(i))} \langle y \rangle^{-(2m-i)}\mathbf{1}_{\{|x|<\frac 12|y|\}},
$$
which implies that $W_{i,j,2}^{L}\in \mathbb{B}(L^p)$ for all  $1\le p<\frac{n}{n-2m+\vartheta(i)}$.

 When $i=2m-\frac n2, \ 2m-\frac n2 +\frac 12$, decompose $ T_{i,j}^\pm(\lambda,x,y)$ as in Lemma \ref{lemma-T-even-1}. Then
$$\eqref{eq-estimates-for-T-2},~\eqref{eq-estimates-for-T-3}~\implies~ |W_{i,j,2}^{L}(x,y)|  \lesssim \langle x \rangle^{-\frac{n}{2}} \langle y \rangle^{-\frac{n}{2}}\mathbf{1}_{\{|x|<\frac 12|y|\}},$$
which implies that $W_{i,j,2}^{L}\in \mathbb{B}(L^p)$ for all  $1\le p<2$.

We next consider the case $2m<n\le 4m$.  

When $i=m-\frac n2<0$,
$$
\eqref{eq-estimates-for-T-2-2-1}~ \text{and}~ \eqref{eq-estimates-for-T-2-2-2}~\implies~|W_{i,j,2}^{L}(x,y)|  \lesssim
 \begin{cases} 
 	\langle y \rangle^{-n}\mathbf{1}_{\{|x|<\frac 12|y|\}}  & \text{if $j\neq 0$ or $n<4m$,}\\[3mm]
    \langle x\rangle^{-2m}\langle y \rangle^{-n+\varepsilon}\mathbf{1}_{\{|x|<\frac 12|y|\}} & \text{if $j=0$ and $n=4m$ },\\[3mm]
 \end{cases}
$$
 for each $\varepsilon \in (0,1)$.
When $ 0\le i< 2m-\frac n2$, 
$$\eqref{eq-estimates-for-T-2-2-3} ~\implies~ |W_{i,j,2}^{L}(x,y)|\lesssim 
\langle x \rangle^{-(n-2m+\vartheta(i))} \langle y \rangle^{-(2m-i)}\mathbf{1}_{\{|x|<\frac 12|y|\}}.$$
We conclude that $W_{i,j,2}^{L}\in \mathbb{B}(L^p)$ for all $p$ satisfying 
$$ 1\le p<\infty ~ \text{if $i=m-\frac n2$ and} ~ 1\le p<\frac{n}{n-2m+\vartheta(i)}  ~\text{if $i\neq m-\frac n2$}.$$ 

Meanwhile, when $i=2m-\frac n2$ or $2m-\frac n2 +\frac 12$, decompose $ T_{i,j}^\pm(\lambda,x,y)$ as in Lemma \ref{lemma-T-even-1}. Then it follows that
$$\eqref{eq-estimates-for-T-2},~\eqref{eq-estimates-for-T-3}~\Rightarrow~ |W_{i,j,2}^{L}(x,y)|  \lesssim \langle x \rangle^{-\frac{n}{2}} \langle y \rangle^{-\frac{n}{2}}\mathbf{1}_{\{|x|<\frac 12|y|\}}.$$
In this case, we conclude that $W_{i,j,2}^{L}\in \mathbb{B}(L^p)$ for all  $1\le p<2$.

Summing up, the $L^p$-boundedness of $W_{i,j,2}^{L}$ in this proposition has been established.
\end{proof}

\subsection{The $L^p$-boundedness of $W_{i,j,3}^L$}\label{subsec:3.2} 
% In this subsection, we study the $L^p$-boundedness of $W_{i,j,3}^L$, as stated in the following proposition. 
\begin{proposition}\label{estforw_3}
Let $1\le n\le 4m$ be even and $0\le \mathbf{k}\le m_n+1$. Then
for any $i,j\in J_\mathbf{k}$, $W_{i,j,3}^L$ is a bounded operator on $L^p(\mathbb{R}^n)$ for all $1<p<\infty$.	
\end{proposition}

To prove the results on $W^L_{i,j,3}$, it suffices to establish the $L^p$-boundedness of both the good parts $I_{i,j,1,t}^{\pm,g}$ (for $t=2,3$) and $I_{i,j,s}^{\pm,g}$ (for $s=2,3$), as well as the bad part $I_{i,j,1,1}^{\pm,b}$, for all $1<p<\infty$. These operators are precisely defined in equations \eqref{eq-def-Iijb} through \eqref{eq-def-Iij-g2}.

\subsubsection{The  good parts $I_{i,j,1,t}^{\pm,g}$ and $I_{i,j,s}^{\pm,g}$ }

\begin{proposition}\label{lem.estforwij3g} For any $i,j\in J_\mathbf{k}$ and $s=2,3$, 
$I_{i,j,1,s}^{\pm,g}$ and $I_{i,j,s}^{\pm,g}$ are  bounded operators on $L^p(\mathbb{R}^n)$ for all $1<p<\infty$.
\end{proposition}
We first prove the following lemma.
\begin{lemma}\label{est for T-3}
Let $1\le n\le 4m$ be an even integer and $0\le \mathbf{k}\le m_n+1$. Then,
for any $i,j\in J_\mathbf{k}$, $\lambda\in(0,\lambda_0)$ and $l=0,\ldots, \frac{n}{2}+1$, we have the following estimates:
\vskip0.2cm
\indent\emph{(\romannumeral1)}~ For any $\varepsilon \in [0,1)$, ${T}_{i,j,2}^\pm(\lambda;x,y)$ and  $T_{i,j,1,2}^\pm(\lambda;x,y)$ given in \eqref{eq4.101012} and \eqref{eq.4.100123} satisfy
\begin{equation} \label{eq.4.10001}
	\big|\partial^l_\lambda T_{i,j,1,2}^\pm(\lambda;x,y)\big| +\big|\partial^l_\lambda {T}_{i,j,2}^\pm(\lambda;x,y)\big|\le C
	\begin{cases}
		\lambda^{-l-\varepsilon}\langle x\rangle^{-(n-1)-\varepsilon}, \quad &\text{if $i< 2m- \frac n2$}, \\[3mm]
		\lambda^{-\frac 12-l-\frac{\varepsilon}{2}}\langle x\rangle^{-n+\frac{1}{2}-\frac{\varepsilon}{2}}, &\text{if $i= 2m- \frac n2,\ 2m- \frac n2+\frac 12$},
	\end{cases}
\end{equation}
where the absolute constant $C$ in the above inequality is independent of $\varepsilon$;
\vskip0.2cm
\indent\emph{(\romannumeral2)}~ In addition, ${T}_{i,j,3}^\pm(\lambda;x,y)$ and  $T_{i,j,1,3}^\pm(\lambda;x,y)$ given in \eqref{eq4.101012} and \eqref{eq.4.100123} satisfy
\begin{equation} \label{eq.4.10002}
	\big|\partial^l_\lambda T_{i,j,1,3}^\pm(\lambda;x,y)\big| +\big|\partial^l_\lambda 
	{T}_{i,j,3}^\pm(\lambda;x,y)\big|\lesssim 
	\lambda^{1-l}\langle x\rangle^{-(n-1)}.
\end{equation}
\end{lemma}
\vskip0.3cm
\begin{proof}
Note that by the definition of $\omega_{i,1}^{\pm}(\lambda,\cdot,x)$ in \eqref{eq3.100001}, for any $l=0,\ldots, \frac{n}{2}+1$,  a direct computation shows that
\begin{equation*}
\begin{aligned}
  \left\|\partial^l_\lambda \omega_{i,1}^{\pm}(\lambda,x,\cdot) \right\|_{L^2}
  \lesssim &
  \lambda^{\frac{n+1}{2}-2m+\vartheta(i)-l} \langle x\rangle^{-\frac{n-1}{2}} \left\|\sum_{|\alpha|=|\beta|=\vartheta(j)}A_{\alpha,\beta}^\pm Q_j(v(z)z^\beta)(y)	\right\|_{L^2_y}\\ \lesssim & \lambda^{\frac{n+1}{2}-2m+\vartheta(i)-l} \langle x\rangle^{-\frac{n-1}{2}},  
\end{aligned}
\end{equation*}
which is consistent with the form of the estimate in equation \eqref{eq-k2-even-1}. Applying Leibniz's rule like \eqref{eq-leibnitz-expression- for-T}, the above estimate and \eqref{eq-k2-even-2} and \eqref{eq3.1001} yield \eqref{eq.4.10001}. Similarly, \eqref{eq-k2-even-2} combined with \eqref{eq3.10011} shows \eqref{eq.4.10002}.
\end{proof}

%We now prove  Proposition \ref{lem.estforwij3g}.
\begin{proof}[\bf{Proof of Proposition \ref{lem.estforwij3g}.}]
We first consider the operator  $I_{i,j,1,s}^{\pm,g}$ and $I_{i,j,s}^{\pm,g}$ for $s=2,\,3$.
It follows from Lemma \ref{oscillatory estimates} and \eqref{eq.4.10001} that for each $\varepsilon\in [0,1)$, in the region where $\frac 12|y|\le |x| \le 2|y|$, we have 
\begin{equation*}
|I_{i,j,2}^{\pm,g}(x,y)|+|I_{i,j,1,2}^{\pm,g}(x,y)|\lesssim \langle|x|-|y|\rangle^{-1+\varepsilon}\langle x\rangle^{-(n-1+\varepsilon)}, \quad i< 2m- \frac n2.
\end{equation*}
and when $i= 2m- \frac n2,\ 2m- \frac n2+\frac 12$,
\begin{equation*}
\begin{aligned}
&|I_{i,j,2}^{\pm,g}(x,y)|+|I_{i,j,1,2}^{\pm,g}(x,y)|\\
\lesssim& \langle|x|-|y|\rangle^{-1/2+\varepsilon/2}\langle x\rangle^{-n+\frac 12-\varepsilon/2}  \lesssim  \langle|x|-|y|\rangle^{-1+\varepsilon}\langle x\rangle^{-(n-1+\varepsilon)}, 	
\end{aligned}
\end{equation*}
where we have used the fact $\langle|x|-|y|\rangle\lesssim \langle x\rangle$ when  $\frac 12|y|\le |x| \le 2|y|$.
Thus, Lemma \ref{lemma-bounedness-2} shows that for all $i,j\in J_{\k}$, the operators $I_{i,j,2}^{\pm,g}$ and  $I_{i,j,1,2}^{\pm,g}$ are bounded on $L^p(\mathbb{R}^n)$ for $1<p<\infty$.

 Similarly, according to Lemma \ref{oscillatory estimates} and \eqref{eq.4.10001}, we have that $I_{i,j,1,3}^{\pm,g}(x,y)$ are supported in the region where $\frac{1}{2}|y|\le |x|\le2|y|$ and satisfy 
\begin{equation*}
|I_{i,j,3}^{\pm,g}(x,y)|+|I_{i,j,1,3}^{\pm,g}(x,y)|\lesssim
   \langle|x|-|y|\rangle^{-2}\langle x\rangle^{-(n-1)}, 
\end{equation*}
which implies that it is admissible. Thus,  the operators $I_{i,j,3}^{\pm,g}$ and $I_{i,j,1,3}^{\pm,g}$  are bounded on $L^p(\mathbb{R}^n)$ for all $1 \le p \le \infty$.
\end{proof}

\subsubsection{The  bad parts $I_{i,j,1,1}^{\pm,b}$}

\begin{proposition}\label{lem.estforwij3b1} Under Assumption \ref{Assumption}, for any $i,j\in J_\mathbf{k}$, 
$I_{i,j,1,1}^{\pm,b}$ are  bounded operators on $L^p(\mathbb{R}^n)$ for all $1<p<\infty$.
\end{proposition}

Before the proof of  Proposition \ref{lem.estforwij3b1}, let us make some reductions for the bad parts $I_{i,j,1,1}^{\pm,b}$.  
By the definition of $T_{i,j,1,s}^{\pm}(\lambda;x,y)$, we decompose it as 
\begin{equation*}
T_{i,j,1,s}^{\pm}(\lambda;x,y)=T_{i,j,1,1}^{\pm,M}(\lambda;x,y)+T_{i,j,1,1}^{\pm,\varGamma}(\lambda;x,y),
\end{equation*}
where
$$T_{i,j,1,1}^{\pm,M}(\lambda;x,y)=\lambda^{4m-n-1}\widetilde{\rho}_i(\lambda)\widetilde{\rho}_j(\lambda)\Big\langle M_{i,j}^+\omega_{j,1}^{\pm}(\lambda,y,\cdot),\, \omega_{i,1}^{-}(\lambda, x,\cdot)\Big\rangle,$$
and 
$$T_{i,j,1,1}^{\pm,\Gamma}(\lambda;x,y)=\lambda^{4m-n-1}\widetilde{\rho}_i(\lambda)\widetilde{\rho}_j(\lambda)\Big\langle \Gamma_{i,j}^{+,r}(\lambda)\omega_{j,1}^{\pm}(\lambda,y,\cdot),\, \omega_{i,1}^{-}(\lambda, x,\cdot)\Big\rangle.$$

We denote the corresponding operators as
\begin{align*}
I_{i,j,1,1}^{\pm,b,M}(x,y)&=\mathbf{1}_{\big\{\frac12|y|\le|x|\le 2|y|\big\}}\int_0^\infty e^{i\lambda(|x|\pm|y|)}T_{i,j,1,1}^{\pm,M}(\lambda;x,y)\chi(\lambda)d\lambda, \\
I_{i,j,1,1}^{\pm,b,\Gamma}(x,y)&=\mathbf{1}_{\big\{\frac12|y|\le|x|\le 2|y|\big\}}\int_0^\infty e^{i\lambda(|x|\pm|y|)}T_{i,j,1,1}^{\pm,\Gamma}(\lambda;x,y)\chi(\lambda)d\lambda.	
\end{align*}

\begin{proposition}\label{pro.est-for-IGamma} Under Assumption \ref{Assumption}, for any $i,j\in J_\mathbf{k}$, 
	the operator $I_{i,j,1,1}^{\pm,\Gamma}$ is  bounded on $L^p(\mathbb{R}^n)$ for all $1\le p\le\infty$.
\end{proposition}
\begin{proof}
	Note that, according to the definition of $\widetilde{\rho}_j(\lambda)$ in \eqref{eq-widetilde-rho-j} and the expression of $\omega_{j,1}^{\pm}(\lambda, x,y)$ in \eqref{eq3.100001}, the function $T_{i,j,1,1}^{\pm,\Gamma}(\lambda;x,y)$ satisfies, for all $i,j\in J_{\k}$, that
	\begin{equation*}
	\big|\partial^l_\lambda T_{i,j,1,1}^{\pm,\Gamma}(\lambda;x,y)\big| \lesssim 
	\lambda^{-l}|\log\lambda|^{-\frac 52}\langle x\rangle^{-\frac{n-1}{2}}\langle x\rangle^{-\frac{n-1}{2}}, \quad l=0,\,1,\,\ldots,\frac n2+1.	
	\end{equation*}

If  $\big||x|-|y|\big|\le 1$,  we immediately have that
$$\left|I_{i,j,1,1}^{\pm,b,\Gamma}(x,y)\right|\lesssim  \langle x\rangle^{-\frac{n-1}{2}} \langle y\rangle^{-\frac{n-1}{2}}\mathbf{1}_{\big\{\frac12|y|\le|x|\le 2|y|\big\}} \lesssim \langle |x|-|y|\rangle^{-2}\langle x\rangle^{-\frac{n-1}{2}} \langle y\rangle^{-\frac{n-1}{2}}\mathbf{1}_{\big\{\frac12|y|\le|x|\le 2|y|\big\}},$$
where the right hand side is admissible. 

If  $\big| |x|-|y|\big| > 1$, we split the integral as a sum of 
\begin{align}\label{eq-part-1-using-Phi}
&\mathbf{1}_{\big\{\frac12|y|\le|x|\le 2|y|\big\}}\int_0^\infty e^{i\lambda(|x|\pm|y|)}T_{i,j,1,1}^{\pm,\Gamma}(\lambda;x,y)\chi(\lambda)\Psi_0\left( \lambda\langle|x|- |y||\rangle \right) d\lambda,\\ \label{eq-part-2-using-Phi}
&\mathbf{1}_{\big\{\frac12|y|\le|x|\le 2|y|\big\}}\int_0^\infty e^{i\lambda(|x|\pm|y|)}T_{i,j,1,1}^{\pm,\Gamma}(\lambda;x,y)\chi(\lambda)\Psi_1\left( \lambda\langle|x|- |y||\rangle\right) d\lambda.
\end{align}
In the above, $\Psi_0$ is smooth function supported in the interval $[-1,1]$ such that $\Psi_0(\lambda)=1$ for $|\lambda|<1/2$, and $\Psi_1=1-\Psi_0$.
Applying integration by parts once to \eqref{eq-part-1-using-Phi}, we obtain
\begin{align*}
\left|\eqref{eq-part-1-using-Phi}\right|\lesssim& \frac{1}{\langle|x|\pm |y||\rangle } \langle x\rangle^{-\frac{n-1}{2}} \langle y\rangle^{-\frac{n-1}{2}}\mathbf{1}_{\big\{\frac12|y|\le|x|\le 2|y|\big\}} \int_0^{\langle|x|- |y||\rangle^{-1}} \lambda^{-1}|\log\lambda|^{-\frac 52} d\lambda \\
\lesssim& \langle|x|\pm |y||\rangle^{-1}\left(\log \langle|x|- |y||\rangle\right)^{-\frac 32 }\min\{ \langle x\rangle^{-(n-1)},\, \langle y\rangle^{-(n-1)} \},
\end{align*}
where the right hand side is admissible.
Applying the integration by parts twice to \eqref{eq-part-1-using-Phi}, then we have 
\begin{align*}
	\left|\eqref{eq-part-2-using-Phi}\right|\lesssim& {\langle|x|\pm |y||\rangle }^{-2} \langle x\rangle^{-\frac{n-1}{2}} \langle y\rangle^{-\frac{n-1}{2}}\mathbf{1}_{\big\{\frac12|y|\le|x|\le 2|y|\big\}} \int_{\langle|x|-|y||\rangle^{-1}}^{-\lambda_0} \lambda^{-2}|\log\lambda|^{-\frac 52} d\lambda \\
	\lesssim& \langle|x|\pm |y||\rangle^{-1} \left(\log \langle|x|- |y||\rangle\right)^{-\frac 32 }\min\{ \langle x\rangle^{-(n-1)},\, \langle y\rangle^{-(n-1)} \}.
\end{align*}
This shows $\eqref{eq-part-2-using-Phi}$ is also admissible. 
In conclusion, the operator $I_{i,j,1,1}^{\pm,\Gamma}$ is admissible.
\end{proof}

Next, we consider the operator $I_{i,j,1,1}^{\pm,b,M}$. 
For simplicity, we denote 
\begin{equation}\label{bad parts}
\Omega_{i,j}^{\pm}(x,y)=\int_0^\infty e^{i\lambda(|x|\pm|y|)}T_{i,j,1,1}^{\pm,b,M}(\lambda,x,y)\chi(\lambda)d\lambda.
\end{equation}
Let $\chi_0={\bf{1}}_{\{x\in\mathbb{R}^n;\,|x|\le 1\}}$ and $\chi_s={\bf{1}}_{\{x\in\mathbb{R}^n;\,2^{s-1}<|x|\le 2^s\}}$ for $s=1,\,2,\cdots$. Then, we have that  
\begin{equation*}
	\psi(x,y)=\mathbf{1}_{\big\{\frac12|y|\le|x|\le 2|y|\big\}}-\sum_{\substack{s,h\in \mathbb{N}_0\\|s-h|<2}}\chi_s(x)\chi_h(y).
\end{equation*}
Thus,  we have the following corresponding division:
\begin{equation}\label{sumofw3b}
I_{i,j,1,1}^{\pm,b}(x,y)=\sum_{\substack{s,h\in \mathbb{N}_0\\|s-h|<2}}\chi_s(x)\Omega_{i,j}^{\pm}(x,y)\chi_h(y)+ \psi(x,y)\Omega_{i,j}^{\pm}(x,y).
\end{equation}
Since both $\frac14|y|\le |x|\le 4|y|$ and $|x|\lesssim||x|-|y||$ hold in the support of $\psi(x,y)$, applying the same analysis as that used for the good parts $I_{i,j,1,s}^{\pm,g}(x,y)$ for $s=2,3$, it is not difficult to conclude that $\psi(x,y)\Omega_{i,j}^{\pm}(x,y)$  are admissible, and are bounded on $L^p(\mathbb{R}^n)$ for all $1 \le p\le \infty$. We omit the details here. 
Therefore, to prove Proposition \ref{lem.estforwij3b1}, it suffices to show that the operators $\sum_{\substack{|s-h|<2}}\chi_s\Omega_{i,j}^{\pm}(x,y)\chi_h$ are bounded on $L^p(\mathbb{R}^n)$ for all $1<p<\infty$.

%Now we present the proof of Proposition \ref{lem.estforwij3b1}.  

\begin{proof}[\bf{Proof of Proposition \ref{lem.estforwij3b1}.}]
We first claim that for any $i,j\in J_\mathbf{k}$, there exists a constant $C_p$ independent of all $s,h\in \mathbb{N}_0$ with $|s-h|<2$ such that the operators \begin{equation}\label{claim}
\|\chi_s\Omega_{i,j}^{\pm}\chi_h\|_{\mathbb{B}(L^p(\mathbb{R}^n))}\le C_p\end{equation}
holds for all $1<p<\infty$. Once this claim is established,  then for $1<p<\infty$,  we have that
\begin{equation}
\begin{aligned}
 \Big\|\sum_{\substack{|s-h|<2}} \chi_s\Omega_{i,j}^{\pm}\chi_h f\Big\|^p_{L^p}&=\int_{\mathbb{R}^n} \Big|\sum_{s\in \mathbb{N}_0}  \sum_{\substack{h:|s-h|<2}}\chi_s\Omega_{i,j}^{\pm}\chi_h f(x) \Big|^p dx\\  
 &=\sum_{s\in \mathbb{N}_0} \int_{\mathbb{R}^n} \Big| \sum_{\substack{h:|s-h|<2}}\chi_s\Omega_{i,j}^{\pm}\chi_h f(x) \Big|^p dx\\  
 &\lesssim \sum_{s\in \mathbb{N}_0} \sum_{\substack{h:|s-h|<2}} \int_{\mathbb{R}^n} \Big| \chi_h f(x) \Big|^p dx 
 \lesssim  \big\| f\big\|^p_{L^p}. 
\end{aligned}
\end{equation}
This, together with Proposition \ref{pro.est-for-IGamma}, shows that $I_{i,j,1,1}^{\pm,b}$ are bounded on $L^p(\mathbb{R}^n))$ for all $1<p<\infty$.

In the following, we turn to prove the claim that $\chi_s\Omega_{i,j}^{\pm}(x,y)\chi_h$ satisfies the estimate \eqref{claim} uniformly for all $s,h\in \mathbb{N}_0$ with $|s-h|<2$. By  the definition of $\omega_{i,1}^{\pm}$ in \eqref{eq3.100001}, it follows that
\begin{equation*}
	T_{i,j,1,1}^{\pm}(\lambda,x,y)=\sum_{|\alpha|=|\beta|=\vartheta(j)}\sum_{|\tilde{\alpha}|=|\tilde{\beta}|=\vartheta(i)}A_{\alpha,\beta}^\pm A_{\tilde{\alpha},\tilde{\beta}}^-\, {T}_{\alpha,\beta,\tilde{\beta},\tilde{\alpha}}(\lambda,x,y) ,
\end{equation*}
where
\begin{equation}\label{Tab0}
\begin{aligned}
{T}_{\alpha,\beta,\tilde{\beta},\tilde{\alpha}}(\lambda,x,y)=C(M_{i,j}^+,\alpha, \beta) \,\lambda^{\vartheta(i)+\vartheta(j)}&\widetilde{\rho}_i(\lambda)\widetilde{\rho}_j(\lambda) \\ 
&\times \mathbf{1}_{\{|x|>1,\,|y|>1\}} \cdot \frac{y^\alpha}{|y|^{\frac{n-1}{2}+\vartheta(j)}}	\frac{x^{\tilde{\alpha}}}{|x|^{\frac{n-1}{2}+\vartheta(i)}},
\end{aligned}
\end{equation}
and $C(M_{i,j}^+,\alpha, \beta)=\Big\langle M_{i,j}^+Q_jv(z^\beta),\, Q_iv(z^{\tilde{\beta}}) \Big \rangle$.
 Thus, $\chi_s(x)\Omega_{i,j}^{\pm}(x,y)\chi_h(y)$ can be written as a linear combination of finite terms as follows: 
\begin{equation}\label{omegaq}
\Omega_{s,h,\alpha,\beta,\tilde{\beta},\tilde{\alpha}}^\pm(x,\, y)=\chi_s(x)\chi_h(y)\int_0^\infty e^{i\lambda(|x|\pm|y|)}T_{\alpha,\beta,\tilde{\beta},\tilde{\alpha}}(\lambda,x,y)\chi(\lambda)d\lambda,
\end{equation}
where $|\alpha|=|\beta|=\vartheta(j)$, $|\tilde{\alpha}|=|\tilde{\beta}|=\vartheta(i)$. Note that by the definitions $\vartheta(i)$ given in \eqref{thetaj}, and $\widetilde{\rho}_i(\lambda)$ given in \eqref{eq-widetilde-rho-j}, when $i,j \in J_\mathbf{k}$, one of the following conditions must holds: \vspace{0.5em}  

 (i) \ \ \ \   $\lambda^{\vartheta(i)+\vartheta(j)}\widetilde{\rho}_i(\lambda)\widetilde{\rho}_j(\lambda)=1$, \vspace{0.5em}  
 
(ii) \ \ \  $\left| \frac{d^l}{d\lambda^l}\Big(\lambda^{\vartheta(i)+\vartheta(j)}\widetilde{\rho}_i(\lambda)\widetilde{\rho}_j(\lambda)\Big)\right| \lesssim_l \lambda^{\frac 12-l}, \quad \lambda\in (0,1/2), \ l=0,\,1,\,\ldots$ 
 \vspace{0.5em} 

If (ii) holds, the same analysis as in the proof of Proposition \ref{pro.est-for-IGamma} yields that
$$\big|\Omega_{s,h,\alpha,\beta,\tilde{\beta},\tilde{\alpha}}^\pm(x,y) \big| \le C \langle|x|-|y|\rangle^{-\frac 32}\min\Big\{ \langle x\rangle^{-(n-1)},\, \langle y\rangle^{-(n-1)} \Big\}$$
holds uniformly for all  $|s-h|<2$. This shows $\left\| \Omega_{s,h,\alpha,\beta,\tilde{\beta},\tilde{\alpha}}^\pm\right\|_{\mathbb{B}(L^p)}\le C_p $ uniformly for $s,h$ with $|s-h|<2$ and $1\le p\le \infty$.

We are left to consider the case that (i) is true. In this case, set
\begin{equation}\label{zeta}
\zeta^\pm(x,y)= \int_0^\infty e^{i\lambda(|x|\pm|y|)}\chi(\lambda)d\lambda.	
\end{equation}
Then
$$\Omega_{s,h,\alpha,\beta,\tilde{\beta},\tilde{\alpha}}^\pm(x,\, y)=C \, \chi_s(x)\chi_h(y)\mathbf{1}_{\{|x|>1,\,|y|>1\}} \cdot \frac{y^\alpha}{|y|^{\frac{n-1}{2}+\vartheta(j)}}	\frac{x^{\tilde{\alpha}}}{|x|^{\frac{n-1}{2}+\vartheta(i)}}\zeta^\pm(x,y).$$
In the following, we further split the situation into $||x|-|y||< 1$ and $||x|-|y||> 1$ to estimate $\Omega_{s,h,\alpha,\beta,\tilde{\beta},\tilde{\alpha}}^\pm(x,\, y)$.

When $\big||x|-|y|\big|\le 1$, a direct computation shows that when $|s-h|< 2$, there is
\begin{equation}\label{zetasmall}
\begin{aligned}
\Big|\Omega_{s,h,\alpha,\beta,\tilde{\beta},\tilde{\alpha}}^\pm(x,\, y)\Big|&\lesssim \langle x\rangle^{-\frac{n-1}{2}}\langle y\rangle^{-\frac{n-1}{2}}\chi_s(x)\chi_h(y),
\end{aligned}	
\end{equation}
which implies that the integral kernels $\Omega_{s,h,\alpha,\beta,\tilde{\beta},\tilde{\alpha}}^\pm(x,\, y)$ are admissible. Therefore, we have that the corresponding operators  are bounded on $L^p(\mathbb{R}^n)$ for $1\le p\le \infty$ with the bound independent of $s,h$ with $|s-h|<2$. 

When $\big||x|-|y|\big|> 1$,  by  integration by parts once for the integral \eqref{zeta},  we get
\begin{equation*}
\begin{aligned}
 \zeta^\pm(x,y)
=&C\frac{1}{|x|\pm|y|}+C'\frac{1}{|x|\pm|y|}\int e^{i\lambda(|x|\pm|y|)}\chi'(\lambda/\lambda_0)d\lambda\\
:=&  \zeta_{0}^{\pm}(x,y)+\zeta_{1}^{\pm}(x,y).
\end{aligned}
\end{equation*}
Using Lemma \ref{oscillatory estimates} again, it is not difficult to obtain that 
\begin{equation*}
\left|\chi_s(x)\chi_h(y)\zeta_{1}^{\pm}(x,y)\right| \lesssim \langle|x|\pm |y|\rangle^{-2},
\end{equation*}
which implies 
\begin{equation*}
\begin{aligned}
\Big|C(M_{i,j}^+,\alpha, \beta)\, \chi_s(x)\chi_h(y)\mathbf{1}_{\{|x|>1,\,|y|>1\}}& \cdot \frac{y^\alpha}{|y|^{\frac{n-1}{2}+\vartheta(j)}}	\frac{x^{\tilde{\alpha}}}{|x|^{\frac{n-1}{2}+\vartheta(i)}}\zeta_{1}^\pm(x,y)\Big|  \\
&\le C \langle|x|\pm |y|\rangle^{-2} \min\{ \langle x\rangle^{-(n-1)},\, \langle y\rangle^{-(n-1)} \},
\end{aligned}
\end{equation*}
with a constant $C$ independent of $s,h$.  Thus, the corresponding operator is admissible with a  $L^p$-bound being independent of $s,h$.

It finally remains to deal with $\zeta_{0}^{\pm}(x,y)$ when $\big||x|-|y|\big|> 1$.  Note that when $s=0$ or $h=0$,
$\chi_s(x)\chi_h(y)\mathbf{1}_{\{|x|>1,\,|y|>1\}}=0.$
Thus, $\Omega_{s,h,\alpha,\beta,\tilde{\beta},\tilde{\alpha}}^\pm(x,\, y)=0$ when $s=0$ or $h=0$ by  the definition of $\zeta_{0}^{\pm}$. Hence, we only focus on the case where $s>1$, $h>1$ and $|s-h|<2$. Since $|\alpha|=\vartheta(j)$ and $|\tilde{\alpha}|=\vartheta(i)$,  the both functions: $$\mathbf{1}_{\{|y|>1\}}\frac{y^\alpha}{|y|^{\vartheta(j)}} \ \ {\rm and } \ \ \mathbf{1}_{\{|x|>1\}}\frac{x^{\tilde{\alpha}}}{|x|^{\vartheta(i)}}$$ 
are bounded on $\R^n$. So we can reduce the part related to  $\zeta_{0}^{\pm}(x,y)$ to $L^p$ boundedness of  the following truncated operator:
$$(\Lambda^\pm f)(x)=\int_{||x|\pm|y||> 1}\frac{\chi_h(|y|)}{|y|^{\frac{n-1}{2}}}\ \frac{\chi_s(|x|)}{|x|^{\frac{n-1}{2}}}\ \frac{1}{|x|\pm|y|}\ f(y)dy,$$
Denote the sphere mean of $f$ by 
$$(M f)(r)=\int_{S^{n-1}} f(r\omega)d\omega, \ \  dy=r^{n-1}drd\omega.$$
Then, for $1<p<\infty$,  by the polar coordinates transform, one has
\begin{equation}\label{eq-trans-to-hilbert}
	\begin{split}
		\left\|\Lambda^\pm f\right\|_{L^p}^p
		=\int_0^\infty&\Big|\int_{|\rho\pm r|>1}\frac{ \chi_h(r)\ M f(r) r^{\frac{n-1}{2}}}{\rho\pm r}dr\Big|^p
		\chi_s(\rho)\rho^{n-1-\frac{(n-1)p}{2}}d\rho\\
		\lesssim& \int_0^\infty\Big|\int_{|\rho\pm r|>1}\frac{\chi_h(r)\ M f(r)r^{\frac{n-1}{2}}}{\rho\pm r}dr\Big|^p
		2^{h(n-1)(1-\frac{p}{2})}\chi_s(\rho)d\rho\\
		\lesssim& 2^{h(n-1)(1-\frac{p}{2})}\int_0^\infty\Big|\int_{|\rho\pm r|>1}\frac{\chi_h(r)\ M f(r)r^{\frac{n-1}{2}}}{\rho\pm r}dr\Big|^p d\rho\\
	\end{split}
\end{equation}
where  we have used $|s-h|<2$, $s,h\ge 1$ and $r\sim\rho\sim 2^s\sim 2^h$ in the support of $\chi_s(\rho)$.   Since the truncated Hilbert transforms 
$$\mathbb{H}: g\mapsto \int_{|\rho\pm r|>1}\frac{g(r)}{\rho\pm r}dr,$$
are $L^p$-boundedness for $1<p<\infty$, so by the left-side of  \eqref{eq-trans-to-hilbert}, it follows that
\begin{equation}\label{eqhilbert}
	\begin{split}
		\left\|\Lambda^\pm f\right\|_{L^p}^p\lesssim &2^{h(n-1)(1-\frac{p}{2})}
		\int_{\mathbb{R}} |\chi_h(r)M f(r) r^{\frac{n-1}{2}}|^p dr\\
		\lesssim &\int_{\mathbb{R}} |\chi_h(r)M f(r)|^p r^{\frac{n-1}{2}}dr
		= \|\chi_h f\|^p_{L^p(\mathbb{R}^n)}\le \|f\|^p_{L^p(\mathbb{R}^n)},
	\end{split}
\end{equation}
Therefore, the inequality ~\eqref{eqhilbert},  together with \eqref{zetasmall},  implies that that for all $1< p<\infty$,
$\|\Omega_{s,h,\alpha,\beta,\tilde{\beta},\tilde{\alpha}}^\pm\|_{\mathbb{B}(L^p)} \le C_p$
hold uniformly for $|s-h|<2$ when $\lambda^{\vartheta(i)+\vartheta(j)}\widetilde{\rho}_i(\lambda)\widetilde{\rho}_j(\lambda)=1$.

In summary, there exists a constant $C_p$, independent of $s,h\in \mathbb{N}_0$ with $|s-h|<2$, such that $\|\chi_sI_{i,j}^{\pm,b}\chi_h\|_{\mathbb{B}(L^p(\mathbb{R}^n))}\le C_p$ for all $1<p<\infty$.
\end{proof}

Finally, we prove Theorem \ref{thm-low energy-regular-even} by combining the results from the previous three subsection.

\begin{proof}[{\bf Proof of Theorem \ref{thm-low energy-regular-even}.}]
Recall that we reduced the $L^p$-boundedness of $W$ to that of $W_{i,j,s}^L$ for $i,j\in J_\mathbf{k}$ and $s=1,2,3$ at the beginning of this section. From Proposition \ref{propo-low energy-even-1} and Proposition \ref{estforw_3}, it follows that $W_{i,j,1}^L$ and $W_{i,j,3}^L$ are bounded on $L^p(\mathbb{R}^n)$ for all $1<p<\infty$ and $i,j\in J_\mathbf{k}$. Thus, it remains to consider $W_{i,j,2}^L$.

Note that by the definition of $J_\mathbf{k}$ in \eqref{eq-Jk-even} and \eqref{eq-Jk-even-2}, the set $J_\mathbf{k}$ is increasing with respect to $\k$. Moreover, when zero energy is a $\mathbf{k}$-th kind resonance, for even $n\le 2m$,
\begin{equation}\label{j_k}
j_{\k}=\max J_\mathbf{k}=
\begin{cases}
m-\frac{n}{2}+\frac 12 , &\quad \text{when}~\mathbf{k}=0,\\
m-\frac{n}{2}+\frac 34 , &\quad \text{when}~\mathbf{k}=1,\\
k_c+\mathbf{k}-2\le 2m-n,	&\quad \text{when}~1<\mathbf{k}\le k_c=m-\frac n2+1,\\
k_c+\mathbf{k}-2,&\quad \text{when}~k_c<\mathbf{k}\le m_n,\\
2m-\frac{n}{2}+\frac 12, &\quad \text{when}~\mathbf{k}= m_n+1,
\end{cases}
\end{equation}
and for an even integer $n$ such that $2m<n\le 4m$, 
\begin{equation}\label{j_k_1}
	j_{\k}=\max J_\mathbf{k}=
	\begin{cases}
		m-\frac{n}{2}, &\quad \text{when}~\mathbf{k}=k_c=0,\\
		
		k_c+\mathbf{k}-1,&\quad \text{when}~0=k_c<\mathbf{k}\le m_n,\\
		2m-\frac{n}{2}+\frac 12, &\quad \text{when}~\mathbf{k}= m_n+1.
	\end{cases}
\end{equation}
Combining this with Proposition \ref{propo-low energy-even-2} and the fact that $W^L\in \mathbb{B}(L^2)$, we conclude that when zero energy is a $\mathbf{k}$-th kind resonance, all $W_{i,j,2}^L$ for $i,j\in J_{\k}$ belong to $\mathbb{B}(L^p)$ in the range stated in Theorem \ref{thm-low energy-regular-even}.
To summarize, we conclude Theorem \ref{thm-low energy-regular-even}.
\end{proof}

%\vspace{1.2em}
\section{High energy estimates}\label{sec:high energy}

In this section, we will prove the $L^p$-boundedness of $W^H$. 
The main result is stated below. 
\begin{theorem}\label{thm-high energy-n<2m}
	Let $m\ge 2$ and $n$ be even such that $1\le n\le 4m$, and Assumption \ref{Assumption} holds. Then the high energy part $W^H$ is bounded on $L^p(\mathbb{R}^n)$ for all $1\le p\le\infty$.
\end{theorem}

\begin{remark}
	Let $H=-\Delta+V$ (i.e. $m=1$). Then Yajima \cite{Yajima-99-CMP} proved that when $n = 2$, the operator $W^H$ is $L^p$-bounded for $1 < p < \infty$,  and \cite{Yajima-1995-even} showed that when $n = 4$, the operator $W^H$ is $L^p$-bounded for $1\leq p\leq\infty$.
\end{remark}

\emph{Since Erdo\u{g}an-Green \cite{EG22} has proved that $W^H$ is bounded on $L^p(\R^n)$ for all $1\le p\le\infty$ when $n>2m,$ we only need to prove the case when $1\le n\le 2m$. Theorem \ref{thm-high energy-n<2m} follows directly from the result for $n=2m$ $($see Theorem \ref{thm-high energy-1}$)$ and the result for $1\le n<2m$ $($see Theorem \ref{thm-high energy-2}$)$.}

Recall that 
\begin{equation*} W^H=\frac{m}{\pi i}\int_0^\infty\lambda^{2m-1}(1-\chi(\lambda))R^+(\lambda^{2m})V\big(R_0^+(\lambda^{2m})-R_0^-(\lambda^{2m})\big)d\lambda, \end{equation*} 
where the even cut-off function $\chi(\lambda)\in C_0^\infty(\mathbb{R}^n)$ satisfies $\chi \equiv 1$ on $(-\lambda_0/2, \lambda_0/2)$ and $\text{supp} \chi \subset [-\lambda_0, \lambda_0]$.

\subsection{The cases for even dimension $n=2m$}
We first consider the dimension $n=2m$, which can be viewed as the critical dimension. 

\begin{theorem}\label{thm-high energy-1}
	Let $n=2m$ and $m\ge 2$. Assume that $H$ has no embedded positive eigenvalues  and $|V(x)|\lesssim\langle x\rangle^{-\beta}$ with some $\beta>n+4.$ Then the operator $W^H$  is bounded on $L^p(\mathbb{R}^n)$ for all $1\le p\le\infty.$
\end{theorem}

By the following resolvent identity 
\begin{equation*}
	R^+(\lambda^{2m})=R_0^+(\lambda^{2m})-R_0^+(\lambda^{2m})VR^+(\lambda^{2m}),
\end{equation*}
we can divide the high energy portion $W^H$ into the following integrals:
\begin{equation*}
	\begin{split}
		W_1^H=&\frac{m}{\pi i}\int_0^\infty\lambda^{2m-1}(1-\chi(\lambda))R_0^+(\lambda^{2m})V(R_0^+(\lambda^{2m})-R_0^-(\lambda^{2m}))d\lambda,\\
		W_2^H=&\frac{m}{\pi i}\int_0^\infty\lambda^{2m-1}(1-\chi(\lambda))R_0^+(\lambda^{2m})VR^+(\lambda^{2m})V(R_0^+(\lambda^{2m})-R_0^-(\lambda^{2m}))d\lambda.
	\end{split}
\end{equation*}
Next, we will prove that all of $W_1^H$ and $W_2^H$ are bounded on $L^p$ for $1\le p\le \infty$ to obtain Theorem \ref{thm-high energy-1} (see Proposition \ref{propo-Lp--W_{1,1}{H}} and \ref{propo-Lp--W_{2}{H}} below).

We will first prove the result for $W_1^H$. The proof can be carried out using the methods from \cite{Yajima-JMSJ-95}, \cite{Yajima-99-CMP} and \cite{EG22}. Applying the identity \eqref{eq-repren-spectral-measure}, we rewrite $W_1^H$ as $$W_1^H=\frac{m}{\pi i}\int_{0}^{\infty}\lambda^{2m+\delta}R_0^+(\lambda^{2m})V(R_0^+(\lambda^{2m})-R_0^-(\lambda^{2m})) d\lambda \circ \frac{(1-\chi(D))}{|D|^{1+\delta}},$$ 
for a small number $0<\delta<1/10$ , where $D=-i\nabla$.
It is clear that the Fourier multiplier $\frac{(1-\chi(D))}{|D|^{1+\delta}}$ is bounded on $L^p$ for all $1\le p\le \infty$. Thus, it suffices to consider the operator 
$$W_{1,1}^H=\frac{m}{\pi i}
\int_{0}^{\infty}\lambda^{2m+\delta}R_0^+(\lambda^{2m})V(R_0^+(\lambda^{2m})-R_0^-(\lambda^{2m}))d\lambda$$
Applying the identity \eqref{eq-repren-spectral-measure} again, for a Schwartz function $f$, we have 
\begin{equation}\label{eq: Fourier tran of high W_11^H}
	\begin{aligned}
		\mathscr{F}(W_{1,1}^H f)(\xi)=&\lim_{\varepsilon\rightarrow 0^+}\int_{\R^n}\frac{\widehat{V}(k)}{|\xi|^{2m}-|\xi-k|^{2m}-i\varepsilon}\cdot |\xi-k|^{1+\delta} \widehat{f}(\xi-k)dk \\
		=&\lim_{\varepsilon\rightarrow 0^+}\int_{\R^n}\frac{\widehat{V}(k)}{|\xi|^{2}-|\xi-k|^{2}-i\varepsilon}\cdot \frac{|\xi-k|^{1+\delta}}{P_{k}(\xi)} \widehat{f}(\xi-k)dk \\
		=&\lim_{\varepsilon\rightarrow 0^+}\int_{\R^n}\frac{\widehat{V}(k)}{2\omega\cdot\xi-|k|-i\varepsilon}\cdot \frac{\big|\xi/|k|-\omega\big|^{1+\delta}}{P_{\omega}(\xi/|k|)} \widehat{f}(\xi-k) |k|^{-2m+2+\delta}dk
	\end{aligned}
\end{equation} 
where $P_{k}(\xi)=\sum_{\ell=0}^{m-1}|\xi|^{2\ell}|\xi-k|^{2m-2-2\ell}$, and $\omega=k/|k|$. 
Note that 
$$\frac{1}{2\omega\cdot\xi-|k|-i\varepsilon}=\int_0^{\infty} e^{-it(2\omega\cdot\xi-|k|-i\varepsilon)} dt.$$
Formally,  $W_{1,1}^H f(x)$ can be represented as 
\begin{equation*}
	\begin{split}
		W_{1,1}^H f(x)=&\int_0^{\infty}\int_{\R^n} e^{-i|k|t}\, \widehat{V}(k) \left( h_{|k|,\omega}\ast (e^{ikz}f(z))\right)  (x-2t\omega) |k|^{-2m+2+\delta}dk dt \\
		=& \int_0^{\infty}\int_{S^{n-1}} \int_{\R^n} \int_0^{\infty} s^{1+\delta} e^{-is(t-\omega\cdot x+\omega\cdot y)}\, \widehat{V}(s \omega) h_{s,\omega}(y) ds \, f (x-2t\omega-y) dy d\omega dt,
	\end{split}
\end{equation*}
where $s=|k|$, $h_{s,\omega}=\mathscr{F}^{-1}(\frac{\big|\xi/|k|-\omega\big|^{1+\delta}}{P_{\omega}(\xi/|k|)})$. Upon substituting $t$ with $t+\omega\cdot x-\omega\cdot y$, then we have
\begin{equation}\label{eq: Fourier tran of high W_11^H-1}
	W_{1,1}^H f(x)= \int_{S^{n-1}} \int_{\R^n}\int_{-\omega\cdot x+\omega\cdot y}^{\infty} \int_0^{\infty} s^{1+\delta} e^{-ist}\, \widehat{V}(s \omega) h_{s,\omega}(y) ds \, f (x-2\omega(\omega\cdot x)-y+2\omega(\omega\cdot y)-2t\omega) dy d\omega dt.
\end{equation}

Now we prove the $L^p$-boundedness of $W_{1}^H$.
\begin{proposition}\label{propo-Lp--W_{1,1}{H}}
	Assume that  $\langle x\rangle V(x)\in L^1(\R^n)\cap L^2(\R^n)$ with $n=2m$ and $m\ge 2$. Then the operator $W_{1}^{H}$  is bounded on $L^p(\mathbb{R}^n)$ for all $1\le p\le\infty.$
\end{proposition} 
\begin{proof}
	Thanks to the previous reduction, we only need to show that $W_{1,1}^H$ is bounded on $L^p$ for all $1\le p\le \infty$.
	Noting that $x\rightarrow x-2\omega(\omega\cdot x)$ is a isometry, by \eqref{eq: Fourier tran of high W_11^H-1} and Minkowski’s integral inequality, it suffices to prove that 
	$$K(t,\omega,y)=\int_0^{\infty}  e^{-ist}s^{1+\delta}\, \widehat{V}(s \omega) h_{s,\omega}(y) ds \in L^{1}(\R \times S^{n-1} \times \R^n ), \quad \text{with $\delta\le 1/10$}.$$ 
	
	We first consider the property of $h_{s,\omega}(y) $. By the definition, we have  
	$$h_{1,\omega}(y)=\frac{1}{(2\pi)^n} \int_{\R^n} e^{i\xi\cdot y} \frac{\big|\xi-\omega\big|^{1+\delta}}{P_{\omega}(\xi)} d\xi,$$
    and $$h_{s,\omega}(y)= s^{n} h_{1,\omega}(sy).$$
	For $m\ge 2$, a simple calculation shows that 
	\begin{equation}\label{eq: est fof h-1}
		\left|\partial_{\xi}^{\alpha} \left( \frac{\big|\xi-\omega\big|^{1+\delta}}{P_{\omega}(\xi)} \right) \right|\lesssim_{\alpha} \min\{ |\xi-\omega|^{1+\delta-|\alpha|}, |\xi-\omega|^{1+\delta-(2m-2) -|\alpha|}\}, 	
	\end{equation}
	for each $\alpha \in \N_0^n$. 
	We decompose the integral above as 
	$$I_{i}(y)=\int_{\R^n} e^{i(\xi-\omega) y} \frac{\big|\xi-\omega\big|^{1+\delta}}{P_{\omega}(\xi)}\phi_i(|y||\xi-\omega|) d\xi,$$ 
	where $\phi_0 \in C^{\infty}(\R)$ such that $\phi_0(\lambda)=1$ when $|\lambda|\leq 1/2$ and $\phi_0(\lambda)=0$ when $|\lambda|\geq 1$, and $\phi_1(\lambda)=1-\phi_0$. Then we have from \eqref{eq: est fof h-1} with $m\ge 2$ that 
	$$\left| I_{0}(y)\right|\lesssim  \min\{ |y|^{-1-n-\delta}, |y|^{1-n-\delta }\}.$$
	Using integration by parts $n+2$ times, it follows from \eqref{eq: est fof h-1} with $m\ge 2$ that 
	$$
	\left| I_{1}(y)\right|\lesssim \min\{ |y|^{-1-n-\delta}, |y|^{1-n-\delta }\}. $$
	     We conclude that 
	\begin{equation}\label{eq: est fof h-2}
		\left| h_{s,\omega}(y)\right|\lesssim s^{\delta}|y|^{-n+\delta}\mathbf{1}_{\{|y|\le 1\}}+ s^{-\delta}|y|^{-n-\delta}\mathbf{1}_{\{|y|>1\}}.
	\end{equation}
	Meanwhile, note that 
	$$\partial_s (h_{s,\omega}(y))= s^{n-1}h_{1,\omega}(sy)+  \frac{s^n i}{(2\pi)^n} \int_{\R^n} e^{i s\xi\cdot y} \,  \frac{\xi\cdot y \big|\xi-\omega\big|^{1+\delta}}{P_{\omega}(\xi)} d\xi.$$
	Similarly, one can conclude that
	\begin{equation}\label{eq: est fof h-3}
		\left| \partial_s (h_{s,\omega}(y)) \right|\lesssim s^{\delta-1}|y|^{-n+\delta}\mathbf{1}_{\{|y|\le 1\}}+ s^{-\delta-1}|y|^{-n-\delta}\mathbf{1}_{\{|y|>1\}}.
	\end{equation}
	
	We now turn our attention back to  $K(t,\omega,y)$. We decompose $\R$ as $|t|\le 1$ and $|t|>1$. In the region where $|t| \le 1$, using the H\"{o}lder inequality and Plancherel theorem,  one has that
	\begin{align*}
		\int_{S^{n-1} \times \R^n}\int_{|t|\le 1}\left|K(t,\omega,y) \right| \lesssim& \frac{1}{2\pi} \left\|K(t,\omega,y)\right\|_{L^1\left(S^{n-1} \times \R^n, L^2(\R)\right)} \\
		=&2\pi  \left\|s^{1+\delta} \, \widehat{V}(s \omega) h_{s,\omega}(y)\right\|_{L^1\left(S^{n-1} \times \R^n, L^2(\R)\right)} \\
		\lesssim &\left\|s^{1+2\delta} \, \widehat{V}(s \omega)\right\|_{L^1\left(S^{n-1}, L^2(\R)\right)}+\left\|s \, \widehat{V}(s \omega)\right\|_{L^1\left(S^{n-1}, L^2(\R)\right)} \\
		\lesssim & \left\| \frac{\widehat{V}(\xi)}{|\xi|^{\frac{n-3}{2}-2\delta}}\right\|_{L^2\left(\R^{n}\right)}+
		\left\| \frac{\widehat{V}(\xi)}{|\xi|^{\frac{n-3}{2}}}\right\|_{L^2\left(\R^{n}\right)} \\
		\lesssim& \left\|\widehat{V}(\xi)\right\|_{L^{\infty}\left(\R^{n}\right)}+ \left\|\widehat{V}(\xi)\right\|_{L^2\left(\R^{n}\right)}\lesssim 1,
	\end{align*}
	where in the third line above, we have used \eqref{eq: est fof h-2}, and in the last line, we used the facts $V \in L^1\cup L^2$ and $n\ge 4$.
	
    Similarly, in the region where $|t|> 1$, using the integration by parts,  the H\"{o}lder inequality and Plancherel theorem, one concludes that 
	\begin{align*}
		\int_{S^{n-1} \times \R^n}\int_{|t|> 1}\left|K(t,\omega,y) \right| \lesssim& \int_{S^{n-1} \times \R^n}\int_{|t|> 1}|t|^{-1}\left| \int_0^{\infty}  e^{-ist} \partial_s \left( s^{1+\delta}\, \widehat{V}(s \omega) h_{s,\omega}(y)\right)  ds \right| \\
		\lesssim &\left\|\int_0^{\infty}  e^{-ist} \partial_s \left( s^{1+\delta}\, \widehat{V}(s \omega) h_{s,\omega}(y)\right) \right\|_{L^1\left(S^{n-1} \times \R^n, L^2(\R)\right)} \\
		=&2\pi  \left\|\partial_s \left( s^{1+\delta}\, \widehat{V}(s \omega) h_{s,\omega}(y)\right)\right\|_{L^1\left(S^{n-1} \times \R^n, L^2(\R)\right)}. 
	\end{align*}
	Since 
	$$ \left| \partial_s \left( s^{1+\delta}\, \widehat{V}(s \omega) h_{s,\omega}(y)\right) \right| \lesssim \left|\partial_s \left(  \widehat{V}(s \omega)\right) s^{1+\delta}\, h_{s,\omega}(y)\right|+ \left| \partial_s \left( s^{1+\delta}\,  h_{s,\omega}(y)\right) \widehat{V}(s \omega)\right|,$$
	combining with the estimate for $\partial_s (h_{s,\omega}(y))$ in \eqref{eq: est fof h-3}, we conclude that 
	\begin{align*}
		\int_{S^{n-1} \times \R^n}\int_{|t|> 1}\left|K(t,\omega,y) \right| dt dy d\omega 
		\lesssim & \left\| \frac{\widehat{|x|V}(\xi)}{|\xi|^{\frac{n-3}{2}-2\delta}}\right\|_{L^2\left(\R^{n}\right)}+
		\left\| \frac{\widehat{|x|V}(\xi)}{|\xi|^{\frac{n-3}{2}}}\right\|_{L^2\left(\R^{n}\right)}+\\
		&\ \ \ \ \  \left\| \frac{\widehat{V}(\xi)}{|\xi|^{\frac{n-1}{2}-2\delta}}\right\|_{L^2\left(\R^{n}\right)}
		+\left\| \frac{\widehat{V}(\xi)}{|\xi|^{\frac{n-1}{2}}}\right\|_{L^2\left(\R^{n}\right)} 
		\lesssim  1.
	\end{align*}
	We have completed the proof.
\end{proof}

It remains to prove the $L^p$-boundedness of $W_2^{H}$.
To prove it, we need the following estimate for the resolvent in high energy part (see \cite[Lemma 3.9]{FSWY20}).
\begin{lemma}\label{lemma-high energy-second part-n<2m}
	Assume that $H$ has no embedded positive eigenvalues and $|V(x)|\lesssim \langle x\rangle^{-2(k+1)}$, then for any $\sigma>k+\frac{1}{2},$ $\partial^k_\lambda (R^\pm(\lambda^{2m}))$ and $\partial^k_\lambda(R_0^\pm(\lambda^{2m}))$ are bounded from $L^{2}_\sigma(\mathbb{R}^n)$ to $L^{2}_{-\sigma}(\mathbb{R}^n).$ Furthermore, we have 
	\begin{equation}\label{eq-resolvent high energy estimates}
		\Big\|\partial^k_\lambda (R_0^\pm(\lambda^{2m}))\Big\|_{L^{2}_\sigma(\mathbb{R}^n)\rightarrow L^{2}_{-\sigma}(\mathbb{R}^n)}+\Big\|\partial^k_\lambda (R^\pm(\lambda^{2m}))\Big\|_{L^{2}_\sigma(\mathbb{R}^n)\rightarrow L^{2}_{-\sigma}(\mathbb{R}^n)}\lesssim \lambda^{-(2m-1)}.
	\end{equation}
\end{lemma}

\begin{proposition}\label{propo-Lp--W_{2}{H}}
	Let $n=2m$ and $m\ge 2$.	Assume that $H$ has no embedded positive eigenvalues and $|V(x)|\lesssim\langle x\rangle^{-\beta}$ with some $\beta>n+4.$ Then the operator $W_{2}^{H}$  is bounded on $L^p(\mathbb{R}^n)$ for all $1\le p\le\infty.$
\end{proposition} 
\begin{proof}
	Recall that the operator $W_2^H$ is defined by
	\[W_2^H=\frac{m}{\pi i}\int_0^\infty\lambda^{2m-1}(1-\chi(\lambda))R_0^+(\lambda^{2m})VR^+(\lambda^{2m})V(R_0^+(\lambda^{2m})-R_0^-(\lambda^{2m}))d\lambda.\]
	Applying \eqref{eq-free kernel-F}, we write its integral kernel as
	\begin{equation*}
		\begin{aligned}
			W_2^{H}(x,y)=&\frac{m}{\pi i}\int_0^\infty(1-\chi(\lambda))\left(  e^{i\lambda(|x|+ |y|)}T^{H,+}(\lambda,x,y)-e^{i\lambda(|x|- |y|)}T^{H,-}(\lambda,x,y)\right)d\lambda \\
			:=&W_2^{H,+}(x,y)-W_2^{H,-}(x,y),
		\end{aligned}
	\end{equation*}	
	where 
	$$T^{H,\pm}(\lambda,x,y)=\frac{m}{\pi i}\Big\langle R^+(\lambda^{2m})V\frac{e^{\pm i(|y-\cdot|-|y|)}F_n^\pm(\lambda|y-\cdot|)}{|y-\cdot|^\frac{n-1}{2}},\, V\frac{e^{- i(|x-\cdot|-|x|)}F_n^-(\lambda|x-\cdot|)}{|x-\cdot|^\frac{n-1}{2}} \Big\rangle.$$	
	Using Lemma \ref{lemma-high energy-second part-n<2m} and \eqref{eq-F decay-odd},  a direct computation shows that for $l=0,1,2,\ldots, \frac{n}{2}+1$, one has 
	$$\left| \partial_\lambda^lT^{H,\pm}(\lambda,x,y)\right| \lesssim\lambda^{-2m+1}\langle x\rangle^{-\frac{n-1}{2}}\langle y\rangle^{-\frac{n-1}{2}}, \quad \lambda\ge\lambda_0/2.$$
	The estimate relies on the boundedness of the $\frac{n}{2}+1$-order derivative of $R^+(\lambda^{2m})$, which requires the decay condition on $V$.
	Note that $1-2m\le -2$ since $m\ge 2$. Thus, when $||x|\pm |y||\le 1$, it follows immediately from the expressions of $W_2^{H,\pm}(x,y)$ that
	\begin{equation*}
		\big| W_2^{H,+}(x,y)\big|+ \big| W_2^{H,-}(x,y)\big|	\lesssim \langle x\rangle^{-\frac{n-1}{2}}\langle y\rangle^{-\frac{n-1}{2}}.	
	\end{equation*}
	When $||x|\pm |y||\ge 1$, applying integration by parts $\frac{n}{2}+1$ times, we have 
	\begin{equation*}
		\big| W_2^{H,+}(x,y)\big|+ \big| W_2^{H,-}(x,y)\big|	\lesssim\left||x|-|y| \right|^{-\frac{n+2}{2}}  \langle x\rangle^{-\frac{n-1}{2}}\langle y\rangle^{-\frac{n-1}{2}}.
	\end{equation*}
	These two estimates show that  $W_2^{H,\pm}(x,y)$ are admissible, thereby ensuring that $W_2^H$ is bounded on $L^p(\mathbb{R}^n)$ for all $1 \leq p \leq \infty$.
\end{proof}

\subsection{ The cases of  even dimension $1\le n< 2m$}
\begin{theorem}\label{thm-high energy-2}
	Let $1\le n<2m$ be even. Assume that $H$ has no embedded positive eigenvalues  and $|V(x)|\lesssim\langle x\rangle^{-\beta}$ with some $\beta>n+4.$ Then the operator $W^H$  is bounded on $L^p(\mathbb{R}^n)$ for all $1\le p\le\infty.$
\end{theorem}
\begin{proof}
We still divide the high energy portion $W^H$ into (up to constants)
\begin{equation*}
	\begin{split}
		W_1^H=&\frac{m}{\pi i}\int_0^\infty\lambda^{2m-1}(1-\chi(\lambda))R_0^+(\lambda^{2m})V(R_0^+(\lambda^{2m})-R_0^-(\lambda^{2m}))d\lambda,\\
		W_2^H=&\frac{m}{\pi i}\int_0^\infty\lambda^{2m-1}(1-\chi(\lambda))R_0^+(\lambda^{2m})VR^+(\lambda^{2m})V(R_0^+(\lambda^{2m})-R_0^-(\lambda^{2m}))d\lambda.\\
	\end{split}
\end{equation*} 
Using \eqref{eq-free kernel-F},  we can write the integral kernel $W_1^H(x,y)$ as follows:
$$	W_1^{H}(x,y)=	W_{1}^{H,+}(x,y)-W_{1}^{H,-}(x,y)$$
where
\begin{equation}\label{eq-int-kernel-for-W1H(x,y)-small-d}
	\begin{split}
		W_1^{H,\pm}(x,y)=\frac{m}{\pi i}\int_{0}^\infty\lambda^{n-2m}(1-\chi(\lambda))\Big(\int_{\mathbb{R}^n}e^{i\lambda|x-z|\pm i\lambda|z-y|}
		\frac{F_n^+(\lambda|x-z|)V(z)F_n^-(\lambda|z-y|)}{|x-z|^{\frac{n-1}{2}}|z-y|^{\frac{n-1}{2}}}dz\Big)d\lambda.
	\end{split}
\end{equation}
Then using integration by parts and \eqref{eq-F decay-odd}, $W_1^{H}(x,y)$ is bounded by 
$$\left|	W_1^{H}(x,y)  \right|\lesssim \int_{\mathbb{R}^n} \frac{|V(z)|}{|x-z|^{\frac{n-1}{2}}\langle|x-z|-|z-y|\rangle^{(n+1)}|z-y|^{\frac{n-1}{2}}}dz,   $$
which is admissible.
The proof of $W_2^H(x,y)$ can be carried out by the same method in the case $n=2m$. We omit it here.
\end{proof}

%\vspace{1.2em}
\section{The $L^p$-unboundedness of $W_\pm(H;(-\Delta)^m)$}\label{sec-unb}
In this section, we will analyze the $L^p$-unboundedness of  wave operators $W_\pm(H;(-\Delta)^m)$  when zero energy is the $\k$-th kind resonance of $H$ with $\k>k_c$, including the case of zero eigenvalue. We establish this result through proof by contradiction.   The key idea is that we reduce the problem to the optimality of time-decay estimates of $e^{itH}$ in the weighted $L^2$ space.

Recall that number $p_{\k}$ defined in \eqref{eq: def of p_k} is
\begin{equation*}
	p_{\mathbf{k}}=
	\begin{cases}
		\frac{2n}{2\mathbf{k}+n-2m-2}, &\quad k_c<\mathbf{k}\le m_n-1,\ 1\le n\le 2m,\\
		\frac{n}{\mathbf{k}+n-2m-1}, &\quad k_c<\mathbf{k}\le m_n-1,\ 2m< n\le 4m,\\
		2, &\quad \mathbf{k}=m_n,\ m_n+1.
	\end{cases}	
\end{equation*}

\begin{proposition}\label{Pro: lower bound of the time decay}
Suppose that the assumption of Theorem \ref{thm-unbound} holds. Let $\sigma>\frac n2+1$. Then, when $k_c<\mathbf{k}<m_n$, there exists  a non-zero constant $C>0$ and a large number $L>0$ such that for $|t|>L$, 
	\begin{equation}\label{eq: time decay lower bound}
		\Big\|e^{i tH}P_{ac}(H)\Big\|_{\mathbb{B}(L^{2}_{\sigma},\, L^{2}_{-\sigma})} \ge C |t|^{-\frac{n}{m}(\frac 12-\frac{1}{p_{\k}})};
	\end{equation}
When  $\mathbf{k}= m_n$ or $\mathbf{k}=m_n+1$, there exists  a non-zero constant $C>0$ and  a large number $L>0$ such that for $|t|>L$, 
\begin{equation}\label{eq: time decay lower bound-1}
	\Big\|e^{i tH}P_{ac}(H)\Big\|_{\mathbb{B}(L^{2}_{\sigma},\, L^{2}_{-\sigma})}\ge C\big(\log |t|\big)^{-1}.
\end{equation}
\end{proposition}
\begin{proof}
	From Stone's formula, the propagator $e^{i tH}P_{ac}(H)$ can be represented as
	$$e^{i tH}P_{ac}(H)=\frac{m}{\pi}\int_0^{\infty} e^{  it\lambda^{2m}}\Big(R^{+}(\lambda^{2m})-R^{-}(\lambda^{2m})\Big) \lambda^{2m-1}d\lambda.$$
	With this identity, we can decompose $e^{i tH}P_{ac}(H)$ into two parts as
	$$\Omega^{L}=\frac{m}{\pi i}\int_0^{\infty} e^{it\lambda^{2m}}\Big(R^{+}(\lambda^{2m})-R^{-}(\lambda^{2m})\Big) \lambda^{2m-1}\chi(\lambda)d\lambda,$$
	and
	$$\Omega^{H}=\frac{m}{\pi i}\int_0^{\infty} e^{it\lambda^{2m}}\Big(R^{+}(\lambda^{2m})-R^{-}(\lambda^{2m})\Big) \lambda^{2m-1}(1-\chi)(\lambda)d\lambda,$$
	where $\chi(\lambda)$ is the same as that in \eqref{eq-wave operator-low-0}. By using the higher energy estimates \eqref{eq-resolvent high energy estimates} of $R^\pm(\lambda^2m)$,  performing $\frac{n}{2}$ times of integration by parts on $\Omega^{H}$, we obtain that the time-decay for $\Omega^{H}$ can be bounded by $|t|^{-\frac{n}{2}}$ from $L^{2}_{\sigma}$ to $L^{2}_{-\sigma}$. Note that $\frac{n}{2}> \frac{n}{m}(\frac 12-\frac{1}{p_{\k}})$, which means that the higher energy part $\Omega^{H}$ is negligible in the estimate \eqref{eq: time decay lower bound} and \eqref{eq: time decay lower bound-1} if $t$ is large.
	
Next to focus on the lower energy part $\Omega^{L}$.	By using the symmetric second resolvent identity:
$$R^\pm(\lambda^{2m})=R_0^\pm(\lambda^{2m})-R_0^{\pm}(\lambda^{2m})v\left(M^{\pm}(\lambda) \right)^{-1}vR_0^{\pm}(\lambda^{2m}),$$
it follows that 
	\begin{align*}
		\Omega&^{L}=e^{it(-\Delta)^{m}}\chi(|D|) \\
		&-\frac{m}{\pi i}\int_0^{\infty} e^{it\lambda^{2m}}\Big[R_0^{+}(\lambda^{2m})v\left(M^{+}(\lambda) \right)^{-1}vR_0^{+}(\lambda^{2m}) -R_0^{-}(\lambda^{2m}) v\left(M^{-}(\lambda) \right)^{-1}vR_0^{-}(\lambda^{2m})\Big] \lambda^{2m-1}\chi(\lambda)d\lambda,
	\end{align*}
	where $D=i\nabla$.
	Since for any $2\le p\le \infty$ (see e.g.,\cite{ZYF})
	\begin{equation*}
		\left\|e^{it(-\Delta)^{m}}\chi(|D|)\right\|_{L^{p'}-L^p}\lesssim |t|^{-\frac{n}{m}(\frac {1}{2}-\frac{1}{p})}, \quad t\neq0,
	\end{equation*} 
     by the H\"{o}lder inequality,   we have that $\mathbb{B}(L^{2}_{\sigma},\, L^{2}_{-\sigma})$-bound of $e^{it(-\Delta)^{m}}\chi(|D|)$ satisfies the same estimate. Hence when we choose $p>p_{\k}$, this implies that the free part $e^{it(-\Delta)^{m}}\chi(|D|)$ of $\Omega^L$ also is negligible to the estimate \eqref{eq: time decay lower bound} and \eqref{eq: time decay lower bound-1}  if $t$ is large.

     Hence we just consider the integral part of $\Omega^L$  associated with
	$$R_0^{+}(\lambda^{2m})v\left(M^{+}(\lambda) \right)^{-1}vR_0^{+}(\lambda^{2m}) -R_0^{-}(\lambda^{2m}) v\left(M^{-}(\lambda) \right)^{-1}vR_0^{-}(\lambda^{2m}),$$
    which are divided into the following two cases.
		\vskip0.2cm
 \textbf{A.\ Proof of \eqref{eq: time decay lower bound} when $k_c<\k< m_n$.} 
	
    When zero energy is the $\k$-th kind resonance of $H$ with $k_c<\k< m_n$ and $1\le n\le 2m$,  using the expansions of $R_0^{\pm}(\lambda^{2m})$ in \eqref{eq-free expansion-even-2} with $\theta=2m-n+1$ and the expansions of $(M^{\pm}(\lambda))^{-1}$ in \eqref{eq-M expansion-even},  the following identity in $\mathbb{B}(L^{2}_{\sigma},\, L^{2}_{-\sigma})$ holds:
	\begin{equation}\label{eq: expan for spectral measure-0}
		\begin{aligned}
			R_0^{+}(\lambda^{2m})v\left(M^{+}(\lambda) \right)^{-1}&vR_0^{+}(\lambda^{2m}) -R_0^{-}(\lambda^{2m}) v\left(M^{-}(\lambda) \right)^{-1}vR_0^{-}(\lambda^{2m}) \\
			&=b_0^2\, \lambda^{2m-n-2j_{\k}} \, G_{2m-n}^0vQ_{j_{\k}}\left( M^{+}_{j_{\k},j_{\k}}-M^{-}_{j_{\k},j_{\k}}\right) Q_{j_{\k}}vG_{2m-n}^0+\Phi(\lambda),	
		\end{aligned}
	\end{equation}
	where  we have used the cancellation property of $Q_j$ in \eqref{eq-cancellation-Q_j}. In the above identity, $j_{\k}=\max J_{\k}$, and $b_0G^0_{2m-n}(x,y)$ is defined in \eqref{eq-G-even}, which is the fundamental solution of $(-\Delta)^{m}$ when $1\le n\le 2m$, and $\Phi(\lambda)$ satisfies 
	\begin{equation}\label{eq: est for Phi}
    \left\|\frac{d^l}{d\lambda^{l}}\Phi(\lambda)\right\|_{\mathbb{B}(L^{2}_{\sigma},\, L^{2}_{-\sigma})}\lesssim |\log \lambda|^{-\frac 12} \lambda^ {2m-n-2j_{\k}-l} , 
	\end{equation}
	for $l=0,1,\ldots, \frac{n}{2}+1$ and $\lambda\in (0,\lambda_0)$. 
	
	Similarly, when $2m<n\le 4m$, $b_0G_{2m-n}(x,y)=b_0|x-y|^{2m-n}$ is  the fundamental solution of $(-\Delta)^{m}$.  By the expansion \eqref{eq-M expansion-even-1}, we have that 
	\begin{equation}\label{eq: expan for spectral measure-1}
		\begin{aligned}
			R_0^{+}(\lambda^{2m})v\left(M^{+}(\lambda) \right)^{-1}&vR_0^{+}(\lambda^{2m}) -R_0^{-}(\lambda^{2m}) v\left(M^{-}(\lambda) \right)^{-1}vR_0^{-}(\lambda^{2m}) \\
			&=b_0^2\, \lambda^{2m-n-2j_{\k}} \, G_{2m-n}vQ_{j_{\k}}\left( M^{+}_{j_{\k},j_{\k}}-M^{-}_{j_{\k},j_{\k}}\right) Q_{j_{\k}}vG_{2m-n}+\Phi(\lambda),	
		\end{aligned}
	\end{equation}
and $\Phi(\lambda)$ also satisfies \eqref{eq: est for Phi}.

We first  consider the part of $\Omega^{L}$ associated to $\Phi(\lambda)$, which is written as 
	$$ \frac{m}{\pi i}\int_0^{\infty} e^{it\lambda^{2m}}\lambda^{2m-1}\Phi(\lambda)\chi(\lambda)d\lambda. $$
	We decompose it as 
	$$\frac{m}{\pi i}\int_0^{\infty} e^{it\lambda^{2m}}\lambda^{2m-1}\Phi(\lambda)\chi(\lambda)\chi(\lambda/|t|^{\frac {1}{2m}})d\lambda+\frac{m}{\pi i}\int_0^{\infty} e^{it\lambda^{2m}}\lambda^{2m-1}\Phi(\lambda)\chi(\lambda)\big(1-\chi(\lambda/|t|^{\frac {1}{2m}})\big)d\lambda. $$
	By a direction calculation on the first integral,  and  using integration by parts $\frac{n}{2}+1$ times on the second integral, we can show that the part of $\Omega^{L}$ related to $\Phi(\lambda)$  contributes the time decay $|t|^{-\frac{4m-n-2j_{\k}}{2m}}(\log|t|)^{-\frac 12}$. Note that from the \eqref{j_k} and \eqref{j_k_1}, we know that when $1\le n\le 2m$,
    $$j_{\k}=k_c+k-2=m-\frac{n}{2}+k-1,$$  and  $j_{\k}=k_c+k-1=k-1$  for $2m<n\le 4m$. Therefore, it follows from the definition of $p_{\k}$ above ( i.e.,\eqref{p_k} )  that 
   $$ \frac{4m-n-2j_{\k}}{2m}=\frac{n}{m}\big(\frac 12-\frac{1}{p_{\k}}\big),$$
which  means that the contribution of the part related to $\Phi(\lambda)$ is  is negligible to the estimate \eqref{eq: time decay lower bound} when $t$ is large due to the additional log-decay factor $(\log|t|)^{-\frac 12}$ .
    
Now it remains to consider the part of $\Omega^{L}$ related to the principal terms 
	$$b_0^2\, \lambda^{2m-n-2j_{\k}} \, G^0_{2m-n}vQ_{j_{\k}}\left( M^{+}_{j_{\k},j_{\k}}-M^{-}_{j_{\k},j_{\k}}\right) Q_{j_{\k}}vG^0_{2m-n},\ \  1\le n\le 2m,$$
    and
    $$b_0^2\, \lambda^{2m-n-2j_{\k}} \, G_{2m-n}vQ_{j_{\k}}\left( M^{+}_{j_{\k},j_{\k}}-M^{-}_{j_{\k},j_{\k}}\right) Q_{j_{\k}}vG_{2m-n},\ \  2m< n\le 4m. $$
	Because  of similarity, we only deal with  the case $1\le n\le 2m$. This corresponding  part can be written as 
	$$
	I_{M}(t)= \frac{m b_0^2}{\pi i}\int_0^{\infty} e^{it\lambda^{2m}} \lambda^{4m-n-2j_{\k}-1} \, \chi(\lambda) d\lambda \cdot\Big( G^0_{2m-n}vQ_{j_{\k}}\left( M^{+}_{j_{\k},j_{\k}}-M^{-}_{j_{\k},j_{\k}}\right)  Q_{j_{\k}}vG^0_{2m-n}\Big)
	$$
	To prove the desired estimate \eqref{eq: time decay lower bound}, it suffices to show that $I_{M}(t)$ satisfies \eqref{eq: time decay lower bound}. 
 Note that from Lemma \ref{lem: est for osci int} below,  we have 
	$$\Big|\int_0^{\infty} e^{it\lambda^{2m}} \lambda^{4m-n-2j_{\k}-1} \, \chi(\lambda) d\lambda\Big| \sim |t|^{\frac{4m-n-2j_{\k}}{2m}}=|t|^{-\frac{n}{m}\big(\frac 12-\frac{1}{p_{\k}}\big)},$$ 
	as $|t|$ is sufficiently large. Hence it is just enough to show that the principal terms in \eqref{eq: expan for spectral measure-0}  don't vanish in $\mathbb{B}(L^{2}_{\sigma},\, L^{2}_{-\sigma})$. That is, there exist non-zero functions $\phi,\,\psi\in L_{\sigma}^2$ such that 
	$$
	\left\langle b_0^2 G_{2m - n}^0vQ_{j_{\k}}\left( M^{+}_{j_{\k},j_{\k}}-M^{-}_{j_{\k},j_{\k}}\right) Q_{j_{\k}}vG_{2m - n}^0\psi,\, \phi\right\rangle\neq 0.
	$$
	%or 
	%$$
	%\left\langle b_0^2 G_{2m - n}vQ_{j_{\k}}\left( M^{+}_{j_{\k},j_{\k}}-M^{-}_{j_{\k},j_{\k}}\right) Q_{j_{\k}}vG_{2m - n}\psi,\, \phi\right\rangle\neq 0.
	%$$
	Firstly, we have $Q_{j_k}\neq 0$ by Proposition \ref{prop-resonance-equivalence}. Meanwhile,  $Q_{j_{\k}}$ is a finite rank projection on $L^2$ by the definition. Thus, it can be written as 
	$$
	Q_{j_{\k}}f= \sum_{s = 1}^{\tau}\langle f,\, \psi_{j_{\k}}^s\rangle\psi_{j_{\k}}^s,  
	$$ 
	where $\{\psi_{j_{\k}}^s~|~s = 1,\cdots,\tau\}$ is the orthonormal basis of $Q_{j_{\k}}L^2$. Since $\k >k_c$, one has that $Q_{j_{\k}}T_0 = 0$, which implies that $Q_{j_{\k}}(U+b_0vG_{2m - n}^0v)U \psi_{j_{\k}}^s=0$. Hence, we have
	$$
	b_0Q_{j_k}vG_{2m - n}^0 vU \psi_{j_{\k}}^s= -\psi_{j_{\k}}^s.
	$$
	Note that $vU \psi_{j_{\k}}^s \in L^{2}_{\sigma}$ by the decay assumption \eqref{decayonV} on $V$. Let $\psi=vU \psi_{j_{\k}}^s \in L^{2}_{\sigma}$, which is not zero since $\psi_{j_{\k}}^s \neq 0$. then $$
	\left\langle b_0^2 G_{2m - n}^0vQ_{j_{\k}}\left( M^{+}_{j_{\k},j_{\k}}-M^{-}_{j_{\k},j_{\k}}\right) Q_{j_{\k}}vG_{2m - n}^0\psi,\, \psi\right\rangle=\left\langle Q_{j_{\k}}\left( M^{+}_{j_{\k},j_{\k}}-M^{-}_{j_{\k},j_{\k}}\right) Q_{j_{\k}} \psi_{j_{\k}}^s,  \psi_{j_{\k}}^s\right\rangle
	$$ which does not vanish from Proposition \ref{pro: M^pm_{jj} dose not vanish}.
	Therefore, we conclude the \eqref{eq: time decay lower bound}.
	\vskip0.2cm
\textbf{B. \ Proof of the  \eqref{eq: time decay lower bound-1} when $\mathbf{k}=m_n$ or $\mathbf{k}=m_n+1$.}  
        
       1. \underline{The case of $\k=m_n$}.  
       
       When $1\le n \le 2m$, by the expansions in \eqref{eq-M expansion-even} and noting $M_{2m-\frac n2, 2m-\frac n2}^{+}=M_{2m-\frac n2, 2m-\frac n2}^{-}$ in Theorem \ref{thm-M inverse-even}  we have that  
	\begin{equation}\label{eq: expan for spectral measure-2}
		\begin{aligned}
			&R_0^{+}(\lambda^{2m})v\left(M^{+}(\lambda) \right)^{-1}vR_0^{+}(\lambda^{2m}) -R_0^{-}(\lambda^{2m}) v\left(M^{-}(\lambda) \right)^{-1}vR_0^{-}(\lambda^{2m}) \\
			=&\lambda^{2m-n} (\rho_{2m-\frac n2}(\lambda))^2 \,b_0^2\,  G_{2m-n}^0vQ_{2m-\frac n2}\left( \Gamma^{+}_{2m-\frac n2, 2m-\frac n2}-\Gamma^{-}_{2m-\frac n2, 2m-\frac n2}\right)(\lambda)  Q_{2m-\frac n2}vG_{2m-n}^0+{\Phi}(\lambda) \\
			=&\lambda^{-2m} (\log\lambda)^{-2} \,b_0^2\,  G_{2m-n}^0vQ_{2m-\frac n2}\mathscr{M}_{2m-\frac n2, 2m-\frac n2} Q_{2m-\frac n2}vG_{2m-n}^0+\widetilde{\Phi}(\lambda),	
		\end{aligned}
	\end{equation}
	where 
	$$\mathscr{M}_{2m-\frac n2, 2m-\frac n2}=-(b_1)^{-2}\left( a_{2m-\frac n2}^{+}-a_{2m-\frac n2}^{-}\right)  \left(Q_{2m-\frac n2}vG_{4m-n}vQ_{2m-\frac n2}\right)^{-1},$$
	and in the last line, we used  \eqref{eq:no vanish for j=2m-fracn2} and $\rho_{2m-\frac n2}(\lambda)=\lambda^{2m-\frac n2} (\log \lambda)^{-\frac12}$.
	Moreover, by \eqref{eq-expansion-for-Gamma-i,j} and \eqref{eq: est for Phi(lambda)}, the operator $\widetilde{\Phi}(\lambda)$ satisfies 
	\begin{equation}\label{eq: est for tildePhi}
		\left\|\frac{d^l}{d\lambda^{l}}\widetilde{\Phi}(\lambda)\right\|_{\mathbb{B}(L^{2}_{\sigma},\, L^{2}_{-\sigma})}\lesssim |\log \lambda|^{-\frac 52} \lambda^ {-l},
	\end{equation}
	for $l=0,1,\cdots,\frac n2+1$. Since $Q_{2m-\frac n2}\neq 0$ according to Proposition \ref{prop-resonance-equivalence}, we have 
	$$\left(Q_{2m-\frac n2}vG_{4m-n}vQ_{2m-\frac n2}\right)^{-1} \neq 0\  \text{on} \  Q_{2m-\frac n2}L^2.$$ 
	Furthermore, it follows from \eqref{eq： def for a_j^pm} that
	$\left( a_{2m-\frac n2}^{+}-a_{2m-\frac n2}^{-}\right)=-i\, \mu_{j+\frac n2 -1} \neq 0$. Thus, by combining these facts with the self-adjointness of the operator $\left(Q_{2m-\frac n2}vG_{4m-n}vQ_{2m-\frac n2}\right)^{-1}$, we can find a function $\psi_{2m-\frac n2}^1$ in $Q_{2m-\frac n2}L^2$ such that 
	\begin{align*}
		\theta_0=&\big\langle b_0^2\,  G_{2m-n}vQ_{j_{\k}}\mathscr{M}_{2m-\frac n2, 2m-\frac n2} Q_{j_{\k}}vG_{2m-n} \psi_{2m-\frac n2}^1,\,  \psi_{2m-\frac n2}^1\big\rangle \\
		&\quad = -(b_1)^{-2}\left( a_{2m-\frac n2}^{+}-a_{2m-\frac n2}^{-}\right) \Big\langle \left(Q_{2m-\frac n2}vG_{4m-n}vQ_{2m-\frac n2}\right)^{-1}\psi_{2m-\frac n2}^1,\,  \psi_{2m-\frac n2}^1\Big\rangle \neq 0. 
	\end{align*} 
	Therefor, when $\lambda$ is sufficiently small, we have 
	\begin{equation}\label{eq: the asym behavior of spectral measure}
		\begin{aligned}
			\big\langle  &\left((\log\lambda)^{-2}b_0^2\,  G_{2m-n}vQ_{2m-\frac n2}\mathscr{M}_{2m-\frac n2, 2m-\frac n2} Q_{2m-\frac n2}vG_{2m-n}+\widetilde{\Phi}(\lambda)\right)  \psi_{2m-\frac n2}^1,\,  \psi_{2m-\frac n2}^1\big\rangle \\
			&\qquad \qquad \qquad\qquad=\theta_0(\log\lambda)^{-2}+ \big\langle \widetilde{\Phi}(\lambda) \psi_{2m-\frac n2}^1,\,  \psi_{2m-\frac n2}^1\big\rangle =:E_{M}(\lambda).
		\end{aligned}
	\end{equation}
	Substituting this term into the expression of $\Omega^{L}$, we consider the integral
	\begin{equation}\label{eq: def for I_M}
		I_{M}(t)=\frac{m}{\pi i}\int_0^{\infty} e^{it\lambda^{2m}}\lambda^{-1} E_{M}(\lambda)\chi(\lambda)d\lambda=\frac{1}{2\pi i}\int_0^{\infty} e^{it\lambda}\lambda^{-1} E_{M}(\lambda^{\frac{1}{2m}})\chi(\lambda^{\frac{1}{2m}})d\lambda.
	\end{equation}
	To prove the inequality \eqref{eq: time decay lower bound-1}.
	it suffices to prove that $I_{M}(t)$ still satisfies \eqref{eq: time decay lower bound-1}, which is actually given in \eqref{eq: est for osci int-2} below.
	
	When $2m< n\le 4m$,  replacing $G_{2m-n}^0$ by $G_{2m-n}$, one also has
	\begin{equation}\label{eq: expan for spectral measure-3}
		\begin{aligned}
			R_0^{+}(\lambda^{2m})v&\left(M^{+}(\lambda) \right)^{-1}vR_0^{+}(\lambda^{2m}) -R_0^{-}(\lambda^{2m}) v\left(M^{-}(\lambda) \right)^{-1}vR_0^{-}(\lambda^{2m}) \\
			&=(\log\lambda)^{-2}\lambda^{-2m} \,b_0^2\,  G_{2m-n}vQ_{2m-\frac n2}\mathscr{M}_{2m-\frac n2, 2m-\frac n2} Q_{2m-\frac n2}vG_{2m-n}+\widetilde{\Phi}(\lambda).
		\end{aligned}
	\end{equation}
	And by the same argument, one can also conclude  \eqref{eq: time decay lower bound-1}.
	
	2. \underline{The case of $\k=m_n+1$}. 
    
   Firstly, noting  that there exists a non-zero solution of \eqref{equ1.2} in $\displaystyle\cap_{s<0}L^2_{s}(\R^n)\setminus{L^2}(\R^n)$, we have $Q_{2m-\frac n2}\neq 0$ according to Proposition \ref{prop-resonance-equivalence}.
	Hence when $\k=m_n+1$ and $1\le n\le 2m$, we can also find a function $\psi_{2m-\frac n2}^1$ in $Q_{2m-\frac n2}L^2$ such that 
	\begin{align*}
		\Big\langle R_0^{+}(\lambda^{2m})v&\left(M^{+}(\lambda) \right)^{-1}vR_0^{+}(\lambda^{2m}) -R_0^{-}(\lambda^{2m}) v\left(M^{-}(\lambda)  \right)^{-1}vR_0^{-}(\lambda^{2m})\psi_{2m-\frac n2}^1,\, \psi_{2m-\frac n2}^1\Big\rangle \\
		=&\big\langle  \left((\log\lambda)^{-2}b_0^2\,  G_{2m-n}vQ_{2m-\frac n2}\mathscr{M}_{2m-\frac n2, 2m-\frac n2} Q_{2m-\frac n2}vG_{2m-n}+\widetilde{\Phi}(\lambda)\right)  \psi_{2m-\frac n2}^1,\,  \psi_{2m-\frac n2}^1\big\rangle \\
		\sim& (\log\lambda)^{-2}.
	\end{align*}
	The same proof applied to $\mathbf{k}=m_n$  can be extended to  the case $\mathbf{k}=m_n+1$. Similarly, for $2m < n \le 4m$ and $\mathbf{k}=m_n+1$, \eqref{eq: time decay lower bound-1} also remains valid. Thereby we conclude the whole proof.
	\end{proof}
	
	We need the lower bound estimates of the oscillatory integral. The first one can be seen in Section 5.1 (d) in Chapter 8 of \cite{Stein-1993}.
	\begin{lemma}\label{lem: est for osci int}
	Suppose that $\psi\in C_0^\infty(\R)$ and $\Re(\mu)>-1$. Then, as $x$ tends to $\infty$, we have 
	$$\int_{0}^{\infty} e^{i\lambda x}\psi (\lambda) \lambda^{\mu}d\lambda\sim \sum_{j=0}^{\infty} a_j x^{-j-1-\mu}, \quad a_j=i^{j+\mu+1} \frac{j!}{\Gamma(j+\mu+1)}\psi^{(j)}(0),$$
	where $\Gamma(x)$ denotes the Gamma function.
	\end{lemma}
	\begin{lemma}\label{lem: est for osci int-1}
	The oscillatory integral $I_{M}(t)$ given in \eqref{eq: def for I_M} satisfies that for $|t|\rightarrow\infty$,  
	\begin{equation}\label{eq: est for osci int-2}
		I_{M}(t)\sim \left( \log|t| \right)^{-1}.
	\end{equation}
	\end{lemma}
	\begin{proof}
	Let $\phi_0(\lambda)$ be a smooth function supported in $[-1/10,1/10]$ satisfying that $\phi_0(\lambda)=1$ when $|\lambda|\le 1/20$, and $\phi_1(\lambda)=1-\phi_0(\lambda)$.
	We decompose the integral $I_{M}(t)$ into 
	\begin{equation}\label{eq: decomposition fo I_M}
	\frac{1}{2\pi i}\int_0^{\infty} e^{it\lambda}\lambda^{-1} E_{M}(\lambda^{\frac{1}{2m}})\chi(\lambda^{\frac{1}{2m}})\phi_0(|t|\lambda)d\lambda+\frac{1}{2\pi i}\int_0^{\infty} e^{it\lambda}\lambda^{-1} E_{M}(\lambda^{\frac{1}{2m}})\chi(\lambda^{\frac{1}{2m}})\phi_1(|t|\lambda)d\lambda.
	\end{equation}
    Note that by \eqref{eq: est for tildePhi}, there is 
	$$\left| \frac{d^\ell}{d\lambda^\ell}E_{M}(\lambda^{\frac{1}{2m}})\right|\lesssim |\log\lambda|^{-2}|\lambda|^{-\ell}, \quad \ell=0,1. $$
	When $|t|$ is large, we perform the integration by parts once to the second integral in \eqref{eq: decomposition fo I_M}. Then,  this integral is bounded by 
	$$|t|^{-1}\int_{\frac {1}{20}|t|^{-1}}^{\lambda_0} \lambda^{-2}(\log \lambda)^{-2} d\lambda \lesssim (\log |t|)^{-2}.$$
	For the first integral in \eqref{eq: decomposition fo I_M}, note that 
	$$E_{M}(\lambda^{\frac{1}{2m}})=4m^2\theta_0(\log\lambda)^{-2}+\big\langle \widetilde{\Phi}(\lambda^{\frac {1}{2m}}) \psi_{2m-\frac n2}^1,\,  \psi_{2m-\frac n2}^1\big\rangle.$$
	On the one hand, by \eqref{eq: est for tildePhi}, one has 
	\begin{align*}
		\Big| \frac{1}{2\pi i}\int_0^{\infty} e^{it\lambda}\lambda^{-1} \big\langle \widetilde{\Phi}(\lambda^{\frac {1}{2m}}) \psi_{2m-\frac n2}^1,\,  \psi_{2m-\frac n2}^1\big\rangle & \chi(\lambda^{\frac{1}{2m}})\phi_0(|t|\lambda)d\lambda\Big| \\
		&\lesssim \int_0^{\frac{1}{10}|t|^{-1}}\lambda^{-1}(\log\lambda)^{-\frac 52} d\lambda \sim (\log|t|)^{-\frac 32}. 
	\end{align*}
	On the other hand, since $\theta_0\neq 0$, we have 
	\begin{align*}
		\left| \frac{4m^2\theta_0}{2\pi i}\int_0^{\infty} e^{it\lambda}\lambda^{-1}  (\log\lambda)^{-2} \chi(\lambda^{\frac{1}{2m}})\phi_0(|t|\lambda)d\lambda\right| 
		\sim& \left|\int_0^{\infty} e^{it\lambda}\lambda^{-1}  (\log\lambda)^{-2} \chi(\lambda^{\frac{1}{2m}})\phi_0(|t|\lambda)d\lambda\right| \\
		\sim& \int_{0}^{\frac{1}{10}|t|^{-1}} \lambda^{-1}  (\log\lambda)^{-2} d\lambda
		\sim \big|\log |t|\big|^{-1},
	\end{align*}
    where in the second equivalence, we used the fact that the real part 
     $$\Re \left(e^{it\lambda}\lambda^{-1}  (\log\lambda)^{-2} \chi(\lambda^{\frac{1}{2m}})\right) \gtrsim \lambda^{-1}  (\log\lambda)^{-2},$$ 
    for $\lambda |t|<1/10.$ At now, we have concluded the \eqref{eq: est for osci int-2}.
    \end{proof}

We now prove Theorem \ref{thm-unbound}.
\begin{proof}[{\bf Proof of Theorem \ref{thm-unbound}}]
	Assume that when zero is the $\k$-th kind of resonance with $k_c<\k$, the wave operators $W_\pm(H; (-\Delta)^m)$ are bounded on $L^{\bf{q}}$ for some $\bf{q}$ such that $p_{\k} <\bf{q} \le \infty$. Then by the intertwining property in \eqref{eq-intertwing identity} and \cite{ZYF} ,  one has 
	\begin{equation*}
		\Big\|e^{i tH}P_{ac}(H)\Big\|_{L^{{\bf{q}}^{\prime}}-L^{\bf{q}}}\lesssim |t|^{-\frac{n}{m}(\frac 12-\frac{1}{\bf{q}})}, \quad t\neq0.
	\end{equation*}
	Through H\"{o}lder inequality, the estimate 
	\begin{equation}\label{eq: disper est in weighted L2}
		\Big\|e^{i tH}P_{ac}(H)\Big\|_{L^{2}_{\sigma}-L^{2}_{-\sigma}}\lesssim |t|^{-\frac{n}{m}(\frac 12-\frac{1}{\bf{q}})}, \quad t\neq0,
	\end{equation}
    holds for any $\sigma> \frac n2+1$. However, since $\frac{n}{m}(\frac 12-\frac{1}{\bf{q}})\ge\frac{n}{m}(\frac 12-\frac{1}{p_{\k}})$ when $k_c<\k<m_n$, the inequality \eqref{eq: disper est in weighted L2} is contradictory to the lower bound \eqref{eq: time decay lower bound} as $|t|\rightarrow\infty$. When $\k=m_n$ or $\k=m_n+1$, noting $p_{\k}=2$ and $\frac{n}{m}(\frac 12-\frac{1}{\bf{q}})>0$,  the inequality \eqref{eq: disper est in weighted L2} contradicts the inequality \eqref{eq: time decay lower bound-1}. From these facts, we conclude that \eqref{eq: disper est in weighted L2} is invalid, which implies that the desired wave operator $W_\pm(H; (-\Delta)^m)$ cannot be bounded $L^{\bf{q}}$ for each $\bf{q}$ such that $p_{\k} <\bf{q} \le \infty$. 
\end{proof}

%\vspace{1.2em}
 \section{Proofs of  Theorems \ref{thm-M inverse-even} and  \ref{thm-M inverse-even-1}}\label{sec-proof-for-sec2.2}

In this section, we prove Theorems \ref{thm-M inverse-even} and \ref{thm-M inverse-even-1} from Section \ref{sec-expansion}, which characterize the asymptotic expansion of $\bigl(M^\pm(\lambda)\bigr)^{-1}$ near zero energy.

We first introduce the symbol class $\mathcal{S}^b(\Omega)$ consisting of all $\mathbb{B}(L^2)$-valued $C^{\frac{n}{2}+1}$ functions $f$ on $\Omega\subset\R$ such that
\begin{equation}\label{eq:symb class}
\bigl\|\tfrac{d^l}{d\lambda^l}f(\lambda)\bigr\|_{\mathbb{B}(L^2)} \lesssim |\lambda|^{b-l}, \qquad l=0,\dots,\tfrac{n}{2}+1,\ \lambda\in\Omega.
\end{equation}
Occasionally, we also use $\mathcal{S}^b(\Omega)$ to denote a function in this class. In the rest of this section, we no longer specify that the dimension $n$
is even.

\subsection{Expansions of the operators $M^\pm(\lambda)$}\,\par
Recall that  
$$
M^\pm(\lambda)=U+vR_0^\pm(\lambda^{2m})v,
$$
where 
$$v(x) = \sqrt{|V(x)|}\quad  \text{and}\quad  U (x) = \sgn V(x).$$
Suppose that zero energy be the $\mathbf{k}$-th kind resonance of $H=(-\Delta)^m+V$. 
Using the expansions \eqref{eq-free expansion-even-1}--\eqref{eq-free expansion-even-5} of $R_0^{\pm}(\lambda^{2m})(x,y)$,
along with the decay assumption \eqref{decayonV} and the Hilbert–Schmidt property \footnote{For $n = 4m$, the operator $v(-\Delta)^{-m}v$ is $L^2$-bounded by Riesz potential theory, cf.\ \cite[Lemma 2.3]{J80}, and thus the Hilbert–Schmidt property is not needed in this case.}, for $\lambda\in(0,1/2)$,
we derive the following asymptotic expansion for the Kato–Birman operators $M^\pm(\lambda)$ on $\mathbb{B}(L^2)$.

\begin{itemize}
	\item[(i)] Let $n\le 2m$. When $0\le \mathbf{k}\le m_n-1$, we have
	\begin{equation}\label{eq-expansion-Mpm-1}
		\begin{aligned}
			M^\pm(\lambda)&=\sum_{j=0}^ {\lfloor\frac{\theta-1}{2}\rfloor} a_j^\pm \lambda^{n-2m+2j} v G_{2j} v 
			 + b_0\log\lambda \, v G_{2m-n} v + T_0 + \lambda^{n-2m} v r_{\theta}^\pm(\lambda) v,
		\end{aligned}
	\end{equation}
	where
	\begin{equation*}
		\theta=
		\begin{cases}
			2m-n+1, & \text{for } \mathbf{k}=0,\\[4pt]
			2m-n+2\mathbf{k}-1, & \text{for } 1\le \mathbf{k}\le m_n-1=m.
		\end{cases}
	\end{equation*}
    And when $\mathbf{k}=m_n$ or $m_n+1$, we set $\theta=4m-n+1$ and have
    \begin{equation}\label{eq-expansion-Mpm-3}
		\begin{aligned}
		M^\pm(\lambda)=&\sum_{j=0}^{2m-\frac n2} a_j^\pm \lambda^{n-2m+2j} v  G_{2j} v + b_0\log\lambda \, v G_{2m-n} v
		+ T_0 \\
         &\quad+ b_1 \lambda^{2m}\log\lambda \, v G_{4m-n} v+ b_1 \lambda^{2m} v G^0_{4m-n} v 
		   + \lambda^{n-2m} v r_{\theta}^\pm(\lambda) v.
		\end{aligned}
	\end{equation}

	\item[(ii)] Let $2m<n\le 4m$. When $0\le \mathbf{k}\le m_n-1$, we have
	\begin{equation}\label{eq-expansion-Mpm-2}
		\begin{aligned}
	M^\pm(\lambda)&=T_0+\sum_{j=0}^ {\lfloor\frac{\theta-1}{2}\rfloor} a_j^\pm \lambda^{n-2m+2j} v G_{2j} v + \lambda^{n-2m} v r_{\theta}^\pm(\lambda) v,
		\end{aligned}
	\end{equation}
	where
	\begin{equation*}
		\theta=
		\begin{cases}
			0, & \text{for } \mathbf{k}=0,\\[4pt]
			2\mathbf{k}-1, & \text{for } 1\le \mathbf{k}\le m_n-1=2m-\frac n2.
		\end{cases}
	\end{equation*} 
    When $\mathbf{k}=m_n$ or $m_n+1$, we set $\theta=4m-n+1$ and have
	\begin{equation}\label{eq-expansion-Mpm-4}
		\begin{aligned}
			M^\pm(\lambda)=&T_0+\sum_{j=0}^{2m-\frac n2} a_j^\pm \lambda^{n-2m+2j} v G_{2j} v 
			+ b_1 \lambda^{2m}\log\lambda \, v G_{4m-n} v  \\
			&\quad+ b_1 \lambda^{2m} v G^0_{4m-n} v + \lambda^{n-2m} v r_{\theta}^\pm(\lambda) v.
		\end{aligned}
	\end{equation}
\end{itemize}
In the above expansions \eqref{eq-expansion-Mpm-1}--\eqref{eq-expansion-Mpm-4}, the operators $r_{\theta}^\pm(\lambda)$ are of the integral kernels $r_{\theta}^\pm(\lambda|x-y|)$ and satisfy
\begin{equation}\label{eq-est-for-rlam}
r_{\theta}^\pm(\lambda)\in \mathcal{S}^\theta\bigl((0,1/2)\bigr).
\end{equation}
We emphasize that the decay assumption on $V$ in \eqref{decayonV} is primarily needed to ensure the $L^2$-boundedness of all terms in the expansion of $M^{\pm}(\lambda)$.  

For convenience, we denote
\begin{equation}\label{eq-def-for-g0-g1}
g_{0}^\pm(\lambda)=b_0\log\lambda + a^\pm_{2m-n},
\end{equation}
and
\begin{equation}\label{eq-def-for-g0-g1-1}
g_{1}^\pm(\lambda)=b_1\log\lambda + a^\pm_{4m-n}.
\end{equation}

\vskip0.5cm
\subsection{Transformation to the block operator matrix }\,\par
By Proposition~\ref{prop-resonance-equivalence}, whenever zero is a resonance of the $\mathbf{k}$-th kind, the family $\{Q_j L^2\}_{j\in J_{\mathbf{k}}}$ furnishes a direct sum decomposition of $L^2$.
For fixed $\lambda\in(0,1/2)$, we define an isomorphism $B_\lambda$ by
\[
B_\lambda := \bigl( \rho_j(\lambda) Q_j \bigr)_{j\in J_{\mathbf{k}}} : 
\bigoplus_{j\in J_{\mathbf{k}}} Q_j L^2 \ni (f_j)_{j\in J_{\mathbf{k}}}
\longmapsto \sum_{j\in J_{\mathbf{k}}} \rho_j(\lambda)f_j \in L^2,
\]
where $\bigoplus_{j\in J_{\mathbf{k}}} Q_j L^2$ denotes the corresponding vector-valued Hilbert space, and the coefficients $\rho_j(\lambda)$ for $j\in J_{\mathbf{k}}$ are given in \eqref{eq-def-for-rhoj-1} and \eqref{eq-def-for-rhoj-2}.
Let $B_{\lambda}^T$ denote the transpose of $B_{\lambda}$. Since $B_{\lambda}$ is bijective, it follows from \cite[Lemma 3.12]{JK} that if $B_{\lambda}^T M^\pm(\lambda) B_{\lambda}$ is invertible on $\bigoplus_{j\in J_{\mathbf{k}}} Q_j L^2$, then
\begin{equation}\label{eq:A to M inverse}
\bigl( M^\pm(\lambda) \bigr)^{-1} = B_{\lambda} \bigl( B_{\lambda}^T M^\pm(\lambda) B_{\lambda} \bigr)^{-1} B_{\lambda}^T.
\end{equation}
Accordingly, to analyze the low-energy asymptotic behavior of $\bigl( M^\pm(\lambda) \bigr)^{-1}$ on $L^2$, it suffices to establish the invertibility of the block operator matrix $B_{\lambda}^T M^\pm(\lambda) B_{\lambda}$ and study the asymptotic properties of its inverse on $\bigoplus_{j\in J_{\mathbf{k}}} Q_j L^2$.

% We first prove Theorem  \ref{thm-M inverse-even}, and the proof of Theorem  \ref{thm-M inverse-even-1} follows a similar approach. 
Denoted by the index sets as follows:
\begin{equation}\label{eq-def-sets}
\begin{aligned}
\tilde{J}_{0}=&\{j\in J_{m_n+1};\, 0\le j \le m-\frac{n}{2}-1\},\\
\tilde{J}_{1}=&\{j\in J_{m_n+1};\, m-\frac{n}{2}\le j \le\ m-\frac{n}{2}+\frac 34\},\\
\tilde{J}_{2}=&\{j\in J_{m_n+1};\, m-\frac{n}{2}+1\le j \le 2m-\frac{n}{2}-1\},\\
\tilde{J}_{3}=&\{2m-\frac{n}{2},\ 2m-\frac{n}{2}+\frac 12\}.	
\end{aligned}
\end{equation}
We note that by the construction, 
$$
\begin{cases}
 \tilde{J}_0 = \emptyset \quad  &\text{if $n\ge 2m$},  \\
 \tilde{J}_1 = \{m-\frac{n}{2},\, m-\frac{n}{2}+\frac 12,\, m-\frac{n}{2}+\frac 34 \} \quad  &\text{if $n\le 2m$},  \\
 \tilde{J}_1 = \{m-\frac{n}{2}\} \quad  &\text{if $n> 2m$}.  
\end{cases}
$$
From the definitions of $S_j$ in \eqref{eq-S_j-even}, $Q_j$ in \eqref{eq-Q_j-even}, along with \eqref{eq-cancellation-Q_j}, we obtain
\begin{equation}\label{eq:cancell of Q_j-1}
\begin{cases}
Q_iT_0Q_j = 0 & \text{for } i,j \in J_{m_n+1}\setminus\tilde{J}_1  \text{ and } i+j > 2m-n, \\
Q_ivG_{2l}vQ_j = 0 & \text{for } i,j \in J_{m_n+1} \text{ and } \vartheta(i) + \vartheta(j) > 2l \geq 0,
\end{cases}
\end{equation}
where the function $\vartheta(j)$ is given by
\begin{equation}\label{thetaj-1}
\vartheta(j)=\max\{0,\lfloor j+1/2 \rfloor\}.
\end{equation}
 
% For the case where zero is a resonance of the $\mathbf{k}$-th kind with $0 \leq \mathbf{k} \leq m_n + 1$, we only consider indices $i,j \in J_{\mathbf{k}}$. 

Let $n\le 2m$ and $\rho_j(\lambda)$ be given in \eqref{eq-def-for-rhoj-1}. If zero is a resonance of the $\mathbf{k}$-th kind, using the cancellation property \eqref{eq:cancell of Q_j-1}, together with the expansions \eqref{eq-expansion-Mpm-1} and \eqref{eq-expansion-Mpm-3}, we obtain via direct computation the decomposition
\begin{equation}\label{eq: def of the elements of BMB}
\rho_i(\lambda)\rho_j(\lambda)\, Q_i M^\pm(\lambda) Q_j = \lambda^{n-2m}\bigl(d_{i,j}^\pm + r_{i,j}^\pm(\lambda)\bigr) \qquad \text{for all } i,j \in J_{\mathbf{k}},
\end{equation}
where $d_{i,j}^\pm$ and $r_{i,j}^\pm(\lambda)$ are defined  by
\begin{equation}\label{eq-d_ij-k<m}
d_{i,j}^\pm=
\begin{cases}
a_{\frac{i+j}{2}}^\pm Q_i v G_{i+j} v Q_j, 
& \text{if } i,j\in \tilde{J}_{0}\cup\tilde{J}_{2},\ i+j \text{ is even, and } i+j\neq 2m-n,\\[6pt]
g_0^\pm(\lambda) Q_i v G_{2m-n} v Q_j + Q_i T_0 Q_j, 
& \text{if } i,j\in \tilde{J}_{0}\cup\tilde{J}_{2},\ i+j = 2m-n,\\[6pt]
b_0 Q_i G_{2m-n} Q_j, 
& \text{if } i=j=m-\frac{n}{2},\\[6pt]
Q_i T_0 Q_j, 
& \text{if } i=j=m-\frac{n}{2}+\frac{1}{2},\\[6pt]
Q_i T_0 Q_j, 
& \text{if } (i,j)=\bigl(m-\tfrac{n}{2},\, m-\tfrac{n}{2}+\tfrac{3}{4}\bigr),\ \bigl(m-\tfrac{n}{2}+\tfrac{3}{4},\, m-\tfrac{n}{2}\bigr),\\[6pt]
b_{1} Q_i v G_{4m-n} v Q_j, 
& \text{if } i=j=2m-\frac{n}{2},\\[6pt]
b_{1} Q_i v G^0_{4m-n} v Q_j, 
& \text{if } i=j=2m-\frac{n}{2}+\frac{1}{2},\\[6pt]
0, & \text{otherwise},
\end{cases}
\end{equation}
and 
\begin{equation}\label{eq: expansion for r_ij^pm}
r_{i,j}^\pm(\lambda)=(\log\lambda)^{-\frac12}\varTheta_{i,j}^{\pm,0}+(\log\lambda)^{-1}\varTheta_{i,j}^{\pm,1}+(\log\lambda)^{-\frac 32}\varTheta_{i,j}^{\pm,2}+(\log\lambda)^{-2}\varTheta_{i,j}^{\pm,3}+\varTheta_{i,j}^{\pm,r}(\lambda).
\end{equation}
Here, $\varTheta_{i,j}^{\pm,s}$ for $s=0,1,2,3$ are operators in $\mathbb{B}(L^2)$ independent of $\lambda$(some of which may be zero), furthermore,
\begin{equation}\label{eq-est-rijlamb}
 \varTheta_{i,j}^{\pm,r}(\lambda)\in (\log\lambda)^{-\frac 52} \mathcal{S}^0\big((0,1/2)\big).
\end{equation}

Let $2m<n\le 4m$ and $\rho_j(\lambda)$ be given in \eqref{eq-def-for-rhoj-2}.
\begin{equation}\label{eq-d_ij-k-nl2m}
d_{i,j}^\pm=
\begin{cases}
a_{\frac{i+j}{2}}^\pm Q_i v G_{i+j} v Q_j, 
& \text{if } i,j\in \tilde{J}_{2} \text{ and } i+j \text{ is even},\\[6pt]
Q_i T_0 Q_j, 
& \text{if } i=j=m-\frac{n}{2},\\[6pt]
b_{1} Q_i v G_{4m-n} v Q_j, 
& \text{if } i=j=2m-\frac{n}{2},\\[6pt]
b_{1} Q_i v G_{4m-n}^0 v Q_j, 
& \text{if } i=j=2m-\frac{n}{2}+\frac{1}{2},\\[6pt]
0, & \text{otherwise}.
\end{cases}
\end{equation}
Moreover, for $i,j\in J_{m_n+1}$ we have
\begin{equation}\label{eq: expansion for r_ij^pm}
r_{i,j}^\pm(\lambda)=(\log\lambda)^{-\frac12}\varTheta_{i,j}^{\pm,0}+(\log\lambda)^{-1}\varTheta_{i,j}^{\pm,1}+(\log\lambda)^{-\frac 32}\varTheta_{i,j}^{\pm,2}+(\log\lambda)^{-2}\varTheta_{i,j}^{\pm,3}+\varTheta_{i,j}^{\pm,r}(\lambda),
\end{equation}
where $\varTheta_{i,j}^{\pm,s}$ for $s=0,1,2,3$ are operators in $\mathbb{B}(L^2)$ independent of $\lambda$ (some of which may be zero). Furthermore,
\begin{equation}\label{eq-est-rijlamb-0}
 \varTheta_{i,j}^{\pm,r}(\lambda)\in (\log\lambda)^{-\frac 52} \mathcal{S}^0\big((0,1/2)\big).
\end{equation}

In particular, if $2m<n\le 4m$, we have
\[
\Theta_{m-\frac{n}{2},i}^{\pm,s}=\Theta_{i, m-\frac{n}{2}}^{\pm,s}=0,\qquad s=0,1,2,3.
\]
Moreover, when $i\neq 2m-\frac n2$ or $2m<n<4m$, we have 
\begin{equation}\label{eq-est-Theta-ij-2-1}
\Theta_{i, m-\frac{n}{2}}^{\pm,r}(\lambda)~\text{and}~ 
\Theta_{ m-\frac{n}{2}, i}^{\pm,r}(\lambda) \in \mathcal{S}^{\frac n2-m}\bigl((0,1/2)\bigr);
\end{equation}
when $i=2m-\frac n2$ and $n=4m$, we have
\begin{equation}\label{eq-est-Theta-ij-2-2}
\Theta_{i, m-\frac{n}{2}}^{\pm,r}(\lambda)~\text{and}~ 
\Theta_{ m-\frac{n}{2}, i}^{\pm,r}(\lambda) \in (\log\lambda)^{\frac12}\,\mathcal{S}^{\frac n2-m}\bigl((0,1/2)\bigr).
\end{equation}
Additionally,
\begin{equation}\label{eq: expansion for rij when i,j=2m-n2}
\begin{aligned}
r_{2m-\frac n2,\,2m-\frac n2}^\pm(\lambda)
&= (\log \lambda)^{-1} a^{\pm}_{4m-n}\, Q_{2m-\frac n2} v G_{4m-n} v Q_{2m-\frac n2} 
  \\
&\quad+ (\log \lambda)^{-1} b_1\, Q_{2m-\frac n2} v G_{4m-n}^0 v Q_{2m-\frac n2} + \Theta_{2m-\frac n2,\,2m-\frac n2}^{\pm,r}(\lambda),
\end{aligned}
\end{equation}
where $\Theta_{2m-\frac n2,\,2m-\frac n2}^{\pm,r}(\lambda)$ satisfies the improved estimate
\begin{equation}\label{eq: expansion for rij when i,j=2m-n2 -1}
\varTheta_{2m-\frac n2,\,2m-\frac n2}^{\pm,r}(\lambda)\in \mathcal{S}^{1}\bigl((0,1/2)\bigr).
\end{equation}

\begin{remark}
The results in\eqref{eq: expansion for rij when i,j=2m-n2} and \eqref{eq: expansion for rij when i,j=2m-n2 -1} are specifically used to prove the unboundedness of the wave operators.   
\end{remark}

Therefore, by the identity \eqref{eq: def of the elements of BMB}, the operator $B^T M^\pm(\lambda) B$ can be written in the following form:
\begin{equation}\label{eq-B*MB}
B^T M^\pm(\lambda) B = \lambda^{n-2m} \bigl( D_{\mathbf{k}}^\pm + S_{\mathbf{k}}^\pm(\lambda) \bigr)
= \lambda^{n-2m} \Bigl( \bigl( d_{i,j}^\pm \bigr)_{i,j\in J_{\mathbf{k}}} + \bigl( r_{i,j}^\pm(\lambda) \bigr)_{i,j\in J_{\mathbf{k}}} \Bigr),
\end{equation}
where 
\[
D_{\mathbf{k}}^\pm = \bigl( d_{i,j}^\pm \bigr)_{i,j\in J_{\mathbf{k}}}, \qquad 
S_{\mathbf{k}}^\pm(\lambda) = \bigl( r_{i,j}^\pm(\lambda) \bigr)_{i,j\in J_{\mathbf{k}}}
\]
are block operator matrices on $L^2$, with the row and column indices arranged in ascending order. By applying perturbation theory, the problem reduces to studying the invertibility of $D_{\mathbf{k}}^\pm$ and characterizing the structure of $(D_{\mathbf{k}}^\pm)^{-1}$, whenever the inverse exists on $\bigoplus_{j\in J_{\mathbf{k}}} Q_j L^2$.

To further clarify, from \eqref{eq-d_ij-k<m} and \eqref{eq-d_ij-k-nl2m},  when $\mathbf{k}=m_n+1$, the matrices 
\[
D_{m_n+1}^\pm=\big( d_{i,j}^\pm\big)_{i,j\in J_{m_n+1}}
\]
possess the following structure: \\
\noindent when $1\le n\le 2m$,
\begin{equation}\label{eq-structure-D-mn+1-1}
D_{m_n+1}^\pm=
	\begin{pmatrix}
		D_{0,0}^{\pm} & 0 & D_{0,2}^{\pm}&0\\
		0 & D_{1,1}^{\pm} & 0&0\\
	    D_{2,0}^{\pm}  & 0 & D_{2,2}^{\pm}&0\\
	     0 & 0 & 0&D_{3,3}^{\pm}\\
	\end{pmatrix},
\end{equation}
when $2m<n\le 4m$, 
\begin{equation}\label{eq-structure-D-mn+1-2}
	D_{m_n+1}^\pm=
	\begin{pmatrix}
	     d^\pm_{m-\frac{n}{2},m-\frac{n}{2}} & 0&0\\
		 0 & D_{2,2}^{\pm}&0\\
	     0 & 0&D_{3,3}^{\pm}\\
	\end{pmatrix},
\end{equation}
where the submatrices $D_{l,h}^\pm$ are  
\begin{equation}\label{eq-structure-submatrices}
D_{l,h}^\pm=\Big( d_{i,j}^\pm\Big)_{i\in \tilde{J}_{l},\ j\in \tilde{J}_{h}}, \quad l,\, h=0,1,2,3,
\end{equation}
and $D_{1,1}^\pm$ in dimensions $1\le n\le 2m$ are given by
\begin{equation}\label{eq-def-D11}
D_{1,1}^\pm=
\begin{pmatrix}
d^\pm_{m-\frac{n}{2},m-\frac{n}{2}} & 0 & d^\pm_{m-\frac{n}{2},m-\frac{n}{2}+\frac{3}{4}} \\
0 & d^\pm_{m-\frac{n}{2}+\frac{1}{2},m-\frac{n}{2}+\frac{1}{2}} & 0 \\
d^\pm_{m-\frac{n}{2}+\frac{3}{4},m-\frac{n}{2}} & 0 & 0 \\ 	
\end{pmatrix},
\end{equation} 
$D_{3,3}^\pm$  are given by
\begin{equation}\label{eq-def-D33}
	D_{3,3}^\pm=
	\begin{pmatrix}
		d^\pm_{2m-\frac{n}{2},2m-\frac{n}{2}} & 0  \\
		0 & d^\pm_{2m-\frac{n}{2}+\frac{1}{2},2m-\frac{n}{2}+\frac{1}{2}}  \\
	\end{pmatrix}.
\end{equation} 
It should be emphasized  that the matrix $D_{\k}^\pm$ is the submatrix of $D_{m_n+1}^\pm$ with the row and 
column indices restricted to $J_\mathbf{k}\subset J_{m_n+1}$. Therefore, once the analysis for $D_{m_n+1}^\pm$ is conducted, the analysis  for $D_{\k}^\pm$ with $\k=0,1,\ldots, m_n$ can be handed in a parallel manner. 

In what follows, we present two key observations that establish the invertibility of $D_{\k}^{\pm}$.  In even dimensions $1\le n\le 2m$, we decompose $Q_j$ for $j\in \tilde{J}_2$ as 
\begin{equation}\label{eq-decom-for-Qj}
	Q_j=Q_j^0+Q_j^1,	
\end{equation} 
where 
$$Q_j^0L^2=Q_jL^2\cap \{x^\alpha v(x):\, |\alpha|\le j\}^{\perp}.$$
Let 
$$X_0=\bigoplus_{j\in \tilde{J}_{0}} Q_jL^2, \quad X_2=\bigoplus_{j\in \tilde{J}_{2}} Q_jL^2,$$
and 
$$ X_2^1=\bigoplus_{j\in \tilde{J}_{2}} Q^1_jL^2.$$
Now we define the
block matrices  $\textbf{D}_{0,0}=(d_{l,h})_{l,h\in \tilde{J}_{0}}$ and $\textbf{D}_{2,2}=(d_{l,h})_{l,h\in \tilde{J}_{2}}$ as follows:
\begin{equation}\label{eq-D00}
	d_{l,h}=\begin{cases}
		(-\i)^{l+h}(-1)^ha_{\frac{l+h}{2}} Q_lvG_{l+h}vQ_h,\quad&\mbox{if}\,\, l+h\, \,\,\mbox{is even},\\
		0,  &\mbox{if}\,\, l+h\,\, \mbox{is odd},
	\end{cases}
\end{equation}
where  $a_{\frac{l+h}{2}}$ are defined by \eqref{eq:def of aj-1} in Lemma \ref{lem-Aalbet}. We remark that $\textbf{D}_{0,0}$ vanishes in even dimensions $2m<n\le 4m$. Here, we consider these two matrices since the following relationship holds by the definition of $a_j$,
\begin{equation}\label{eq-rotation-D00}
D_{0,0}^{\pm}=U_0^\mp\, \textbf{D}_{0,0}\, U_0^\mp \text{diag}\left((-1)^hQ_h \right)_{h\in\tilde{J}_{0} }, \quad D_{2,2}^{\pm}=U_2^\mp\, \textbf{D}_{2,2}\, U_2^\mp \text{diag}\left((-1)^hQ_h \right)_{h\in\tilde{J}_{2} },
\end{equation}
where the $U_0^\pm$ and $U_2^\pm$ are diagonal operator-valued matrices defined by 
\begin{equation}\label{eq-def-U0-U2}
\begin{aligned}
 U_0^\pm=&\text{diag}\left(i^l \,  \exp\left( \pm \frac{i(l+\frac n2-m)\pi}{2m}\right)Q_l \right)_{l\in\tilde{J}_{0}}, \\
U_2^\pm=&\text{diag}\left(i^l \,  \exp\left( \pm \frac{i(l+\frac n2-m)\pi}{2m}\right)Q_l \right)_{l\in\tilde{J}_{2}}.   
\end{aligned}
\end{equation}
Here, the notation $\text{diag}\left((-1)^hQ_h \right)_{h\in\tilde{J}_{0}}$ represents a diagonal matrix with entries $(-1)^hQ_h$ for each $h\in\tilde{J}_{0}$, and similarly for $\text{diag}\left((-1)^hQ_h \right)_{h\in\tilde{J}_{2}}$ and the ohters.

We make the following crucial observations about $\mathbf{D}_{0,0}$ and $\mathbf{D}_{2,2}$, which are essential for establishing the invertibility of $D_{\k}^\pm$:

\begin{proposition}\label{pro-proerty-of-D00-D22}
Assume that zero energy is the $m_n+1$-th kind of resonance (i.e. eigenvalue) of $H$, and $V$ satisfies \eqref{decayonV}. Then the block operator matrix $\textbf{D}_{0,0}$ is strictly positive on the space $X_0$, while, $\textbf{D}_{2,2}$ is strictly negative on $ X_2^1$ when $n\le 2m$, and on $ X_2$ when $2m<n\le 4m$ .
\end{proposition}
\begin{proof}
	We firs prove the statements for even dimensions $1\le n\le 2m$, and the proof for the dimension $2m<n\le 4m$ follows by the same arguments.
	It follows from Lemma \ref{lem-Aalbet} and the definitions of $G_{\frac{l+h}{2}}$ in \eqref{eq-G-even} that for any $f\in  L^2$,
	\begin{equation*}
		a_{\frac{l+h}{2}}Q_lvG_{l+h}vQ_hf=\sum_{|\alpha|=l, |\beta|=h} A_{\alpha,\beta}Q_l(x^\alpha v)\int_{\R^n}y^\beta v(y)Q_hf(y)\d y, \,\,\,\mbox{if}\,\,l, h\in \tilde{J}_{0},
	\end{equation*}
and
 \begin{equation*}
 	a_{\frac{l+h}{2}}Q^1_lvG_{l+h}vQ^1_hf=\sum_{|\alpha|=l, |\beta|=h} A_{\alpha,\beta}Q^1_l(x^\alpha v)\int_{\R^n}y^\beta v(y)Q^1_hf(y)\d y, \,\,\,\mbox{if}\,\,l, h\in \tilde{J}_{2}.
 \end{equation*}
	Let $\{u_\alpha\}_{|\alpha|=l}$ be an orthonormal basis of $Q_lL^2$ for $l\in \tilde{J}_{0}$, and $\{u_\alpha\}_{\alpha\in F_l}$ be an orthonormal basis of $Q^1_lL^2$ for $l\in \tilde{J}_{2}$, where $F_l$ is a set of indices associated with $Q^1_lL^2$. For $|\alpha|, |\beta|\in \tilde{J}_{0}$ ,  we denote
\begin{equation}\label{eq:Lambda-1}
 \begin{aligned}
		\tilde{A}_{\alpha,\beta}&=  \langle d_{|\alpha|, |\beta|}\ u_{\alpha}, u_{\beta}\rangle=(-\i)^{|\alpha|+|\beta|}(-1)^{|\beta|}\sum_{|\alpha'|=|\alpha|, |\beta'|=|\beta|} A_{\alpha,\beta} \langle u_{\alpha}, y^{\alpha'}v\rangle \langle x^{\beta'}v, u_{\beta}\rangle.\\
		&=(-\i)^{|\alpha|+|\beta|}(-1)^{|\beta|}\sum_{|\alpha'|=|\alpha|, |\beta'|=|\beta|}\Lambda_{\alpha,\alpha'}A_{\alpha',\beta'}\bar{\Lambda}_{\beta',\beta},
	\end{aligned}   
\end{equation}
where  $\Lambda_{\alpha,\beta}=\langle Q_{|\alpha|}u_{\alpha},\, Q_{|\beta|}(x^\beta v)\rangle$ for $|\alpha|, |\beta|\in \tilde{J}_{0}$. 
Let $ l,\, h\in \tilde{J}_{2}$, $ |\alpha|\in F_l$, $|\beta|\in F_h$. We  define 
\begin{equation}\label{eq:Lambda-2}
\begin{aligned}
	\tilde{A}_{\alpha,\beta}&=  \langle d_{|\alpha|, |\beta|}\ u_{\alpha}, u_{\beta}\rangle=(-\i)^{|\alpha|+|\beta|}(-1)^{|\beta|}\sum_{|\alpha'|=l,\,  |\beta'|=h}\Lambda_{\alpha,\alpha'}A_{\alpha',\beta'}\bar{\Lambda}_{\beta',\beta},
\end{aligned}
\end{equation}
where $\Lambda_{\alpha,\alpha^\prime}=\langle Q^1_{l}u_{\alpha},\, Q^1_{l}(x^{\alpha^\prime} v)\rangle$ and $\Lambda_{\alpha^\prime,\alpha}=\Lambda_{\alpha,\alpha^\prime}$.
	
  Denote $$\tilde{E}_{0}=(\tilde{A}_{\alpha,\beta})_{|\alpha|,|\beta|\in \tilde{J}_{0}},\quad   
	 \Lambda_0=(\Lambda_{\alpha,\beta})_{|\alpha|,|\beta|\in \tilde{J}_{0}}$$ 
and $$\tilde{E}_{2}=(\tilde{A}_{\alpha,\beta})_{|\alpha|,|\beta|\in \tilde{J}_2}, \quad   
	 \Lambda_2=(\Lambda_{\alpha,\beta})_{\alpha\in \cup_{l\in \tilde{J}_{2}}F_l,\, |\beta|\in \tilde{J}_{2}}.$$ 
Thus, we know that the positivity of $\textbf{D}_{0,0}$ and $\tilde{E}_0$ is equivalent,  so is the negativity of $\textbf{D}_{2,2}$ and  $\tilde{E}_2$.  
 
 Now we define two block matrices
 \begin{equation*}
 	E_0=\left((-1)^{|\beta|}(-\i)^{|\alpha|+|\beta|}(A_{\alpha,\beta})_{|\alpha|=l, |\beta|=h}\right)_{l,h \in \tilde{J}_{0}},\qquad E_2=\left((-1)^{|\beta|}(-\i)^{|\alpha|+|\beta|}(A_{\alpha,\beta})_{|\alpha|=l, |\beta|=h}\right)_{l,h \in \tilde{J}_{2}}.
 \end{equation*}
 It follows from the identities \eqref{eq:Lambda-1} and \eqref{eq:Lambda-2} above that
	 \begin{equation*}
	 \tilde{E}_{0}=\Lambda_0{E}_{0}(\overline{\Lambda}_0)^{T}, \quad \tilde{E}_{2}=\Lambda_2{E}_{2}(\overline{\Lambda}_2)^{T},
	 \end{equation*}
 where 
 the $(\overline{\Lambda}_0)^{T}$ and $(\overline{\Lambda}_0)^{T}$ are the conjugate transpose matrices of $\Lambda_0$ and $\Lambda_2$ respectively.  
 Note that by the definition of $Q_l$ for $l\in \tilde{J}_{0}$ and $Q_l^1$ for $l\in \tilde{J}_{2}$, it is not difficult to conclude the block diagonal matrix $\Lambda_0$ and  $\Lambda_2$ is nonsingular as well as $\Lambda_2$ is of row full rank.  Therefore, the positivity of $\textbf{D}_{0,0}$ and $E_0$ are equivalent,  and $\textbf{D}_{2,2}$ is negative if $E_2$ is negative. Rigorously, if  $\Lambda_0$ is singular, then the row vectors of this matrix are linearly dependent. In other words, there exists a sequence ofnon-zerornumberss $z_{\alpha},\ |\alpha|\in   \tilde{J}_{0}$  such that for each $|\beta|\in  \tilde{J}_{0}$, the following equalities hold: 
 $$0=\sum_{|\alpha|\in\tilde{J}_{0}}  z_{\alpha}\Lambda_{\alpha,\beta}=\Big \langle \sum_{|\alpha|\in  \tilde{J}_{0}} z_{\alpha} u_{\alpha},\, x^\beta v\Big\rangle.$$
However, the identity above combined with the definition of $Q_j$ in \eqref{eq-Q_j-even} implies $\sum_{|\alpha|\in \tilde{J}_0} z_\alpha u_\alpha = 0$. Consequently, $z_\alpha = 0$ for all $\alpha$, yielding a contradiction. The same argument also applies to $\Lambda_2$.

It suffices to prove that $E_0$ is strictly positive definite and $E_2$ is strictly negative definite. Note that the $L^2$ vectors
$$\left\{\frac{\xi^\alpha}{(2\pi)^{\frac{n}{2}}\alpha!(1+|\xi|^{2m})^\frac12}\right\}_{|\alpha|\in \tilde{J}_{0}}$$
are linearly independent, and we also see by \eqref{eq-def-Aalbet-2} that $E_0$ is the Gram matrix with respect to such vectors. This implies that $E_0$ is strictly positive definite.
From equation \eqref{eq-def-Aalbet-2}, we observe that $-E_2$ constitutes a Gram matrix associated with the linearly independent vector set
$$
\left\{\frac{\xi^\alpha}{(2\pi)^{\frac{n}{2}}\alpha!(1+|\xi|^{2m})^{\frac{1}{2}}|\xi|^m}\right\}_{|\alpha|\in \tilde{J}_{2}}.
$$
Therefore $E_2$ is strictly negative, and the proof is completed.
\end{proof}

\begin{remark}\label{rem-for-D22}
Note that $(d_{l,h})_{l,h\in \tilde{J}_{2}\cap \tilde{J}_{\k}}$ is a diagonal submatrix  of $\textbf{D}_{2,2}$.
\end{remark}
%\vskip0.5cm
\subsection{Asymptotic expansion for $(M^\pm(\lambda))^{-1}$}
The following lemma contains the abstract form of Feshbach formula (see \cite[Lemma 2.3]{JN01}) which will be used to study the inverse of the operator matrices $\mathbb{A^\pm}$.

\begin{lemma}\label{lemma-Feshbach formula}
	Let $\mathbb{A^\pm}$ be an operator matrix on $\mathcal{H}=\mathcal{H}_1\oplus\mathcal{H}_2:$
	\begin{equation}
		\mathbb{A^\pm}=\begin{pmatrix}
			a_{11} & a_{12}\\
			a_{21} & a_{22}
		\end{pmatrix},\ \ \ a_{ij}:\ \mathcal{H}_j\rightarrow\mathcal{H}_i,
	\end{equation}
	where $a_{11},a_{22}$ are closed and $a_{12},a_{21}$ are bounded. Assume that $a_{11}$ has a bounded inverse, then $\mathbb{A^\pm}$ has a bounded inverse if and only if the operator
	\begin{equation}
		d:= a_{22}-a_{21}a_{11}^{-1}a_{12}
	\end{equation}
	has a bounded inverse.  Furthermore, if $d$ has a bounded inverse, we then have 
	\begin{equation}\label{eq-inverse-feshbach}
		\left(\mathbb{A^\pm} \right)^{-1}=
		\begin{pmatrix}
			a_{11}^{-1}a_{12}d^{-1}a_{21}a_{11}^{-1}+a_{11}^{-1} & -a_{11}^{-1}a_{12}d^{-1}\\
			-d^{-1}a_{21}a_{11}^{-1} & d^{-1}
		\end{pmatrix}.
	\end{equation}
\end{lemma}

We are now in the position to prove Theorem \ref{thm-M inverse-even} and Theorem \ref{thm-M inverse-even-1}. We first prove Theorem  \ref{thm-M inverse-even}, and the proof of Theorem  \ref{thm-M inverse-even} follows similarly.

\begin{proof}[{\bf Proof of Theorem \ref{thm-M inverse-even}.}]
To prove this theorem, we focus on the case $\k = m_n + 1$. The results for $\k < m_n + 1$ can be obtained from the $\k = m_n + 1$ case by setting $Q_j = 0$ for $j \in J_{m_n+1} \setminus J_{\k}$.

As mentioned before, since $B$ is bijective, by Lemma 3.12 in \cite{JK}, if we have that there exists a small positive $\lambda_0$ such that the block operator matrix $B^TM^\pm(\lambda)B$ is invertible on $\bigoplus_{j\in J_{\k}}Q_jL^2$ when $0<\lambda<\lambda_0$, then it is obtained that 
\begin{equation}\label{eq-inverse-A-to-M}
\left( M^\pm(\lambda)\right)^{-1}=	B(B^TM^\pm(\lambda)B)^{-1}B^T.
\end{equation} 
On the other hand, we have
$$B^TM^\pm(\lambda)B=\lambda^{n-2m}\left( D_\k^\pm+S_\k^\pm(\lambda)\right) =\lambda^{n-2m}\Big(\big( d_{i,j}^\pm\big)_{i,j\in J_{\k}}+\big( r_{i,j}^\pm(\lambda)\big)_{i,j\in J_{\k}}\Big).$$	
Using the argument of perturbation, when zero is the $\k$-th kind of resonance,   to study the invertibility of $B^TM^\pm(\lambda)B$ when $0<\lambda<\lambda_0$ for a sufficiently small $\lambda_0>0$, it suffices to analyze the invertibility  of $D_\k^\pm$ on $\bigoplus_{j\in J_{\k}}Q_jL^2$. 

When $\k=m_n+1$, observe the property of structure listed in \eqref{eq-structure-D-mn+1-1} -- \eqref{eq-def-D33}. It suffices to prove that in even dimensions $1\le n\le 2m,$ 
 \begin{equation*}
 A_0^\pm(\lambda)=
 \begin{pmatrix}
 	D_{0,0}^{\pm} & D_{0,2}^{\pm}\\
 	D_{2,0}^{\pm} & D_{2,2}^{\pm}\\
 \end{pmatrix}, \  \  \  \  \ D_{1,1}^{\pm},\  \  \  \  \ D^\pm_{3,3} \end{equation*}
are invertible on $\bigoplus_{j\in \tilde{J}_{0}\cup\tilde{J}_{2}}Q_jL^2$ , $\bigoplus_{j\in \tilde{J}_{1}}Q_jL^2$ and $\bigoplus_{j\in \tilde{J}_{3}}Q_jL^2$ respectively. 

\underline{Part I: the invertibility of $A_0^\pm(\lambda)$.}\\
\indent
Note that by the definition in \eqref{eq-d_ij-k<m},
\[
d_{i,j}^{\pm}=g_0^\pm(\lambda)Q_ivG_{2m-n}vQ_j+Q_iT_0Q_j \qquad \text{for } i+j=2m-n,\ i,j\in \tilde{J}_0\cup\tilde{J}_2.
\]
Then we can split
\begin{equation}\label{eq-Dk-hat}
\begin{split}
A_0^\pm(\lambda) &= g_0^\pm(\lambda)
\begin{pmatrix}
0 & \mathbf{D}_{0,2}^{s} \\
\mathbf{D}_{2,0}^{s} &  
\end{pmatrix}
+ \begin{pmatrix}
D_{00}^\pm & D_{0,2}^{\pm,r} \\
D_{2,0}^{\pm,r} & D_{2,2}^\pm 
\end{pmatrix} \\
&=: g_0^\pm(\lambda)A_0^{s}+A_0^{\pm,r},
\end{split}
\end{equation}
where 
\begin{equation*}
\begin{split}
\mathbf{D}_{0,2}^{s} &= \bigl( d_{i,j}^s \bigr)_{i\in\tilde{J}_0,\,j\in\tilde{J}_2}, \qquad 
\mathbf{D}_{2,0}^{s} = \bigl( \mathbf{D}_{0,2}^{s} \bigr)^{\!T}, \\
D_{0,2}^{\pm,r} &= \bigl( d_{i,j}^{\pm,r} \bigr)_{i\in\tilde{J}_0,\,j\in\tilde{J}_2}, \qquad 
D_{2,0}^{\pm,r} = \bigl( D_{0,2}^{\pm,r} \bigr)^{\!T}.
\end{split}
\end{equation*}
In the above, the elements are given by
\begin{equation*}
\begin{split}
d_{i,j}^s &=
\begin{cases}
Q_i v G_{2m-n} v Q_j, & i+j=2m-n,\\[2pt]
0, & \text{otherwise},
\end{cases} \\[4pt]
d_{i,j}^{\pm,r} &=
\begin{cases}
Q_i T_0 Q_j, & i+j=2m-n,\\[2pt]
a_{\frac{i+j}{2}}^\pm Q_i v G_{i+j} v Q_j, & \text{if } i+j \text{ is even},\\[2pt]
0, & \text{otherwise}.
\end{cases}
\end{split}
\end{equation*}
Using the relations in \eqref{eq:def of aj-1} and \eqref{eq-def-U0-U2}, we then obtain
\begin{equation}\label{eq: A0}
\begin{aligned}
\mathbb{A}_0^\pm(\lambda) &:=
\operatorname{diag}\bigl((-1)^i Q_i \bigr)_{i\in\tilde{J}_{0}\cup\tilde{J}_{2}} \,
\operatorname{diag}\bigl( U_l^\pm \bigr)_{l=0,2} \cdot A_0^{\pm}(\lambda) \cdot 
\operatorname{diag}\bigl( U_l^\pm \bigr)_{l=0,2} \\
&= g_0^\pm(\lambda) A_0^{s} + A_0^{r},
\end{aligned}
\end{equation}
where 
\[
A_0^{r}= \begin{pmatrix}
\mathbf{D}_{00} & \mathbf{D}_{0,2}^{r} \\
\mathbf{D}_{2,0}^{r} & \mathbf{D}_{2,2}
\end{pmatrix},
\]
and the elements $d_{i,j}^{r}$ of $\mathbf{D}_{0,2}^{r}$ and $\mathbf{D}_{2,0}^{r}$ are given by 
\begin{equation*}
d_{i,j}=
\begin{cases}
(-1)^{m-\frac n2+i} Q_i T_0 Q_j, & i+j=2m-n,\\[2pt]
a_{\frac{i+j}{2}} Q_i v G_{i+j} v Q_j, & i+j \text{ is even and } i+j\neq 2m-n,\\[2pt]
0, & \text{otherwise}.
\end{cases}
\end{equation*}
The blocks $\mathbf{D}_{00}$ and $\mathbf{D}_{2,2}$ are defined in \eqref{eq-D00} and satisfy \eqref{eq-rotation-D00}.
Based on this identity, we reduce the study of $A_0^\pm(\lambda)$ to that of $\mathbb{A}_0^\pm(\lambda)$. The key point of the proof is that $\mathbf{D}_{0,0}$ is strictly positive definite on $X_0=\bigoplus_{j\in \tilde{J}_0}Q_j$ and $\mathbf{D}_{2,2}$ is strictly negative definite on $X_2^1=\bigoplus_{j\in \tilde{J}_2}Q_{j}^1$, which are established in Proposition \ref{pro-proerty-of-D00-D22}.

We next study the invertibility and asymptotic behavior of $\mathbb{A}_0^\pm(\lambda)$. To this end, note that $\mathbf{D}_{0,0}$ is strictly positive definite on $X_0$, hence invertible. Then, by the abstract Feshbach formula in Lemma \ref{lemma-Feshbach formula}, the operator $\mathbb{A}_0^\pm(\lambda)$ is invertible on $X_0\oplus X_2=\bigoplus_{j\in \tilde{J}_0\cup\tilde{J}_2} Q_j L^2$ if and only if the operator
\begin{align*}
\mathbf{d}^\pm = \mathbf{D}_{2,2} - \bigl(g_0^\pm(\lambda)\mathbf{D}_{2,0}^{s}+\mathbf{D}_{2,0}^{r}\bigr) \bigl(\mathbf{D}_{0,0}\bigr)^{-1} \bigl(g_0^\pm(\lambda)\mathbf{D}_{0,2}^{s}+\mathbf{D}_{0,2}^{r}\bigr)
\end{align*}
is invertible on $X_2$.  
Moreover, if $\mathbf{d}^\pm$ is invertible, then 
\begin{equation}\label{eq-inverse-feshbach-1}
\bigl(\mathbb{A}_0^\pm \bigr)^{-1}=
\begin{pmatrix}
a_{11}^{-1}a_{12}(\mathbf{d}^\pm)^{-1}a_{21}a_{11}^{-1}+a_{11}^{-1} & -a_{11}^{-1}a_{12}(\mathbf{d}^\pm)^{-1}\\[4pt]
-(\mathbf{d}^\pm)^{-1}a_{21}a_{11}^{-1} & (\mathbf{d}^\pm)^{-1}
\end{pmatrix},
\end{equation}
where we set $a_{11}=\mathbf{D}_{0,0}$, $a_{12}=g_0^\pm(\lambda)\mathbf{D}_{0,2}^{s}+\mathbf{D}_{0,2}^{r}$ and $a_{21}=g_0^\pm(\lambda)\mathbf{D}_{2,0}^{s}+\mathbf{D}_{2,0}^{r}$.

We now expand $\mathbf{d}^\pm$ as 
\begin{align*}
\mathbf{d}^\pm =& \bigl(g_0^\pm(\lambda)\bigr)^2 \mathbf{D}_{2,0}^{s} \bigl(\mathbf{D}_{0,0}\bigr)^{-1}\mathbf{D}_{0,2}^{s} 
+ g_0^\pm(\lambda)\Bigl( \mathbf{D}_{2,0}^{r}\bigl(\mathbf{D}_{0,0}\bigr)^{-1}\mathbf{D}_{0,2}^{s}
+ \mathbf{D}_{2,0}^{s}\bigl(\mathbf{D}_{0,0}\bigr)^{-1}\mathbf{D}_{0,2}^{r}\Bigr) \\
&+ \mathbf{D}_{2,2} - \mathbf{D}_{2,0}^{r}\bigl(\mathbf{D}_{0,0}\bigr)^{-1}\mathbf{D}_{0,2}^{r}.
\end{align*}

The analysis now proceeds depending on whether $\mathbf{D}_{2,0}^{s} \bigl(\mathbf{D}_{0,0}\bigr)^{-1}\mathbf{D}_{0,2}^{s}$ is invertible. 

Note that $g_0^\pm(\lambda)=b_0\log\lambda+a_{2m-n}^\pm$ with $a_{2m-n}^\pm\in\mathbb{C}\setminus \mathbb{R}$.
When $\mathbf{D}_{2,0}^{s} \bigl(\mathbf{D}_{0,0}\bigr)^{-1}\mathbf{D}_{0,2}^{s}$ is invertible, perturbation theory shows that $\mathbf{d}^\pm$ is invertible and admits the expansion
\begin{equation}\label{eq:expan for dpm-1}
\bigl(\mathbf{d}^\pm(\lambda)\bigr)^{-1}= \bigl(g_0^\pm(\lambda)\bigr)^{-2}\Bigl(\mathbf{D}_{2,0}^{s} \bigl(\mathbf{D}_{0,0}\bigr)^{-1}\mathbf{D}_{0,2}^{s}\Bigr)^{-1}
+ \bigl(g_0^\pm(\lambda)\bigr)^{-3}\mathcal{S}^{0}\bigl((0,\lambda_0)\bigr)
\end{equation}
for a sufficiently small $\lambda_0>0$. In this case, using \eqref{eq-inverse-feshbach-1}, we can expand $\bigl(\mathbb{A}_0^\pm \bigr)^{-1}$ as
\begin{equation}\label{eq:expan of A0pm-1}
\bigl(\mathbb{A}_0^\pm(\lambda) \bigr)^{-1}= \mathbf{U}_0 + \bigl(g_0^\pm(\lambda)\bigr)^{-1}\mathbf{U}_1 + \bigl(g_0^\pm(\lambda)\bigr)^{-2}\mathbf{U}_2 + \bigl(g_0^\pm(\lambda)\bigr)^{-3}\mathcal{S}^0, \quad \lambda\in(0,\lambda_0),
\end{equation}
where $\mathbf{U}_1, \mathbf{U}_2$ are bounded operators on $X_0\oplus X_2$ independent of $\lambda$.

If $\mathbf{D}_{2,0}^{s} \bigl(\mathbf{D}_{0,0}\bigr)^{-1}\mathbf{D}_{0,2}^{s}$ is not invertible, let $\mathcal{P}_1$ be the projection onto its kernel on $X_0$, and set $\mathcal{P}_0=I-\mathcal{P}_1$. Since $\bigl(\mathbf{D}_{0,0}\bigr)^{-1}$ is stictly positive and $\mathbf{D}_{2,0}^{s}$ is the adjoint of $\mathbf{D}_{0,2}^{s}$, we have
\begin{equation}\label{eq:pro pf P1}
\mathbf{D}_{2,0}^{s}\mathcal{P}_1 = \mathbf{D}_{0,2}^{s}\mathcal{P}_1 = 0 .    
\end{equation}
Then, with the isomorphism $\sigma_\lambda = \bigl((g_0^\pm(\lambda))^{-1} \mathcal{P}_0,\; \mathcal{P}_1 \bigr)$, we obtain
\begin{equation}\label{eq-inverse-of-A0-1}
\begin{split}
\sigma_\lambda^{T} \mathbf{d}^\pm \sigma_\lambda
&= \begin{pmatrix}
\mathcal{P}_0\mathbf{D}_{2,0}^{s} \bigl(\mathbf{D}_{0,0}\bigr)^{-1}\mathbf{D}_{0,2}^{s}\mathcal{P}_0 
& \mathcal{P}_0\mathbf{D}_{2,0}^{s}\bigl(\mathbf{D}_{0,0}\bigr)^{-1}\mathbf{D}_{0,2}^{r}\mathcal{P}_1 \\[4pt]
\mathcal{P}_1\mathbf{D}_{2,0}^{r}\bigl(\mathbf{D}_{0,0}\bigr)^{-1}\mathbf{D}_{0,2}^{s}\mathcal{P}_0
& \mathcal{P}_1\Bigl(\mathbf{D}_{2,2} - \mathbf{D}_{2,0}^{r}\bigl(\mathbf{D}_{0,0}\bigr)^{-1}\mathbf{D}_{0,2}^{r}\Bigr)\mathcal{P}_1
\end{pmatrix} \\
&\qquad + \begin{pmatrix}
g_0^\pm(\lambda)^{-1}\mathcal{S}^{0} & g_0^\pm(\lambda)^{-1}\mathcal{S}^{0} \\[4pt]
g_0^\pm(\lambda)^{-1}\mathcal{S}^{0} & 0
\end{pmatrix} \\
&=: \tilde{\mathbf{d}}_0 + \tilde{\mathbf{d}}_r,
\end{split}
\end{equation} 
and 
\begin{equation}\label{eq:inverse of dpm}
\bigl(\mathbf{d}^\pm\bigr)^{-1} = \sigma_\lambda \bigl(\sigma_\lambda^{T} \mathbf{d}^\pm \sigma_\lambda\bigr)^{-1} \sigma_\lambda^{T}
\end{equation}
whenever $\bigl(\sigma_\lambda^{T} \mathbf{d}^\pm \sigma_\lambda\bigr)^{-1}$ exists.
By perturbation theory, the invertibility of $\sigma_\lambda^{T} \mathbf{d}^\pm \sigma_\lambda$ reduces to that of $\tilde{\mathbf{d}}_0$. Moreover, by the abstract Feshbach formula, $\tilde{\mathbf{d}}_0$ is invertible if and only if 
\begin{equation}\label{eq:reduction to an operator}
\begin{aligned}
\mathbf{e} = &\ \mathcal{P}_1\Bigl(\mathbf{D}_{2,2} - \mathbf{D}_{2,0}^{r}\bigl(\mathbf{D}_{0,0}\bigr)^{-1}\mathbf{D}_{0,2}^{r}\Bigr)\mathcal{P}_1 \\
&- \mathcal{P}_1\mathbf{D}_{2,0}^{r}\bigl(\mathbf{D}_{0,0}\bigr)^{-1}\mathbf{D}_{0,2}^{s}\mathcal{P}_0
\Bigl(\mathcal{P}_0\mathbf{D}_{2,0}^{s} \bigl(\mathbf{D}_{0,0}\bigr)^{-1}\mathbf{D}_{0,2}^{s}\mathcal{P}_0\Bigr)^{-1}
\mathcal{P}_0\mathbf{D}_{2,0}^{s}\bigl(\mathbf{D}_{0,0}\bigr)^{-1}\mathbf{D}_{0,2}^{r}\mathcal{P}_1    
\end{aligned}
\end{equation}
is invertible on $\mathcal{P}_1 X_2$.

We next show that the operator $\mathbf{e}_0 := \mathbf{D}_{2,2} - \mathbf{D}_{2,0}^{r}\bigl(\mathbf{D}_{0,0}\bigr)^{-1}\mathbf{D}_{0,2}^{r}$ is strictly negative on $X_2$. To this end, decompose 
\[
Q_j = Q_j^0 + Q_j^1 \qquad (j\in \tilde{J}_2)
\] 
as in \eqref{eq-decom-for-Qj}.
Since $Q_j^0 L^2 \subset \{x^\alpha v : |\alpha|\le j \}^\perp$, we have 
\begin{equation}\label{eq-for-dijr}
Q_i^0 d_{i,j}^{r} Q_j^0 = 0 \quad \text{when } i+j \neq 2m-n,
\end{equation}
and consequently
\[
\operatorname{diag}\bigl(Q_j^0 \bigr)_{j\in \tilde{J}_2} \mathbf{D}_{2,2} = \mathbf{D}_{2,2}\, \operatorname{diag}\bigl(Q_j^0 \bigr)_{j\in \tilde{J}_2}=0.
\]
Now one has
\begin{equation}\label{eq:repre of d}
\mathbf{e}_0 = \operatorname{diag}\bigl(Q_j^1 \bigr)_{j\in \tilde{J}_2} \cdot \mathbf{D}_{2,2} \cdot \operatorname{diag}\bigl(Q_j^1 \bigr)_{j\in \tilde{J}_2}
- \mathbf{D}_{2,0}^{r}\bigl(\mathbf{D}_{0,0}\bigr)^{-1}\mathbf{D}_{0,2}^{r}.   
\end{equation}
By Proposition \ref{pro-proerty-of-D00-D22}, $\mathbf{D}_{2,2}$ is strictly negative on $X_2^1 = \operatorname{diag}(Q_j^1) X_2$. Meanwhile, $\mathbf{D}_{2,0}^{r}\bigl(\mathbf{D}_{0,0}\bigr)^{-1}\mathbf{D}_{0,2}^{r}$ is non‑negative because $\bigl(\mathbf{D}_{0,0}\bigr)^{-1}$ is strictly positive. Hence $\mathbf{e}_0$ is non‑positive. To prove that $\mathbf{e}_0$ is strictly negative, it suffices to show 
\[
\mathbf{e}_0 f = 0 \;\Longrightarrow\; f = 0.
\]
To this end, take $f = (f_{m-\frac{n}{2}+1},\dots,f_{2m-\frac{n}{2}-1})\in X_2$ with $\mathbf{e}_0 f = 0$.
From \eqref{eq:repre of d} we obtain
\[
\operatorname{diag}\bigl(Q_j^1 \bigr)_{j\in \tilde{J}_2} f = 0, \qquad \mathbf{D}_{0,2}^{r} f = 0.
\]
Together with \eqref{eq-for-dijr}, these imply $f_j \in Q_j^0 L^2$.
This fact, combined with $S_{2m-n-j} T_0 Q_j = 0$, yields
\[
f_j \in \{x^\alpha v : |\alpha|\le j\}^\perp \cap \ker\bigl(S_{2m-n-j-1} T_0 S_{m-\frac n2}\bigr) = S_j L^2 
\]
for each  $j\in \tilde{J}_2$. By the definition of $Q_j$, we have $Q_j L^2 \cap S_j L^2 = \{0\}$, so $f_j = 0$ for all $j\in \tilde{J}_2$. Thus $\mathbf{e}_0$ is strictly negative and invertible.

We have shown that $\mathbf{D}_{2,2} - \mathbf{D}_{2,0}^{r}\bigl(\mathbf{D}_{0,0}\bigr)^{-1}\mathbf{D}_{0,2}^{r}$ is strictly negative.
Furthermore, the operator $$\bigl(\mathcal{P}_0\mathbf{D}_{2,0}^{s} \bigl(\mathbf{D}_{0,0}\bigr)^{-1}\mathbf{D}_{0,2}^{s}\mathcal{P}_0\bigr)^{-1}$$ is strictly positive on $\mathcal{P}_0 X_2$. Consequently, the operator $\mathbf{e}$ defined in \eqref{eq:reduction to an operator} is strictly negative and invertible on $\mathcal{P}_1 X_2$, which finally implies that $\tilde{\mathbf{d}}_0$ and hence $\sigma_\lambda^{T} \mathbf{d}^\pm \sigma_\lambda$ are invertible on $X_0\oplus X_2$.

At this point, combining the Feshbach formula and \eqref{eq:inverse of dpm}, we obtain
\begin{equation}\label{eq:expan for dpm-1-2}
\bigl(\mathbf{d}^\pm(\lambda)\bigr)^{-1}= \mathcal{P}_1 \mathbf{e}^{-1} \mathcal{P}_1
+ g_0^\pm(\lambda)^{-1}\bigl(\mathbf{F}_1\mathcal{P}_1+\mathcal{P}_1\mathbf{F}_1^*\bigr)
+ g_0^\pm(\lambda)^{-2}\mathbf{F}_2
+ g_0^\pm(\lambda)^{-3}\mathcal{S}^{0}
\end{equation}
for $\lambda\in (0,\lambda_0)$ with a sufficiently small $\lambda_0$, where $\mathbf{F}_1$ and $\mathbf{F}_2$ are independent of $\lambda$. Using the property \eqref{eq:pro pf P1} and setting $a_{11}=\mathbf{D}_{0,0}$, $a_{12}=g_0^\pm(\lambda)\mathbf{D}_{0,2}^{s}+\mathbf{D}_{0,2}^{r}$, $a_{21}=g_0^\pm(\lambda)\mathbf{D}_{2,0}^{s}+\mathbf{D}_{2,0}^{r}$, we have 
\begin{align*}
a_{11}^{-1}a_{12}\bigl(\mathbf{d}^\pm(\lambda)\bigr)^{-1}a_{21}a_{11}^{-1}+a_{11}^{-1}
&= \mathbf{u}_{0,0} + \bigl(g_0^\pm(\lambda)\bigr)^{-1}\mathbf{u}_{1,0} + \bigl(g_0^\pm(\lambda)\bigr)^{-2}\mathbf{u}_{2,0} + \bigl(g_0^\pm(\lambda)\bigr)^{-3}\mathcal{S}^0, \\ 
-a_{11}^{-1}a_{12}\bigl(\mathbf{d}^\pm(\lambda)\bigr)^{-1}
&= \mathbf{u}_{0,1} + \bigl(g_0^\pm(\lambda)\bigr)^{-1}\mathbf{u}_{1,1} + \bigl(g_0^\pm(\lambda)\bigr)^{-2}\mathbf{u}_{2,1} + \bigl(g_0^\pm(\lambda)\bigr)^{-3}\mathcal{S}^0,
\end{align*}
and
\[
-\bigl(\mathbf{d}^\pm(\lambda)\bigr)^{-1}a_{21}a_{11}^{-1}
= \mathbf{u}_{0,2} + \bigl(g_0^\pm(\lambda)\bigr)^{-1}\mathbf{u}_{1,2} + \bigl(g_0^\pm(\lambda)\bigr)^{-2}\mathbf{u}_{2,2} + \bigl(g_0^\pm(\lambda)\bigr)^{-3}\mathcal{S}^0,
\]
for $\lambda\in(0,\lambda_0)$. 
Based on \eqref{eq-inverse-feshbach-1}, we can expand $\bigl(\mathbb{A}_0^\pm (\lambda)\bigr)^{-1}$ as
\begin{equation}\label{eq:expan for A0pm-2}
\bigl(\mathbb{A}_0^\pm(\lambda) \bigr)^{-1}= \mathbf{U}_0 + \bigl(g_0^\pm(\lambda)\bigr)^{-1}\mathbf{U}_1 + \bigl(g_0^\pm(\lambda)\bigr)^{-2}\mathbf{U}_2 + \bigl(g_0^\pm(\lambda)\bigr)^{-3}\mathcal{S}^0, \quad \lambda\in(0,\lambda_0),
\end{equation}
for a sufficiently small $\lambda_0>0$, where $\mathbf{U}_0,\mathbf{U}_1,\mathbf{U}_2$ are bounded operators on $X_0\oplus X_2$ independent of $\lambda$. Moreover,
\[
\mathbf{U}_0=
\begin{pmatrix}
\mathbf{u}_{0,0} & \mathbf{u}_{0,1}\\
\mathbf{u}_{0,2} & \mathcal{P}_1\mathbf{e}^{-1}\mathcal{P}_1
\end{pmatrix}.
\]

Finally, combining \eqref{eq:expan of A0pm-1} with \eqref{eq:expan for A0pm-2} and using the identity
\begin{equation}\label{eq-the-inverese-of-A0pm}
\bigl(A_0^{\pm}(\lambda)\bigr)^{-1}
= \operatorname{diag}\bigl( U_l^\pm\bigr)_{l=0,2} \,
\bigl(\mathbb{A}_0^\pm(\lambda) \bigr)^{-1} \,
\operatorname{diag}\bigl((-1)^i Q_i \bigr)_{i\in\tilde{J}_{0}\cup\tilde{J}_{2}} \,
\operatorname{diag}\bigl( U_l^\pm\bigr)_{l=0,2},
\end{equation}
we obtain 
\begin{equation}\label{eq-the-inverese-of-A0pm-1}
\bigl(A_0^{\pm}(\lambda)\bigr)^{-1}
= \Pi_0^\pm + \bigl(g_0^\pm(\lambda)\bigr)^{-1}\Pi_1^\pm + \bigl(g_0^\pm(\lambda)\bigr)^{-2}\Pi_2^\pm + \bigl(g_0^\pm(\lambda)\bigr)^{-3}\Pi_r^\pm(\lambda),
\end{equation}
for $\lambda\in (0,\lambda_0)$, where $\Pi_0^\pm,\Pi_1^\pm,\Pi_2^\pm$ are bounded operators on $X_0\oplus X_2$ independent of $\lambda$, and $\Pi_r^\pm(\lambda) \in \mathcal{S}^0\bigl((0,\lambda_0)\bigr)$. Moreover,
one has 
\[
\Pi_k^\pm=\operatorname{diag}\bigl(U_l^\pm\bigr)_{l=0,2} \,\mathbf{U}_k
 \,
\operatorname{diag}\bigl((-1)^i Q_i \bigr)_{i\in\tilde{J}_{0}\cup\tilde{J}_{2}} \,
\operatorname{diag}\bigl( U_l^\pm\bigr)_{l=0,2}.
\]

\underline{Part II: the invertibility of $D_{1,1}^\pm$.}\\
\indent 
Recall that 
\begin{equation*}
D_{1,1}^\pm=
\begin{pmatrix}
d^\pm_{m-\frac{n}{2},m-\frac{n}{2}} & 0 & d^\pm_{m-\frac{n}{2},m-\frac{n}{2}+\frac{3}{4}} \\
0 & d^\pm_{m-\frac{n}{2}+\frac{1}{2},m-\frac{n}{2}+\frac{1}{2}} & 0 \\
d^\pm_{m-\frac{n}{2}+\frac{3}{4},m-\frac{n}{2}} & 0 & 0 \\ 	
\end{pmatrix},
\end{equation*}
as in \eqref{eq-def-D11}.

First note that $d^\pm_{m-\frac{n}{2}+\frac{1}{2},m-\frac{n}{2}+\frac{1}{2}}=Q_{m-\frac{n}{2}+\frac{1}{2}}T_0Q_{m-\frac{n}{2}+\frac{1}{2}}$ is invertible  on $Q_{m-\frac{n}{2}+\frac{1}{2}}L^2$ by the definition of $Q_{m-\frac{n}{2}+\frac{1}{2}}$ and the Riesz-Fredholm theory. Note that for any $\psi \in Q_{m-\frac{n}{2}}L^2$, we have
\begin{align*}
\langle vG_{2m-n}v\psi, \psi \rangle &= \int_{\mathbb{R}^{2n}} |x-y|^{2m-n}\psi(x) \overline{\psi}(y) dx dy \\
&= \int_{\mathbb{R}^{2n}} (|x|^2 + |y|^2 - 2xy)^{m-\frac{n}{2}} \psi(x) \overline{\psi}(y) dx dy \\
&= (-1)^{m-\frac{n}{2}} \sum_{0 \leq 2l \le m-\frac{n}{2}} \frac{(m-\frac{n}{2})!}{2 l! (m-\frac{n}{2}-2l)!} \int_{\mathbb{R}^{2n}} |x|^{2l} (2xy)^{m-\frac{n}{2}-2l} |y|^{2l} \psi(x) \overline{\psi}(y) dx dy \\
&= (-1)^{m-\frac{n}{2}} \sum_{ 0 \leq 2l  \leq m-\frac{n}{2}} \sum_{|\alpha| \leq m-\frac{n}{2}-2l} \frac{(m-\frac{n}{2})! 2^{m-\frac{n}{2}-2l}}{2 l! \alpha!} \int_{\mathbb{R}^{2n}} |x|^{2l} x^\alpha y^\alpha |y|^{2l} \psi(x) \overline{\psi}(y) dx dy,
\end{align*}
where $\alpha! = \alpha_1! \alpha_2! \cdots \alpha_n!$, and we have used the property $Q_{m-\frac{n}{2}}(x^\beta v) = 0$ for $|\beta| \leq m-\frac{n}{2}-1$. 
This identity shows that $Q_{m-\frac{n}{2}}vG_{2m-n}vQ_{m-\frac{n}{2}}$ is strictly positive on $Q_{m-\frac{n}{2}}L^2$ when $m-\frac{n}{2}$ is even, and strictly negative when $m-\frac{n}{2}$ is odd. Consequently, $d^\pm_{m-\frac{n}{2},m-\frac{n}{2}} = b_0 Q_{m-\frac{n}{2}}vG_{2m-n}vQ_{m-\frac{n}{2}}$ is invertible on $Q_{m-\frac{n}{2}}L^2$ since $b_0 \neq 0$.

Then by the Feshbach formula, we know that $D_{1,1}^\pm$ is invertible if and only if 
\begin{equation*}
	Q_{m-\frac{n}{2}+\frac{3}{4}}T_0Q_{m-\frac{n}{2}}\left(b_0 Q_{m-\frac{n}{2}}vG_{2m-n}vQ_{m-\frac{n}{2}}\right) ^{-1}Q_{m-\frac{n}{2}}T_0Q_{m-\frac{n}{2}+\frac{3}{4}}
\end{equation*}
is invertible on $Q_{m-\frac{n}{2}+\frac{3}{4}}L^2$. Since $Q_{m-\frac{n}{2}}vG_{2m-n}vQ_{m-\frac{n}{2}}$ is strictly positive or negative definite on $Q_{m-\frac{n}{2}}L^2$, the problem reduces to proving that for any $\psi \in L^2$,
\begin{equation*}
Q_{m-\frac{n}{2}}T_0Q_{m-\frac{n}{2}+\frac{3}{4}}\psi = 0 \implies \psi = 0.
\end{equation*}
This implication follows immediately from the definitions 
\begin{equation*}
Q_{m-\frac{n}{2}+\frac{3}{4}} = S_{m-\frac{n}{2}+\frac{1}{2}} - S_{m-\frac{n}{2}+\frac{3}{4}}, \quad S_{m-\frac{n}{2}+\frac{3}{4}}L^2=\Ker(Q_{m-\frac{n}{2}} T_0)^{\perp}\cap S_{m-\frac{n}{2}+\frac{1}{2}}L^2. 
\end{equation*}

\underline{Part III: the invertibility of $D_{3,3}^\pm$.}\\
\indent On the one hand, similar to the analysis of $d^\pm_{m-\frac{n}{2},m-\frac{n}{2}}$, the operator 
$$d^\pm_{2m-\frac{n}{2},2m-\frac{n}{2}}=b_1Q_{2m-\frac{n}{2}}vG_{4m-n}vQ_{2m-\frac{n}{2}}$$ 
is invertible on $Q_{2m-\frac{n}{2}}L^2$ by the definition of $Q_{2m-\frac{n}{2}}$. On the other hand, 
$$d^\pm_{2m-\frac{n}{2}+\frac{1}{2},2m-\frac{n}{2}+\frac{1}{2}}=b_1Q_{2m-\frac{n}{2}+\frac12}vG_{4m-n}^0vQ_{2m-\frac{n}{2}+\frac12}$$
is invertible since 
\begin{equation*}
\left\langle d^\pm_{2m-\frac{n}{2},2m-\frac{n}{2}}f,\, f \right\rangle=\int_{\R^n} |(-\Delta)^{-m}f(x)|^2dx, \quad f\in Q_{2m-\frac{n}{2}+\frac12}L^2.  	
\end{equation*} 
These two statements show that $D_{3,3}^\pm$ is invertible on $\bigoplus_{j\in \tilde{J}_3}Q_jL^2$.

To sum up, we have obtained that $D_{m_n+1}^\pm$ is invertible on $\bigoplus_{j\in{J}_{m_n+1}}Q_jL^2$. Thus, there exists a small $\lambda_0$ such that for $0<\lambda<\lambda_0$, $M^\pm(\lambda)$ is invertible on $L^2$.  
	
Next,  we will \underline{provide the asymptotic} \underline{expansion of the inverse of $M^\pm(\lambda)$}  by using \eqref{eq-inverse-A-to-M}. For a sufficiently small $\lambda_0$ and $0<\lambda<\lambda_0$, we denote that 
$$\left( A_0^\pm\right) ^{-1}= \begin{pmatrix}
	\mathcal{A}_{0,0}^\pm & 	\mathcal{A}_{0,2}^\pm \\
	\mathcal{A}_{2,0}^\pm & \mathcal{A}_{2,2}^\pm
\end{pmatrix}
\ +\ 
\begin{pmatrix}
	\Upsilon_{0,0}^\pm & 	\Upsilon_{0,2}^\pm \\
	\Upsilon_{2,0}^\pm & \Upsilon_{2,2}^\pm
\end{pmatrix}. 
$$
where $A_{l,h}^\pm$ is  bounded  from $ \bigoplus_{j\in \tilde{J}_l}Q_jL^2$ to $ \bigoplus_{j\in \tilde{J}_ h}Q_jL^2$ and the second matrix in the right hand satisfy \eqref{eq: expansion for r_ij^pm}. Thus, when $\k=m_n+1$, in even dimensions $1\le n\le 2m$, $(D_{m_n+1}^\pm)^{-1}$ has the  structure as 
\begin{equation}\label{eq-structure-D-mn+1-1-invese}
\begin{aligned}
	(D_{m_n+1}^\pm)^{-1}=&
\begin{pmatrix}
	\mathcal{A}_{0,0}^{\pm} & 0 & \mathcal{A}_{0,2}^{\pm}&0\\
	0 & (D_{1,1}^{\pm})^{-1} & 0&0\\
	\mathcal{A}_{2,0}^{\pm}  & 0 & \mathcal{A}_{2,2}^{\pm}&0\\
	0 & 0 & 0&(D_{3,3}^{\pm})^{-1}\\
\end{pmatrix}
\ +\ 
\begin{pmatrix}
	\Upsilon_{0,0}^{\pm} & 0 & \Upsilon_{0,2}^{\pm}&0\\
	0 & 0 & 0&0\\
	\Upsilon_{2,0}^{\pm}  & 0 & \Upsilon_{2,2}^{\pm}&0\\
	0 & 0 & 0&0\\
\end{pmatrix} \\
=\mathcal{M}^{\pm}+\Upsilon_{M}^{\pm}(\lambda).  
\end{aligned}
\end{equation}
Denote the elements in the first matrix $\mathcal{M}^{\pm}$ in \eqref{eq-structure-D-mn+1-1-invese} by $M_{i,j}^{\pm}$ for $i,j \in J_{m_n+1}$. When $\k=m_n+1$, it follows from \eqref{eq-inverse-A-to-M} and Neumann series that 
\begin{equation}\label{eq：Neumann series to main result} 
\begin{split}
\left( M^\pm(\lambda)\right)^{-1}&= \lambda^{2m-n}B\left(D_{m_n+1}^\pm+S_{m_n+1}^\pm(\lambda)\right)^{-1}B^  *\\
&=\lambda^{2m-n}B\left(\sum_{l=0}^{\infty}(D_{m_n+1}^\pm)^{-1}\left( -S_{m_n+1}^\pm(\lambda)\,(D_{m_n+1}^\pm) ^{-1}\right)^l \right)B^T, 
\end{split}	
\end{equation}	
which, together with \eqref{eq-est-rijlamb} and \eqref{eq-the-inverese-of-A0pm}--\eqref{eq-structure-D-mn+1-1-invese}, implies the results in \eqref{eq-M expansion-even}--\eqref{eq-est-gamma-ij-2}. Specially, \eqref{eq-Mij=0} is derived from the structure of $D_{m_n+1}^\pm$ listed in \eqref{eq-structure-D-mn+1-1} -- \eqref{eq-def-D33}. 
Thus, the proof is completed.
\end{proof}
%\vspace{0.4 em}
\begin{proof}[{\bf Proof of Theorem \ref{thm-M inverse-even-1}.}]
We only sketch the proof in the even dimensions $2m <n\le 4m$. As before, because of \eqref{eq-B*MB}, the key is to prove the the invertibility  of $D_{m_n+1}^\pm$ ( or $D_{\k}^\pm$ ) on $\bigoplus_{j\in J_{m_n+1}}Q_jL^2$ ( or $\bigoplus_{j\in J_{\k}}Q_jL^2$). By observing the  structure listed in \eqref{eq-structure-D-mn+1-2} and \eqref{eq-def-D33}, it suffices to prove that $A_0^\pm=D_{2,2}^{\pm}$ and $D_{3,3}$  
are invertible on 	$\bigoplus_{j\in \tilde{J}_{2}}Q_jL^2$ and 	$\bigoplus_{j\in \tilde{J}_{3}}Q_jL^2$ respectively. 

We know that 
\begin{equation}\label{eq: D22 for n>2m}
	\textbf{D}_{2,2}=\text{diag}\left((-1)^iQ_i \right)_{i\in\tilde{J}_{2} }\, U_2^\pm D_{2,2}^\pm U_2^\pm.
\end{equation}
It follows from  Proposition \ref{pro-proerty-of-D00-D22} that $\textbf{D}_{2,2}$ is strictly negative definite on $\bigoplus_{j\in \tilde{J}_2}Q_j L^2$, and thus $D_{2,2}^\pm$ is invertible. Meanwhile, by using the same analysis as in even dimensions  $n\le 2m$, $D_{3,3}^\pm$ is invertible.

Thus, in even dimensions $2m<n\le 4m$, for a sufficiently small $\lambda_0$ and $0<\lambda<\lambda_0$, we have  
\begin{equation}\label{eq-structure-D-mn+1--invrse-2}
	\left( D_{m_n+1}^\pm\right)^{-1}=
	\begin{pmatrix}
		\left( d^\pm_{m-\frac{n}{2},m-\frac{n}{2}}\right) ^{-1} & 0&0\\
		0 & \left( D_{2,2}^{\pm}\right) ^{-1}&0\\
		0 & 0&\left( D_{3,3}^{\pm}\right) ^{-1}\\
	\end{pmatrix}.
\end{equation}
 Denote the elements in $\left( D_{m_n+1}^\pm\right)^{-1}$ by $M_{i,j}^{\pm}$ for $i,j \in J_{m_n+1}$. Then, using \eqref{eq-inverse-A-to-M} and the Neumann series expansion,  the results in this theorem are deduced. Specifically, \eqref{eq-est-gamma-ij-2-1} and \eqref{eq-est-gamma-ij-2-2} are derived from \eqref{eq-est-Theta-ij-2-1}, \eqref{eq-est-Theta-ij-2-2} and the structure of $\left( D_{m_n+1}^\pm\right)^{-1}$ in \eqref{eq-structure-D-mn+1--invrse-2}.
\end{proof}

%\vskip0.5cm
\subsection{Two propositions for the unboundedness of wave operators}
Finally, we present two propositions that serve to prove the unboundedness of the wave operators. These two results play key roles in studying the asymptotic behavior of, especially, the lower bound of $R^{+}(\lambda^2)-R^{-}(\lambda^2)$ as $\lambda$ tends to $0$. The proof of these two propositions hinges upon more details from the proofs of Theorem \ref{thm-M inverse-even} and Theorem \ref{thm-M inverse-even-1}.
\begin{proposition}\label{pro: M^pm_{jj} dose not vanish}
The operator $Q_j(M^{+}_{j,j}-M^{-}_{j,j})Q_j$ doesn't vanish on $Q_jL^2$ if $Q_j\neq 0$ for $j\in J_{\k}$ and $2m-n+1\le j< 2m-\frac n2$. Specifically, for any non-zero $\psi\in Q_jL^2$, we have $\langle Q_j(M^{+}_{j,j}-M^{+}_{j,j})Q_j\psi,\, \psi \rangle \neq 0$.
\end{proposition}
\begin{proof}
 When $2m < n \leq 4m$, for $2m - n + 1 \leq j \leq 2m - \frac{n}{2} - 1$, the operator $M^{\pm}_{j,j}$ lies on the diagonal of $\left( D_{2,2}^{\pm}\right) ^{-1}$. Moreover, according to \eqref{eq-rotation-D00}, we have
 $$
 \left(D_{2,2}^{\pm}\right) ^{-1}=  U_2^\pm \, \mathbf{D}_{2,2}\, U_2^\pm  \text{diag}\left((-1)^iQ_i \right)_{i\in\tilde{J}_{2}}.
 $$
By the definition of $U_2^\pm$ given in \eqref{eq-def-U0-U2}, we obtain
 \begin{equation}\label{eq: Mjj^+-Mjj^-}
 Q_j(M^{+}_{j,j}-M^{-}_{j,j})Q_j=\left(e^{i\frac{(j+\frac n2-m)\pi}{2m}}-e^{- i\frac{(j+\frac n2-m)\pi}{2m}}\right)Q_jM_{j,j}Q_j,
 \end{equation}
 where $M_{i,j}$ is an element of $\mathbf{D}_{2,2}$. This equation, combined with the fact that $\mathbf{D}_{2,2}$ is strictly negative definite on $\bigoplus_{j\in \tilde{J}_2}Q_j L^2$, yields our claim for the case of $2m < n \leq 4m$.

In the dimensions $1\le n \leq 2m$, for $2m - n + 1 \leq j \leq 2m - \frac{n}{2} - 1$, the operator $M^{\pm}_{j,j}$ lies on the diagonal of $\mathcal{A}_{2,2}$ as presented in \eqref{eq-structure-D-mn+1-1-invese}. Moreover, from \eqref{eq-Dk-hat}, and \eqref{eq-the-inverese-of-A0pm} and the Feshbach formula, one has 
$$\mathcal{A}_{2,2}= U_2^\pm\cdot K^{1}_{2,2} \mathbf{d}^{-1} K^{1}_{2,2}\cdot 
U_2^\pm \, \text{diag}\left((-1)^iQ_i \right)_{i\in \tilde{J}_{2} },$$
where $ K^{1}_{2,2}=(k_{i,j})_{i,j\in \tilde{J}_2}$ is a sub-matrix of the projection matrix $K_1=(k_{i,j})_{i,j\in\tilde{J}_0\cup \tilde{J}_2}$ onto the kernel of $A_{0}^{s}$. Here, $A_{0}^{s}$ is defined in \eqref{eq-Dk-hat}. Notice that, due to the cancellation of $Q_j$ for $2m - n + 1 \leq j \leq 2m - \frac{n}{2} - 1$, we have $G_{2m - n}Q_j = 0$. This fact implies that
$$A_{0}^{s} \cdot \begin{pmatrix}
	0&0\\
	0& \text{diag}\left( Q_j\right) _{j\in \tilde{J}_2}\\
\end{pmatrix}= 0.$$
  Hence, since $K_1$ is the projection onto the kernel of $A_{0}^{s}$, by the definition of $K^{1}_{2,2}$, one has 
$$K^{1}_{2,2}\, \text{diag}\left( Q_j\right) _{j\in \tilde{J}_2}=\text{diag}\left( Q_j\right) _{j\in \tilde{J}_2}\, K^{1}_{2,2}=\text{diag}\left( Q_j\right) _{j\in \tilde{J}_2}.$$
Then, we have 
$$\mathcal{A}_{2,2}= U_2^\pm\cdot  \mathbf{d}^{-1} \cdot 
U_2^\pm \, \text{diag}\left((-1)^iQ_i \right)_{i\in \tilde{J}_{2} }$$
Combining this identity with the strictly-negativity of $\mathbf{d}$, we can conclude that $M^{\pm}_{j,j}$ still satisfies \eqref{eq: Mjj^+-Mjj^-} with a negative operator $M_{j,j}$ on $Q_jL^2$, and thus we obtain the desired result.
\end{proof}

\begin{proposition}\label{pro: dose not vanish for j=2m-n2}
When $\mathbf{k}= m_n$ or $m_n+1$,  we have that 
\begin{equation}\label{eq:no vanish for j=2m-fracn2}
\begin{split}
	\Gamma^{+}&_{2m-\frac n2, 2m-\frac n2}(\lambda)-\Gamma^{-}_{2m-\frac n2, 2m-\frac n2}(\lambda)\\
	&=-(\log\lambda)^{-1}(b_1)^{-2}\left( a_{2m-\frac n2}^{+}-a_{2m-\frac n2}^{-}\right)  \left(Q_{2m-\frac n2}vG_{4m-n}vQ_{2m-\frac n2}\right)^{-1}+\Psi(\lambda), 	
\end{split}
\end{equation}
where $\Gamma^{\pm}_{2m-\frac n2, 2m-\frac n2}(\lambda)$ are given in Theorem \ref{thm-M inverse-even} and Theorem \ref{thm-M inverse-even-1}. Moreover, the operator $\Psi(\lambda)$ satisfies 
\begin{equation}\label{eq: est for Phi(lambda)}
\left\|\frac{d^l}{d\lambda^l}\Psi(\lambda)\right\|_{\mathbb{B}(L^2)}\lesssim |\log \lambda|^{-\frac 32} \lambda^{-l}, 
\end{equation}
for $0\le l \le \frac n2+1$ and $0<\lambda<\lambda_0$.
\end{proposition}
\begin{proof}
 We claim that 
\begin{equation}\label{eq: claim for Gamma when i=j=2m-fracn2}
	\Gamma_{2m-\frac n2, 2m-\frac n2}^\pm(\lambda)=-(\log\lambda)^{-1}\left(\frac{a_{2m-\frac n2}^{\pm}}{(b_1)^2}\, Q_{2m-\frac n2}\left(Q_{2m-\frac n2}G_{4m-n}Q_{2m-\frac n2}\right)^{-1}Q_{2m-\frac n2}+Z\right)+\Psi^{\pm}(\lambda), 
\end{equation}
where $Z$ is an operator on $L^2$ independent of the sign $\pm$ . Moreover, $\Psi^{\pm}(\lambda)$ satisfies \eqref{eq: est for Phi(lambda)}. If the claim \eqref{eq: claim for Gamma when i=j=2m-fracn2} holds, then \eqref{eq:no vanish for j=2m-fracn2} readily follows. We only consider the case in the spatial dimension $n\le 2m$. The same analysis applies to the even dimensions $2m<n\le 4m$, and thus we omit the repetitive details here for brevity.  Furthermore, we only prove the case $\k=m_n+1$ since the result is similar for $\k=m_n$ by letting $Q_{2m-\frac n2 +\frac 12}=0$ in the subsequent proof, due to Proposition \ref{prop-resonance-equivalence}.

To prove the  claim \eqref{eq: claim for Gamma when i=j=2m-fracn2}, from \eqref{eq-structure-D-mn+1-1-invese} and the Neumann series given in \eqref{eq：Neumann series to main result}, we have that
\begin{align*}
&\left(\Gamma^{\pm}_{i,j}(\lambda)\right)_{i,j\in J_{m_n+1}}\\
&=\text{diag}(Q_i)_{i\in J_{m_n+1}}\left(\Upsilon_{M}^{\pm}(\lambda)+ \sum_{l=1}^{\infty}(D_{m_n+1}^\pm)^{-1}
\left( -S_{m_n+1}^\pm(\lambda)\,(D_{m_n+1}^\pm) ^{-1}\right)^l\right) \text{diag}(Q_j)_{j\in J_{m_n+1}},
\end{align*}
where $\big(\Gamma^{\pm}_{i,j}(\lambda)\big)_{i,j\in J_{m_n+1}}$ denotes the matrix of the elements $\Gamma^{\pm}_{i,j}(\lambda)$.  It is clear that $\Upsilon_{M}^{\pm}(\lambda)$ makes no contribution to $\Gamma^{\pm}_{2m-\frac n2, 2m-\frac n2}(\lambda)$  since the element of the subscript $(2m-\frac n2, 2m-\frac n2)$ in $\Upsilon_{M}^{\pm}(\lambda)$ is zero by noting \eqref{eq-structure-D-mn+1-1-invese}. 

Noting that the expansion \eqref{eq: expansion for r_ij^pm} and $S_{m_n+1}^\pm(\lambda)=\big(r_{i,j}^{\pm}(\lambda)\big)_{i,j\in J_{m_n+1}}$, the series $$\sum_{l=3}^{\infty}(D_{m_n+1}^\pm)^{-1}
\left( -S_{m_n+1}^\pm(\lambda)\,(D_{m_n+1}^\pm) ^{-1}\right)^l$$ only contributes to the part $\Psi^{\pm}(\lambda)$. Therefore, we now only need to consider the contribution of 
\begin{equation}\label{eq: the main term on Gamma 2m-frac n2}
\sum_{l=1}^{2}(D_{m_n+1}^\pm)^{-1}
\big( -S_{m_n+1}^\pm(\lambda)\,(D_{m_n+1}^\pm)^{-1}\big)^l.	
\end{equation}

\textbf{I. The contribution of $(D_{m_n+1}^\pm)^{-1}\big( -S_{m_n+1}^\pm(\lambda)\,(D_{m_n+1}^\pm)^{-1}\big)$.}

One first has that
$$(D_{m_n+1}^\pm)^{-1}   \big( -S_{m_n+1}^\pm(\lambda)\,(D_{m_n+1}^\pm)^{-1}\big)=-\left( \mathcal{M}^{\pm}+\Upsilon_{M}^{\pm}(\lambda)\right)S_{m_n+1}^\pm(\lambda)\left( \mathcal{M}^{\pm}+\Upsilon_{M}^{\pm}(\lambda)\right).$$
By \eqref{eq-d_ij-k<m} and \eqref{eq-def-D33},  the block $\left( D_{3,3}^{\pm}\right)^{-1}$ has the following structure:
$$ \left( D_{3,3}^{\pm}\right)^{-1}= \begin{pmatrix}
	\left(	b_{1} Q_{m-\frac n2}vG_{4m-n}vQ_{m-\frac n2}\right)^{-1}  & 0  \\
	0 &  	\left(	b_{1} Q_{m-\frac n2+\frac 12}vG_{4m-n}^0vQ_{m-\frac n2+\frac 12}\right)^{-1} \\
\end{pmatrix}.$$
Therefore, via the matrix multiplication and the structure of $\mathcal{M}^{\pm}$ (which has only one nonzero element $\left(b_{1} Q_{2m-\frac n2}vG_{4m-n}vQ_{2m-\frac n2}\right)^{-1}$ in the $(2m-\frac{n}{2})$-column and $(2m-\frac{n}{2})$-row, i.e., $M_{i,2m-\frac n2}^{\pm}=M_{2m-\frac n2, i}^{\pm}=0$ when $i\neq 2m-\frac n2$), the contribution of the matrix $-\mathcal{M}^{\pm}S_{m_n+1}^\pm(\lambda)\mathcal{M}^{\pm}$ to  $\Gamma_{2m-\frac n2, 2m-\frac n2}^\pm(\lambda)$ is given by
\begin{equation}\label{eq: contribution of the main-0}
 -\left(b_{1} Q_{2m-\frac n2}vG_{4m-n}vQ_{2m-\frac n2}\right)^{-1} r_{2m-\frac n2, 2m-\frac n2}^{\pm}(\lambda)\left(b_{1} Q_{2m-\frac n2}vG_{4m-n}vQ_{2m-\frac n2}\right)^{-1}.
\end{equation}
According to the expansion of $r_{2m-\frac n2, 2m-\frac n2}^{\pm}(\lambda)$ in  \eqref{eq: expansion for rij when i,j=2m-n2}, the operator in \eqref{eq: contribution of the main-0} can be written as
\begin{equation}\label{eq: contribution of the main-1}
\begin{split}
-(\log\lambda)^{-1}\left( {(b_1)^{-2}}\, a_{2m-\frac n2}^{\pm} Q_{2m-\frac n2}\left(Q_{2m-\frac n2}vG_{4m-n}vQ_{2m-\frac n2}\right)^{-1}Q_{2m-\frac n2}+Z_0\right)+\Psi_0^{\pm}(\lambda).
\end{split}
\end{equation}
Here, 
$$Z_0 =(b_1)^{-1}\,\left(Q_{2m-\frac n2}vG_{4m-n}vQ_{2m-\frac n2}\right)^{-1}Q_{2m-\frac n2}vG_{4m-n}^0vQ_{2m-\frac n2}\left(Q_{2m-\frac n2}vG_{4m-n}vQ_{2m-\frac n2}\right)^{-1}$$ 
is independent of the sign $\pm$, and $\Psi_0^{\pm}(\lambda)$  also satisfies \eqref{eq: est for Phi(lambda)}. Meanwhile, by the structure of $\Upsilon_{M}^{\pm}(\lambda)$ (the elements in $(2m-\frac{n}{2})$-column and $(2m-\frac{n}{2})$-row of $\Upsilon_{M}^{\pm}(\lambda)$ are zero), the matrix $$
-\Upsilon_{M}^{\pm}(\lambda)S_{m_n+1}^\pm(\lambda)\mathcal{M}^{\pm}-\mathcal{M}^{\pm}S_{m_n+1}^\pm(\lambda)\Upsilon_{M}^{\pm}(\lambda)-\Upsilon_{M}^{\pm}(\lambda)S_{m_n+1}^\pm(\lambda)\Upsilon_{M}^{\pm}(\lambda)
$$ 
makes no contribute to $\Gamma^{\pm}_{2m-\frac n2, 2m-\frac n2}(\lambda)$. 

\textbf{II. The contribution of $(D_{m_n+1}^\pm)^{-1}\big( -S_{m_n+1}^\pm(\lambda)\,(D_{m_n+1}^\pm)^{-1}\big)^2$.}

Next, we consider the contribution of 
$$(D_{m_n+1}^\pm)^{-1}
\big( -S_{m_n+1}^\pm(\lambda)\,(D_{m_n+1}^\pm)^{-1}\big)^2$$
 in \eqref{eq: the main term on Gamma 2m-frac n2}. Noting that $(D_{m_n+1}^\pm)^{-1}=\mathcal{M}^{\pm}+\Upsilon_{M}^{\pm}(\lambda)$,  this operator can be decomposed into three parts as follows
 $$
 O_1^{\pm}=\big((D_{m_n+1}^\pm)^{-1}\, S_{m_n+1}^\pm(\lambda)\big)^2\Upsilon_{M}^{\pm}(\lambda) +\Upsilon_{M}^{\pm}(\lambda) 
 \big( S_{m_n+1}^\pm(\lambda)\,(D_{m_n+1}^\pm)^{-1}\big)^2,$$
 
$$O_2^{\pm}=(D_{m_n+1}^\pm)^{-1}\, S_{m_n+1}^\pm(\lambda)\Upsilon_{M}^{\pm}(\lambda) S_{m_n+1}^\pm(\lambda)\,(D_{m_n+1}^\pm)^{-1},$$
and 
$$O_3^{\pm}=\mathcal{M}^{\pm}
\big( S_{m_n+1}^\pm(\lambda)\mathcal{M}^{\pm}\big)^2.$$
Since the elements in $(2m-\frac{n}{2})$-column and $(2m-\frac{n}{2})$-row of $\Upsilon_{M}^{\pm}(\lambda)$ are zero,  $O_1^{\pm}$ makes no contribution to  $\Gamma^{\pm}_{2m-\frac n2, 2m-\frac n2}(\lambda)$. Meanwhile, the second operator $O_2^{\pm}$ only contributes to $\Psi^{\pm}(\lambda)$. Hence, we are left to observe the contribution of $O_3^{\pm}$.   

By \eqref{eq: def of the elements of BMB} and the property in \eqref{eq-cancellation-Q_j}, for $j\in J_{m_n+1}$, a careful calculation shows that $r_{2m-\frac n2, j}^{\pm}(\lambda)$ are equal to
\begin{align*}
\begin{cases}
(\log \lambda)^{-\frac 12}a_{m-\frac n4+j/2}^{\pm}Q_{2m-\frac n2}vG_{2m-\frac n2+j}v Q_{j}+\varTheta_{2m-\frac n2, j}^{\pm, l}(\lambda),\    \text{ if $2m-\frac n2+j$ is even and $j\in \tilde{J}_0\cup \tilde{J}_2$},& \\
(\log \lambda)^{-\frac 12}b_1 Q_{2m-\frac n2}vG_{4m-n}^{0}v Q_{2m-\frac{n}{2}+\frac 12}
+\varTheta_{2m-\frac n2, j}^{\pm, l}(\lambda) ,\quad\ \  \text{ if $j=2m-\frac n2+\frac 12$},& \\ 
\varTheta_{2m-\frac n2, j}^{\pm, l}(\lambda), \ \quad\qquad\qquad\qquad\qquad\qquad\qquad\qquad\qquad\qquad\,\,\text{else},& 
\end{cases}
\end{align*}
where $\tilde{J}_0$ and $\tilde{J}_2$ are defined in \eqref{eq-def-sets}. Meanwhile, for $i\in J_{m_n=1}$,  $r_{i,2m-\frac n2}^{\pm}(\lambda)$ are equal to
\begin{align*}
	\begin{cases}
		(\log \lambda)^{-\frac 12}a_{m-\frac n4+i/2}^{\pm}Q_{i}vG_{2m-\frac n2+i}vQ_{2m-\frac n2} +\varTheta_{i,2m-\frac n2}^{\pm, l}(\lambda),\,   \text{ if $2m-\frac n2+i$ is even and $i\in \tilde{J}_0\cup \tilde{J}_2$},& \\
		(\log \lambda)^{-\frac 12}b_1 Q_{2m-\frac{n}{2}+\frac 12}vG_{4m-n}^{0}vQ_{2m-\frac n2} +\varTheta_{i,2m-\frac n2}^{\pm, l}(\lambda),\ \ \ \  \text{ if $i=2m-\frac n2+\frac 12$}& \\ 
	\varTheta_{i,2m-\frac n2}^{\pm, l}(\lambda),\quad\qquad\qquad\qquad\qquad\qquad\qquad\qquad\qquad\qquad\,\,\text{else}.& 
	\end{cases}
\end{align*}
In the above two identities, $\varTheta_{i,j}^{\pm, l}(\lambda)$ satisfy 
\begin{equation}\label{eq: est for vartheta,l}
\left\|\frac{d^l}{d\lambda^l}\varTheta_{2m-\frac n2, 2m-\frac{n}{2}+\frac12}^{\pm, l}(\lambda) \right\|_{\mathbb{B}(L^2)}\lesssim |\log\lambda|^{-1}\lambda^{-l},\quad \text{for $0<\lambda<\frac12$}.
\end{equation}
Thus, we decompose $S_{m_n+1}^\pm(\lambda)$ as a sum of 
{\small
$$S_{m_n+1}^{\pm,0}(\lambda)=\begin{pmatrix}
	0 &\left( (\log \lambda)^{-\frac 12}\varTheta_{i,2m-\frac n2}^{\pm, 0}\right) _{ i\in \tilde{J}_0\cup\tilde{J}_1\tilde{J}_2}&0 \\
	\left((\log \lambda)^{-\frac 12}\varTheta_{2m-\frac n2,j}^{\pm, 0} \right)_{ i\in \tilde{J}_0\cup\tilde{J}_1\tilde{J}_2} &0 &(\log \lambda)^{-\frac 12}\varTheta_{2m-\frac n2,2m-\frac n2+\frac 12}^{\pm,0}\\
    0 &(\log \lambda)^{-\frac 12}\varTheta_{2m-\frac n2+\frac 12, 2m-\frac n2}^{\pm,0} &0\\
\end{pmatrix}$$
}
and 
$${\small S_{m_n+1}^{\pm,1}(\lambda)=\begin{pmatrix}
	 \big(r_{i,j}^{\pm}(\lambda)\big)_{i,j\in \tilde{J}_0\cup\tilde{J}_1\tilde{J}_2}& \left(\varTheta_{i,2m-\frac n2}^{\pm, l}(\lambda) \right)_{i\in \tilde{J}_0\cup\tilde{J}_1\tilde{J}_2} &0 \\
	\left(\varTheta_{2m-\frac n2,j}^{\pm, l}(\lambda) \right)_{j\in \tilde{J}_0\cup\tilde{J}_1\tilde{J}_2} & \varTheta_{2m-\frac n2,2m-\frac n2}^{\pm,l}(\lambda)&\varTheta_{2m-\frac n2,2m-\frac n2+\frac 12}^{\pm,l}(\lambda)\\
	0  & \varTheta_{2m-\frac n2+\frac12,2m-\frac n2}^{\pm,l}(\lambda) &0\\
	\end{pmatrix},}$$ 
\noindent where if $i\in \tilde{J}_0\cup\tilde{J}_2$ and $i+2m-\frac n2$ is even,  
$$\varTheta_{i,2m-\frac n2}^{\pm, 0}=a_{m-\frac n4+\frac i2}^{\pm}Q_{i}vG_{2m-\frac n2+i}vQ_{2m-\frac n2},\quad \varTheta_{2m-\frac n2,i}^{\pm, 0}=a_{m-\frac n4+\frac i2}^{\pm}Q_{2m-\frac n2}vG_{2m-\frac n2+i}v Q_{i};$$
if $i\in \tilde{J}_0\cup\tilde{J}_2$ and $i+2m-\frac n2$ is old,
$$\varTheta_{i,2m-\frac n2}^{\pm, 0},\, \varTheta_{2m-\frac n2,i}^{\pm, 0}=0;$$
if $i,j \in \tilde{J}_1$, 
$$\varTheta_{i,2m-\frac n2}^{\pm, 0},\, \quad \varTheta_{2m-\frac n2,j}^{\pm, 0}=0;$$
and
$$\varTheta_{2m-\frac n2,2m-\frac n2+\frac 12}^{\pm,0}=b_1Q_{2m-\frac n2}vG_{4m-n}^0vQ_{2m-\frac n2+\frac12},\ \varTheta_{2m-\frac n2+\frac 12,2m-\frac n2}^{\pm, 0}=b_1Q_{2m-\frac n2+\frac12}vG_{4m-n}^0vQ_{2m-\frac n2}.$$
 At this time, $O_3^\pm$ can be written as 
$$O_3^{\pm}=\mathcal{M}^{\pm}
\big( S_{m_n+1}^{\pm,0}(\lambda)\mathcal{M}^{\pm}\big)^2+O_3^{\pm, r}.$$

We now analyze the contribution of $O_3^{\pm}$. Noting \eqref{eq: est for vartheta,l}, it is straightforward to show that $O_3^{\pm,r}$ contributes only to $\Psi^{\pm}(\lambda)$. Hence it remains to study the contribution of $\mathcal{M}^{\pm}\big( S_{m_n+1}^{\pm,0}(\lambda)\mathcal{M}^{\pm}\big)^2$. Let $$S_{m_n+1}^{\pm,0}(\lambda)\mathcal{M}^{\pm}S_{m_n+1}^{\pm,0}(\lambda)=\left(s_{i,j}^{\pm}(\lambda) \right)_{i,j\in J_{m_n+1}}.$$
By the matrix multiplication along with the structures of $S_{m_n+1}^{\pm,0}(\lambda)$ and $\mathcal{M}^{\pm}$, a simple computation shows that 
\begin{align*}
s_{2m-\frac{n}{2}, 2m-\frac{n}{2}}^{\pm}(\lambda)&=(\log\lambda)^{-1}\sum_{i,j \in J_{m_n+1}}\varTheta_{2m-\frac n2,i}^{\pm, 0}\, M_{i,j}^{\pm}\varTheta_{j,2m-\frac n2}^{\pm, 0} \\
&=(\log\lambda)^{-1}\sum_{i,j \in \tilde{J}_{0}\cup \tilde{J}_{2}}\varTheta_{2m-\frac n2,i}^{\pm, 0}\, M_{i,j}^{\pm}\varTheta_{j,2m-\frac n2}^{\pm, 0}+(\log\lambda)^{-1}L,
\end{align*}
where by the definition, the operator $L$ is independent of the sign $\pm$ and   
\begin{align*}
L=&(b_1)^{2} Q_{2m-\frac n2}vG_{4m-n}^0vQ_{2m-\frac n2+\frac12}\left(d^\pm_{2m-\frac{n}{2}+\frac{1}{2},2m-\frac{n}{2}+\frac{1}{2}}\right)^{-1} Q_{2m-\frac n2+\frac12}vG_{4m-n}^0vQ_{2m-\frac n2}\\
=& (b_1)^{-1}Q_{2m-\frac n2}vG_{4m-n}^0vQ_{2m-\frac n2+\frac12}\left(Q_{2m-\frac n2+\frac12}vG_{4m-n}^0vQ_{2m-\frac n2+\frac12}\right)^{-1} Q_{2m-\frac n2+\frac12}vG_{4m-n}^0vQ_{2m-\frac n2}.
\end{align*}
In the following, we shall show the part $(\log\lambda)^{-1}\sum_{i,j \in \tilde{J}_{0}\cup \tilde{J}_{2}}\varTheta_{2m-\frac n2,i}^{\pm, 0}\, M_{i,j}^{\pm}\varTheta_{j,2m-\frac n2}^{\pm, 0}$ is also independent of the sign $\pm$. By the definition of $\varTheta_{2m-\frac n2,s}^{\pm, 0}$ and $\varTheta_{h,2m-\frac n2}^{\pm, 0}$, if $2m-\frac n2+s$ and $2m-\frac n2+h$ are even, then 
$$\varTheta_{2m-\frac n2,s}^{\pm, 0}\, M_{s,h}^{\pm}\,\varTheta_{h,2m-\frac n2}^{\pm, 0}=a_{m-\frac n4+\frac s2}^{\pm}a_{m-\frac n2+\frac h2}^{\pm}Q_{2m-\frac n2}vG_{2m-\frac n2+s}v Q_{s}\cdot M_{s,h}^{\pm}\cdot Q_{h}vG_{2m-\frac n2+h}vQ_{2m-\frac n2},$$
otherwise, 
$$\varTheta_{2m-\frac n2,s}^{\pm, 0}\, M_{s,h}^{\pm}\,\varTheta_{h,2m-\frac n2}^{\pm, 0}=0.$$
According to \eqref{eq-def-U0-U2} and the relation in \eqref{eq-the-inverese-of-A0pm}, for $s,h\in \tilde{J}_0\cup \tilde{J}_2$, we can $M_{s,h}^{\pm}$ as
$$ M_{s,h}^{\pm}=(-1)^{s}(-i)^{s+h}e^{\pm i\frac{s+h+n-2m}{2m}}M_{s,h},$$
where the operator $M_{s,h}$ is independent of $\pm$. Note that by \eqref{eq:def of aj-1},
$$a_{m-\frac n4+\frac s2}^{\pm}=a_{m-\frac n4+ \frac s2} \cdot e^{\mp i\frac{n-2m+2(m-\frac n4+\frac s2)}{2m}}=a_{m-\frac n4+ \frac s2}e^{\mp i\frac{\frac n2+s}{2m}}$$ 
Thus, if $2m-\frac n2+s$ and $2m-\frac n2+h$ are even, then
$$\varTheta_{2m-\frac n2,s}^{\pm, 0}\, M_{s,h}^{\pm}\,\varTheta_{h,2m-\frac n2}^{\pm, 0}=a_{m-\frac n4+\frac s2}a_{m-\frac n4+\frac h2}Q_{2m-\frac n2}vG_{2m-\frac n2+s}v Q_{s}\cdot M_{s,h}\cdot Q_{h}vG_{2m-\frac n2+h}vQ_{2m-\frac n2},$$
where we have used the fact $(-1)^s (-i)^{s+h}=1$. We have shown that $s_{2m-\frac{n}{2}, 2m-\frac{n}{2}}^{\pm}(\lambda)$ is independent of the sign $\pm$. By the structure of $\mathcal{M}^{\pm}$, the contribution of $\mathcal{M}^{\pm}
\big(S_{m_n+1}^{\pm,0}(\lambda)\mathcal{M}^{\pm}\big)^2$ to $\Gamma_{2m-\frac n2, 2m-\frac n2}^\pm(\lambda)$ is 
$$\left(	b_{1} Q_ivG_{4m-n}vQ_j\right)^{-1}s_{2m-\frac{n}{2}, 2m-\frac{n}{2}}^{\pm}(\lambda)\left(	b_{1} Q_ivG_{4m-n}vQ_j\right)^{-1},$$
which is also independent of the sign $\pm$. This fact shows $\mathcal{M}^{\pm}
\big(S_{m_n+1}^{\pm,0}(\lambda)\mathcal{M}^{\pm}\big)^2$ only contributes to the part $Z$ in \eqref{eq: claim for Gamma when i=j=2m-fracn2}. To sum up, we obtain the claim in \eqref{eq: claim for Gamma when i=j=2m-fracn2} for the dimensions $n\le 2m$. 
The proof is completed. 
\end{proof}

%\vspace{1.2em}   
\section{Proof of  Lemmas \ref{lemma-QjvR-even-1}--\ref{lemma-QjvR-even-2}}\label{sec:proof}

In this section, we will give the proof of Lemmas \ref{lemma-QjvR-even-1}--\ref{lemma-QjvR-even-2}. Before the proof, we first give some lemmas on the integral estimate which will be used in the rest of this section. The following lemma appears in \cite[Lemma 3.8]{GV}. 
\begin{lemma}\label{lem-est-integral-1}
	Let $n\ge 1$. Then there is some absolute constant $C>0$ such that
	\begin{equation*}\label{eq-integral}
		\int_{\mathbb{R}^n}|x-y|^{-k}\langle y\rangle^{-l}\,\d y\le C \langle x\rangle^{-\min\{k,\, k+l-n\}}, \quad x\in\mathbb{R}^n,
	\end{equation*}
	provided  $l\ge 0$, $0\le k<n$ and $k+l>n$.
\end{lemma}

The following two lemmas  deals with functions satisfying certain vanishing moments.
\begin{lemma}\label{lemma-chhz-vanish-1} \textnormal{(\cite[Lemma 4.2]{CHHZ-2024})}
Suppose that  $\frac{1-n}{2}\le k\in\mathbb{Z}$ and $j\in\mathbb{N}_0$.
Let $f\in L_\sigma^2(\mathbb{R}^n)$ with $\sigma >\max\{j,\, k\}+\frac{n}{2}$ and satisfy that $\int_{\R^n}x^\alpha f(x)dx=0$ for all~ $|\alpha|\le j-1$. Then for any $x\in\mathbb{R}^n$, one has
\begin{equation}
\Big|\int_{\R^{n}}|x-y|^kf(y)dy\Big|\lesssim\|f\|_{L^2_\sigma}\langle x\rangle^{k-j}, \quad x\in \mathbb{R}^n.
\end{equation}
\end{lemma}

\begin{lemma}\label{lemma-kernel estimate-for-vanish-2}
Under the same hypothesis of Lemma \ref{lem-est-integral-1}, we have that
	\begin{equation}\label{eq-kernel estimate-log}
		\left| \int_{\mathbb{R}^n} |x-y|^k\log|x-y|f(y)dy\right| \lesssim\|f\|_{L^2_\sigma}\langle x\rangle^{p-j}(1+\log\langle x\rangle)^{\delta_{k,j}},\ \ x\in\mathbb{R}^n.
	\end{equation}
	where 
	$$\delta_{k,j}=\begin{cases}
		0, & \mbox{if }\ p\ \text{is even and}\ k\le j-1,\\
		1, & \mbox{otherwise}.
	\end{cases}$$
\end{lemma}
\begin{proof}
When $|x|\le 1$, note that 
$$\Big| |x-y|^k\log|x-y| \Big|\lesssim |x-y|^{k-\varepsilon}+|x-y|^{k+\varepsilon} \quad \text{for a sufficiently small $\varepsilon>0$.} $$
Thus, in this case, it follows from  the H\"older inequality that 
$$\left| \int_{\mathbb{R}^n} |x-y|^k\log|x-y|f(y)dy\right| \lesssim 1,$$ 
which shows \eqref{eq-kernel estimate-log}.

Now we set $|x|>1$. And we have 
$$|x-y|^k\log|x-y|=|x-y|^k\log\frac{|x-y|}{|x|}+|x-y|^k\log|x|.$$
First, Lemma \ref{lemma-chhz-vanish-1} derives that 
\begin{equation}\label{eq-est-log-2}
\left|\int_{\mathbb{R}^n} |x-y|^k\log|x|f(y) dy\right|\lesssim \|f\|_{L^2_\sigma}\, \langle x\rangle^{k-j}(1+\log\langle x\rangle).
\end{equation}
Specially, when $k$ is even and $k\le j-1$, one has 
\begin{equation}\label{eq-est-log-3}
\int_{\mathbb{R}^n} |x-y|^k\log|x|f(y) dy=0,
\end{equation}
since $$|x-y|^k=\sum_{|\alpha|+|\beta|=k}C_{\alpha, \beta}(-1)^{|\beta|}x^{\alpha}y^{\beta}.$$
Next, using the Talyor expansion, then we have 
\begin{align*}
&|x-y|^k\log\frac{|x-y|}{|x|}
=|x|^k\left( 1+\frac{|x-y|-|x|}{|x|}\right)^k \log\left( 1+\frac{|x-y|-|x|}{|x|}\right) \\
=& |x|^k\sum_{l=0}^{j-1}C_l\left(\frac{|x-y|-|x|}{|x|} \right)^l \\
&+ \frac{|x|^k}{(j-1)!}\left(\frac{|x-y|-|x|}{|x|} \right)^j \int_0^1 g^{(j)}\left(1+s\frac{|x-y|-|x|}{|x|} \right)(1-s)^{j-1}ds\\
=:& I(x,y)+II(x,y),	
\end{align*}
where the function $g(\lambda)=\lambda^k\log\lambda$ and $g^{(j)}(\lambda)=\frac{d^j}{d\lambda^j}g(\lambda)$. 
By  Lemma \ref{lemma-chhz-vanish-1}, we have 
$$\left| \int_{\mathbb{R}^n} I(x,y) f(y) dy\right|\lesssim \|f\|_{L^2_\sigma}\, \langle x\rangle^{k-j}.$$
Note that a direct computation shows that 
$$|g^{(j)}(\lambda)|\lesssim
\begin{cases}
|\lambda|^{k-j}(1+|\log\lambda|), \quad & j\le k, \\
|\lambda|^{k-j}, \quad & j> k. \\
\end{cases}	  
$$
A simple computation yields that if  $l\ge 0$,
\begin{align*}
	\left(1+s\frac{|x-y|-|x|}{|x|} \right)^l=\left( \frac{(1-s)|x|+s|x-y|}{|x|}\right) ^{l}&\lesssim \left( 1+\frac{|y|}{|x|}\right) ^{l},
\end{align*}
and if $i\ge 1$
\begin{align*}
\left( \frac{(1-s)|x|+s|x-y|}{|x|}\right) ^{-i}&\lesssim \left( \frac{(1-s)|x|+s|x-y|}{|x|}\right) ^{-i+\frac 12}\left( \frac{(1-s)|x|+s|x-y|}{|x|}\right) ^{-\frac 12}\\
&\lesssim (1-s)^{-i+\frac 12}\left( 1+\left( \frac{|x|}{|x-y|}\right)^{-\frac 12}\right).
\end{align*}
And $|\log\lambda|\lesssim |\lambda|^{\varepsilon}+|\lambda|^{-\varepsilon} $ for some positive small $\varepsilon$ such that $\varepsilon<\sigma-k-\frac n2$.
Thus,
$$\left| II(x,y)\right|\lesssim
\begin{cases}
	|x|^{k-j}|y|^j\left( 1+\frac{|y|}{|x|}\right) ^{k-j+\varepsilon}+|x|^{k-j}|y|^j\left( 1+\frac{|y|}{|x|}\right) ^{k-j}\left(1+(\frac{|x|}{|x-y|})^{\varepsilon}\right), \quad & j\le k, \\
	|x|^{k-j}|y|^j\left(1+(\frac{|x|}{|x-y|})^{\frac 12}\right), \quad & j> k. \\
\end{cases}	  
$$ 
There, applying Lemma \ref{lem-est-integral-1}, it follows that 
$$\left| \int_{\mathbb{R}^n} II(x,y) f(y) dy \right| \lesssim \|f\|_{L^2_\sigma}\, \langle x\rangle^{k-j}.$$ 
This, together with \eqref{eq-est-log-2} and \eqref{eq-est-log-3} , yields \eqref{eq-kernel estimate-log}. 
\end{proof}

We then have an immediate Proposition. 
\begin{proposition}\label{corcancell}
Assume that $V$ satisfies \eqref{decayonV} in Assumption \ref{Assumption},  $j\in J_{\mathbf{k}}$, $Q_j$ is the projection given by \eqref{eq-Q_j-even},  and $l\in\mathbb{Z}$ such that $\frac{1-n}{2}\le l\le j+1$. Then we have 
\begin{equation}\label{eq-kernel estimate-even-1}
\Big\|Q_j(v(\cdot)|\cdot-y|^l)(x) \Big\|_{L_x^2(\mathbb{R}^n)}\lesssim\langle y\rangle^{l-\vartheta(j)}, \quad x\in \mathbb{R}^n,
\end{equation}
where $\vartheta(j)$ given in \eqref{thetaj}. Similarly, when $\vartheta(j)>2m-n\ge 0$, we have 
\begin{equation}\label{eq-kernel estimate-even-fundamental-solution}
	\Big\|Q_j\left( v(\cdot)|\cdot-y|^{2m-n}\log|\cdot-y|\right) (x) \Big\|_{L_x^2(\mathbb{R}^n)}\lesssim\langle y\rangle^{l-\vartheta(j)}, \quad x\in \mathbb{R}^n,
\end{equation}
\end{proposition}
\begin{proof}
By the definition of $Q_j$, we known that $Q_j(y^{\alpha}v)=0$ for $|\alpha|\le \vartheta(j)-1$. Thus, one has 
\[\int_{\R^n}x^\alpha v(x)Q_j(v(\cdot)|\cdot-y|^l)(x)dx=0.\]
Then by Lemma \ref{lemma-chhz-vanish-1}, we have
\begin{align*}
	\left\|Q_j(v(\cdot)|\cdot-y|^l)(x)\right\|_{L_x^2}^2 
  = &\Big\langle v(x)|x-y|^l,\,  Q_j(v(\cdot)|\cdot-y|^l)(x)\Big\rangle=\int_{\R^n} |x-y|^lv(x)Q_j(v(\cdot)|\cdot-y|^l)(x)dx\\
	\lesssim& \left\|v(x)Q_j(v(\cdot)|\cdot-y|^l)(x)\right\|_{L_{\sigma}^2} \left\langle  y\right\rangle ^{l-\vartheta(j)}
	\lesssim \Big\|Q_j(v(\cdot)|\cdot-y|^l)(x)\Big\|_{L_x^2}\left\langle  y\right\rangle ^{l-\vartheta(j)},
\end{align*}
where $\sigma=\max\{j,l\}+\frac{n}{2}+\varepsilon$ for some sufficiently small $\varepsilon>0$. In the last inequality above, we have used the fact that $|v|\lesssim \langle x\rangle^{-\beta/2} \lesssim \langle x\rangle^{-\sigma}$ from \eqref{decayonV} in Assumption \ref{Assumption} and the range of $j$ given in \eqref{eq-max-Jk}. Hence, the estimate \eqref{eq-kernel estimate-even-1} holds.  And \eqref{eq-kernel estimate-even-fundamental-solution} is obtained in parallel.
\end{proof}

\subsection{Proof of Lemmas \ref{lemma-QjvR-even-1}, \ref{lemma-QjvR-even-1-1} and  \ref{lemma-QjvR+-R-even-1} }\label{secA.1}
 In this subsection, we will prove the conclusions in Lemmas \ref{lemma-QjvR-even-1}--\ref{lemma-QjvR+-R-even-1} by applying the property of free resolvents and Proposition \ref{corcancell}.
\begin{proof}[{\bf Proof of Lemma \ref{lemma-QjvR-even-1}.}]
We focus on the proof of \eqref{eq-QjR kernel-even}--\eqref{eq-wj-estimate-even-4}. We first
define
\begin{equation*}
	\omega_{j}^{\pm}(\lambda,x,y)=e^{\mp i\lambda|y|}\Big(Q_jvR_0^\pm(\lambda^{2m})(\cdot,\,y)\Big)(x).
\end{equation*}
Then, the identity \eqref{eq-QjR kernel-even} follows.

We next establish \eqref{eq-QjR kernel-even} and the estimates \eqref{eq-wj-estimate-even-1}--\eqref{eq-wj-estimate-even-4} \underline{in even dimensions $n\le 4m-2$} by considering two cases: \(\lambda \langle y \rangle \le 1\) and \(\lambda \langle y \rangle \ge 1\).

\vskip0.2cm
\textbf{\emph{Case 1: $\lambda\langle y\rangle\le 1$.}}
In this region, it is obvious that 
\begin{equation}\label{eq-est-for-exp}
\left| \partial_\lambda^l\left(e^{\mp i\lambda|y|} \right)\right|\lesssim \lambda^{-l}.  	
\end{equation}

Thus, to obtain \eqref{eq-wj-estimate-even-1}-- \eqref{eq-wj-estimate-even-3}, it suffices to consider the estimates on $\Big(Q_jvR_0^\pm(\lambda^{2m})(\cdot,\,y)\Big)(x)$.

For each $j\in J_\mathbf{k}$, note that $Q_j(x^{\alpha}v(x))=0$ for all $|\alpha|\le\vartheta(j)-1$ where $\alpha\in \mathbb{N}_0^n$ . Thus, $Q_j(v(\cdot)|\cdot-y|^{2i})=0$  for $2i\le \vartheta(j)-1$, $i\in\mathbb{N}_0$ and a fixed  $y$. Also, note that  $\vartheta(j)\neq 2m-n$ if $j\neq 2m-n$ by the definition $\vartheta(j)$ in \eqref{thetaj}. 

In the even dimensions $1\le n\le 2m$, if $\vartheta(j)< 2m-n$, by applying the expansion \eqref{eq-free expansion-even-1} in Lemma \ref{lemma-free resolvent expansio-smalln} with $\theta=\vartheta(j)$, one has 
\begin{equation}\label{eq-expression-vartheta(j)-1}
\Big(Q_jvR_0^\pm(\lambda^{2m})(\cdot,\,y)\Big)(x)=\lambda^{n-2m}\Big(Q_jvr_{\vartheta(j)}^\pm(\lambda|\cdot-y|)\Big)(x),
\end{equation}
and if $\vartheta(j)\ge 2m-n$, we by applying the expansion \eqref{eq-free expansion-even-2} with $\theta=2m-n+1$, one has that
\begin{equation}\label{eq-expression-vartheta(j)-2}
\begin{aligned}
\Big(Q_j&vR_0^\pm(\lambda^{2m})(\cdot,\,y)\Big)(x)= \\
&\lambda^{n-2m}\Big(Q_jv\left( a_{m-\frac n2}^\pm\lambda^{2j}|\cdot-y|^{2j}+b_0(\lambda|\cdot-y|)^{2m-n}\log(\lambda|\cdot-y|)+r_{2m-n+1}^\pm(\lambda|\cdot-y|)\right) \Big)(x).
\end{aligned}
\end{equation} 
In the even dimensions $2m<n< 4m$, by applying the expansion \eqref{eq-free expansion-even-4} in Lemma \ref{lemma-free resolvent expansio-largen} with $\theta=0$, one has 
\begin{equation}\label{eq-expression-vartheta(j)-3}
	\Big(Q_jvR_0^\pm(\lambda^{2m})(\cdot,\,y)\Big)(x)=\Big(Q_jv\left( b_0|\cdot-y|^{2m-n}+\lambda^{n-2m}r_{0}^\pm(\lambda|\cdot-y|)\right) \Big)(x)
\end{equation}
In the even dimensions $n=4m$,  by applying the expansion \eqref{eq-free expansion-even-4} in Lemma \ref{lemma-free resolvent expansio-largen} with $\theta=0$, one has 
\begin{equation}\label{eq-expression-vartheta(j)-4}
	\Big(Q_jvR_0^\pm(\lambda^{2m})(\cdot,\,y)\Big)(x)=\Big(Q_jv\left( b_0|\cdot-y|^{2m-n}+b_1\lambda^{n-2m}\log(\lambda|\cdot-y|)\lambda^{n-2m}+r_{1}^\pm(\lambda|\cdot-y|)\right) \Big)(x)
\end{equation}

\underline{We first give the proof for estimates \eqref{eq-wj-estimate-even-1}--\eqref{eq-wj-estimate-even-4} when $\vartheta(j)< 2m-n$ and $1\le n\le 2m$.} 

Applying the Taylor expansion for $r_\theta^\pm$  yields
\begin{equation*}
	r_\theta^\pm(\lambda|x-y|)=\sum_{s=0}^{\vartheta(j)-1}\big(r_{\vartheta(j)}^\pm\big)^{(s)}(\lambda|y|)\frac{(\lambda|x-y|-\lambda|y|)^s}{s!}+\tilde{r}_{\vartheta(j)}^\pm(\lambda,|x-y|,|y|),
\end{equation*}
where
\begin{equation}\label{eq-remainder-r}
\begin{split}
\tilde{r}_\theta^\pm&(\lambda,|x-y|,|y|) \\
&=\frac{(\lambda|x-y|-\lambda|y|)^{\vartheta(j) }}{(\vartheta(j)-1)!}\int_0^1(1-t)^{\vartheta(j)-1}\big(r_{\vartheta(j)}^\pm\big)^{(\vartheta(j))}\big(\lambda|y|+t\lambda(|x-y|-|y|)\big)dt.    
\end{split}
\end{equation}
Here, $\big(r_\theta^\pm\big)^{(s)}(z)=\frac{d^s}{dz^s}r_\theta^\pm(z)$ for $z\in \R$.
Then, we have that $\omega_{j}^{\pm}(\lambda,x,y)$ can be written as a sum of 
\begin{equation}\label{eq-case1-t1}
e^{\mp i\lambda|y|}\sum_{s=0}^{\vartheta(j)-1}\frac{\lambda^{n-2m+s}}{s!}\big(r_\theta^\pm\big)^{(s)}(\lambda|y|)\, Q_jv\big(|\cdot-y|-|y|\big)^s(x),
\end{equation}
and
\begin{equation}\label{eq-case1-t2}
\lambda^{n-2m}e^{\mp i\lambda|y|} Q_jv\big(\tilde{r}_{\vartheta(j) }^\pm(\lambda,|\cdot-y|,|y|)\big)(x).	
\end{equation}

We first provide the estimate for the term \eqref{eq-case1-t2}. 

Based on $\eqref{eq-r estimate-even}$ and the fact that $ \theta=  \vartheta(j)$ when $\vartheta(j)< 2m-n$, for $ \lambda \in (0,1) $ and $ l = 0, \ldots, \frac{n}{2}+1 $, the following estimates hold uniformly for $ t \in [0,1]$:\\
when $ \lambda|y|+t\lambda(|x-y|-|y|)\le 1$, then
\begin{equation*}
\Big|\partial_\lambda^l \Big[\big(r_\theta^\pm\big)^{(\vartheta(j))} (\lambda|y|+t\lambda(|x-y|-|y|))\Big]\Big|\lesssim \lambda^{-l};	
\end{equation*} 
when $\lambda|y|+t\lambda(|x-y|-|y|)\ge 1$, then
\begin{equation*}
\begin{aligned}
	\Big|\partial_\lambda^l \Big[\big(r_\theta^\pm\big)^{(\vartheta(j))} (\lambda|y|+t\lambda(|x-y|-|y|))\Big]\Big|&\lesssim \lambda^{-l}(\lambda|y|+t\lambda(|x-y|-|y|))^{\max\{0,\, l-\frac{n-1}{2}\}}\\
	&\lesssim \lambda^{-l}\langle x\rangle^2,
\end{aligned}
\end{equation*}
where we utilize the facts $\big||x-y|-|y|\big|\le |x|$ and $\lambda|y|\le 1$.
Thus, a direct computation  shows that for $l=0,\ldots,\frac{n}{2}+1$, the remainders  $\tilde{r}_{\vartheta(j)}^\pm(\lambda,|x-y|,|y|)$ given by \eqref{eq-remainder-r} satisfy
\begin{equation*}
	\Big|\partial^l_\lambda\Big[\tilde{r}_{\vartheta(j)}^\pm(\lambda,|x-y|,|y|)\Big]\Big| \lesssim \lambda^{\vartheta(j)-l}\langle x\rangle^{\vartheta(j)+2}.\ \ \ 
\end{equation*}
 In addition, note that $j\in J_{\mathbf{k}}$, so that $\vartheta(j)\le \max\{m-\frac n2, 0\}+\mathbf{k}+1$ by \eqref{eq-max-Jk}.
 It follows immediately that $\frac{\beta}{2}>\vartheta(j)+2+\frac{n}{2}$ where  $\beta$ is given in \eqref{decayonV}. Since $|v(x)|\lesssim \langle x\rangle^{-\frac \beta 2}$, therefore, for $\lambda\langle y\rangle\le 1$,  $ 0<\lambda<1$ and $l=0, \ldots, \frac{n}{2}+1$, by the Leibniz rule,  we have that  
\begin{equation}\label{eq3.901}
\begin{aligned}
&\Big\|\partial^l_\lambda\Big[\lambda^{n-2m}e^{\mp i\lambda|y|} Q_j\big(v\tilde{r}_{\vartheta(j) }^\pm(\lambda,|\cdot-y|,|y|)\big)(x)\Big]\Big\|_{L_x^2}\\
\lesssim&\sum_{s=0}^l\big|\partial^s_\lambda \big(\lambda^{n-2m}e^{\mp i\lambda|y|}\big)\big| \,\Big\|Q_j\big(v \partial^{l-s}_\lambda\tilde{r}_{\vartheta(j)}^\pm(\lambda,|\cdot-y|,|y|)\big)(x) \Big\|_{L_x^2}\\
\lesssim&\sum_{s=0}^l\big|\partial^s_\lambda \big(\lambda^{n-2m}e^{\mp i\lambda|y|}\big)\big| \,\Big\|v(x) \partial^{l-s}_\lambda\tilde{r}_{\vartheta(j)}^\pm(\lambda,|x-y|,|y|) \Big\|_{L_x^2}\\
\lesssim &\lambda^{n-2m+\vartheta(j)-l} \Big\|v(x)\langle x\rangle^{\vartheta(j)+2}\Big\|_{L_x^2}\\
\lesssim &\lambda^{n-2m+\vartheta(j)-l}.
\end{aligned}
\end{equation}

We now estimate the term \eqref{eq-case1-t1}.
Applying \eqref{eq-r estimate-even} again, we have 
\begin{equation*}
\Big|\partial^l_\lambda\Big[({r}_\theta^\pm)^{(s)}(\lambda|y|)\Big]\Big| \lesssim \lambda^{\theta-s-l}\langle y\rangle^{\vartheta-s}, \quad l=0,\ldots,\frac{n}{2}+1.
\end{equation*}
This, together with Proposition \ref{corcancell}, shows that  for all  $0\le s\le\vartheta(j)-1$, $\lambda\langle y\rangle\le 1$ and $ 0<\lambda<1$, the following holds:
\begin{equation}\label{eq3.902}
\begin{aligned}
&\Big\|\partial^l_\lambda\Big[\lambda^{n-2m+s}e^{\mp i\lambda|y|}\big(r_\theta^\pm\big)^{(s)}(\lambda|y|)\, Q_jv\big(|\cdot-y|-|y|\big)^s(x)\Big]\Big\|_{L_x^2}\\
\lesssim& \sum_{k=0}^l \left|\partial^k_\lambda \left( \lambda^{n-2m+s}e^{\mp i\lambda|y|}\right)\right|\cdot \left|\partial^{l-k}_\lambda \big(r_\theta^\pm\big)^{(s)}(\lambda|y|) \right|\cdot \Big\|Q_jv(|\cdot-y|-|y|)^s(x)\Big\|_{L_x^2}\\
\lesssim&  \lambda^{n-2m+\theta-l}\langle y\rangle^{\theta-s} \sum_{\rho=0}^{s}|y|^{s-\rho}\Big\| Q_jv(|\cdot-y|^\rho)(x)\Big\|_{L^2_x}\\
\lesssim& \lambda^{n-2m+\theta-l}\langle y\rangle^{\theta-\vartheta(j)}.	
\end{aligned}
\end{equation}
In the above inequality, $\rho\le \vartheta(j)-1\le j$ satisfies the condition in Proposition \ref{corcancell}.  

Therefore, combining \eqref{eq3.901} with \eqref{eq3.902}, if $\vartheta(j)< 2m-n$ and $n\le 2m$,  we have that in the region where $\lambda\langle y\rangle\le 1$,  
\begin{equation}\label{eq3.9031}
	\|\partial^l_\lambda \omega_{j}^{\pm}(\lambda,x,y)\|_{L_x^2}\lesssim \lambda^{n-2m+\theta-l}\langle y\rangle^{\theta-\vartheta(j)}=\lambda^{n-2m+\vartheta(j)-l}, \quad l=0,\ldots,\frac{n}{2}+1.
\end{equation}
To sum up, when $1\le n\le 2m$, the estimats \eqref{eq-wj-estimate-even-1} and \eqref{eq-wj-estimate-even-3} hold in the region where $\lambda\langle y\rangle\le 1$.

\underline{We next consider the case of $\vartheta(j)>  2m-n$ and $1\le n\le 2m$, in the region where $\lambda\langle y\rangle\le 1$.}  

In the same vein, from \eqref{eq-expression-vartheta(j)-2}, we have that for $l=0,\ldots,\frac{n}{2}+1$, there is
\begin{equation}\label{eq3.9032}
	\left\|\partial^l_\lambda \left( Q_jv\left(r_{2m-n+1}^\pm(\lambda|\cdot-y|)\right)(x) \right) \right\|_{L_x^2}\lesssim \lambda^{1-l}\langle y\rangle^{2m-n+1-\vartheta(j)}\lesssim \lambda^{-l}\langle y\rangle^{2m-n-\vartheta(j)}.
\end{equation}
In addition, note that  
$$Q_jv\left( a_{m-\frac n2}^\pm|\cdot-y|^{2j}+b_0|\cdot-y|^{2m-n}\log(\lambda|\cdot-y|)\right)=Q_jv\left( b_0|\cdot-y|^{2m-n}\log|\cdot-y|\right) .$$
From Corollary \ref{corcancell},  there is 
\begin{align*}
\left\| \partial^l_\lambda \Big(Q_jv\left( a_{m-\frac n2}^\pm|\cdot-y|^{2j}+b_0|\cdot-y|^{2m-n}\log(\lambda|\cdot-y|)\right)(x)\Big)\right\|_{L_x^2}
\lesssim  \lambda^{-l} \langle y\rangle^{2m-n-\vartheta(j)}.  
\end{align*}

This combined with \eqref{eq3.9032} shows that If $\vartheta(j)> 2m-n$ and $1\le n\le 2m$, in the region where $\lambda\langle y\rangle\le 1$, it follows that 
\begin{equation}\label{eq3.9033}
	\|\partial^l_\lambda \omega_{j}^{\pm}(\lambda,x,y)\|_{L_x^2}\lesssim \lambda^{-l}\langle y\rangle^{2m-n-\vartheta(j)}, \quad l=0,\ldots,\frac{n}{2}+1.
\end{equation}

Therefore, combining this with \eqref{eq3.9031},  when $1\le n\le 2m$, the estimates \eqref{eq-wj-estimate-even-1}, \eqref{eq-wj-estimate-even-3} and \eqref{eq-wj-estimate-even-4} hold in the region where $\lambda\langle y\rangle\le 1$.

\underline{We proceed to consider the case of $\vartheta(j)=2m-n$, $1\le n\le 2m$, and $\lambda\langle y\rangle\le 1$.} 

In this case, we apply the expansion \eqref{eq-expression-vartheta(j)-2}. By \eqref{eq-kernel estimate-even-1} in Proposition \ref{corcancell},  it is not difficult to see that 
\begin{align*}
\left\| \partial^l_\lambda \Big(Q_jv\left( a_{m-\frac n2}^\pm|\cdot-y|^{2j}\right)(x)\right\|_{L_x^2}
	\lesssim  \lambda^{2m-n-l}.  
\end{align*}
Meanwhile, one has that 
$$b_0|\cdot-y|^{2m-n}\log(\lambda|\cdot-y|)=b_0\log(\lambda|y|)|\cdot-y|^{2m-n}+b_0|\cdot-y|^{2m-n}\log\frac{|\cdot-y|}{|y|}.$$
Using Proposition \ref{corcancell} and the fact $\lambda\langle y\rangle\le 1$, for $l=0,\ldots,\frac{n}{2}+1$, it follows that 
\begin{align*}
	\left\| \partial^l_\lambda \Big(Q_jv\left(b_0\log(\lambda|y|)|\cdot-y|^{2m-n} \right)(x)\right\|_{L_x^2}
	\lesssim  \lambda^{-\varepsilon-l}\langle y\rangle^{-\varepsilon},\, \ \text{for all} \   \varepsilon \in (0,1/2).  
\end{align*}
On the other hand, the same trick used in the proof of Lemma \ref{lemma-kernel estimate-for-vanish-2} leads to  
\begin{align*}
	\left\| \partial^l_\lambda \Big(Q_jv\left(b_0|\cdot-y|^{2m-n}\log\frac{|\cdot-y|}{|y|} \right)(x)\right\|_{L_x^2}
	\lesssim  1.  
\end{align*}
Then, combining with \eqref{eq3.9032}, when $1\le n\le 2m$,  we obtain \eqref{eq-wj-estimate-even-2} in the region $\lambda\langle y\rangle\le 1$.

When $2m<n< 4m$, by using the expansions \eqref{eq-expression-vartheta(j)-3} and \eqref{eq-expression-vartheta(j)-4} instead and in a parallel manner, one can obtain \eqref{eq-wj-estimate-even-1}, \eqref{eq-wj-estimate-even-3} and \eqref{eq-wj-estimate-even-4} as well. The details are omitted for simplicity.

\vskip0.2cm
\textbf{\emph{Case 2 :  $\lambda\langle y\rangle\ge 1$.}}
Noting \eqref{eq-free kernel-F}, we have that
\begin{equation*}
\Big(Q_jvR_0^\pm(\lambda^{2m}(\cdot, y))\Big)(x)=Q_jv\Big(e^{\pm i\lambda|\cdot-y|}\frac{\lambda^{\frac{n+1}{2}-2m}}{|\cdot-y|^{\frac{n-1}{2}}}F_n^\pm(\lambda|\cdot-y|)\Big)(x).
\end{equation*}
Now, $\omega_{j}^{\pm}(\lambda,x,y)$ can be written as follows
\begin{equation*}
	\omega_{j}^{\pm}(\lambda,x,y)= Q_jv\Big(e^{\pm i\lambda(|\cdot-y|-|y|)}\frac{\lambda^{\frac{n+1}{2}-2m}}{|\cdot-y|^{\frac{n-1}{2}}}F_n^\pm(\lambda|\cdot-y|)\Big)(x).
\end{equation*}
Note that 
$$\left|\partial^k_\lambda \Big( e^{\pm i\lambda(|\cdot-y|-|y|)}\Big)\right|\leq \langle\, \cdot\,\rangle^k, \quad \text{for}~ k\in \mathbb{N}_0.$$ 
Then, using \eqref{eq-F decay-odd} and Lemma \ref{lem-est-integral-1}, it is not difficult to find that when $|y|\le 1$,  for any $l=0,\ldots,\frac{n}{2}+1$ and each $N\in \mathbb{N}_0$, we have
\begin{equation}\label{eq3.903}
	\Big\|\partial^l_\lambda \omega_{j}^{\pm}(\lambda,x,y)\Big\|_{L_x^2}\lesssim \lambda^{\frac{n+1}{2}-2m-l}\left\|\frac{\langle y\rangle^lv(x)}{|x-y|^{\frac{n-1}{2}}} \right\|_{L_x^2}  \lesssim_N \lambda^{\frac{n+1}{2}-2m-l}\langle y\rangle^{-N}.
\end{equation}
Hence,  it follows that
 estimates \eqref{eq-wj-estimate-even-1}--\eqref{eq-wj-estimate-even-4} hold in the case when $\lambda\langle y\rangle\ge 1$ and $|y|\le1.$ 

In what follows, we consider the case when $\lambda\langle y\rangle\ge 1$ and $|y|\ge1$. Using the Taylor expansion for $G^\pm_n(z)=e^{\pm iz}F_n^\pm(z)$ of $\big(\vartheta(j)-1\big)$-th order for $z\in \R$, we have
\begin{equation*}
	G^\pm_n(\lambda|x-y|) =\sum_{s=0}^{\vartheta(j)-1}\big(G^\pm_n\big)^{(s)}(\lambda|y|)\frac{(\lambda|x-y|-\lambda|y|)^s}{s!}+g_{\vartheta(j)}^\pm(\lambda,|x-y|,|y|),
\end{equation*}
where
\begin{align*}
	g_{\vartheta(j)}^\pm(&\lambda,|x-y|,|y|)\\
 &=\frac{(\lambda|x-y|-\lambda|y|)^{\vartheta(j) }}{\big(\vartheta(j)-1\big)!}\int_0^1(1-t)^{\vartheta(j)-1}\big(G^\pm_n\big)^{({\vartheta(j)})}(\lambda|y|+t\lambda(|x-y|-|y|))dt.
\end{align*}
Now, it is obtained that $\omega_{j}^{\pm}(\lambda,x,y)$ is a sum of the following two terms:
\begin{equation}\label{eq-case2-t1}
e^{\mp i\lambda|y|}\lambda^{\frac{n+1}{2}-2m} \sum_{s=0}^{\vartheta(j)-1}\lambda^s\big(G^\pm_n\big)^{(s)}(\lambda|y|) Q_jv\Big( \frac{(|\cdot-y|-|y|)^s}{s!|\cdot-y|^{\frac{n-1}{2}}}\Big)(x),
\end{equation}
and
\begin{equation}\label{eq-case2-t2}
e^{\mp i\lambda|y|}\lambda^{\frac{n+1}{2}-2m}Q_jv\Big(\frac{g_{\vartheta(j)}^\pm(\lambda,|\cdot-y|,|y|)}{|\cdot-y|^{\frac{n-1}{2}}}\Big)(x).
\end{equation}
And, the definition of $G^\pm_n$ and  the Leibniz rule show that 
\begin{equation}\label{eq3.904}
	(G^\pm_n)^{(s)}(z)=\sum_{k=0}^s(\pm i)^k C_s^k e^{\pm iz}(F_n^\pm)^{(s-k)}(z), \quad z\in \R,~ s=0,\ldots, \vartheta(j),
\end{equation}
where the constant $C_s^k=\frac{s!}{k!(s-k!)}$.
To estimate $\omega_{j}^{\pm}(\lambda,x,y)$,  it suffices to estimate \eqref{eq-case2-t1} and \eqref{eq-case2-t2}.

We first estimate the term \eqref{eq-case2-t2} in the region of $\lambda\langle y\rangle>1$ and $|y|\ge1$. Noting the expression \eqref{eq3.904}, one has 
\begin{align*}
 e^{\mp i\lambda|y|}(G^\pm_n)^{(\vartheta(j))}\big(\lambda|y|+&t\lambda(|x-y|-|y|)\big)\\
 &=e^{\pm it\lambda(|x-y|-|y|)}\sum_{s=0}^{\vartheta(j)} (\pm i)^s C_{\vartheta(j)}^s (F_n^\pm)^{(\vartheta(j)-s)}\big(\lambda|y|+t\lambda(|x-y|-|y|)\big).	
\end{align*}
Hence, by the Leibniz rule, for any $l=0,\ldots,\frac{n}{2}+1$, $0<\lambda<1$ and $|y|\ge 1$, the following estimate
\begin{equation}\label{eq3.9010}
	\left| \partial_\lambda^l\Big[e^{\mp i\lambda|y|}(G^\pm_n)^{(\vartheta(j))}\big(\lambda|y|+t\lambda(|x-y|-|y|)\big)\Big]\right| \lesssim
	\lambda^{-l}\langle x\rangle^{\frac{n}{2}+1},
\end{equation}
holds uniformly for $t\in[0,\,1]$.
This implies that
$$\bigg| \partial_\lambda^l\Big[\lambda^{\frac{n+1}{2}-2m}e^{\mp i\lambda|y|} g_{\vartheta(j)}^\pm(\lambda,|x-y|,|y|)\Big]\bigg|\lesssim\lambda^{\frac{n+1}{2}-2m+(\vartheta(j)-l}\langle x\rangle^{\vartheta(j)+\frac{n}{2}+1}.$$
Therefore, for $l=0,\ldots,\frac{n}{2}+1$ and $\lambda\in(0,1)$, we have from Lemma \ref{lem-est-integral-1} and the decay assumption \eqref{decayonV} that
\begin{equation}\label{eq3.905}
	\bigg\|\partial^l_\lambda\Big[e^{\mp i\lambda|y|}\lambda^{\frac{n+1}{2}-2m}\big(Q_jv\Big(\frac{g_{\vartheta(j)}^\pm(\lambda,|\cdot-y|,|y|)}{|\cdot-y|^{-\frac{n-1}{2}}}\Big)(x)\Big]\bigg\|_{L_x^2}\lesssim \lambda^{\frac{n+1}{2}-2m+\vartheta(j)-l}\langle y\rangle^{-\frac{n-1}{2}}.
\end{equation}

We next estimate the term \eqref{eq-case2-t1}.
Note that by \eqref{eq3.904}, one has
$$e^{\mp i\lambda|y|}(G^\pm_n)^{(s)}(\lambda|y|)=\sum_{k=0}^s (\pm i)^k C_{s}^k(F_n^\pm)^{(s-k)}(\lambda|y|).$$
And by the estimate \eqref{eq-F decay-odd} of $F^\pm_n$, for  $l=0,\ldots,\frac{n}{2}+1$, there exists
$$
\left| \partial_\lambda^l\Big[e^{\mp i\lambda|y|}(G^\pm_n)^{(s)}(\lambda|y|)\Big]\right| \lesssim \lambda^{-l}.
$$
Hence, for $0\le s \le \vartheta(j)-1$,  $|y|\ge 1$ , $\lambda\langle y\rangle\ge 1$ and $l=0,\ldots,\frac{n}{2}+1$,  Proposition \ref{corcancell} derives that
\begin{equation*}
\begin{aligned}
&\left\| \partial_\lambda^l\Big[\lambda^{\frac{n+1}{2}-2m+s} e^{\mp i\lambda|y|}\big(G^\pm_n\big)^{(s)}(\lambda|y|)\cdot Q_jv\Big(\frac{(|\cdot-y|-|y|)^s}{|\cdot-y|^{\frac{n-1}{2}}s!}\Big)(x)\Big]\right\|_{L_x^2}\\
\lesssim& \left|\partial_\lambda^l\left( \lambda^{\frac{n+1}{2}-2m+s} e^{\mp i\lambda|y|}\big(G^\pm_n\big)^{(s)}(\lambda|y|)\right)\right|\cdot \left\|Q_jv\Big(\frac{(|\cdot-y|-|y|)^s}{|\cdot-y|^{\frac{n-1}{2}}}\Big)(x) \right\|_{L_x^2}\\
\lesssim& \lambda^{\frac{n+1}{2}-2m+s-l}\langle y\rangle^{s-\frac{n-1}{2}-\vartheta(j)}\lesssim \lambda^{\frac{n+1}{2}-2m+\vartheta(j)-l}\langle y\rangle^{-\frac{n-1}{2}},
\end{aligned}
\end{equation*}
where in the last line, we apply the facts $\lambda \langle y\rangle \ge 1$ and $s<\vartheta(j)$.
Combining this with\eqref{eq3.903} and \eqref{eq3.905}, one can find that for $\lambda\langle y\rangle \ge 1$ and $\lambda\in(0,1)$,
\begin{equation}\label{eq3.906}
	\big\|\partial^l_\lambda \omega_{j}^{\pm}(\lambda,x,y)\big\|_{L_x^2}\lesssim \lambda^{\frac{n+1}{2}-2m+\vartheta(j)-l}\langle y\rangle^{-\frac{n-1}{2}}, \quad l=0,\ldots,\frac{n}{2}+1.
\end{equation}
Thus, we obtain \eqref{eq-wj-estimate-even-1}--\eqref{eq-wj-estimate-even-4} in the region $\lambda\langle y\rangle\ge 1$. 

Thus we complete the proof.
\end{proof}

\begin{proof}[{\bf Proof of Lemma \ref{lemma-QjvR-even-1-1}.}]
When $d=4m$,  define 
$$\omega_{j}^{\pm}(\lambda,x,y)= e^{\pm i\lambda|y|}Q_jvR_0^\pm(\lambda^{2m})V\big[R_0^\pm(\lambda^{2m})(\cdot,\,y)\big](x), $$ 
which implies \eqref{eq-QjR kernel-even-1-1}.
According to the expression \eqref{eq-free kernel-F-01},  the integral kernel $R_0^\pm(\lambda^{2m})V\big[R_0^\pm(\lambda^{2m})(\cdot,\,y)\big](x)$ can be expressed as a sum 
\begin{align*}
&e^{\pm i\lambda|y|} R_0^\pm(\lambda^{2m})V\big[R_0^\pm(\lambda^{2m})(\cdot,\,y)\big](x)\\
=&\int_{\R^{n}} e^{\pm i\lambda(|x-z|+|z-y|-|y|)}\left(\frac{ \mathbf{F}_{0,0}^\pm(\lambda,x,z,y)V(z)}{|x-z|^{\frac n2}|z-y|^{\frac n2}}+ \lambda^{\frac{1}{2}}\cdot \frac{\mathbf{F}_{0,1}^\pm(\lambda,x,z,y)\, V(z)\, }{|x-z|^{\frac n2}|z-y|^{\frac{n-1}{2}}} \right)dz \\
&+\int_{\R^{n}} e^{\pm i\lambda(|x-z|+|z-y|-|y|)}\left(\lambda^{\frac{1}{2}}\cdot \frac{ \mathbf{F}_{1,0}^\pm(\lambda,x,z,y)V(z)}{|x-z|^{\frac{n-1}{2}}|z-y|^{\frac n2}}+ \lambda\cdot  \frac{\mathbf{F}_{1,1}^\pm(\lambda,x,z,y)\, V(z)\, }{|x-z|^{\frac{n-1}{2}}|z-y|^{\frac{n-1}{2}}} \right)dz,
\end{align*}
where $\mathbf{F}_{i,j}^\pm(\lambda,x,z,y)$ is defined as
$$\mathbf{F}_{i,j}^\pm(\lambda,x,z,y)=F_{n,i}^\pm(\lambda|x-z|)F_{n,j}^\pm(\lambda|z-y|),\quad \text{for $i,\,j \in \{0,1\}$}. $$
Since  $F_{n,0}^\pm(z)$ and $F_{n,1}^\pm(z)$ satisfy \eqref{eq-F decay-odd}, one has that 
$\mathbf{F}_{i,j}^\pm(\lambda,x,z,y)$ obey the estimates as
$$\left|\partial_{\lambda}^l\mathbf{F}_{i,j}^\pm(\lambda,x,z,y)  \right|\lesssim \lambda^{-l},\quad l\in \N_0.$$
Thus, a direct computation shows that 
\begin{align*}
 \left|\partial_{\lambda}^l\left(e^{\pm i\lambda|y|} R_0^\pm(\lambda^{2m})V\big[R_0^\pm(\lambda^{2m})(\cdot,\,y)\big](x) \right)  \right| 
\lesssim& \lambda^{-l}\int_{\R^{n}}\sum_{i,j=0}^{1} \frac{\langle|x-z|+|z|\rangle^{\frac n2+1} |V(z)|}{|x-z|^{\frac {n-i}{2}}|z-y|^{\frac {n-j}{2}}} dz \\
\lesssim&\lambda^{-l}\int_{\R^{n}} \sum_{j=0}^{1} \frac{\langle x-z \rangle^{\frac{3}{2}}
\langle z\rangle^{\frac n2+1} |V(z)|}{|z-y|^{\frac {n-j}{2}}} dz \\
\lesssim& \langle x\rangle^{\frac{3}{2}} \langle y\rangle^{-\frac{n-1}{2}},
\end{align*} 
for $0<\lambda<1$ and $l=0,1,\ldots,\frac n2+1$.
Since $Q_j$ is $L^2$-bounded and $v(x)\langle x\rangle^{\frac{3}{2}}\in L^2$, we obtain the estimate \eqref{eq-wj-estimate-even-1-2}. 
\end{proof}

\begin{proof}[{\bf Proof of Lemma \ref{lemma-QjvR+-R-even-1}.}]
	We begin by introducing the cut-off function $\phi_0 \in C^{\infty}(\R)$ such that $\phi_0(\lambda)=1$ when $|\lambda|\leq 1/2$ and $\phi_0(\lambda)=0$ when $|\lambda|\geq 1$. Meanwhile, let $\phi_1(\lambda)=1-\phi_0$.
	Note that by the expansion \eqref{eq-free expansion-even-4} when $\theta=\vartheta(j)$ and the cancellation property of $Q_j$, one has 
	\begin{equation*}
		\Big(Q_jv(R_0^+(\lambda^{2m})(\cdot,\,y)-R_0^-(\lambda^{2m})(\cdot,\,y))\Big)(x)=\lambda^{n-2m}
		\Big(Q_jv(r_{\vartheta(j)}^+(\lambda|\cdot-y|)-r_{\vartheta(j)}^-(\lambda|\cdot-y|))\Big)(x).
	\end{equation*}
	Then, we rewrite 
	$$\phi_0(\lambda \langle y\rangle)Q_jv(r_{\vartheta(j)}^\pm (\lambda|\cdot-y|)=\pm \phi_0(\lambda \langle y\rangle) e^{\pm i\lambda|y|}\widetilde{\omega}_{j,0}^{\pm}(\lambda,x,y).$$
	Additionally, in even dimensions $n\le 4m-2$, we define
	$$\phi_1(\lambda \langle y\rangle)Q_jv(R_0^+(\lambda^{2m})(\cdot,\,y))(x)=\pm e^{\pm i\lambda|y|}\phi_1(\lambda \langle y\rangle){\omega}_{j}^{\pm}(\lambda,x,y)=e^{\pm i\lambda|y|}\widetilde{\omega}_{j,1}^{\pm}(\lambda,x,y),$$
	where ${\omega}_{j}^{\pm}(\lambda,x,y)$ is defined in Lemma \ref{lemma-QjvR-even-1}.
	And when $n=4m$, from \eqref{eq-free kernel-F-02}, we have 
	$$
	R_0^\pm(\lambda^{2m})(x,y)-R_0^\pm(\lambda^{2m})(x,y)=\lambda^{\frac{1}{2}}\frac{e^{ i\lambda|x-y|}}{|x-y|^{\frac{n-1}{2}}}\widetilde{F}_{n}^+(\lambda|x-y|)-\lambda^{\frac{1}{2}}\frac{e^{-i\lambda|x-y|}}{|x-y|^{\frac{n-1}{2}}}\widetilde{F}_{n}^-(\lambda|x-y|).
	$$
	When when $n=4m$, define 
	$$\widetilde{\omega}_{j,1}^{\pm}(\lambda,x,y)=\lambda^{\frac{1}{2}}Q_jv\left( \frac{e^{\pm i\lambda(|\cdot-y|-|y|)}}{|\cdot-y|^{\frac{n-1}{2}}}
	\widetilde{F}_{n}^-(\lambda|\cdot-y|)\right)(x). $$
    Let 
	$$\widetilde{\omega}_{j}^{\pm}(\lambda,x,y)=\widetilde{\omega}_{j,0}^{\pm}(\lambda,x,y)+\widetilde{\omega}_{j,1}^{\pm}(\lambda,x,y). $$
	Then, we have \eqref{eq3.900}. 
	By the proof of Lemma \ref{lemma-QjvR-even-1}, we only need to consider estimate of $\widetilde{\omega}_{j,0}^{\pm}(\lambda,x,y)$. Applying the same approach, one has that the inequalities \eqref{eq-k2-even-1} and \eqref{eq-k2-even-2} hold.
\end{proof}

\begin{proof}[{\bf Proof of Lemma \ref{lemma-QjvR-even-2}.}]
This lemma primarily aims to study the refined structure of  ${\omega}_{j}^{\pm}(\lambda,x,y)$, as defined in Lemma \ref{lemma-QjvR-even-1}, in the region where $\lambda\langle y\rangle\ge 1$. The proof is identical to that in odd dimensions, as detailed in Lemma 2.10 of \cite{CSWY24}, and is omitted here for brevity.
\end{proof}
\vskip0.2cm
\noindent{\bf Acknowledgments:} 
Avy Soffer is partially supported by NSF-DMS (No. 2205931) and the Simon’s Foundation (No. 395767).  
 Xiaohua Yao is partially supported by NSFC (No. 12171182 and  12531005).  Zhao Wu is partially supported by NSFC (No. 42450275).

%\noindent{\bf Acknowledgments:} 

%\nocite{*}
%\bibliographystyle{amsalpha}
%\bibliography{referencelib}

\begin{thebibliography}{BAKS00}

\bibitem{Agmon}
S.~Agmon, \emph{Spectral properties of {S}chr\"{o}dinger operators and
  scattering theory}, Ann. Scuola Norm. Sup. Pisa Cl. Sci. (4) \textbf{2}
  (1975), no.~2, 151--218. 
\bibitem{AS}  M. Abramowltz and I. Stegun, \emph{Handbook of Mathematical Functions}, National Bureau of Standards, Applied Mathematics Series .  

\bibitem{Beceanu}
M.~Beceanu, \emph{Structure of wave operators for a scaling-critiacal class of potentials}, Amer. J. Math. \textbf{136}(2014), no.~2 255-308.  

\bibitem{Beceanu-Schlag-20}
M.~Beceanu and W. Schlag, \emph{Structure formulas for wave operators}, Amer. J. Math. \textbf{142}(2020), no.~3, 751-807.  
  
\bibitem{CHHZ-2024}
{\sc H. Cheng, S. Huang, T. Huang and Q. Zheng}, {\em Pointwise estimates for the fundamental solutions of higher order Schrödinger equations in odd dimensions I: low dimensional case}, arXiv e-prints, (2024), arXiv:2401.04969v10.

\bibitem{CSWY24}
{\sc H. Cheng,  A. Soffer, Z. Wu and X. Yao}, \emph{The $L^p$-boundedness of wave operators for  Schr\"odinger type operator with zero singularities in low odd dimensions,} http://arxiv.org/abs/2505.07009.
 
\bibitem{AF06}
P.~D'Ancona and L.~Fanelli, \emph{$L^p$-boundedness of the wave operator for the one dimensional Schr\"odinger operator}, Commun. Math. Phys. \textbf{268}(2006), no.~2, 415--438.

\bibitem{Erdogan-Goldberg-Green-JFA-2018}
M.B.~Erdo\u{g}an, M.~Goldberg and W.R.~Green, \emph{On the $L^p$ boundedness of wave operators for two-dimensional Schr\"{o}dinger operators with threshold obstructions}, J. Funct. Anal. \textbf{274}(2018),  2139--2161.
  

\bibitem{EGG-JFA-2023}
M.B.~Erdo\u{g}an, M.~Goldberg and W.R.~Green, \emph{Counterexamples to $L^p$ boundedness of wave operators for classical and higher order Schr\"{o}dinger operators}, J. Funct. Anal. \textbf{285} (2023), no.~5, Paper No. 110008. 

\bibitem{EGG_ArXiv} M. B. Erdo\u{g}an, M. Goldberg and W. R. Green, {\it Dispersive estimates for higher order Schr\"odinger operators with scaling-critical potentials}. Math. Ann. (2025). https://doi.org/10.1007/s00208-025-03146-1.

\bibitem{EG22}
M.B.~Erdo\u{g}an and W.R.~Green \emph{The {$L^p$}-continuity of wave operators for higher order Schr\"{o}dinger operators}, Adv. Math. \textbf{404} (2022), Paper No.108450, 41. %\MR{4418888}

\bibitem{EG23}
M.B.~Erdo\u{g}an and W.R.~Green, \emph{The {$L^p$}-continuity of wave operators for higher order Schr\"{o}dinger operators}, J. Differential Equations \textbf{355} (2023), 144--161. %\MR{4542551}

\bibitem{EGL24}
M.B.~Erdo\u{g}an, W.R.~Green and K. Lamaster, \emph{{$L^p$}-continuity of wave operators for higher order Schr\"{o}dinger operators with threshold eigenvalues in high dimensions}, Discrete and Continuous Dynamical Systems,
 \textbf{45} (2025), no. 9,  3258--3274.

\bibitem{EGL25}
M.B.~Erdo\u{g}an, W.R.~Green and K. Lamaster, \emph{$L^p$-boundedness of wave operators for higher order schr\"{o}dinger operators with threshold eigenvalues}, arXiv e-prints, (2025), arXiv:2506.16378v1.

\bibitem{EGT21}
M.B.~Erdo\u{g}an, W.R.~Green and E.~Toprak, \emph{On the fourth order Schr\"{o}dinger equation in three dimensions: dispersive estimates and zero energy resonances}, J. Differential Equations \textbf{271} (2021), 152--185.

%\bibitem{Feng-2021}
%H.~Feng, \emph{Dispersive estimates for inhomogeneous fourth-order Schr\"{o}dinger operator in 3D with zero energy obstructions}, Nonlinear Anal. \textbf{207} (2021), Paper No. 112269, 33. 

\bibitem{FSWY20}
H.~Feng, A.~Soffer, Z.~Wu and X. Yao, \emph{Decay estimates for higher-order elliptic operators}, Trans. Amer. Math. Soc. \textbf{373} (2020), no.~4, 2805--2859. 

\bibitem{FSY}
H.~Feng, A.~Soffer and X.~Yao, \emph{Decay estimates and Strichartz estimates of fourth-order Schr\"{o}dinger operator}, J. Funct. Anal. \textbf{274} (2018), no.~2, 605--658. 

\bibitem{Finco-Yajima-2006-even}
D. Finco and K.~Yajima, \emph{The $L^p$ boundedness of wave operators for Schr\"{o}dinger operators with threshold singularities II. the even dimensional case}, J. Math. Sci. Univ. Tokyo  \textbf{13} (2006), no.~3, 277--346.

\bibitem{Gal-Yaj-00}
A.~Galtbayar and K.~Yajima, \emph{The $L^p$-continuity of wave operators for one dimensional Sch\"{o}dinger operators}, J. Math. Sci. Univ. Tokyo. \textbf{7} (2000), 221-240.


\bibitem{GY23}
A.~Galtbayar and K.~Yajima, \emph{The $L^p$-boundedness of wave operators for fourth order Sch\"{o}dinger operators on $\mathbf{R}^4$}, J. Spectr. Theory  \text{14} (2024), no.~1,  271–354.

\bibitem{GY25}
A.~Galtbayar and K.~Yajima, \emph{The $L^p$-boundedness of wave operators for 4-order Sch\"{o}dinger operators on $\mathbf{R}^2$, I}, arXiv e-prints, (2025), arXiv:2504.11753


\bibitem{Goldberg-Green-2016}
M.~Goldberg and W.R.~Green, \emph{The $L^p$ boundedness of wave operators for Schr\"odinger operators with threshold singularities}, Adv. Math. \textbf{303} (2016), 360--389.

\bibitem{Goldberg-Green-2017}
M. Goldberg and W.R.~Green, \emph{On the $L^p$ boundedness of wave operators for four-dimensional Schr\"odinger operators with a threshold eigenvalue}, Ann. Henri Poincar\'{e}, \textbf{18} (2017), 1269--1288.

\bibitem{Goldberg-Green-high-even}
M. Goldberg, W.R.~Green, \emph{Dispersive estimates for higher dimensional Schrödinger operators with threshold eigenvalues. II: The even dimensional case}. J. Spectr. Theory, \textbf{7} (2017),  pp. 33–86

\bibitem{GG21}
M. Goldberg and W.R.~Green, \emph{ $L^p$ boundedness of wave operators for fourth order Schr\"odinger operators}, Trans. Amer. Math. Soc. \textbf{374} (2021), 4075--4092.


\bibitem{Goldberg-Schlag}
M.~Goldberg and W.~Schlag, \emph{Dispersive estimates for {S}chr\"{o}dinger
  operators in dimensions one and three}, Comm. Math. Phys. \textbf{251}
  (2004), no.~1, 157--178.

\bibitem{GV}  M. Goldberg and M. Visan, \emph{A counterexample to dispersive estimates for Schr\"{o}dinger operators in higher dimensions,} Comm. Math. Phys. \textbf{266} (2006), 211-238.

\bibitem{gra} L. Grafakos, Classical Fourier Analysis.  Graduate Texts in Mathematics, second ed., vol. 249, New York: Springer, 2008.

\bibitem{Green-Toprak-4D}
W.~R.~Green and E.~Toprak, \emph{On the fourth order Schr\"{o}dinger equation in four dimensions: dispersive estimates and zero energy
resonances}, J. Differential Equations \textbf{267} (2019), no.~3, 1899--1954. 


%\bibitem {hormaner}  L. H\"ormander, The Analysis of Linear Partial Differential Operators, vol. II: Differential Operators 
%with Constant Coefficients, second edition, Springer, Berlin, 2005.

\bibitem{JK} A. Jensen, T. Kato, \emph{Spectral properties of Schr\"{o}dinger operators and timedecay of the wave functions,}  Duke Math. J. \textbf{46} (1979), no. 3, 583-611.

\bibitem{J80} A. Jensen, \emph{Spectral properties of Schr\"{o}dinger operators and time-decay of the wave functions results in $L^2(\mathbb{R}^m)$, $m\ge 5$}, Duke Math. J. \textbf{47} (1980), no. 1, 57-80.

\bibitem{JN01} A. Jensen, G. Nenciu, \emph{A unified approach to resolvent expansions at thresholds,} Rev. Math. Phys. \textbf{13} (2001), no. 6, 717-754.

\bibitem{Jensen-Yajima-02-CMP}
A.~Jensen and K.~Yajima, \emph{ A remark on $L^p-$boundedness of wave operators for two-dimensional Schr\"{o}dinger operators}, Comm. Math. Phys. \textbf{225} (2002), no.~3, 633-637.

%\bibitem{Jensen-Yajima-08-PLMS}
%A.~Jensen and K.~Yajima, \emph{On $L^p$ boundedness of wave operators for 4-dimensional Schr\"{o}dinger operators with threshold singularities}, Proc. Lond. Math. Soc. \textbf{96} (2008), no.~1, 136-162.

  \bibitem{JSS}
J.-L. Journ\'{e}, A.~Soffer, and C.~D. Sogge, \emph{Decay estimates for
  {S}chr\"{o}dinger operators}, Comm. Pure Appl. Math. \textbf{44} (1991),
  no.~5, 573--604.

\bibitem{Kuroda}
S.~T.~Kuroda, \emph{Scattering theory for differetial operators. I. Operator theory}, J. Math. Soc. Japan \textbf{25} (1973), no.~1,
  75-104.
  
%\bibitem {kenig-R} C. Kenig, A. Ruiz, C. D. Sogge, \emph{Uniform Sobolev inequalities and unique continuation for second
% order constant coe cient di erential operators,} Duke Math. J. \textbf{55} (1987), 329-347. 

\bibitem{LSY23}
P.~Li, A.~Soffer and X.~Yao, \emph{{Decay estimates for fourth-order Schr{\"o}dinger operators in dimension two}}, J. Funct. Anal. \textbf{284} (2023), no.~6, Paper No. 109816, 83 pp. 
 
\bibitem{MWY22}
H.~Mizutani, Z.~Wan and X.~Yao, \emph{$L^p$ boundedness of wave operators for fourth-order Schr\"odinger operators on the line}, Adv. Math. \textbf{451} (2024), Paper No.109806. 
\bibitem{MWY23-1}
H.~Mizutani, Z.~Wan and X.~Yao, \emph{$L^p$ boundedness of wave operators for fourth-order Schr\"odinger operators with zero resonance on $\mathbf{R}^3$},  J. Funct. Anal. \textbf{289} (2025), no.~8, 111013, 68 pp.

\bibitem{MWY23-2}
H.~Mizutani, Z.~Wan and X.~Yao, \emph{Counterexamples and weak $(1,1)$ estimates of wave operators for fourth-order Schr\"odinger operators in dimension three},  J. Spectral theory,  {\bf14} (2024), no. 4, 1409-1450.
  
%\bibitem{Miya} A. Miyachi, \emph{On some estimates for the wave equation in $L^p$ and $H^p$},  J. Fac. Sci. Univ. Tokyo \textbf{27} (1980), 331-354. 

\bibitem{Reed-Simon-III}
M.~Reed and B.~Simon, \emph{Methods of Modern Mathematical Physics III: Scatteing Theoty}, Academic Press, New York, NY, 1972.

\bibitem{RodSchl}
I. Rodnianski and W. Schlag, \emph{Time decay for solutions of
  {S}chr\"{o}dinger equations with rough and time-dependent potentials},
  Invent. Math. \textbf{155} (2004), no.~3, 451--513.

\bibitem{Schlag-CMP}
W.~Schlag, \emph{Dispersive estimates for {S}chr\"{o}dinger operators in
  dimension two}, Comm. Math. Phys. \textbf{257} (2005), no.~1, 87--117.


\bibitem{Schlag}
W.~Schlag, \emph{Dispersive estimates for {S}chr\"{o}dinger operators: a survey},
  Mathematical aspects of nonlinear dispersive equations, Ann. of Math. Stud.,
  vol. 163, Princeton Univ. Press, Princeton, NJ, 2007, pp.~255--285.


\bibitem{SWY22}
A.~Soffer, Z.~Wu and X.~Yao, \emph{Decay Estimates for Bi-Schr\"{o}dinger Operators in Dimension One}, Ann. Henri
Poincar\'{e}. \textbf{23} (2022), no.~8, 2683--2744.

\bibitem{Stein-1993}
E. M.~Stein and T. S.~Murphy, \emph{Harmonic Analysis: Real-Variable Methods, Orthogonality, and Oscillatory Integrals}, Princeton University Press, 1993.
  
\bibitem{Weder-99-CMP}
R.~Weder, \emph{The {$W_{k, p}$}-continuity of the Schr\"odinger wave operator on the line}, Comm. Math. Phys.  \textbf{208} (1999), no.~2, 507--520.

\bibitem{Weder1}R. A. Weder,
 \emph{$L^p-L^{p^\prime}$ estimates for the  Schr\"odinger equation
on the line and inverse scattering for the non linear Schr\"odinger equation
with a potential}, J. Funct. Analysis  \textbf{170}(2000), 37--68.


\bibitem{Yajima-JMSJ-95}
K.~Yajima, \emph{The {$W^{k,p}$}-continuity of wave operators for
  {S}chr\"{o}dinger operators}, J. Math. Soc. Japan \textbf{47} (1995), no.~3,
  551--581. 
  
\bibitem{Yajima-1995-even}
K.~Yajima, \emph{The {$W^{k,p}$}-continuity of wave operators for Schr\"{o}dinger operators.III.Even-dimensional cases $m\ge4$}, J. Math. Sci. Univ. Tokyo.  \textbf{2} (1995), no.~2, 311--346.
551--581.

\bibitem{Yajima-99-CMP}
K.~Yajima, \emph{$L^p-$boundedness of wave operators for two-dimensional Schr\"{o}dinger operators}, Comm. Math. Phys. \textbf{208} (1999), no.~3, 125-152.

  
\bibitem{Yajima-JMSJ-2006}
K.~Yajima, \emph{The $L^p$ boundedness of wave operators for Schr\"{o}dinger operators with threshold singularities I. the odd dimensional case}, J. Math. Soc. Japan. \textbf{13} (2006),
43--93. 

\bibitem{Yajima-2016}
K.~Yajima, \emph{Remark on the $L^p$-boundedness of wave operators for Schr\"{o}dinger operators with threshold singularities }, Doc. Math. \textbf{21} (2016),
391--443.

\bibitem{Yajima-2016-3d}
K.~Yajima, \emph{$L^1$ and $L^\infty$-boundedness of wave operators for three dimensional Schr\"odinger operators with threshold singularities}, Tokyo J. Math. \textbf{41} (2018), no. 2, 385-406.

\bibitem{Yajima-2022}
K.~Yajima, \emph{The $L^p$-boundedness of wave operators for two dimensional Schr\"{o}dinger operators with threshold singularities }, J. Math. Soc. Japan. \textbf{74} (2022), no.~4,
1169--1217. 

\bibitem{Yajima-2022-4}
K.~Yajima, \emph{The $L^p$-boundedness of wave operators for four dimensional Schr\"{o}dinger operators with threshold resonances}, The physics and mathematics of Elliott Lieb—the 90th anniversary. Vol. II, 517–563.
EMS Press, Berlin, 2022.

\bibitem{ZYF} Q. Zheng, X. Yao and D. Fan, \emph{Convex hypersurfaces and $L^p$ estimates for Schr\"{o}dinger equations}, J. Funct. Anal. \textbf{208} (2004), no. 1, 122-139.

\end{thebibliography}

\end{document}